\documentclass{amsart}
\usepackage[doi=false, isbn=false, url=false, style=alphabetic, backend=bibtex]{biblatex}
\AtEveryBibitem{\clearfield{MRCLASS}} 
\AtEveryBibitem{\clearfield{MRNUMBER}} 
\AtEveryBibitem{\clearfield{MRREVIEWER}}  
\usepackage{accents}
\usepackage{graphicx}
\usepackage{amssymb}
\usepackage{amsmath}
\usepackage{mathtools}
\makeatletter
\providecommand*{\twoheadrightarrowfill@}{%
  \arrowfill@\relbar\relbar\twoheadrightarrow
}
\providecommand*{\twoheadleftarrowfill@}{%
  \arrowfill@\twoheadleftarrow\relbar\relbar
}
\providecommand*{\xtwoheadrightarrow}[2][]{%
  \ext@arrow 0579\twoheadrightarrowfill@{#1}{#2}%
}
\providecommand*{\xtwoheadleftarrow}[2][]{%
  \ext@arrow 5097\twoheadleftarrowfill@{#1}{#2}%
}
\makeatother
\usepackage{amsfonts}
\usepackage{amscd}
\usepackage{url}
\usepackage{tikz}
\usepackage{calligra}
\usetikzlibrary{shapes}
\usetikzlibrary{decorations.markings}
\usepackage{tikz-cd}
\usetikzlibrary{arrows}
\usepackage{multirow}
\usepackage[savepos]{zref}
\usepackage{hyperref}
\newtheorem{thm}{Theorem}[section]
\newtheorem{corollary}[thm]{Corollary}
\newtheorem{lemma}[thm]{Lemma}
\newtheorem{proposition}[thm]{Proposition}

\newtheorem{example}[thm]{Example}

\theoremstyle{definition}
\newtheorem{definition}[thm]{Definition}
\newtheorem{conj}[thm]{Conjecture}
\theoremstyle{remark}
\newtheorem{remark}[thm]{Remark}
\newtheorem{assumption}[thm]{Assumption}
\numberwithin{equation}{section}
\numberwithin{figure}{section}

\usepackage[mathscr]{euscript}
\usepackage[T2A,T1]{fontenc}

\newenvironment{notedescription}%
  {\begin{description}%
    \setlength{\itemsep}{2.5pt}%
    \setlength{\parskip}{5pt}}%
  {\end{description}}

\newcommand{\form}{\upsilon}
\DeclareMathOperator{\GCD}{GCD}
\DeclareMathOperator{\Fr}{{Fr}}
\DeclareMathOperator{\can}{can}

\DeclareMathOperator{\OP}{op}

\DeclareMathOperator{\temp}{temp}
\DeclareMathOperator{\fin}{fin}

\DeclareMathOperator{\cusp}{cusp}
\DeclareMathOperator{\dis}{disc}
\DeclareMathOperator{\uni}{unit}

\DeclareMathOperator{\prof}{prof}
\DeclareMathOperator{\mult}{m}
\DeclareMathOperator{\HW}{HW}

\DeclareMathOperator{\Tw}{Tw}

\DeclareMathOperator{\SC}{sc}
\DeclareMathOperator{\hyp}{hyp}
\DeclareMathOperator{\der}{der}
\DeclareMathOperator{\reg}{reg}
\DeclareMathOperator{\creg}{creg}
\DeclareMathOperator{\greg}{greg}
\DeclareMathOperator{\Comm}{Comm}
\DeclareMathOperator{\Spl}{Spl}
\DeclareMathOperator{\Br}{Br}

\DeclareMathOperator{\cent}{cent}

\DeclareMathOperator{\Hom}{Hom}

\DeclareMathOperator{\Par}{Par}

\DeclareMathOperator{\Aut}{Aut}
\DeclareMathOperator{\Int}{Int}

\DeclareMathOperator{\ab}{ab}

\DeclareMathOperator{\Gal}{Gal}

\DeclareMathOperator{\Ker}{Ker}

\DeclareMathOperator{\val}{val}

\DeclareMathOperator{\Hilb}{Hilb}
\DeclareMathOperator{\Char}{char}
\DeclareMathOperator{\rec}{rec}

\DeclareMathOperator{\Spec}{Spec}

\DeclareMathOperator{\Id}{Id}
\DeclareMathOperator{\Res}{Res}

\DeclareMathOperator{\sgn}{sgn}

\DeclareMathOperator{\Tr}{Tr}
\DeclareMathOperator{\Norm}{N}
\DeclareMathOperator{\Kum}{\varkappa}
\DeclareMathOperator{\Obj}{Ob}

\DeclareMathOperator{\et}{\acute{e}t}

\newcommand{\From}{\colon}
\newcommand{\inar}{\ar@{^{(}->}}
\newcommand{\onar}{\ar@{->>}}

\renewcommand{\Im}{\mathrm{Im}}

\newcommand{\defined}[1]{\underline{{#1}}}
\usepackage{accents}
\newlength{\dtildeheight}
\newcommand{\dtilde}[1]{%
    \settoheight{\dtildeheight}{\ensuremath{\tilde{#1}}}%
    \addtolength{\dtildeheight}{-0.15ex}%
    \tilde{\vphantom{\rule{1pt}{\dtildeheight}}%
    \smash{\tilde{#1}}}}

\newcommand{\Irr}{{\boldsymbol{\Pi}}}
\newcommand{\WP}{{\boldsymbol{\Phi}}}
\newcommand{\Baer}{\dotplus}

\newcommand{\Weil}{\mathcal{W}}
\newcommand{\weil}{\mathbf{w}}

\newcommand{\Lang}{\mathcal{L}}
\newcommand{\Del}{\mathcal{D}}
\newcommand{\AB}{\mathcal{A}}
\renewcommand{\P}{\mathcal{P}}
\newcommand{\Inertia}{\mathcal{I}}
\newcommand{\Frob}{\Fr}

\newcommand{\limdir}{\varinjlim}

\newcommand{\Basis}{{\mathcal B}}

\makeatletter
\newcommand{\raisemath}[1]{\mathpalette{\raisem@th{#1}}}
\newcommand{\raisem@th}[3]{\raisebox{#1}{$#2#3$}}
\makeatother
\newcommand{\Vee}{ { \raisemath{-3pt}{\vee} } }

\newcommand{\EL}{{\mathsf L}}
\newcommand{\DEE}{{\mathsf D}}

\newcommand{\Cat}[1]{ {\mathsf{#1}} }
\newcommand{\Fun}[1]{ {\mathsf{#1}} }

\newcommand{\Lie}[1]{ {\mathfrak{#1}} }

\newcommand{\sch}[1]{\underline{\boldsymbol{ \mathrm{#1}}}}
\newcommand{\alg}[1]{\boldsymbol{\mathrm{#1}}}

\newcommand{\sheaf}[1]{{\mathscr{#1}}}
\newcommand{\shom}{\mathscr{H}\mathit{om}}
\newcommand{\stors}{\boldsymbol{\Cat{Tors}}}
\newcommand{\ssym}{\mathscr{S}\mathit{ym}}

\newcommand{\mGal}{\widetilde{\Gal}}
\newcommand{\amu}{\sch{\mu}}
\newcommand{\sAut}{\mathscr{A}\mathit{ut}}

\newcommand{\sspl}{\mathscr{S}\mathit{pl}}
\newcommand{\Whit}{\mathscr{W}\mathit{hit}}
\newcommand{\whit}{\boldsymbol{w}}
\newcommand{\dgp}[1]{\boldsymbol{\mathscr{#1}}}
\newcommand{\gerb}[1]{\boldsymbol{\Cat{#1}}}
\newcommand{\gerbob}[1]{{\boldsymbol{#1}}}

\newcommand{\Ext}{ {\mbox{Ext}} }

\newcommand{\ZZ}{\mathbb Z}
\newcommand{\QQ}{\mathbb Q}
\newcommand{\RR}{\mathbb R}

\newcommand{\CC}{\mathbb C}

\renewcommand{\AA}{\mathbb A}
\newcommand{\VV}{\mathcal V}

\newcommand{\kk}{\mathfrak{f}}

\renewcommand{\SS}{\mathcal S}

\newcommand{\hecke}{\mathcal H}
\newcommand{\ident}{\equiv}

\newcommand{\Sat}{\mathcal S}

\newcommand{\FF}{\mathbb F}
\newcommand{\OO}{\mathcal{O}}

\newcommand{\Into}{\hookrightarrow}
\newcommand{\Onto}{\twoheadrightarrow}
\newcommand{\To}{\rightarrow}

\newcommand{\A}{\mathsf{A}}

\newcommand{\isom}{\cong}
\newcommand{\half}{\tfrac{1}{2}}
\newcommand{\AF}{\mathcal{AF}}

\renewcommand{\th}{\text{th}}

\newcommand{\inarrow}{\arrow[hook]}
\newcommand{\onarrow}{\arrow[two heads]}

\makeatletter
\newcommand\@biprod[1]{%
  \vcenter{\hbox{\ooalign{$#1\prod$\cr$#1\coprod$\cr}}}}
\newcommand\biprod{\mathop{\mathpalette\@biprod\relax}\displaylimits}
\makeatother

\newcommand{\defeq}{:=}
\DeclareMathAlphabet{\mathcalligra}{T1}{calligra}{m}{n}
\DeclareMathOperator{\Zar}{Zar}

\begin{document}

\title{L-groups and parameters for covering groups}%
\author{Martin H. Weissman}%
\date{\today}

\address{Yale-NUS College, 6 College Ave East, \#B1-01, Singapore 138614}
\email{marty.weissman@yale-nus.edu.sg}%


\begin{abstract}
We incorporate covers of quasisplit reductive groups into the Langlands program, defining an L-group associated to such a cover.  We work with all covers that arise from extensions of quasisplit reductive groups by $\alg{K}_2$ -- the class studied by Brylinski and Deligne.  We use this L-group to parameterize genuine irreducible representations in many contexts, including covers of split tori, unramified representations, and discrete series for double covers of semisimple groups over $\RR$.  An appendix surveys torsors and gerbes on the \'etale site, as they are used in the construction of the L-group.
\end{abstract}

\maketitle

\tableofcontents

\section*{Introduction}

\subsection*{Constructions and conjectures}

Let $\alg{G}$ be a quasisplit reductive group over a local or global field $F$.  In \cite{B-D}, Brylinski and Deligne introduce objects called {\em central extensions of $\alg{G}$ by $\alg{K}_2$}, and they express hope that for ``a global field this will prove useful in the study of `metaplectic' automorphic forms.''  We pursue their vision in this paper, and elaborate below.

Let $n$ be a positive integer and let $\mu_n$ denote the group of $n^{\th}$ roots of unity in $F$.  Assume that $\mu_n$ has order $n$.  Let $\alg{G}'$ be a central extension of $\alg{G}$ by $\alg{K}_2$, in the sense of \cite{B-D}.  We call the pair $\alg{\tilde G} \defeq (\alg{G}', n)$ a ``degree $n$ cover'' of $\alg{G}$.  Fix a separable closure $\bar F / F$ and write $\Gal_F = \Gal(\bar F / F)$.  Fix an injective character $\epsilon \From \mu_n \Into \CC^\times$.  

Associated to $\alg{G}'$, Brylinski and Deligne associate three invariants, which we call $Q$, $\sheaf{D}$, and $f$.  The first is a Weyl- and Galois-invariant quadratic form on the cocharacter lattice of a maximal torus in $\alg{G}$.  The second is a central extension of this cocharacter lattice by $\sheaf{G}_m$ (in the category of sheaves of groups on $F_{\et}$).  The third is difficult to describe here, but it reflects the rigidity of central extensions of simply-connected semisimple groups.  

\textbf{Part 1} of this article defines the L-group ${}^\EL \tilde G$ of such a cover $\alg{\tilde G}$, from the three invariants $Q$, $\sheaf{D}$, and $f$ (as well as $n$, $\epsilon$, and $\bar F$).  It is an extension of $\Gal_F$ by $\tilde G^\vee$, where $\tilde G^\vee$ is a pinned complex reductive group, on which $\Gal_F$ acts by pinned automorphisms.  Unlike Langlands' L-group, ours does not come equipped with a distinguished splitting -- although a noncanonical splitting often exists.  Throughout the construction of ${}^\EL \tilde G$, the arithmetically-inclined reader may replace $\CC$ by any $\ZZ[1/n]$-algebra $\Omega$ endowed with $\epsilon \From \mu_n \Into \Omega^\times$, without running into much difficulty.  The dual group $\tilde G^\vee$ has been considered by other authors, and it appears in various forms in \cite{FinkelbergLysenko}, \cite{McNamara}, \cite{Reich}, and \cite{ABPTV}.  We tabulate the dual groups; each comes equipped with a 2-torsion element $\tau_Q(-1)$ in its center $\tilde Z^\vee$.

In our previous article \cite{MWCrelle}, we limited our attention to split reductive groups, and constructed an L-group by bludgening Hopf algebras with two ``twists.''  The construction here is more delicate, and more general.  The ``first twist'' of \cite{MWCrelle} is encoded here in the following way.  The quadratic Hilbert symbol may be used to define a canonical 2-cocycle, yielding an extension $\mu_2 \Into \mGal_F \Onto \Gal_F$ which we call the metaGalois group.  The metaGalois group may be of independent interest -- one might look for its representations in nature, e.g., in the \'etale cohomology of a variety over $\QQ(i)$ equipped with twisted descent data to $\QQ$.    Pushing out the metaGalois group via the central 2-torsion element in $\tilde G^\vee$ yields the first twist,
\begin{equation}
\label{FirstTwist}
\tilde Z^\vee \Into (\tau_Q)_\ast \mGal_F \Onto \Gal_F.
\end{equation}

The ``second twist'' of \cite{MWCrelle} provided the greatest challenge there and here.  There, it was defined by twisting the multiplication in a Hopf algebra.  After attempting many reformulations (e.g., a Tannakian approach), we found the gerbe $\gerb{E}_\epsilon(\alg{\tilde G})$ on $F_{\et}$, which applies to covers of quasisplit groups and serves as the second twist here.  It is a bit different from gerbes that typically arise in the Langlands program, and so we include an appendix with relevant background on torsors and gerbes.  The \'etale fundamental group of this gerbe provides the second twist,
\begin{equation}
\label{SecondTwist}
\tilde Z^\vee \Into \pi_1^{\et}(\gerb{E}_\epsilon(\alg{\tilde G})) \Onto \Gal_F.
\end{equation}

The Baer sum of \eqref{FirstTwist} and \eqref{SecondTwist} gives an extension of $\Gal_F$ by $\tilde Z^\vee$.  Pushing out to $\tilde G^\vee$ (respecting the $\Gal_F$-action throughout) yields the L-group
$$\tilde G^\vee \Into {}^\EL \tilde G \Onto \Gal_F,$$
of Part 1.  In addition to its construction, we verify that this L-group behaves well with respect to Levi subgroups, passage between global fields, local fields, and rings of integers therein, and that the L-group construction is functorial for a class of ``well-aligned'' homomorphisms.

The construction of the L-group allows us to consider Weil parameters.  When $F$ is a local or global field, we consider the set $\WP_\epsilon(\alg{\tilde G}/F)$ of $\tilde G^\vee$-orbits of Weil parameters $\Weil_F \To {}^\EL \tilde G$.  When $\OO$ is the ring of integers in a nonarchimedean local field, we consider the set $\WP_\epsilon(\alg{\tilde G} / \OO)$ of $\tilde G^\vee$-orbits of unramified Weil parameters $\Weil_F \To {}^\EL \tilde G$.  One could similarly define Weil-Deligne parameters, (conjectural) global Langlands parameters, etc.

In a set of unpublished notes, we constructed an L-group for split groups without using the gerbe discussed above.  This ``$E_1 + E_2$'' construction has been studied further by Wee Teck Gan and Gao Fan in \cite{GanGao} and \cite{GaoThesis}.  The construction of this paper agrees with the $E_1 + E_2$ construction for split tori, but we have not checked that the constructions agree for split reductive groups.  

With the construction of the L-group complete, we turn our attention to representation theory in \textbf{Part 2}.  The cover $\alg{\tilde G}$ and character $\epsilon$ allow us to define $\epsilon$-genuine irreducible representations of various sorts.  The set $\Irr_\epsilon(\alg{\tilde G} / \bullet)$ is defined in three contexts.
\begin{description}
\item[$F$ a local field]  Brylinski and Deligne construct \cite[\S 10.3]{B-D} a central extension $\mu_n \Into \tilde G \Onto G = \alg{G}(F)$, and we consider the set $\Irr_\epsilon(\alg{\tilde G}/F)$ of equivalence classes of irreducible admissible $\epsilon$-genuine representations of $\tilde G$.  
\item[$F$ a global field] Brylinski and Deligne construct \cite[\S 10.4]{B-D}  a central extension $\mu_n \Into \tilde G_\AA \Onto G_\AA = \alg{G}(\AA)$, canonically split over $\alg{G}(F)$, and we consider the set $\Irr_\epsilon(\alg{\tilde G}/F)$ of equivalence classes of $\epsilon$-genuine automorphic representations of $\tilde G_\AA$.
\item[$\OO$ the integers in a nonarchimedean local field $F$] Brylinski and Deligne construct \cite[\S 10.3, 10.7]{B-D} a central extension $\mu_n \Into \tilde G \Onto G = \alg{G}(F)$, canonically split over $G^\circ = \alg{G}(\OO)$.  We consider the set $\Irr_\epsilon(\alg{\tilde G}/\OO)$ of equivalence classes of irreducible $G^\circ$-spherical representations of $\tilde G$.
\end{description}

Part 2 introduces these classes of representations, reviewing or adapting foundational results as needed.  These include old results, such as the basic theory of admissible, unitary, and tempered representations, and results which are recent for covering groups, such as the Satake isomorphism and Langlands classification.  Some new features arise for covering groups:  we introduce the notion of ``central core character'' which is a bit coarser than ``central character.''  We place irreducible representations into ``pouches'' which should be subsets of L-packets in what follows later.    

The remainder of the paper is devoted to supporting the following ``Local Langlands Conjecture for Covers'' (LLCC), an analogue of the local Langlands conjectures (LLC).  For the LLC, we refer to the excellent survey by Cogdell \cite{Cogdell}.\begin{conj}[LLCC]
\label{ConjectureLang}
When $F$ is a local field, there is a \textbf{natural} finite-to-one parameterization,
$$
\Lang_\epsilon \From \Irr_\epsilon(\alg{\tilde G} / F) \To \WP_\epsilon(\alg{\tilde G} / F).
$$
\end{conj}
Such a conjecture is nearly meaningless, without defining the adjective ``natural.''  Naturality in the (traditional) local Langlands conjectures (LLC) includes compatibility with the bijective parameterizations for split tori (class field theory) and unramified representations (the Satake isomorphism), and with parabolic induction (the Langlands classification).  In \cite[\S 10]{BorelCorvallis}, it is suggested that naturality includes desiderata which specify how central characters and twisting by characters should correspond to various properties of and operations on Weil parameters.  One could add to these desiderata today, specifying for example how the formal degree (of discrete series) should correspond to adjoint $\gamma$-factors for Weil parameters (see \cite{HII}), or how the contragredient operation should correspond to the Chevalley involution (see \cite{AdamsVogan}).

For covers, we can make a similar list of desiderata for the LLCC.  We expect a finite-to-one parameterization for covering groups, compatible with our results for split tori (described in Part 3), with the unramified case (described in Part 4), with parabolic induction (via the Langlands classification for covers) and central core characters and twisting by characters (described in Part 2), and with formal degrees and adjoint $\gamma$-factors.

Unlike the LLC, we have not attempted to characterize the image of our conjectural parameterization in this paper, i.e., we have not identified the ``relevant'' parameters for covering groups.  The cases of split tori and discrete series for real groups  should suggest a characterization in the future.

\textbf{Part 3} focuses on the case of ``sharp'' covers of split tori.  For such a sharp cover (over a local or global field, or in the unramified setting), we define a {\em bijective} parameterization
$$\Lang_\epsilon  \From \Irr_\epsilon(\alg{\tilde T} / \bullet) \To \WP_\epsilon(\alg{\tilde T} / \bullet).$$  
This parameterization is natural for pullbacks of covers via isomorphisms, for isomorphisms of covers of a given split torus, and for Baer sums of covers.  This case constrains and guides many others, and it occupies the largest part of this article.    The ``sharp'' case quickly leads to the general case of split tori, where the parameterization may no longer be surjective.

\textbf{Part 4} includes three more cases where a precise parameterization is possible.  First is the spherical/unramified case.  When $\alg{\tilde G}$ is a cover of a quasisplit group over $\OO$, we define a {\em bijective} parameterization
$$\Lang_\epsilon \From \Irr_\epsilon(\alg{\tilde G} / \OO) \To \WP_\epsilon(\alg{\tilde G} / \OO).$$  
This parameterization follows from the Satake isomorphism (for covering groups), the parameterization for sharp covers of split tori above, and careful tracking of the Weyl group actions.

In \cite{GanGao}, Gan and Gao have verified the LLCC for split tori and unramified representations of split reductive groups, in the context of the L-group given by the $E_1 + E_2$ construction.

Second, we consider an anisotropic torus $\alg{T}$ over $\RR$ (i.e., $T = \alg{T}(\RR)$ is compact) and a sharp double cover $\alg{\tilde T}$ of $\alg{T}$.  In this case, we define a {\em bijective} parameterization
$$\Lang_\epsilon \From \Irr_\epsilon(\alg{\tilde T} / \RR) \To \WP_\epsilon(\alg{\tilde T} / \RR).$$ 
This case -- when $\tilde T$ is a torus in the oldest-fashioned sense and the representation theory is dead-simple -- is interesting from the standpoint of parameterization.  It validates our choice of gerbe, since the equivalence class of the gerbe must be exactly right to correspond to the correct lattice of characters.

Third, when $\alg{\tilde G}$ is a double cover of a semisimple quasisplit group over $\RR$, we define a finite-to-one parameterization of discrete series, 
$$\Lang_\epsilon \From \Irr_\epsilon^{\dis}(\alg{\tilde G} / \RR) \To \WP_\epsilon^{\dis}(\alg{\tilde G} / \RR).$$
This is based on the Harish-Chandra parameterization of discrete series, our parameterization for anisotropic tori, and careful tracking of involutions in the Weyl group.

\subsection*{Further questions}

As in the original Langlands conjectures, Conjecture \ref{ConjectureLang} suggests directions for further investigation.  Here are a few examples, within reach.
\begin{enumerate}
\item
(Inspired by recent work of Wee Teck Gan and Fan Gao \cite{GanGao} on $PGL_2$).  When $\alg{\tilde G}$ is a cover of a split group with trivial first invariant ($Q = 0$), one can find a z-extension (in the sense of Kottwitz \cite{Kot}) $\alg{H} \To \alg{G}$ for which the pullback cover $\alg{\tilde H}$ is isomorphic to the trivial cover.  This identifies genuine representations of $\tilde G$ with ordinary representations of $H$ satisfying a constraint on the central character.  This identification should be reflected on the side of Weil parameters, and the LLCC for such covers $\alg{\tilde G}$ should relate to the LLC for $\alg{H}$.
\item
(Inspired by a conversation with Wee Teck Gan).  When $\alg{G}$ is a group over $\CC$, there are nontrivial covers $\alg{\tilde G}$, but the resulting covers of complex Lie groups split canonically.  In this way, the genuine representations of $\tilde G$ correspond to ordinary representations of $G$.  This should be reflected on the side of Weil parameters, and the LLCC for $\alg{\tilde G}$ should relate to the LLC for $\alg{G}$.
\item
(Inspired by the work of Adams and Vogan \cite{AdamsVogan}).  When $\alg{\tilde G}$ is a cover of $\alg{G}$, a group over a local field $F$, one may define an inverse cover $\alg{\tilde G^{\OP}}$ with respect to the Baer sum.  If $\pi$ is an $\epsilon$-genuine representation of $\tilde G$, then the contragredient representation of $\pi$ is naturally an $\epsilon$-genuine representation of the inverse cover $\tilde G^{\OP}$.  The dual groups of $\alg{\tilde G}$ and $\alg{\tilde G^{\OP}}$ are the same, but the L-groups may not be.  The contragredient should be reflected on the side of Weil parameters, extending the Chevalley involution of the dual group to a map connecting the L-group of $\alg{\tilde G}$ with the L-group of the inverse cover $\alg{\tilde G^{\OP}}$.
\end{enumerate}

Much broader investigations are possible as well, given the L-group constructed in this paper.  Even restricting to the local case, examples include:  the study of pure inner forms (strong inner forms for real groups) and ``stability'' for covering groups; endoscopy for covering groups; the completion of the local Langlands conjectures for covers of reductive groups over $\RR$; base change for covering groups; Langlands-Shahidi L-functions for covering groups (initiated by Fan Gao in \cite{GaoThesis} and Szpruch \cite{Szpruch}); the parameterization of Iwahori-spherical and depth-zero representations; etc.

\subsection*{Philosophies}

A few principles are helpful when considering any putative Langlands program for covering groups.

\begin{enumerate}
\item
There is no $\epsilon$-genuine trivial (or Steinberg) representation for general covers, and so one should not expect a single distinguished splitting of the L-group.
\item
If some set of things is parameterized by cohomology in degree $2$, then that set of things should be viewed as the set of objects in a 2-category.
\item
Things which ``are trivial'' (e.g., extensions, gerbes) can be isomorphic to trivial things in interesting ways.
\end{enumerate}

\subsection*{Acknowledgments}

The ideas of this paper have evolved over the past few years, and I am very grateful for the numerous mathematicians who discussed the constructions in various stages of completeness and correctness.  The American Institute of Mathematics hosted a conference at which I spoke with Wee Teck Gan, Gordan Savin, Jeffrey Adams, Sergey Lysenko, Tamotsu Ikeda, Kaoru Hiraga, Tasho Kaletha and others.  Their previous work, and our discussions at the conference and elsewhere, have been very helpful.  I also appreciate the support of Harvard University during a short visit, where I gained insight from discussions with Dennis Gaitsgory, Dick Gross, and John Tate.  During a visit to the University of Michigan, I gained from feedback from Stephen DeBacker, and learned about gerbes from James Milne.  I appreciate the hospitality of the University of California, Berkeley, where I finished this paper during a summer visit.

In \cite{GanGao} and \cite{GaoThesis}, Wee Teck Gan and Fan Gao have tested some of the conjectures of this paper, and they have gone further in developing the Langlands program for covering groups.  I have greatly appreciated our frequent conversations.  Their results provided constraints which kept the constructions of this paper on track. 

Pierre Deligne has kindly corresponded with me in a series of letters, and his generosity on these occasions has been incredibly helpful.  His ideas led me to a deeper understanding of the crucial questions, and his correspondence motivated me to pursue this project further.  

\newpage
\section*{Notation}
\begin{notedescription}
\item[$F$]  A field, typically local or global.
\item[$\bar F$]  A separable closure of $F$.
\item[$\OO$] The ring of integers in $F$, in the nonarchimedean local case.
\item[$\AA$]  The ring of adeles of $F$, in the global case.
\item[${\Fr}$]  The geometric Frobenius automorphism.
\item[$S$]  A connected scheme, typically $\Spec(F)$ or $\Spec(\OO)$.
\item[$\bar s$]  The geometric point of $S$ corresponding to $\bar F$.
\item[$\Gal_S$]  The absolute Galois group $\pi_1^{\et}(S, \bar s)$.

\item[$\alg{X}$]  An algebraic variety over $S$, or sheaf on $S_{\Zar}$.
\item[$\alg{G}_m$]  The multiplicative group.
\item[$X$ or $X_F$]  The $F$-points $\alg{X}(F)$ for such a variety.

\item[$\sch{X}$]  A scheme over $\ZZ$.
\item[$\amu_n$] The group scheme over $\ZZ$ of $n^{\th}$ roots of unity.

\item[$\mu_n$]  The group $\amu_n(S)$, assumed to be cyclic of order $n$.

\item[$\sheaf{S}$]  A sheaf on $S_{\et}$.
\item[$\sheaf{S}{[U]}$] The sections of $\sheaf{S}$ over $U$ ($U \To S$ \'etale).
\item[$\dgp{G}$]  A local system on $S_{\et}$, of group schemes over $\ZZ$.
\item[$\Cat{C}$]  A category, with objects $\Obj(\Cat{C})$.
\item[$\gerb{E}$] A gerbe on $S_{\et}$.
\item[$\gerb{E}{[U]}$] The groupoid of sections of $\gerb{E}$ over $U$.

\item[$A$]  An abelian group.
\item[$A_{[n]}$]  Its $n$-torsion subgroup.
\item[$A_{/n}$]  The quotient $A / n A$.

\end{notedescription}

\newpage


\part{Covering groups and their L-groups}

\section{Covering groups}

Throughout this article, $S$ will be a scheme in one of the following two classes:  $S = \Spec(F)$ for a field $F$, or $S = \Spec(\OO)$ for a discrete valuation ring $\OO$ with fraction field $F$.  In the latter case, we assume that $\OO$ contains a field, or that $\OO$ has finite residue field.  We will often fix a positive integer $n$, and we will assume that $\mu_n = \amu_n(S)$ is a cyclic group of order $n$.  In Section \ref{MetaGaloisSection}, we will place further restrictions on $S$.

\subsection{Reductive groups}

Let $\alg{G}$ be a reductive group over $S$.  We follow \cite{SGA3} in our conventions, so this means that $\alg{G}$ is a smooth group scheme over $S$ such that $\alg{G}_{\bar s}$ is a connected reductive group for all geometric points $\bar s$ of $S$.  Assume moreover that $\alg{G}$ is \textbf{quasisplit} over $S$.  

Let $\alg{A}$ be a maximal $S$-split torus in $\alg{G}$, and let $\alg{T}$ be the centralizer of $\alg{A}$ in $\alg{G}$.  Then $\alg{T}$ is a maximal torus in $\alg{G}$, and we say that $\alg{T}$ is a \defined{maximally split} maximal torus.  Let $\sheaf{X}$ and $\sheaf{Y}$ be the local systems (on $S_{\et}$) of characters and cocharacters of $\alg{T}$.  

Let $\alg{N}$ be the normalizer of $\alg{T}$ in $\alg{G}$.  Let $\sheaf{W}$ denote the Weyl group of the pair $(\alg{G}, \alg{T})$, viewed as a sheaf on $S_{\et}$ of finite groups.  Then $\sheaf{W}[S] = \alg{N}(S) / \alg{T}(S)$ (see \cite[Expos\'e XXVI, 7.1]{SGA3}).  Let $\alg{B}$ be a Borel subgroup of $\alg{G}$ containing $\alg{T}$, defined over $S$.  Let $\alg{U}$ be the unipotent radical of $\alg{B}$.
\begin{proposition}
\label{BTConj}
Assume as above that $\alg{G}$ is quasisplit, and $S$ is the spectrum of a field or of a DVR.  The group $\alg{G}(S)$ acts transitively, by conjugation, on the set of pairs $(\alg{B}, \alg{T})$ consisting of a Borel subgroup (defined over $S$) and a maximally split maximal torus therein.
\end{proposition}
\proof
As we work over a local base scheme $S$, \cite[Expos\'e XXVI, Proposition 6.16]{SGA3} states that the group $\alg{G}(S)$ acts transitively on the set of maximal split subtori of $\alg{G}$ (defined over $S$).

Every maximally split maximal torus of $\alg{G}$ is the centralizer of such a maximal split torus, and thus $\alg{G}(S)$ acts transitively on the set of maximally split maximal tori in $ \alg{G}$.  The stabilizer of such a maximally split maximal torus $\alg{T}$ is the normalizer $\alg{N}(S)$.  The Weyl group $\sheaf{W}[S] = \alg{N}(S) / \alg{T}(S)$ acts simply-transitively on the minimal parabolic subgroups containing $\alg{T}$ by \cite[Expos\'e XXVI, Proposition 7.2]{SGA3}. This proves the proposition.  
\qed

The roots and coroots (for the adjoint action of $\alg{T}$ on the Lie algebra of $\alg{G}$) form local systems $\Phi$ and $\Phi^\vee$ on $S_{\et}$, contained in $\sheaf{X}$ and $\sheaf{Y}$, respectively.  The simple roots (with respect to the Borel subgroup $\alg{B}$) and their coroots form local systems of subsets $\Delta \subset \Phi$ and $\Delta^\vee \subset \Phi^\vee$, respectively.  In this way we find a local system on $S_{\et}$ of based root data (cf. \cite[\S 1.2]{BorelCorvallis}),
$$\Psi = \left( \sheaf{X}, \Phi, \Delta, \sheaf{Y}, \Phi^\vee, \Delta^\vee \right).$$
Write $\sheaf{Y}^{\SC}$ for the subgroup of $\sheaf{Y}$ spanned by the coroots.  

\subsection{Covers}

In \cite{B-D}, Brylinski and Deligne study central extensions of $\alg{G}$ by $\alg{K}_2$,  where $\alg{G}$ and $\alg{K}_2$ are viewed as sheaves of groups on the big Zariski site $S_{\Zar}$.  These extensions form a category we call $\Cat{CExt}_S(\alg{G}, \alg{K}_2)$.  Such a central extension will be written $\alg{K}_2 \Into \alg{G}' \Onto \alg{G}$ in what follows.  We add one more piece of data in the definition below.
\begin{definition}
\label{DefCover}
A degree $n$ \defined{cover} of $\alg{G}$ over $S$ is a pair $\alg{\tilde G} = (\alg{G}', n)$, where
\begin{enumerate}
\item
$\alg{K}_2 \Into \alg{G}' \Onto \alg{G}$ is a central extension of $\alg{G}$ by $\alg{K}_2$ on $S_{\Zar}$;
\item
$n$ is a positive integer;
\item
For all scheme-theoretic points $s \in S$, with residue field $\kk(s)$, $\# \amu_n( \kk(s) ) = n$.
\end{enumerate}
\end{definition}

Define $\Cat{Cov}_n(\alg{G})$ (or $\Cat{Cov}_{n/S}(\alg{G})$ to avoid confusion) to be the category of degree $n$ covers of $\alg{G}$ over $S$.  The objects are pairs $\alg{\tilde G} = (\alg{G}', n)$ as above, and morphisms are those from $\Cat{CExt}_S(\alg{G}, \alg{K}_2)$ (with $n$ fixed).

If $\gamma \From S_0 \To S$ is a morphism of schemes, then pulling back gives a functor $\gamma^\ast \From \Cat{Cov}_{n/S}(\alg{G}) \To \Cat{Cov}_{n/S_0}(\alg{G}_{S_0})$.  Indeed, a morphism of schemes gives inclusions of residue fields (in the opposite direction) and so Condition (3) is satisfied by the scheme $S_0$ when it is satisfied by the scheme $S$.

Central extensions $\alg{K}_2 \Into \alg{G}' \Onto \alg{G}$ are classified by a triple of invariants $(Q, \sheaf{D}, f)$.  For fields, this is carried out in \cite{B-D}, and the extension to DVRs (with finite residue field or containing a field) is found in \cite{MWIntegral}.  The first invariant $Q \From \sheaf{Y} \To \ZZ$ is a Galois-invariant Weyl-invariant quadratic form, i.e., $Q \in H_{\et}^0(S, \ssym^2(\sheaf{X})^{\sheaf{W}})$.  The second invariant $\sheaf{D}$ is a central extension of sheaves of groups on $S_{\et}$, $\sheaf{G}_m \Into \sheaf{D} \Onto \sheaf{Y}$.  The third invariant $f$ will be discussed later.  

A cover $\alg{\tilde G}$ yields a symmetric $\ZZ$-bilinear form $\beta_Q \From \sheaf{Y} \otimes_\ZZ \sheaf{Y} \To n^{-1} \ZZ$,
$$\beta_Q(y_1, y_2) \defeq n^{-1} \cdot \left( Q(y_1+y_2) - Q(y_1) - Q(y_2) \right).$$
This defines a local system $\sheaf{Y}_{Q,n} \subset \sheaf{Y}$,
$$\sheaf{Y}_{Q,n} = \{ y \in \sheaf{Y} : \beta_Q(y, y') \in \ZZ \text{ for all } y' \in \sheaf{Y} \}.$$

The category of covers is equipped with the structure of a Picard groupoid; one may ``add'' covers via the Baer sum.  If $\alg{\tilde G}_1, \alg{\tilde G}_2$ are two covers of $\alg{G}$ of degree $n$, one obtains a cover $\alg{\tilde G}_1 \Baer \alg{\tilde G}_2 = (\alg{G}_1' \Baer \alg{G}_2', n)$.

When $\alg{\tilde G} = (\alg{G}', n)$ is a degree $n$ cover of $\alg{G}$, and $\alg{H} \To \alg{G}$ is a homomorphism of groups over $S$, write $\alg{\tilde H} = (\alg{H}', n)$ for the cover of $\alg{H}$ resulting from pulling back extensions by $\alg{K}_2$.

In three arithmetic contexts, a cover $\alg{\tilde G}$ yields a central extension of topological groups according to \cite[\S 10.3, 10.4]{B-D}.
\begin{description}
\item[Global] If $S = \Spec(F)$ for a global field $F$, then $\alg{\tilde G}$ yields a central extension $\mu_n \Into \tilde G_\AA \Onto G_\AA$, endowed with a splitting $\sigma_F \From G_F \Into \tilde G_\AA$.
\item[Local] If $S = \Spec(F)$ for a local field $F$, then $\alg{\tilde G}$ yields a central extension $\mu_n \Into \tilde G \Onto G$, where $G = \alg{G}(F)$.
\item[Local integral]  If $S = \Spec(\OO)$, with $\OO$ the ring of integers in a nonarchimedean local field $F$, then $\alg{\tilde G}$ yields a central extension $\mu_n \Into \tilde G \Onto G$, where $G = \alg{G}(F)$, endowed with a splitting $\sigma^\circ \From G^\circ \Into G$.
\end{description} 

Fix an injective character $\epsilon \From \mu_n \Into \CC^\times$.  This allows one to define $\epsilon$-genuine automorphic representations of $\tilde G_\AA$ in the global context, $\epsilon$-genuine admissible representations of $\tilde G$ in the local context, and $\epsilon$-genuine $G^\circ$-spherical representations of $\tilde G$ in the local integral context.

The purpose of this article is the construction an \defined{L-group} associated to such a $\alg{\tilde G}$ and $\epsilon$.  We believe that this L-group will provide a parameterization of irreducible $\epsilon$-genuine representations in the three contexts above.

\subsection{Well-aligned homomorphisms}

Let $\alg{G}_1 \supset \alg{B}_1 \supset \alg{T}_1$ and $\alg{G}_2 \supset \alg{B}_2 \supset \alg{T}_2$ be quasisplit groups over $S$, endowed with Borel subgroups and maximally split maximal tori.    Let $\alg{\tilde G}_1 = (\alg{G}_1', n)$ and $\alg{\tilde G}_2 = (\alg{G}_2', n)$ be covers (of the same degree) of $\alg{G}_1$ and $\alg{G}_2$, respectively.  Write $\sheaf{Y}_1$ and $\sheaf{Y}_2$ for the cocharacter lattices of $\alg{T}_1$ and $\alg{T}_2$, and $Q_1, Q_2$ for the quadratic forms arising from the covers.  These quadratic forms yield sublattices $\sheaf{Y}_{1,Q_1,n}$ and $\sheaf{Y}_{2,Q_2,n}$.  

\begin{definition}
\label{DefWA}
A \defined{well-aligned homomorphism} $\tilde \iota$ from $(\alg{\tilde G}_1, \alg{B}_1, \alg{T}_1)$ to $(\alg{\tilde G}_2, \alg{B}_2, \alg{T}_2)$ is a pair $\tilde \iota = (\iota, \iota')$ of homomorphisms of sheaves of groups on $S_{\Zar}$, making the following diagram commute,
\begin{equation}
\label{CDiota}
\begin{tikzcd}
\alg{K}_2 \inarrow{r} \arrow{d}{=} & \alg{G}_1' \onarrow{r} \arrow{d}{\iota'} & \alg{G}_1 \arrow{d}{\iota} \\
\alg{K}_2 \inarrow{r} & \alg{G}_2' \onarrow{r} & \alg{G}_2
\end{tikzcd}
\end{equation}
and satisfying the following additional axioms:
\begin{enumerate}
\item
$\iota$ has normal image and smooth central kernel;
\item
$\iota(\alg{B}_1) \subset \alg{B}_2$ and $\iota(\alg{T}_1) \subset \alg{T}_2$.  Thus $\iota$ induces a map $\iota \From \sheaf{Y}_1 \To \sheaf{Y}_2$;
\item
$(\iota, \iota')$ realizes $\alg{G}_1'$ as the pullback of $\alg{G}_2'$ via $\iota$;
\item
The homomorphism $\iota$ satisfies $\iota(\sheaf{Y}_{1,Q_1,n}) \subset \sheaf{Y}_{2,Q_2,n}$.
\end{enumerate}
\end{definition}

\begin{remark}
Conditions (1) and (2) are inspired by \cite[\S 1.4, 2.1,2.5]{BorelCorvallis}, though more restrictive.  By ``normal image,'' we mean that for any geometric point $\bar s \To S$, the homomorphism $\iota \From \alg{G}_{1,\bar s} \To \alg{G}_{2, \bar s}$ has normal image.  Condition (3) implies that for all $y \in \sheaf{Y}_1$, $Q_1(y) = Q_2(\iota(y))$.  In other words, $Q_1$ is the image of $Q_2$ via the map
$$\iota^\ast \From H_{\et}^0(S, \ssym^2 (\sheaf{X}_2) ) \To H_{\et}^0(S, \ssym^2 (\sheaf{X}_1) ).$$
But Condition (3) does not imply Condition (4); one may cook up an example with $\alg{G}_1 = \alg{G}_{\mult}$ and $\alg{G}_2 = \alg{G}_{\mult}^2$ which satisfies (3) but not (4).
\end{remark}

\begin{proposition}
\label{ComposeWellaligned}
The composition of well-aligned homomorphisms is well-aligned.
\end{proposition}
\proof
Suppose that $(\iota_1, \iota_1')$ and $(\iota_2, \iota_2')$ are well-aligned homomorphisms, with $\iota_1 \From \alg{G}_1 \To \alg{G}_2$ and $\iota_2 \From \alg{G}_2 \To \alg{G}_3$.  Conditions (2), (3), and (4) are obviously satisfied by the composition $(\iota_2 \circ \iota_1, \iota_2' \circ \iota_1')$.  For condition (1), notice that the kernel of $\iota_2 \circ \iota_1$ is contained in the kernel of $\iota_1$, and hence is central.  The only thing left is to verify that $\iota_2 \circ \iota_1$ has normal image.  This may be checked by looking at geometric fibres; it seems well-known (cf. \cite[\S 1.8]{KotSTF}).
\qed

Inner automorphisms are well-aligned homomorphisms.
\begin{example}
Suppose that $\alg{\tilde G}$ is a degree $n$ cover of a quasisplit group $\alg{G}$.  Suppose that $\alg{B}_0 \supset \alg{T}_0$ and $\alg{B} \supset \alg{T}$ are two Borel subgroups containing maximally split maximal tori.  Suppose that $g \in \alg{G}(S)$, and write $\Int(g)$ for the resulting inner automorphism of $\alg{G}$.  As noted in \cite[0.N.4]{B-D}, $\Int(g)$ lifts canonically to an automorphism $\Int(g)' \in \Aut(\alg{G}')$.  If $\alg{B} = \Int(g) \alg{B}_0$ and $\alg{T} = \Int(g) \alg{T}_0$, then the pair $\left( \Int(g), \Int(g)' \right)$ is a well-aligned homomorphism from $(\alg{\tilde G}, \alg{B}_0, \alg{T}_0)$ to $(\alg{\tilde G}, \alg{B}, \alg{T})$.
\end{example}

While we focus on quasisplit groups in this article, the lifting of inner automorphisms allows one to consider ``pure inner forms'' of covers over a field.  
\begin{definition}
Let $\alg{\tilde G} = (\alg{G}', n)$ be a degree $n$ cover of a quasisplit group $\alg{G}$, over a field $F$.  Let $\xi \in Z_{\et}^1(F, \alg{G})$ be a 1-cocycle.  The image $\Int(\xi)$ in $Z_{\et}^1(F, \alg{Aut}(\alg{G}))$ defines an inner form $\alg{G}_\xi$ of $\alg{G}$.  These are called the \defined{pure inner forms} of $\alg{G}$.  On the other hand, we may consider the image $\Int(\xi)'$ in $Z_{\et}^1(F, \alg{Aut}(\alg{G}'))$, which by \cite[\S 7.1, 7.2]{B-D} defines a central extension $\alg{G}_\xi'$ of $\alg{G}_\xi$ by $\alg{K}_2$.  The cover $\alg{\tilde G}_\xi = (\alg{G}_\xi', n)$ of $\alg{G}_\xi$ will be called a \defined{pure inner form} of the cover $\alg{\tilde G}$.  
\end{definition}  

We have not attempted to go further in the study of inner forms for covers, but presumably one should study something like strong real forms as in \cite[Definition 1.12]{AdamsBarbaschVogan}, and more general rigid forms as in \cite{Kaletha}, if one wishes to assemble L-packets for covering groups.

The next example of a well-aligned homomorphism is relevant for the study of central characters of genuine representations.
\begin{example}
Let $\alg{\tilde G}$ be a degree $n$ cover of a quasisplit group $\alg{G} \supset \alg{B} \supset \alg{T}$.  Let $\alg{H}$ be the maximal torus in the center of $\alg{G}$, with cocharacter lattice $\sheaf{Y}_H \subset \sheaf{Y}$.  Let $\alg{C}$ be the algebraic torus with cocharacter lattice $\sheaf{Y}_H \cap \sheaf{Y}_{Q,n}$, and $\iota \From \alg{C} \To \alg{G}$ the resulting homomorphism (with central image).  Let $\alg{\tilde C}$ denote the pullback of the cover $\alg{\tilde G}$ via $\iota$.  Then $\iota$ lifts to a well-aligned homomorphism from $\alg{\tilde C}$ to $\alg{\tilde G}$.
\end{example}

The final example of a well-aligned homomorphism is relevant for the study of twisting genuine representations by one-dimensional representations of $\alg{G}$.
\begin{example}
Let $\alg{H}$ denote the maximal toral quotient of $\alg{G}$.  In other words, $\alg{H}$ is the torus whose character lattice equals $\Hom(\alg{G}, \alg{G}_m)$.  Let $p \From \alg{G} \To \alg{H}$ denote the canonical homomorphism, and write $\iota \From \alg{G} \To \alg{G} \times \alg{H}$ for the homomorphism $\Id \times p$.  Write $\alg{\tilde G} \times \alg{H}$ for the cover $(\alg{G}' \times \alg{H}, n)$.  

The homomorphism $\iota$ realizes $\alg{\tilde G}$ as the pullback via $\iota$ of the cover $\alg{\tilde G} \times \alg{H}$.  A Borel subgroup and torus in $\alg{G}$ determines a Borel subgroup and torus in $\alg{G} \times \alg{H}$.  In this way, $\iota$ lifts to a well-aligned homomorphism of covers from $\alg{\tilde G}$ to $\alg{\tilde G} \times \alg{H}$.  
\end{example}

\section{The dual group}

In this section, fix a degree $n$ cover $\alg{\tilde G}$ of a quasisplit group $\alg{G}$ over $S$.  Associated to $\alg{\tilde G}$, we define the ``dual group,'' a local system on $S_{\et}$ of affine group schemes over $\ZZ$.  We refer to Appendix \S \ref{GroupsTorsors}, for background on such local systems.  We begin by reviewing the Langlands dual group of $\alg{G}$ in a framework suggested by Deligne (personal communication).

\subsection{The Langlands dual group}

Choose, for now, a Borel subgroup $\alg{B} \subset \alg{G}$ containing a maximally split maximal torus $\alg{T}$.  The based root datum of $(\alg{G}, \alg{B}, \alg{T})$ was denoted $\Psi$, and the dual root datum,
$$\Psi^\vee = \left( \sheaf{Y}, \Phi^\vee, \Delta^\vee, \sheaf{X}, \Phi, \Delta \right),$$
is a local system of root data on $S_{\et}$.  

This defines a unique (up to unique isomorphism) local system $\dgp{G^\vee}$ on $S_{\et}$ of pinned reductive groups over $\ZZ$, called the \defined{Langlands dual group} of $\alg{G}$.  The center of $\dgp{G^\vee}$ is a local system on $S_{\et}$ of groups of multiplicative type over $\ZZ$, given by
$$\dgp{Z^\vee} = \Spec \left( \ZZ [ \sheaf{Y} / \sheaf{Y}^{\SC} ] \right).$$
See Example \ref{LocSpec} for more on local systems and $\Spec$ in this context.

\subsection{The dual group of a cover}

Now we adapt the definition of the dual group to covers.  The ideas here are the same as those of \cite{MWCrelle}.  The ideas for modifying root data originate in \cite[\S 2.2]{LusztigQuantumGroups} in the simply-connected case, in \cite[Theorem 2.9]{FinkelbergLysenko} for the almost simple case, in \cite[\S 11]{McNamara} and \cite{Reich} in the reductive case.  This dual group is also compatible with \cite{ABPTV} and the Hecke algebra isomorphisms of Savin \cite{SavinUnramified}, and the most recent work of Lysenko \cite{LysenkoSatake}.

Associated to the cover $\alg{\tilde G}$ of degree $n$, recall that $Q \From \sheaf{Y} \To \ZZ$ is the first Brylinski-Deligne invariant, and $\beta_Q \From \sheaf{Y} \otimes \sheaf{Y} \To n^{-1} \ZZ$ a symmetric bilinear form, and
$$\sheaf{Y}_{Q,n} = \{ y \in \sheaf{Y} : \beta_Q(y, y') \in \ZZ \text{ for all } y' \in \sheaf{Y} \} \subset \sheaf{Y}.$$
Define $\sheaf{X}_{Q,n} = \{ x \in n^{-1} \sheaf{X} : \langle x, y \rangle \in \ZZ \text{ for all } y \in \sheaf{Y}_{Q,n} \} \subset n^{-1} \sheaf{X}$.  For each root $\phi \in \Phi$, define constants $n_\phi$ and $m_\phi$,
\begin{equation}
\label{nm}
n_\phi = \frac{n}{\GCD(n, Q(\phi^\vee))} , \quad m_\phi = \frac{Q(\phi^\vee)}{\GCD(n, Q(\phi^\vee))}.
\end{equation}
Define \defined{modified roots} and \defined{modified coroots} by
$$\tilde \phi = n_\phi^{-1} \phi, \quad \tilde \phi^{\vee} = n_\phi \phi^\vee.$$
These define subsets $\tilde \Phi = \{ \tilde \phi : \phi \in \Phi \} \subset \sheaf{X}_{Q,n}$ and $\tilde \Phi^\vee = \{ \tilde \phi^\vee : \phi^\vee \in \Phi^\vee \} \subset \sheaf{Y}_{Q,n}$, as in \cite{MWCrelle}.  Modifying the simple roots and their coroots, we have subsets $\tilde \Delta \subset \tilde \Phi$ and $\tilde \Delta^\vee \subset \tilde \Phi$.

By \cite[Construction 1.3]{MWCrelle}, this defines a local system of based root data $\tilde \Psi$ on $S_{\et}$.  Write $\tilde \Psi^\vee$ for its dual,
$$\tilde \Psi^\vee = (\sheaf{Y}_{Q,n}, \tilde \Phi^\vee, \tilde \Delta^\vee, \sheaf{X}_{Q,n}, \tilde \Phi, \tilde \Delta).$$
Write $\sheaf{Y}_{Q,n}^{\SC}$ for the subgroup of $\sheaf{Y}_{Q,n}$ spanned by the modified coroots $\tilde \Phi^\vee$.  

Define $\dgp{\tilde G^\vee}$ to be the (unique up to unique isomorphism) local system on $S_{\et}$ of pinned reductive groups over $\ZZ$, associated to the local system of based root data $\tilde \Psi^\vee$.  Its maximal torus is a local system on $S_{\et}$ of split tori over $\ZZ$,
$$\dgp{\tilde T^\vee} = \Spec \left( \ZZ[ \sheaf{Y}_{Q,n}] \right).$$
The center of $\dgp{\tilde G^\vee}$ is a local system on $S_{\et}$ of groups of multiplicative type over $\ZZ$,
$$\dgp{\tilde Z^\vee} = \Spec \left( \ZZ \left[ \sheaf{Y}_{Q,n} /  \sheaf{Y}_{Q,n}^{\SC} \right] \right).$$
We call $\dgp{\tilde G^\vee}$ (endowed with its pinning) the \defined{dual group} of the cover $\alg{\tilde G}$.

\begin{proposition}
\label{DualGroupModn}
Suppose that $\alg{\tilde G}_0$ is another cover of $\alg{G}$ of degree $n$, with first Brylinski-Deligne invariant $Q_0$.  If $Q \equiv Q_0$ modulo $n$, i.e. $Q(y) - Q_0(y) \in n \ZZ$ for all $y \in \sheaf{Y}$, then the resulting modified root data are equal:  $\tilde \Psi^\vee = \tilde \Psi_0^\vee$.  Thus the dual groups are equal, $\dgp{\tilde G}^\Vee = \dgp{\tilde G}_0^\Vee$.
\end{proposition} 
\proof
One checks directly that $\beta_Q \equiv \beta_{Q_0}$ modulo $\ZZ$, from which it follows that
$$\sheaf{Y}_{Q,n} = \sheaf{Y}_{Q_0,n}, \quad \sheaf{X}_{Q,n} = \sheaf{X}_{Q_0, n}.$$
Similarly, one checks that the constants $n_\phi$ are equal,
$$\frac{n}{\GCD(n, Q_0(\phi^\vee))} = \frac{n}{\GCD(n, Q(\phi^\vee))}.$$
The result follows.
\qed

The Weyl group of $\dgp{\tilde G^\vee}$ with respect to $\dgp{\tilde T^\vee}$ forms a local system $\sheaf{\tilde W}$ on $S_{\et}$ of finite groups, generated (locally on $S_{\et}$) by reflections $s_{\tilde \phi}$ for every $\tilde \phi \in \tilde \Phi$.  The action of $\sheaf{\tilde W}$ on $\sheaf{Y}_{Q,n}$ is given by the formula
$$s_{\tilde \phi}(y) = y - \langle \tilde \phi, y \rangle \tilde \phi^\vee = y - \langle \phi, y \rangle \phi^\vee.$$
This identifies the root reflections $s_{\tilde \phi}$ with the root reflections $s_{\phi}$, and hence identifies the Weyl group $\sheaf{\tilde W}$ with the Weyl group $\sheaf{W}$ of $\alg{G}$ with respect to $\alg{T}$ (where both are viewed as local systems on $S_{\et}$ of finite groups).

The dual group $\dgp{\tilde G^\vee}$ comes with a distinguished 2-torsion element in its center, described here.  From the quadratic form $Q \From \sheaf{Y} \To \ZZ$, observe that $2 Q(y) = n \beta_Q(y,y) \in n \ZZ$ for all $y \in \sheaf{Y}_{Q,n}$.  Moreover, we have 
$$Q(\tilde \phi^\vee) = n_\phi^2 Q(\phi^\vee) = n_\phi m_\phi n \in n \ZZ,$$
for all $\phi \in \Phi$.  Of course, $Q( n y) \in n \ZZ$ as well, for all $y \in \sheaf{Y}_{Q,n}$.  We find a homomorphism of local systems of abelian groups on $S_{\et}$,
$$\overline{n^{-1} Q} \From \frac{\sheaf{Y}_{Q,n}}{\sheaf{Y}_{Q,n}^{\SC} + n \sheaf{Y}_{Q,n}} \To \half \ZZ / \ZZ, \quad y \mapsto n^{-1} Q(y) \text{ mod } \ZZ.$$

Applying $\Spec$ yields a homomorphism of local systems on $S_{\et}$ of diagonalizable group schemes over $\ZZ$,
$$\tau_Q \in \Hom ( \amu_2, \dgp{\tilde Z}_{[n]}^\Vee ).$$
Thus $\tau_Q(-1)$ is a Galois-invariant 2-torsion element in the center of $\dgp{\tilde G^\vee}$.  If $n$ is odd, then $\tau_Q(-1) = 1$.

\subsection{Well-aligned functoriality}
\label{WAFDualGroup}

Consider a well-aligned homomorphism $\tilde \iota \From \alg{\tilde G}_1 \To \alg{\tilde G}_2$ of covers, each endowed with Borel subgroup and maximally split maximal torus.  Here we construct a corresponding homomorphism of dual groups,
$$\iota^\vee \From \dgp{\tilde G}_2^\Vee \To \dgp{\tilde G}_1^\Vee.$$
These dual groups are constructed, locally on $S_{\et}$, from root data:
\begin{align*}
\tilde \Psi_1^\vee &= \left( \sheaf{Y}_{1,Q_1,n}, \tilde \Phi_1^\vee, \tilde \Delta_1^\vee, \sheaf{X}_{1,Q_1,n}, \tilde \Phi_1, \tilde \Delta_1 \right); \\
\tilde \Psi_2^\vee &= \left( \sheaf{Y}_{2,Q_2,n}, \tilde \Phi_2^\vee, \tilde \Delta_2^\vee, \sheaf{X}_{2,Q_2,n}, \tilde \Phi_2, \tilde \Delta_2 \right).
\end{align*}

For the construction of $\iota^\vee$, it suffices to work locally on $S_{\et}$, on a finite \'etale cover over which $\alg{G}_1$ and $\alg{G}_2$ split.  Condition (4) of Definition \ref{DefWA} gives a homomorphism $\iota \From \sheaf{Y}_{1,Q_1,n} \To \sheaf{Y}_{2,Q_2,n}$, and its dual homomorphism $\iota^\ast \From \sheaf{X}_{2,Q_2,n} \To \sheaf{X}_{1,Q_1,n}$.  As $\iota$ has normal image, the coroots from $\Phi_1^\vee$ map to coroots from $\Phi_2^\vee$.  Condition (3) implies that the scaled coroots in $\tilde \Phi_1^\vee \subset \sheaf{Y}_{1,Q_1,n}$ map to scaled coroots in $\tilde \Phi_2^\vee \subset \sheaf{Y}_{2,Q_2,n}$.  Since the Borel subgroups are aligned, the simple scaled coroots in $\tilde \Delta_1^\vee$ map to simple scaled coroots in $\tilde \Delta_2^\vee$.  Dually, the map $\iota^\ast \From \sheaf{X}_{2,Q_2,n} \To \sheaf{X}_{1,Q_1,n}$ sends $\tilde \Phi_2$ to $\tilde \Phi_1$ and $\tilde \Delta_2$ to $\tilde \Delta_1$.  

This allows us to assemble a homomorphism $\iota^\vee \From \dgp{\tilde G}_2^\Vee \To \dgp{\tilde G}_1^\Vee,$ (cf.~\cite[\S 2.1, 2.5]{BorelCorvallis}).  On tori, let $\iota^\vee \From \dgp{\tilde T}_2^\Vee \To \dgp{\tilde T}_1^\Vee$ be the homomorphism dual to the map of character lattices $\iota \From \sheaf{Y}_{1,Q_1,n} \To \sheaf{Y}_{2,Q_2,n}$.  Using the pinnings on $\dgp{\tilde G}_2^\Vee$ and $\dgp{\tilde G}_1^\Vee$, and the map of roots to roots, we obtain a homomorphism from the simply-connected cover $\dgp{\tilde G}_{2,\SC}^\Vee$ of the derived subgroup, $\dgp{\tilde G}_{2, \der}^\Vee$,
$$\iota_{\SC}^\vee \From \dgp{\tilde G}_{2,\SC}^\Vee \To \dgp{\tilde G}_1^\Vee.$$
Let $\dgp{\tilde T}_{2,\SC}^\Vee$ be the pullback of $\dgp{\tilde T}_2^\Vee$.  The following diagram commutes.
$$\begin{tikzcd}
\dgp{\tilde T}_{2, \SC}^\Vee \arrow{r} \inarrow{d} & \dgp{\tilde T}_2^\Vee \arrow{d}{\iota^\vee} \\
\dgp{\tilde G}_{2, \SC}^\Vee \arrow{r}{\iota_{\SC}^\Vee} & \dgp{\tilde G}_1^\Vee
\end{tikzcd}$$

The homomorphism $\iota_{\SC}^\Vee$ descends to the derived subgroup $\dgp{\tilde G}_{2,\der}^\Vee$, since it is trivial on the kernel of $\dgp{\tilde T}_{2,\SC}^\Vee \To \dgp{\tilde T}_2^\Vee$.  In this way, we have a pair of homomorphisms of groups over $\ZZ$,
$$\iota_{\der}^\Vee \From \dgp{\tilde G}_{2,\der}^\Vee \To \dgp{\tilde G}_1^\Vee, \quad \iota^\vee \From \dgp{\tilde T}_{2}^\Vee \To \dgp{\tilde T}_1^\Vee \subset \dgp{\tilde G}_1^\Vee.$$
These homomorphisms agree on their intersection, giving the desired homomorphism $\iota^\vee \From \dgp{\tilde G}_2^\Vee \To \dgp{\tilde G}_1^\Vee$.

Since modified coroots are sent to modified coroots, we find that $\iota^\vee$ sends $\dgp{\tilde Z}_2^\Vee$ to $\dgp{\tilde Z}_1^\Vee$.  The quadratic forms $Q_1$ and $Q_2$ induce two group homomorphisms:
$$\overline{n^{-1} Q_1} \From \frac{\sheaf{Y}_{1,Q_1,n}}{\sheaf{Y}_{1,Q_1,n}^{\SC} + n \sheaf{Y}_{1,Q_1,n}} \To \half \ZZ / \ZZ, \quad \overline{n^{-1} Q_2} \From \frac{\sheaf{Y}_{2,Q_2,n}}{\sheaf{Y}_{2,Q_2,n}^{\SC} + n \sheaf{Y}_{2,Q_2,n}} \To \half \ZZ / \ZZ.$$
These define two homomorphisms of group schemes over $\ZZ$,
$$\tau_{Q_1} \From \amu_2 \To \dgp{\tilde Z}_{1, [n]}^\Vee, \quad \tau_{Q_2} \From \amu_2 \To \dgp{\tilde Z}_{2, [n]}^\Vee.$$
As $Q_2(\iota(y)) = Q_1(y)$ for all $y \in \sheaf{Y}_1$, $\iota^\vee \circ \tau_{Q_1} = \tau_{Q_2}$.  In other words, the homomorphism $\iota^\vee \From \dgp{\tilde G}_2^\Vee \To \dgp{\tilde G}_1^\Vee$ sends the center to the center, and respects the 2-torsion elements therein, $\iota^\vee(\tau_{Q_2}(-1)) = \tau_{Q_1}(-1)$.

Given a pair of well-aligned homomorphisms,
$$\alg{\tilde G}_1 \xrightarrow{(\iota_1, \iota_1')} \alg{\tilde G}_2 \xrightarrow{(\iota_2, \iota_2')} \alg{\tilde G}_3,$$
their composition is a well-aligned homomrphism $(\iota_3, \iota_3') = (\iota_2 \iota_1, \iota_2' \iota_1')$ from $\alg{\tilde G}_1$ to $\alg{\tilde G}_3$ by Proposition \ref{ComposeWellaligned}.  

This gives a commutative diagram of sheaves of abelian groups.
$$\begin{tikzcd}
\sheaf{Y}_{1,Q_1,n} \arrow{r}[swap]{\iota_1} \arrow[bend left=20]{rr}{\iota_3} & \sheaf{Y}_{2,Q_2,n} \arrow{r}[swap]{\iota_2} & \sheaf{Y}_{3,Q_3,n}.
\end{tikzcd}$$
We find such a commutative diagram for dual groups, in the opposite direction.
$$\begin{tikzcd}
\dgp{\tilde G}_3^\Vee \arrow{r}[swap]{\iota_2^\vee} \arrow[bend left=20]{rr}{\iota_3^\vee} & \dgp{\tilde G}_2^\Vee \arrow{r}[swap]{\iota_1^\vee} & \dgp{\tilde G}_1^\Vee.
\end{tikzcd}$$

Let $\Cat{DGp}_S^\ast$ denote the category whose objects are local systems $\dgp{G}^\vee$ on $S_{\et}$ of group schemes over $\ZZ$, endowed with central morphisms $\amu_2 \To \dgp{G}^\vee$ (where $\amu_2$ is viewed as the constant local system of group schemes).  Morphisms in $\Cat{DGp}_S^\ast$ are morphisms of local systems of group schemes over $\ZZ$, compatible with the central morphisms from $\amu_2$.  

Let $\Cat{WAC}_S$ (Well-Aligned-Covers) denote the category whose objects are triples $(\alg{\tilde G}, \alg{B}, \alg{T})$ where $\alg{\tilde G}$ is a cover of a reductive group $\alg{G}$ over $S$, $\alg{B}$ is a Borel subgroup of $\alg{G}$, and $\alg{T}$ is a maximally split maximal torus of $\alg{G}$ contained in $\alg{B}$.  Morphisms in $\Cat{WAC}_S$ are well-aligned homomorphisms of covers.

We have proven the following result.
\begin{proposition}
\label{FuncDualGroup}
The construction of the dual group defines a contravariant functor from $\Cat{WAC}_S$ to $\Cat{DGp}_S^\ast$.
\end{proposition}

\begin{thm}
\label{WellDefinedDualGroup}
Let $\alg{\tilde G}$ be a degree $n$ cover of a quasisplit group $\alg{G}$ over $S$.  The dual group $\dgp{\tilde G^\vee}$ is well-defined, up to unique isomorphism, by $\alg{\tilde G}$ alone.
\end{thm}
\proof
The construction of the dual group depends on the choice of Borel and torus $\alg{B} \supset \alg{T}$.  So it suffices to construct a canonical isomorphism of dual groups $\dgp{\tilde G}^\Vee \To \dgp{\tilde G}_0^\vee$ for any pair of choices $\alg{B}_0 \supset \alg{T}_0$ and $\alg{B} \supset \alg{T}$.  Such ``well-definedness'' is discussed in more detail in \cite[\S 1.1]{DeligneLusztig}.

By Proposition \ref{BTConj}, there exists $g \in \alg{G}(S)$ such that $\Int(g) \alg{T}_0 = \alg{T}$ and $\Int(g) \alg{B}_0 = \alg{B}$.  This automorphism $\Int(g)$ lifts to an automorphism of $\alg{G}'$ and defines a well-aligned isomorphism of covers, 
$$\begin{tikzcd}
\alg{K}_{2} \inarrow{r} \arrow{d}{=} & \alg{G}' \onarrow{r} \arrow{d}{\Int(g)'} & \alg{G} \arrow{d}{\Int(g)} & & \alg{B}_{0} \arrow{d} & \alg{T}_{0} \arrow{d}  \\
\alg{K}_{2} \inarrow{r} & \alg{G}' \onarrow{r} & \alg{G} & & \alg{B} & \alg{T}
\end{tikzcd}$$
As a well-aligned isomorphism of covers, this yields an isomorphism (of local systems on $S_{\et}$ of reductive groups over $\ZZ$)
$$\Int(g)^\vee \From \dgp{\tilde G^\vee} \xrightarrow{\sim} \dgp{\tilde G}_0^\Vee.$$

If $g' \in \alg{G}(S)$ also satisfies $\Int(g') \alg{T}_0 = \alg{T}$ and $\Int(g') \alg{B}_0 = \alg{B}$, then $g' g^{-1} \in \alg{N}(S) \cap \alg{B}(S) = \alg{T}(S)$.  Thus $g' = t g$ for some $t \in \alg{T}(S)$.  Hence $\Int(g') = \Int(t) \Int(g)$, and so by Proposition \ref{FuncDualGroup},
$$\Int(g')^\vee = \Int(g)^\vee \Int(t)^\vee \From \dgp{\tilde G^\vee} \To \dgp{\tilde G}_0^\Vee.$$
But $\Int(t)^\vee = \Id$, since $\Int(t)$ leaves all relevant data unchanged.  Thus $\Int(g')^\vee = \Int(g)^\vee$.  Hence we find a canonical isomorphism $\dgp{\tilde G}^\Vee \xrightarrow{\sim} \dgp{\tilde G}_0^\Vee$.
\qed

\subsection{Change of base scheme}
\label{BaseChangeDualGroup}

Let $\alg{\tilde G}$ be a degree $n$ cover of a quasisplit group $\alg{G} \supset \alg{B} \supset \alg{T}$ over $S$, as before.  Let $\gamma \From S_0 \To S$ be a morphism of schemes, with $S_0 = \Spec(F_0)$ for some field $F_0$ or $S_0 = \Spec(\OO_0)$ for some DVR $\OO_0$ (with finite residue field or containing a field, as usual).  Then $\gamma$ gives rise to a pullback functor $\gamma^\ast$ from sheaves on $S_{\et}$ to sheaves on $S_{0,\et}$.  Pullback via $\gamma$ defines a degree $n$ cover $\alg{\tilde G}_0$ of a quasisplit group $\alg{G}_0 \supset \alg{B}_0 \supset \alg{T}_0$ over $S_0$. 

The cocharacter lattice $\sheaf{Y}_0$ of $\alg{T}_0$ is a sheaf on $S_{0,\et}$.  It is naturally isomorphic to the pullback $\gamma^\ast \sheaf{Y}$, with $\sheaf{Y}$ the cocharacter lattice of $\alg{T}$.  Write $N \From \gamma^\ast \sheaf{Y} \To \sheaf{Y}_0$ for the natural isomorphism.  The quadratic form $Q \From \sheaf{Y} \To \ZZ$ pulls back to a quadratic form $\gamma^\ast Q \From \gamma^\ast \sheaf{Y} \To \ZZ$.  The compatibility of Brylinski-Deligne invariants with pullbacks implies that $\gamma^\ast Q = N^\ast Q_0$, with $Q_0$ the first Brylinski-Deligne invariant of $\alg{\tilde G}_0$.  

\begin{remark}
This compatibility follows straightforwardly from \cite[\S 3.10]{B-D}; if $\alg{\tilde T}$ arises (after a finite \'etale $U \To S$) from the cocycle attached to $C \in \sheaf{X} \otimes \sheaf{X}$, then $\alg{\tilde T}_0$ arises from the pullback of this cocycle, i.e., from an element $C_0 \in \sheaf{X}_0 \otimes \sheaf{X}_0$ with $N^\ast C_0 = \gamma^\ast C$.  The quadratic form $Q$ is given by $Q(y) = C(y,y)$ and similarly $Q_0(y) = C_0(y,y)$.  Since $N^\ast C_0 = \gamma^\ast C$, we find that $N^\ast Q_0 = \gamma^\ast Q$.
\end{remark}

In this way, we find that $N$ restricts to an isomorphism from $\gamma^\ast \sheaf{Y}_{Q,n}$ to $\sheaf{Y}_{0,Q_0,n}$, sending roots and coroots (the sheaves of sets $\gamma^\ast \tilde \Phi$, $\gamma^\ast \tilde \Phi^\vee$ on $S_{0,\et}$) to the corresponding roots and coroots in $\sheaf{Y}_{0,Q_0,n}$.  As $\alg{B}_0 = \gamma^\ast \alg{B}$, simple roots and coroots are identified as well.  

By our construction of the dual group, we find that $N$ gives an isomorphism of local systems on $S_{0,\et}$ of pinned reductive groups over $\ZZ$,
$$N^\vee \From \gamma^\ast \dgp{\tilde G^\vee} \xrightarrow{\sim} \dgp{\tilde G}_0^\Vee.$$

\subsection{Parabolic subgroups}
\label{ParabolicDualGroup}

We keep the degree $n$ cover $\alg{\tilde G}$ of the quasisplit group $\alg{G} \supset \alg{B} \supset \alg{T}$ over $S$.  Let $\alg{P} \subset \alg{G}$ be a parabolic subgroup defined over $S$, containing $\alg{B}$.  Suppose that $\alg{P} = \alg{M} \alg{N}$ is a Levi decomposition defined over $S$, with $\alg{N}$ the unipotent radical of $\alg{P}$ and $\alg{M}$ a Levi factor containing $\alg{T}$.  Let $\alg{B}_{\alg{M}}$ denote the Borel subgroup $\alg{B} \cap \alg{M}$ of $\alg{M}$.  Write $\alg{\tilde M} = (\alg{M}', n)$ for the cover of $\alg{M}$ arising from pulling back $\alg{\tilde G}$.

The first Brylinski-Deligne invariant $Q$ is the same for $\alg{\tilde G}$ as for $\alg{\tilde M}$, as it depends only upon the cover $\alg{\tilde T}$ of their common maximal torus.  Write $\Phi_{\alg{M}}$ and $\Phi_{\alg{M}}^\vee$ for the roots and coroots of $\alg{M}$; these are subsets of $\Phi$ and $\Phi^\vee$, respectively, and by agreement of the first Brylinski-Deligne invariant,
$$\tilde \Phi_{\alg{M}} \subset \tilde \Phi, \quad \tilde \Phi_{\alg{M}}^\vee \subset \tilde \Phi^\vee.$$
As $\alg{B}_{\alg{M}} = \alg{B} \cap \alg{M}$, we have $\tilde \Delta_{\alg{M}} \subset \tilde \Delta$ and $\tilde \Delta_{\alg{M}}^\vee \subset \tilde \Delta^\vee$.  We find a pair of local systems on $S_{\et}$ of based root data,
$$\left( \sheaf{Y}_{Q,n}, \tilde \Phi^\vee, \tilde \Delta^\vee, \sheaf{X}_{Q,n}, \tilde \Phi, \tilde \Delta \right), \quad \left( \sheaf{Y}_{Q,n}, \tilde \Phi_{\alg{M}}^\vee, \tilde \Delta_{\alg{M}}^\vee, \sheaf{X}_{Q,n}, \tilde \Phi_{\alg{M}}, \tilde \Delta_{\alg{M}} \right).$$
The first root datum defines a local system $\dgp{\tilde G^\vee}$ on $S_{\et}$ of pinned reductive groups over $\ZZ$.  The second root datum defines a local system $\dgp{\tilde M^\vee}$ on $S_{\et}$ of pinned reductive Levi subgroups of $\dgp{\tilde G^\vee}$.

Thus the dual group of $\alg{\tilde M}$ is naturally a Levi subgroup of the dual group $\dgp{\tilde G^\vee}$.  Moreover, by agreement of the first Brylinski-Deligne invariants, the central 2-torsion element of $\dgp{\tilde M^\vee}$ coincides with the central 2-torsion element $\tau_Q(-1) \in \dgp{\tilde G^\vee}$.

\subsection{Weyl action on the dual torus}
\label{WeylDualTorus}

As before, keep the degree $n$ cover $\alg{\tilde G}$ of $\alg{G}$.  Write $\alg{\tilde T} = (\alg{T}', n)$ for the resulting cover of a maximally split maximal torus $\alg{T}$.  Assume here that $\alg{T}$ splits over a cyclic Galois cover of $S$.  Suppose that $w \in \sheaf{W}[S]$ is an element of the Weyl group, represented by an element $\dot w \in \sheaf{N}[S]$ (here $\sheaf{N}$ is the sheaf on $S_{\et}$ represented by the normalizer of $\alg{T}$).

Then $\Int(\dot w) \From \alg{\tilde G} \To \alg{\tilde G}$ defines a well-aligned homomorphism from $\alg{\tilde T}$ to itself.
$$\begin{tikzcd}
\alg{K}_2 \inarrow{r} \arrow{d}{=} & \alg{T}' \onarrow{r} \arrow{d}{\Int(\dot w)} & \alg{T} \arrow{d}{\Int(w)} \\
\alg{K}_2 \inarrow{r} & \alg{T}' \onarrow{r} & \alg{T}
\end{tikzcd}$$ 

As such, $\Int(\dot w)$ defines a map of dual groups,
$$\Int(\dot w)^\vee \From \dgp{\tilde T}^\Vee \To \dgp{\tilde T}^\Vee.$$
As $\Int(t)^\vee = \Id$ for all $t \in \alg{T}(S)$, this homomorphism of dual groups depends only on the element $w$ of the Weyl group, not on the chosen representative $\dot w$.  Thus we write
$$\Int(w)^\vee \From \dgp{\tilde T}^\Vee \To \dgp{\tilde T}^\Vee.$$

On the other hand the element $w \in \sheaf{W}[S]$ corresponds to a Galois-invariant element of the Weyl group of the dual group.  From \cite[Lemma 6.2]{BorelCorvallis}, there exists an element $n^\vee \in \dgp{\tilde G}^\Vee[S]$ such that $n^\vee$ normalizes $\dgp{\tilde T}$ and the resulting action of $n^\vee$ on the character lattice $\sheaf{Y}_{Q,n}$ of $\dgp{\tilde T}$ coincides with the action of $w$ on $\sheaf{Y}_{Q,n}$.  Such an element $n^\vee$ is unique up to multiplication by elements of $\dgp{\tilde T}^\Vee[S]$.  Since $\Int(w)^\vee$ and $\Int(n^\vee)$ define the same action on $\sheaf{Y}_{Q,n}$, we find that 
$$\Int(w)^\vee = \Int(n^\vee), \text{ as elements of } \Aut(\dgp{\tilde T}^\Vee).$$

\subsection{Specific cases}

Here we give many specific cases of covers, associated dual groups, and 2-torsion elements in their centers corresponding to $\tau_Q(-1)$.  

\subsubsection{Method for simply-connected groups}

Let $\alg{G}$ be a simply-connected semisimple group over $S$, with Borel subgroup $\alg{B}$ containing a maximally split maximal torus $\alg{T}$ over $S$.  Let $S' / S$ be a Galois cover over which $\alg{T}$ splits.  If $t$ is an integer, there exists a unique Weyl-invariant quadratic form $Q_t$ with value $t$ on all short coroots.  This $Q_t$ is a multiple of the Killing form.

Corresponding to $Q_t$, there is a unique, up to unique isomorphism, object $\alg{G}^{(t)} \in \Cat{CExt}_S(\alg{G}, \alg{K}_2)$ with first Brylinski-Deligne invariant $Q_t$ (by \cite[\S 7.3(i)]{B-D} when working over a field or \cite[\S 3.3.3]{MWIntegral} over a DVR).  Write $\beta_t$ for the resulting $n^{-1} \ZZ$-valued bilinear form.  Weyl-invariance of $Q_t$ implies that
\begin{equation}
\beta_t(\phi^\vee, y) = n^{-1} Q_t(\phi^\vee) \langle \phi, y \rangle, \text{ for all } \phi \in \Phi, \text{ and all } y \in \sheaf{Y}.
\end{equation}
(See the proof of Lemma \ref{DeltaM} for a derivation.)  It follows that, for any $y \in \sheaf{Y}$,
\begin{equation}
\label{BetaObs}
\beta_t(\phi^\vee, y) \in \ZZ \text{ if and only if } \langle \phi,y \rangle \in n_\phi \ZZ.
\end{equation}

Let $Y = \sheaf{Y}[S']$, $X = \sheaf{X}[S']$, and $\Delta = \Delta[S'] = \{ \alpha_1, \ldots, \alpha_\ell \}$ a basis of simple roots corresponding to $\alg{B}$.  Since we assume $\alg{G}$ is simply connected, $Y = \bigoplus_{i = 1}^\ell \alpha_i^\vee \ZZ$.  Write $n_i = n_{\alpha_i}$.  From \eqref{BetaObs}, we find a characterization of $Y_{Q,n} = \sheaf{Y}_{Q,n}[S']$:
\begin{equation}
\label{YQnCharacterization}
y \in Y_{Q,n} \text{ if and only if } \langle \alpha_i, y \rangle \in n_i \ZZ \text{ for all } 1 \leq i \leq \ell.
\end{equation}
This, in turn, can be used to tabulate dual groups.  We provide tables here reference, noting that such information can also be found in the examples of \cite[\S 2.4]{FinkelbergLysenko}.    But we also include data on the central 2-torsion elements that we have not found in the literature.  Our tabulation was greatly assisted by using SAGE \cite{sage}, especially the recently updated package which deftly handles root data.

The following tables only include split groups, and we write $\alg{\tilde G}^\Vee = \dgp{\tilde G}^\Vee[S]$ for the dual group over $S$ (a pinned reductive group scheme over $\ZZ$).  For quasisplit groups, one can view the dual group as a pinned reductive group scheme over $\ZZ$ endowed with Galois action by pinned automorphisms.

\subsubsection{$\alg{SL}_{\ell + 1}$}
For $\alg{G} = \alg{SL}_{\ell + 1}$, the standard Borel subgroup and maximal torus, and system of roots $\alpha_1, \ldots, \alpha_\ell$, the Dynkin diagram is
\begin{center}
\begin{tikzpicture}
\draw (1,0) -- (3.7,0);
\draw (5.3,0) -- (6,0);
\draw[dotted] (3.7,0) -- (5.3,0);
\foreach \i/\l in {1/1,2/2,3/3,6/\ell}
{ 
\filldraw[fill=gray, draw=black] (\i,0) circle (0.1) node[below=2pt] {$\alpha_{\l}$};
}
\end{tikzpicture}
\end{center}
Consider the cover $\alg{\tilde G}$ of degree $n$ arising from the quadratic form satisfying $Q(\phi^\vee) = 1$ for all $\phi \in \Phi$.  Then $n_i  = n$ and $\tilde \alpha_i^\vee = n \alpha_i^\vee$ for all $1 \leq i \leq \ell$.  Hence
$$Y_{Q,n}^{\SC} = n Y = n \alpha_1^\vee \ZZ + \cdots + n \alpha_\ell^\vee \ZZ.$$

The Cartan matrix of $\sch{\tilde G}^\Vee$ has entries $\tilde C_{ij} = \langle \tilde \alpha_i, \tilde \alpha_j^\vee \rangle = \langle \alpha_i, \alpha_j^\vee \rangle$, and so $\sch{\tilde G}^\Vee$ is isogenous to $\sch{SL}_{\ell + 1}$.  To determine the dual group up to isomorphism, it suffices to compute the order of the center, since the center of $\sch{SL}_{\ell + 1}$ is $\amu_{\ell+1}$.  The order of the center is equal to the index $\# Y_{Q,n} / Y_{Q,n}^{\SC}$, and this is computable from \eqref{YQnCharacterization}.  The results are given in Table \ref{DualATable}.

\begin{table}[!htbp]
\begin{tabular}{l|ccccc}
 & \multicolumn{5}{c}{Group $\alg{G}$} \\
$n$ \phantom{M} & $\alg{SL}_2$ & $\alg{SL}_3$ & $\alg{SL}_4$ & $\alg{SL}_5$ & $\alg{SL}_6$ \\ \hline \hline
1 & $\sch{PGL}_2$ & $\sch{PGL}_3$ & $\sch{PGL}_4$ & $\sch{PGL}_5$ & $\sch{PGL}_6$ \\
2 & ${}^{\ast} \sch{SL}_2$  & $\sch{PGL}_3$ & $\sch{SL}_4 / \amu_2$ & $\sch{PGL}_5$ & ${}^{\ast}\sch{SL}_6 / \amu_3$ \\  
3 & $\sch{PGL}_2$ & $\sch{SL}_3$ & $\sch{PGL}_4$ & $\sch{PGL}_5$ & $\sch{SL}_6 / \amu_2$ \\
4 & $\sch{SL}_2$ & $\sch{PGL}_3$ & ${}^\ast \sch{SL}_4$ & $\sch{PGL}_5$ & $\sch{SL}_6 / \amu_3$ \\
5 & $\sch{PGL}_2$ & $\sch{PGL}_3$ & $\sch{PGL}_4$ & $\sch{SL}_5$ & $\sch{PGL}_6$ \\
6 & ${}^\ast \sch{SL}_2$ & $\sch{SL}_3$ & $\sch{SL}_4 / \amu_2$ & $\sch{PGL}_5$ & ${}^\ast \sch{SL}_6$ \\

\end{tabular}
\caption{Table of dual groups for degree $n$ covers of $\alg{SL}_{\ell+1}$.  Groups marked with $\ast$ have $\tau_Q(-1)$ nontrivial.}
\label{DualATable}
\end{table}

The dual groups $\sch{\tilde G}^\Vee$ are consistent with the Iwahori-Hecke algebra isomorphisms found by Savin in \cite[Theorem 7.8]{SavinUnramified}.  In other words, the dual group $\sch{\tilde G}^\Vee$ coincides with the Langlands dual group of $\alg{SL}_{\ell + 1} / \alg{Z}_{[n]}$, where $\alg{Z}_{[n]}$ is the $n$-torsion subgroup of the center of $\alg{SL}_{\ell+1}$.  

The central 2-torsion elements $\tau_Q(-1)$ follow a somewhat predictable pattern.  For covering degree $2$, $\tau_Q(-1)$ is nontrivial for $\alg{SL}_2$, $\alg{SL}_6$, $\alg{SL}_{10}$, $\alg{SL}_{14}$, etc..  In covering degree $4$, $\tau_Q(-1)$ is nontrivial for $\alg{SL}_4$, $\alg{SL}_{12}$, $\alg{SL}_{20}$, etc..  In covering degree $6$, $\tau_Q(-1)$ is nontrivial for $\alg{SL}_2$, $\alg{SL}_6$, $\alg{SL}_{10}$, $\alg{SL}_{14}$, etc..  In covering degree $8$, $\tau_Q(-1)$ is nontrivial for $\alg{SL}_8$, $\alg{SL}_{24}$, etc..  In general, we suspect the following:
\begin{center}
$\tau_Q(-1)$ is nontrivial for a degree $2^e \cdot k$ ($k$ odd) cover of $\alg{SL}_m$ if and only if $m = 2^e \cdot j$ for $j$ odd.
\end{center} 

To illustrate the quasisplit case, consider $\alg{G} = \alg{SU}_3$, a quasisplit special unitary group associated to a degree $2$ Galois cover $S' / S$.  There is a unique degree $n$ cover of $\alg{G}$ arising from the quadratic form taking values $1$ at all coroots.  The dual group of the degree $2$ cover of $\alg{G} = \alg{SU}_3$ is identified with $\sch{PGL}_3$, a pinned reductive group over $\ZZ$, endowed with $\Gal(S' / S)$-action by outer automorphism corresponding to the nontrivial automorphism of the Dynkin diagram \tikz{\draw (1,0) -- (2,0); \filldraw[fill=gray, draw=black] (1,0) circle (0.1);  \filldraw[fill=gray, draw=black] (2,0) circle (0.1);}.  Thus the dual group of the double-cover of $\alg{SU}_3$ coincides with the Langlands dual group of the linear group $\alg{SU}_3$.

\subsubsection{$\alg{Spin}_{2 \ell + 1}$}
For $\alg{G} = \alg{Spin}_{2 \ell + 1}$, the Dynkin diagram has type B.
\begin{center}
\begin{tikzpicture}
\draw (1,0) -- (3.7,0);
\draw (5.3,0) -- (6,0);
\draw[dotted] (3.7,0) -- (5.3,0);
\draw (6,0.05) -- (7,0.05);
\draw (6,-0.05) -- (7,-0.05);
\draw (6.45,0.125) -- (6.55,0) -- (6.45,-0.125);
\foreach \i/\l in {1/1,2/2,3/3,6/{\ell-1}, 7/\ell}
{ 
\filldraw[fill=gray, draw=black] (\i,0) circle (0.1) node[below=2pt] {$\alpha_{\l}$};
}
\end{tikzpicture}
\end{center}
Let $\alg{\tilde G}$ be the cover of degree $n$, associated to the quadratic form taking the value $1$ at all short coroots.  Thus $Q(\alpha_i^\vee) = 2$ for all long coroots $1 \leq i \leq \ell - 1$.  If $n$ is odd, then $n_\alpha = n$ for all coroots $\alpha^\vee$.  If $n$ is even, then $n_i = n/2$ for $1 \leq i \leq \ell-1$ and $n_\ell = n$.  When $n$ is even, short coroots become long and long become short, after modification.  We find that the dual group is isogenous to $\sch{Sp}_{2 \ell}$ if $n$ is odd, and is isogenous to $\sch{Spin}_{2 \ell + 1}$ if $n$ is even.

The centers of $\sch{Sp}_{2 \ell}$ and $\sch{Spin}_{2 \ell + 1}$ are cyclic of order two.  Thus the dual group can be identified by the order of its center.

\begin{table}[!htbp]
\begin{tabular}{l|cccccc}
 & \multicolumn{6}{c}{Group $\alg{G}$} \\
$n$ \phantom{M} & $\sch{Spin}_7$ & $\sch{Spin}_9$ & $\sch{Spin}_{11}$ & $\sch{Spin}_{13}$ & $\sch{Spin}_{15}$ & $\sch{Spin}_{17}$  \\ \hline \hline
1 & $\sch{PGSp}_{6}$ & $\sch{PGSp}_{8}$ & $\sch{PGSp}_{10}$ & $\sch{PGSp}_{12}$ & $\sch{PGSp}_{14}$ & $\sch{PGSp}_{16}$ \\
2 &  $\sch{SO}_{7}$ & $\sch{Spin}_{9}$ & $\sch{SO}_{11}$ & ${}^\ast \sch{Spin}_{13}$ &$ \sch{SO}_{15}$ & $\sch{Spin}_{17}$ \\
3 & $\sch{PGSp}_{6}$ & $\sch{PGSp}_{8}$ & $\sch{PGSp}_{10}$ & $\sch{PGSp}_{12}$ & $\sch{PGSp}_{14}$ & $\sch{PGSp}_{16}$ \\
4 & ${}^\ast  \sch{Spin}_{7}$ & $\sch{Spin}_{9}$ & ${}^\ast \sch{Spin}_{11}$ & $\sch{Spin}_{13}$ & ${}^\ast \sch{Spin}_{15}$ & $\sch{Spin}_{17}$ \\
5 & $\sch{PGSp}_{6} $ & $\sch{PGSp}_{8} $ & $\sch{PGSp}_{10} $ & $\sch{PGSp}_{12}$ &$ \sch{PGSp}_{14}$ & $\sch{PGSp}_{16}$ \\
6 &  $\sch{SO}_{7}$ & $\sch{Spin}_{9}$ & $\sch{SO}_{11}$ & ${}^\ast \sch{Spin}_{13}$ & $\sch{SO}_{15}$ & $\sch{Spin}_{17}$  \\
\end{tabular}
\caption{Table of dual groups for degree $n$ covers of $\alg{Spin}_{2 \ell+1}$.  Groups marked with $\ast$ have $\tau_Q(-1)$ nontrivial.}
\label{DualBTable}
\end{table}
Table \ref{DualBTable} describes the dual groups.  Note that, in this case, the isogeny class of the dual group depends on the covering degree modulo $4$.  In covering degree $4k + 2$, we find that $\tau_Q(-1)$ is nontrivial for $\alg{Spin}_{8 j + 5}$ (corresponding to rank $4j + 2$) for all positive integers $j$.  In covering degree $4k$, we find that $\tau_Q(-1)$ is nontrivial for $\alg{Spin}_{4j + 3}$ for all positive integers $j$.

\subsubsection{$\alg{Sp}_{2 \ell}$}
For $\alg{G} = \alg{Sp}_{2 \ell}$, the Dynkin diagram has type C.
\begin{center}
\begin{tikzpicture}
\draw (1,0) -- (3.7,0);
\draw (5.3,0) -- (6,0);
\draw[dotted] (3.7,0) -- (5.3,0);
\draw (6,0.05) -- (7,0.05);
\draw (6,-0.05) -- (7,-0.05);
\draw (6.55,0.125) -- (6.45,0) -- (6.55,-0.125);
\foreach \i/\l in {1/1,2/2,3/3,6/{\ell-1}, 7/\ell}
{ 
\filldraw[fill=gray, draw=black] (\i,0) circle (0.1) node[below=2pt] {$\alpha_{\l}$};
}
\end{tikzpicture}
\end{center}
Let $\alg{\tilde G}$ be the cover of degree $n$, associated to the quadratic form taking the value $1$ at all short coroots.  As in type B, we find that short coroots become long, and long become short, after modification when $n$ is even.  We find that the dual group is isogenous to $\sch{Sp}_{2 \ell}$ if $n$ is even, and is isogenous to $\sch{Spin}_{2 \ell + 1}$ if $n$ is odd.  As before, the dual group can be identified by the order of its center.

\begin{table}[!htbp]
\begin{tabular}{l|ccc}
 & \multicolumn{3}{c}{Group $\alg{G}$} \\
$n$ \phantom{M} & $\alg{Sp}_6$ & $\alg{Sp}_8$ & $\alg{Sp}_{10}$ \\ \hline \hline
1 & $\sch{SO}_{7}$ & $\sch{SO}_{9}$ & $\sch{SO}_{11}$ \\
2 & ${}^\ast \sch{Sp}_{6}$ & ${}^\ast \sch{Sp}_{8}$ & ${}^\ast \sch{Sp}_{10}$ \\
3 &  $\sch{SO}_{7}$ & $\sch{SO}_{9}$ & $\sch{SO}_{11}$ \\
4 & $\sch{Sp}_{6}$ & $\sch{Sp}_{8}$ & $\sch{Sp}_{10}$ \\
5 &  $\sch{SO}_{7}$ & $\sch{SO}_{9}$ & $\sch{SO}_{11}$ \\
6 & ${}^\ast \sch{Sp}_{6}$ & ${}^\ast \sch{Sp}_{8}$ & ${}^\ast \sch{Sp}_{10}$ \\
\end{tabular}
\caption{Table of dual groups for degree $n$ covers of $\alg{Sp}_{2 \ell}$.  Groups marked with $\ast$ have $\tau_Q(-1)$ nontrivial.}
\label{DualCTable}
\end{table}
As Table \ref{DualCTable} illustrates, the dual group of the degree $n$ cover of $\alg{Sp}_{2 \ell}$ is the simply-connected Chevalley group $\sch{Sp}_{2 \ell}$ when $n$ is even, and the dual group is the adjoint group $\sch{SO}_{2 \ell + 1} = \sch{Spin}_{2 \ell + 1} / \amu_2$ when $n$ is odd.  The central 2-torsion element $\tau_Q(-1)$ is nontrivial when the covering degree is $4k + 2$ for some non-negative integer $k$.  This is consistent (in covering degree $2$) with expectations from the classical theta correspondence for metaplectic groups.

\subsubsection{$\alg{Spin}_{2 \ell}$}
For $\alg{G} = \alg{Spin}_{2 \ell}$, $\ell \geq 4$, the Dynkin diagram has type D.
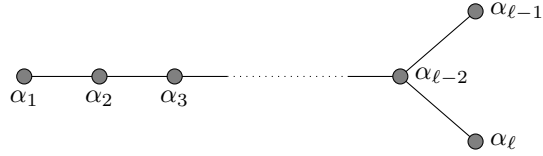
\begin{center}
\begin{tikzpicture}
\draw (1,0) -- (3.7,0);
\draw (5.3,0) -- (6,0);
\draw[dotted] (3.7,0) -- (5.3,0);
\draw (6,0) -- (7,0.866);
\draw (6,0) -- (7,-0.866);
\filldraw[fill=gray, draw=black] (6,0) circle (0.1) node[right=2pt] {$\alpha_{\ell-2}$};
\filldraw[fill=gray, draw=black] (7,0.866) circle (0.1) node[right=2pt] {$\alpha_{\ell-1}$};
\filldraw[fill=gray, draw=black] (7,-0.866) circle (0.1) node[right=2pt] {$\alpha_{\ell}$};
\foreach \i/\j/\l in {1/1,2/2,3/3}
{ 
\filldraw[fill=gray, draw=black] (\i,0) circle (0.1) node[below=2pt] {$\alpha_{\l}$};
}
\end{tikzpicture}
\end{center}
Let $\alg{\tilde G}$ be the cover of degree $n$, associated to the quadratic form taking the value $1$ at all coroots.  By the same methods as in type A, we find that the dual group is isogenous to $\sch{Spin}_{2 \ell}$.  

If $\ell$ is odd, then the center of $\sch{Spin}_{2 \ell}$ is a cyclic group of order $4$.  In this case, the dual group is determined by the order of its center.

If $\ell$ is even, then the center of $\sch{Spin}_{2 \ell}$ is isomorphic to $\amu_2 \times \amu_2$, and so the dual group is not a priori determined by the order of its center.  But fortunately, the order of the center of the dual group always equals $1$ or $4$ when $\ell$ is even, and this suffices to identify the dual group.

\begin{table}[!htbp]
\begin{tabular}{l|cccccc}
 & \multicolumn{6}{c}{Group $\alg{G}$} \\
$n$ \phantom{M} & $\alg{Spin}_8$ & $\alg{Spin}_{10}$ & $\alg{Spin}_{12}$ & $\alg{Spin}_{14}$ & $\alg{Spin}_{16}$ & $\alg{Spin}_{18}$ \\ \hline \hline
1 & $\sch{PGO}_8$ & $\sch{PGO}_{10}$ & $\sch{PGO}_{12}$ & $\sch{PGO}_{14}$ & $\sch{PGO}_{16}$ & $\sch{PGO}_{18}$  \\
2 & $\sch{Spin}_8$ & $\sch{SO}_{10}$ & ${}^\ast \sch{Spin}_{12}$ & $\sch{SO}_{14}$ & $\sch{Spin}_{16}$ & $\sch{SO}_{18}$ \\
3 & $\sch{PGO}_8$ & $\sch{PGO}_{10}$ & $\sch{PGO}_{12}$ & $\sch{PGO}_{14}$ & $\sch{PGO}_{16}$ & $\sch{PGO}_{18}$  \\
4 & $\sch{Spin}_8$ & ${}^\ast \sch{Spin}_{10}$ & $\sch{Spin}_{12}$ & ${}^\ast \sch{Spin}_{14}$ & $\sch{Spin}_{16}$ & ${}^\ast \sch{Spin}_{18}$\\
5 & $\sch{PGO}_8$ & $\sch{PGO}_{10}$ & $\sch{PGO}_{12}$ & $\sch{PGO}_{14}$ & $\sch{PGO}_{16}$ & $\sch{PGO}_{18}$ \\
6 & $\sch{Spin}_8$ & $\sch{SO}_{10}$ & ${}^\ast \sch{Spin}_{12}$ & $\sch{SO}_{14}$ & $\sch{Spin}_{16}$ & $\sch{SO}_{18}$ \\
\end{tabular}
\caption{Table of dual groups for degree $n$ covers of $\alg{Spin}_{2 \ell}$.  Groups marked with $\ast$ have $\tau_Q(-1)$ nontrivial.}
\label{DualDTable}
\end{table}
As Table \ref{DualDTable} illustrates, the dual group of an odd-degree cover of $\alg{Spin}_{2 \ell}$ coincides with the Langlands dual group of the linear group $\alg{Spin}_{2 \ell}$; this dual group is the adjoint group $\sch{PGO}_{2 \ell}$.  But the dual group of an even-degree cover of $\alg{Spin}_{2 \ell}$ depends on the parity of $\ell$ and the covering degree modulo $4$.  As in type A, these dual groups agree with expectations from the Hecke algebra isomorphisms of Savin \cite{SavinUnramified}.

When the covering degree is a multiple of $4$, the element $\tau_Q(-1)$ is nontrivial for $\alg{Spin}_{4j + 2}$ for all $j \geq 2$.  Since $\alg{Spin}_{4j + 2}$ has a unique central element of order two, this suffices to describe $\tau_Q$.  When the covering degree has the form $4k +2$, the element $\tau_Q(-1)$ is nontrivial for $\alg{Spin}_{8j + 4}$ for all $j \geq 1$.  The center of the group $\sch{Spin}_{8j + 4}$ is isomorphic to $\amu_2 \times \amu_2$, which has three distinct 2-torsion elements.  However, only one of these is invariant under the nontrivial outer automorphism of the pinned Chevalley group $\sch{Spin}_{8j + 4}$.  This one must be $\tau_Q(-1)$, since $Q$ is invariant under this outer automorphism.

\subsubsection{Exceptional groups}
Let $\alg{G}$ be a simply-connected split simple group of type $\text{E}_\ell$ (with $\ell = 6, 7, 8$), $\text{F}_4$, or $\text{G}_2$.  Let $\alg{\tilde G}$ be the cover of degree $n$, associated to the quadratic form taking the value $1$ at all short coroots (all coroots in type E).  As in types A and D, we find that the dual group is semisimple and isogenous to the Chevalley group of the same type as $\alg{G}$.  In types $\text{E}_8$, $\text{F}_4$, and $\text{G}_2$, the simply-connected group is centerless, and so $\sch{\tilde G}^\Vee$ coincides with the simply-connected Chevalley group of the same type.

The center of $\sch{E}_6$ has order $3$, and the center of $\sch{E}_7$ has order $2$.  Hence the dual group $\sch{\tilde G}^\Vee$ is determined by the order of its center.  The dual groups are listed in Table \ref{DualETable}.
\begin{table}[!htbp]
\begin{tabular}{l|ccccc}
 & \multicolumn{5}{c}{Group $\alg{G}$} \\
$n$ \phantom{M} & $\alg{E}_6$ & $\alg{E}_7$ & $\alg{E}_8$ & $\alg{F}_4$ & $\alg{G}_2$ \\ \hline
1 & $\sch{E}_6 / \amu_3$ & $\sch{E}_7 / \amu_2$ & $\sch{E}_8$ & $\sch{F}_4$ & $\sch{G}_2$ \\
2 & $\sch{E}_6 / \amu_3$ & ${}^\ast \sch{E}_7$ & $\sch{E}_8$ & $\sch{F}_4$ & $\sch{G}_2$ \\
3 & $\sch{E}_6$ & $\sch{E}_7 / \amu_2$ & $\sch{E}_8$ & $\sch{F}_4$ & $\sch{G}_2$ \\
4 & $\sch{E}_6 / \amu_3$ & $\sch{E}_7$ & $\sch{E}_8$ & $\sch{F}_4$ & $\sch{G}_2$  \\
5 & $\sch{E}_6 / \amu_3$ & $\sch{E}_7 / \amu_2$ & $\sch{E}_8$ & $\sch{F}_4$ & $\sch{G}_2$ \\
6 & $\sch{E}_6$ & ${}^\ast \sch{E}_7$ & $\sch{E}_8$ & $\sch{F}_4$ & $\sch{G}_2$ \\
\hline
\end{tabular}
\caption{Table of dual groups for degree $n$ covers of exceptional groups.  Groups marked with $\ast$ have $\tau_Q(-1)$ nontrivial.}
\label{DualETable}
\end{table}
In type E, these dual groups agree with expectations from \cite{SavinUnramified}.  The central 2-torsion element is nontrivial for $\alg{E}_7$, when the covering degree equals $4j + 2$ for some $j \geq 0$.

\subsubsection{$\alg{GL}_r$}

Suppose that $\alg{G}$ is split reductive, and the derived subgroup of $\alg{G}$ is simply-connected.  Let $\alg{T}$ be a split maximal torus in $\alg{G}$ with cocharacter lattice $Y$.  Then, for any Weyl-invariant quadratic form $Q \From Y \To \ZZ$, there exists a cover $\alg{\tilde G}$ with first Brylinski-Deligne invariant $Q$.

For example, when $\alg{G} = \alg{GL}_r$, there is a two-parameter family of such Weyl-invariant quadratic forms.  Write $\alg{T}$ for the standard maximal torus of diagonal matrices, and identify $Y = \ZZ^r$ in the usual way.  For any pair of integers $q,c$, there exists a unique Weyl-invariant quadratic form $Q_{q,c}$ satisfying
$$Q(1,-1,0,\ldots, 0) = q \text{ and } Q(1,0,\ldots, 0) = 1+c.$$
The $n$-fold covers $\widetilde{GL}_r^{(c)}$ studied by Kazhdan and Patterson \cite[\S 0.1]{KazhdanPatterson} can be constructed from Brylinski-Deligne extensions with first invariant $Q_{1,c}$.  The proof of following result is left to the reader.
\begin{proposition}
Let $\alg{\tilde G}$ be a degree $n$ cover of $\alg{G} = \alg{GL}_r$ with first Brylinski-Deligne invariant $Q_{1,c}$.  If $\GCD(n, 1+r + 2rc) = 1$, then $\sch{\tilde G}^\Vee$ is isomorphic to $\sch{GL}_r$.  If $\GCD(n,r) = 1$, then the derived subgroup of $\sch{\tilde G}^\Vee$ is isomorphic to $\sch{SL}_r$ and thus there exists an isogeny $\sch{\tilde G}^\Vee \To \sch{GL}_r$.
\end{proposition}
This may place the work of Kazhdan and Flicker \cite{KazhdanFlicker} in a functorial context.

\subsubsection{$\alg{GSp}_{2 r}$}

For $\alg{G} = \alg{GSp}_{2r}$, and a standard choice of split maximal torus and Borel subgroup, we write $e_0, \ldots, e_r$ for a basis of $Y$, $f_0, \ldots, f_r$ for the dual basis of $X$, and the simple roots and coroots are
$$\alpha_1 = f_1 - f_2, \ldots, \alpha_{r-1} = f_{r-1} - f_r, \quad \alpha_r = 2 f_r - f_0;$$
$$\alpha_1^\vee = e_1 - e_2, \ldots, \alpha_{r-1}^\vee = e_{r-1} - e_r, \quad \alpha_r^\vee = e_r.$$
The Weyl group is $S_r \ltimes \mu_2^r$, with $S_r$ acting by permutation of indices $1, \ldots, r$ (fixing $e_0$ and $f_0$), and elements $w_j$ (for $1 \leq j \leq r$) of order two which satisfy 
$$w_j(e_j) = - e_j, \quad w_j(e_i) = e_i \text{ for } i \neq j,0, \quad w_j(e_0) = e_0 + e_1.$$
Weyl-invariant quadratic forms on $Y$ are in bijection with pairs $(\kappa, \nu)$ of integers.  For any such pair, there is a unique Weyl-invariant quadratic form $Q_{\kappa, \nu}$ satisfying
$$Q_{\kappa, \nu}(e_0) = \kappa, \quad Q_{\kappa, \nu}(e_i) = \nu \text{ for } 1 \leq i \leq r.$$
The proof of following result is left to the reader.
\begin{proposition}
Let $\alg{\tilde G}$ be a degree $2$ cover of $\alg{G} = \alg{GSp}_{2r}$, with first Brylinski-Deligne invariant $Q_{0,1}$.  Then the dual group $\sch{\tilde G}^\Vee$ is isomorphic to $\sch{GSp}_{2r}$ if $r$ is odd, and to $\sch{PGSp}_{2r} \times \sch{G}_m$ if $r$ is even.
\end{proposition}
We find that that double-covers of $\alg{GSp}_{2r}$ behave differently depending on the parity of $r$; this phenomenon is consistent with the work of Szpruch \cite{Szpruch} on principal series.

\section{The gerbe associated to a cover}

In this section, we construct a gerbe $\gerb{E}_\epsilon(\alg{\tilde G})$ on $S_{\et}$ associated to a degree $n$ cover $\alg{\tilde G}$ of a quasisplit group $\alg{G}$ and an injective character $\epsilon \From \mu_n \Into \CC^\times$.  Fix $\alg{\tilde G}$, $\alg{G}$, and $\epsilon$ throughout.  Also, choose a Borel subgroup containing a maximally split maximal torus $\alg{B} \supset \alg{T}$; we will see that our construction is independent of this choice (in a 2-categorical sense).

We make one assumption about our cover $\alg{\tilde G}$, which enables our construction and is essentially nonrestrictive.
\begin{assumption}[Odd $n$ implies even $Q$]
\label{OddnEvenQ}
If $n$ is odd, then we assume $Q \From \sheaf{Y} \To \ZZ$ takes only even values.
\end{assumption}

If $\alg{\tilde G}$ does not satisfy this assumption, i.e., $n$ is odd and $Q(y)$ is odd for some $y \in \sheaf{Y}$, then replace $\alg{\tilde G}$ by $(n+1) \dot \times \alg{\tilde G}$ (its Baer sum with itself $n+1$ times).  The first Brylinski-Deligne invariant becomes $(n+1) Q$, which is even-valued.  By Proposition \ref{DualGroupModn}, the dual group $\dgp{\tilde G^\vee}$ does not change since $Q \equiv (n+1) Q$ modulo $n$.  Moreover, the resulting extensions of groups over local or global fields, e.g., $\mu_n \Into \tilde G \Onto G$, remain the same (up to natural isomorphism).  Indeed, the Baer sum of $\tilde G$ with itself $n+1$ times is naturally isomorphic to the pushout via the $(n+1)^{\th}$ power map $\mu_n \To \mu_n$, which equals the identity map.   

We work with sheaves of abelian groups on $S_{\et}$, and great care is required to avoid confusion between those in the left column and the right column below.  Define
\begin{align*}
\sheaf{\hat T} = \shom(\sheaf{Y}_{Q,n}, \sheaf{G}_m), & \quad \sheaf{\tilde T^\vee} = \shom(\sheaf{Y}_{Q,n}, \CC^\times); \\
\sheaf{\hat T}_{\SC} = \shom(\sheaf{Y}_{Q,n}^{\SC}, \sheaf{G}_m), & \quad \sheaf{\tilde T}_{\SC}^\Vee = \shom(\sheaf{Y}_{Q,n}^{\SC}, \CC^\times); \\
\sheaf{\hat Z} = \shom(\sheaf{Y}_{Q,n} / \sheaf{Y}_{Q,n}^{\SC}, \sheaf{G}_m), & \quad \sheaf{\tilde Z^\vee} = \shom(\sheaf{Y}_{Q,n} / \sheaf{Y}_{Q,n}^{\SC}, \CC^\times).
\end{align*}
Here $\CC^\times$ denotes the constant sheaf on $S_{\et}$.  Thus, in the right column, we find the complex points of the dual groups,
$$\sheaf{\tilde T^\vee} = \dgp{\tilde T}^\Vee(\CC), \quad \sheaf{\tilde T}_{\SC}^\Vee = \dgp{\tilde T}_{\SC}^\Vee(\CC), \quad \sheaf{\tilde Z^\vee} = \dgp{\tilde Z}^\Vee(\CC).$$

Composing with $\epsilon$ defines homomorphisms of local systems of abelian groups,
\begin{align*}
\sheaf{\hat T}_{[n]} = \shom(\sheaf{Y}_{Q,n}, \mu_n) & \xrightarrow{\epsilon} \sheaf{\tilde T^\vee}; \\
\sheaf{\hat T}_{\SC, {[n]}} = \shom(\sheaf{Y}_{Q,n}^{\SC}, \mu_n) & \xrightarrow{\epsilon} \sheaf{\tilde T}_{\SC}^\Vee; \\
\sheaf{\hat Z}_{[n]} = \shom(\sheaf{Y}_{Q,n} / \sheaf{Y}_{Q,n}^{\SC}, \mu_n) & \xrightarrow{\epsilon} \sheaf{\tilde Z^\vee}.
\end{align*}

\subsection{The gerbe associated to a cover of a torus}

Associated to the cover $\alg{\tilde T} = (\alg{T}', n)$, the second Brylinski-Deligne invariant is a central extension of sheaves of groups on $S_{\et}$,
$$\sheaf{G}_m \Into \sheaf{D} \Onto \sheaf{Y}.$$
The commutator of this extension is given in \cite[Proposition 3.11]{B-D},
\begin{equation}
\label{CommForm}
\Comm(y_1, y_2) = (-1)^{n \beta_Q(y_1, y_2)}, \text{ for all } y_1, y_2 \in \sheaf{Y}.
\end{equation}
Pulling back via $\sheaf{Y}_{Q,n} \Into \sheaf{Y}$, we find an extension of sheaves of groups,
\begin{equation}
\label{BD2Qn}
\sheaf{G}_m \Into \sheaf{D}_{Q,n} \Onto \sheaf{Y}_{Q,n}.
\end{equation}

\begin{proposition}
$\sheaf{D}_{Q,n}$ is a commutative extension.
\end{proposition}
\proof
If $n$ is even and $y_1, y_2 \in \sheaf{Y}_{Q,n}$, then $\beta_Q(y_1, y_2) \in \ZZ$ and $n \beta_Q(y_1, y_2) \in 2 \ZZ$.  On the other hand, if $n$ is odd, Assumption \ref{OddnEvenQ} implies that 
$$n \beta_Q(y_1, y_2) = Q(y_1 + y_2) - Q(y_1) - Q(y_2) \in 2 \ZZ.$$
The commutator formula \eqref{CommForm} finishes the proof.
\qed

Let $\sspl(\sheaf{D}_{Q,n})$ denote the sheaf of {\em splittings} of the commutative extension \eqref{BD2Qn}.  In other words, $\sspl(\sheaf{D}_{Q,n})$ is the subsheaf of $\shom(\sheaf{Y}_{Q,n}, \sheaf{D}_{Q,n})$ consisting of homomorphisms which split \eqref{BD2Qn}.

Over any finite \'etale $U \To S$ splitting $\alg{T}$, $\sheaf{Y}_{Q,n}$ restricts to a constant sheaf of free abelian groups.  Thus $\sspl(\sheaf{D}_{Q,n})$ is a $\sheaf{\hat T}$-torsor on $S_{\et}$, which obtains a point over any such $U$.  The equivalence class of this torsor is determined by its cohomology class $\left[ \sspl(\sheaf{D}_{Q,n}) \right] \in H_{\et}^1(S, \sheaf{\hat T})$.

Consider the Kummer sequence, $\sheaf{\hat T}_{[n]} \Into \sheaf{\hat T}\xtwoheadrightarrow{n} \sheaf{\hat T}$.  Write $\Kum$ (for Kummer) for the connecting map in cohomology, $\Kum \From H_{\et}^1(S, \sheaf{\hat T}) \To H_{\et}^2(S, \sheaf{\hat T}_{[n]})$.  This map in cohomology corresponds to the functor which sends a $\sheaf{\hat T}$-torsor to its gerbe of $n^{\th}$ roots (see \ref{AppendixGerbeRoots} for details).  We write
$\sqrt[n]{\sspl(\sheaf{D}_{Q,n})}$ for the gerbe of $n^{\th}$ roots of the $\sheaf{\hat T}$-torsor $\sspl(\sheaf{D}_{Q,n})$.  It is banded by the local system $\sheaf{\hat T}_{[n]}$ and its equivalence class satisfies
$$\left[ \sqrt[n]{\sspl(\sheaf{D}_{Q,n})} \right] = \Kum [\sspl(\sheaf{D}_{Q,n})].$$
Finally we push out via the homomorphism of local systems,
$$\epsilon \From \sheaf{\hat T}_{[n]} = \shom(\sheaf{Y}_{Q,n}, \mu_n) \To \shom(\sheaf{Y}_{Q,n}, \CC^\times) = \sheaf{\tilde T^\vee}.$$

\begin{definition}
The \defined{gerbe associated to the cover} $\alg{\tilde T}$ is defined by
$$\gerb{E}_\epsilon(\alg{\tilde T}) \defeq \epsilon_\ast \sqrt[n]{\sspl(\sheaf{D}_{Q,n})}.$$
It is a gerbe on $S_{\et}$ banded by the local system of abelian groups $\sheaf{\tilde T^\vee}$.
\end{definition}

\begin{example}
\label{SplitTorusGerbeNeutral}
Suppose that $\alg{T}$ is a split torus.  Then the exact sequence of sheaves $\sheaf{G}_m \Into \sheaf{D}_{Q,n} \Onto \sheaf{Y}_{Q,n}$ splits.  Indeed, $\sheaf{Y}_{Q,n}$ is a constant sheaf of free abelian groups, and Hilbert's Theorem 90 gives a short exact sequence
\begin{equation}
\label{SPointsD}
\sheaf{G}_m[S] \Into \sheaf{D}_{Q,n}[S] \Onto \sheaf{Y}_{Q,n}[S].
\end{equation}
Since $\sheaf{Y}_{Q,n}[S]$ is a free abelian group, this exact sequence splits, and any such splitting defines an $S$-point of the torsor $\sspl(\sheaf{D}_{Q,n})$.  An $S$-point of $\sspl(\sheaf{D}_{Q,n})$, in turn, neutralizes of the gerbe $\sqrt[n]{\sspl(\sheaf{D}_{Q,n})}$.

Thus when $\alg{T}$ is a split torus, the gerbe $\gerb{E}_\epsilon(\alg{\tilde T})$ is trivial.  Any splitting of the sequence \eqref{SPointsD} defines a neutralization of $\gerb{E}_\epsilon(\alg{\tilde T})$.
\end{example}  

\subsection{The gerbe of liftings}

Recall that $\sheaf{Y}_{Q,n}^{\SC}$ denotes the subgroup of $\sheaf{Y}_{Q,n}$ spanned by the modified coroots $\tilde \Phi^\vee$, and $\sheaf{\hat T}_{\SC} = \shom(\sheaf{Y}_{Q,n}^{\SC}, \sheaf{G}_m)$.  The inclusion $\sheaf{Y}_{Q,n}^{\SC} \Into \sheaf{Y}_{Q,n}$ corresponds to a surjective homomorphism,
$$p \From \sheaf{\hat T} \To \sheaf{\hat T}_{\SC}.$$
The extension $\sheaf{G}_m \Into \sheaf{D}_{Q,n} \Onto \sheaf{Y}_{Q,n}$ pulls back via $\sheaf{Y}_{Q,n}^{\SC} \Into \sheaf{Y}_{Q,n}$ to an extension,
$$\sheaf{G}_m \Into \sheaf{D}_{Q,n}^{\SC} \Onto \sheaf{Y}_{Q,n}^{\SC}.$$
A splitting of $\sheaf{D}_{Q,n}$ pulls back to a splitting of $\sheaf{D}_{Q,n}^{\SC}$, providing a map of torsors,
$$p^\ast \From \sspl(\sheaf{D}_{Q,n}) \To \sspl(\sheaf{D}_{Q,n}^{\SC}),$$
lying over $p \From \sheaf{\hat T} \To \sheaf{\hat T}_{\SC}$.  Taking $n^{\th}$ roots of torsors gives a functor of gerbes,
$$\sqrt[n]{p^\ast}\From \sqrt[n]{\sspl(\sheaf{D}_{Q,n})} \To \sqrt[n]{\sspl(\sheaf{D}_{Q,n}^{\SC})},$$
lying over $p \From \sheaf{\hat T}_{[n]} \To \sheaf{\hat T}_{\SC, [n]}$ (see Appendix \ref{AppendixGerbeRoots}).

Recall that $\sheaf{\tilde T^\vee} = \shom(\sheaf{Y}_{Q,n}, \CC^\times)$ and $\sheaf{\tilde T}_{\SC}^\Vee = \shom(\sheaf{Y}_{Q,n}^{\SC}, \CC^\times)$.  Define $\gerb{E}_\epsilon^{\SC}(\alg{\tilde T}) \defeq \epsilon_\ast \sqrt[n]{\sspl(\sheaf{D}_{Q,n}^{\SC})}$ by analogy to $\gerb{E}_\epsilon(\alg{\tilde T}) = \epsilon_\ast \sqrt[n]{\sspl(\sheaf{D}_{Q,n})}$.  Pushing out via $\epsilon$, the functor $\sqrt[n]{p^\ast}$ yields a functor of gerbes
$$\gerb{p} = \epsilon_\ast \sqrt[n]{p^\ast} \From \gerb{E}_\epsilon(\alg{\tilde T}) \To \gerb{E}_\epsilon^{\SC}(\alg{\tilde T}),$$
lying over the homomorphism $p \From \sheaf{\tilde T^\vee} \To \sheaf{\tilde T}_{\SC}^\Vee$.

In the next section, we define the {\em Whittaker torsor}, which gives an object $\whit$ neutralizing the gerbe $\gerb{E}_\epsilon^{\SC}(\alg{\tilde T})$.  We take this construction of $\whit$ for granted at the moment.

\begin{definition}
Define $\gerb{E}_\epsilon(\alg{\tilde G})$ to be the gerbe $\gerb{p}^{-1}(\whit)$ of liftings of $\whit$ via $\gerb{p}$ (see \ref{AppendixGerbeLiftings}).  In other words, $\gerb{E}_\epsilon(\alg{\tilde G})$ is the 
category of pairs $(\gerb{e}, j)$ where $\gerb{e}$ is an object of $\gerb{E}_\epsilon(\alg{\tilde T})$ and $j \From \gerb{p}(\gerb{e}) \To \whit$ is an isomorphism in $\gerb{E}_\epsilon^{\SC}(\alg{\tilde T})$. This is a gerbe on $S_{\et}$ banded by $\sheaf{\tilde Z^\vee} = \Ker(\sheaf{\tilde T^\vee} \xrightarrow{p} \sheaf{\tilde T}_{\SC}^\Vee)$.
\end{definition}

The cohomology classes of our gerbes now fit into a sequence
$$\begin{tikzcd}[row sep = 0em]
\left[ \gerb{E}_\epsilon(\alg{\tilde G}) \right] \arrow[mapsto]{r} & \left[ \gerb{E}_\epsilon(\alg{\tilde T}) \right] \arrow[mapsto]{r} & \left[ \gerb{E}_\epsilon^{\SC}(\alg{\tilde T}) \right] = 0 \\
\rotatebox{-90}{$\in$} & \rotatebox{-90}{$\in$}& \rotatebox{-90}{$\in$} \\
H_{\et}^2(S, \sheaf{\tilde Z^\vee}) \arrow{r} & H_{\et}^2(S, \sheaf{\tilde T^\vee}) \arrow{r} &H_{\et}^2(S, \sheaf{\tilde T_{\SC} ^\vee}) 
\end{tikzcd}$$

\begin{remark}
The construction of this gerbe relies on the (soon-to-be-defined) Whittaker torsor in a crucial way.  We view this as a good thing, since any putative Langlands correspondence should also connect the existence of Whittaker models to properties of the Langlands parameter (cf. \cite{VoganLLC}).
\end{remark}

\subsection{The Whittaker torsor}

Now we construct the object $\whit$ neutralizing the gerbe $\gerb{E}_\epsilon^{\SC}(\alg{\tilde T})$ over $S$.  Let $\alg{U}$ denote the unipotent radical of the Borel subgroup $\alg{B} \subset \alg{G}$, and let $\sheaf{U}$ be the sheaf of groups on $S_{\et}$ that it represents.  Let $\alg{G}_a$ denote the additive group scheme over $S$, and $\sheaf{G}_a$ the sheaf of groups on $S_{\et}$ that it represents.  Recall that $\Delta \subset \Phi$ denotes the subset of simple roots corresponding to the Borel subgroup $\alg{B}$.  

For $S' \To S$ finite \'etale and splitting $\alg{T}$, and $\alpha \in \Delta[S']$, write $\alg{U}_\alpha$ for the one-dimensional root subgroup of $\alg{U}_{S'}$ associated to $\phi$.  Let $\sheaf{U}_\alpha$ be the associated sheaf of abelian groups on $S_{\et}'$.  Write $\shom^\ast(\sheaf{U}_\alpha, \sheaf{G}_{a})$ for the sheaf (on $S_{\et}'$) of isomorphisms from $\sheaf{U}_\alpha$ to $\sheaf{G}_{a,S'}$.  The sheaf $\shom^\ast(\sheaf{U}_\alpha, \sheaf{G}_{a})$ naturally forms a $\sheaf{G}_m$-torsor on $S_{\et}'$, by the formula
$$[h \ast \xi](u) = h^{-1} \cdot \xi(u) \text{ for all } h \in \sheaf{G}_m, \xi \in \shom^\ast(\sheaf{U}_\alpha, \sheaf{G}_a).$$
\begin{definition}
The \defined{Whittaker torsor} is the subsheaf $\Whit \subset \shom(\sheaf{U}, \sheaf{G}_a)$ consisting of those homomorphisms which (locally on $S_{\et}$)  restrict to an isomorphism on every simple root subgroup.  The sheaf $\Whit$ is given the structure of a $\sheaf{\hat T}_{\SC}$-torsor as follows:  for a Galois cover $S' \To S$ splitting $\alg{T}$, we have
$$\sheaf{\hat T}_{\SC}[S'] = \shom(\sheaf{Y}_{Q,n}^{\SC}, \sheaf{G}_m)[S'] = \shom \left( \bigoplus_{\alpha \in \Delta[S']} \ZZ {\tilde \alpha}^\vee, \sheaf{G}_{m} \right) [S'] \ident \prod_{\alpha \in \Delta[S']} \sheaf{G}_{m}[S'].$$
Similarly, we can decompose the Whittaker sheaf
$$\Whit[S'] \ident \bigoplus_{\alpha \in \Delta[S']} \shom^\ast(\sheaf{U}_\alpha, \sheaf{G}_{a})[S'].$$
The $\sheaf{G}_m$-torsor structure on $\shom^\ast(\sheaf{U}_\alpha, \sheaf{G}_{a})$ yields (simple root by simple root) a $\sheaf{\hat T}_{\SC}$-torsor structure on $\Whit$.  Although we have defined the torsor structure locally on $S_{\et}$, the action descends since the $\Gal(S'/S)$-actions are compatible throughout.
\end{definition}

The third Brylinski-Deligne invariant of $\alg{\tilde G}$ is a homomorphism $f \From \sheaf{D}_Q \To \sheaf{D}$ of groups on $S_{\et}$.
$$\begin{tikzcd}
\sheaf{G}_m \inarrow{r} \arrow{d}{=} &  \sheaf{D}_Q \onarrow{r} \inarrow{d}{f} & \sheaf{Y}^{\SC} \inarrow{d} \\
\sheaf{G}_m \inarrow{r} & \sheaf{D} \onarrow{r} & \sheaf{Y}
\end{tikzcd}$$
Here $\sheaf{D}_Q$ is a sheaf on $S_{\et}$ which depends (up to unique isomorphism) only on the Weyl- and Galois-invariant quadratic form $Q \From \sheaf{Y}^{\SC} \To \ZZ$.  This is reviewed in \cite[\S 1.3, 3.3]{MWIntegral}, and characterized in \cite[\S 11]{B-D} when working over a field.

Consider a Galois cover $S' \To S$ splitting $\alg{T}$ as before.  For any $\eta \in \Whit[S']$, and any simple root $\alpha \in \Delta[S']$, there exists a unique element $e_{\eta, \alpha} \in \sheaf{U}_\alpha[S']$ such that $\eta(e_{\eta,\alpha}) = 1$.  From these, \cite[\S 11.2]{B-D} gives elements $[e_{\eta,\alpha}] \in \sheaf{D}_Q[S']$ lying over the simple coroots $\alpha^\vee \in \sheaf{Y}^{\SC}[S']$.
\begin{remark}
When $S = \Spec(F)$ this follows directly from \cite[\S 11.2]{B-D}.  When $S = \Spec(\OO)$, $S' = \Spec(\OO')$, $\eta \in \Whit[S']$, and $F'$ is the fraction field of $\OO'$, we find elements $e_{\eta,\alpha} \in \sheaf{U}_\alpha[F']$; as $\eta$ gives an isomorphism from $\sheaf{U}_\alpha$ to $\sheaf{G}_a$ (as sheaves of groups on $\OO_{\et}'$), it follows that $e_{\eta,\alpha} \in \sheaf{U}_\alpha[\OO']$ as well.  The map $e \mapsto [e]$ of \cite[\S 11.1]{B-D} similarly makes sense over $\OO'$ as well as it does over a field; since we assume $\alg{G}$ is a reductive group over $\OO$, split over $\OO'$, every root $SL_2$ over $F'$ arises from one over $\OO'$.  Thus the results of \cite[\S 11.2]{B-D} apply in the setting of $S = \Spec(\OO)$ as well as in the setting of a field.
\end{remark}

Using the elements $[e_{\eta,\alpha}] \in \sheaf{D}_Q[S']$ lying over the simple coroots $\alpha^\vee$, define
$$\omega(\eta)(\tilde \alpha^\vee) \defeq r_\alpha f([e_{\eta,\alpha}])^{n_\alpha} \in \sheaf{D}_{Q,n}^{\SC}[S'], \text{ lying over } \tilde \alpha^\vee = n_\alpha \alpha^\vee \in \sheaf{Y}_{Q,n}^{\SC}[S'],$$
where the sign $r_\alpha$ is defined by
$$r_\alpha \defeq (-1)^{\frac{Q(\alpha^\vee) n_\alpha (n_\alpha - 1)}{2}}.$$
The map $\omega(\eta) \From \tilde \alpha^\vee \mapsto r_\alpha f([e_{\eta, \alpha}])^{n_\alpha}$ extends uniquely to a splitting of the sequence
\begin{equation}
\begin{tikzcd}
\sheaf{G}_m[S'] \inarrow{r} & \sheaf{D}_{Q,n}^{\SC}[S'] \onarrow{r} & \sheaf{Y}_{Q,n}^{\SC}[S'] \arrow[bend right=20]{l}[swap]{\omega(\eta)}.
\end{tikzcd}
\end{equation}
As $(\sheaf{Y}_{Q,n}^{\SC})_{S'}$ is a constant sheaf, this gives an element $\omega(\eta) \in \sspl(\sheaf{D}_{Q,n}^{\SC})[S']$.  Allowing $\eta$ to vary, and appyling Galois descent (cf.~\cite[Proposition 11.7]{B-D}), we find a map of sheaves on $S_{\et}$,
$$\omega \From \Whit \To \sspl(\sheaf{D}_{Q,n}^{\SC}).$$
To summarize, $\omega$ is the map that sends a nondegenerate character $\eta$ of $\alg{U}$ to the splitting $\omega(\eta)$, which (locally on $S_{\et}$) sends each modified simple coroot $\tilde \alpha^\vee$ to the element $r_\alpha f([e_{\eta, \alpha}])^{n_\alpha}$ of $\sheaf{D}_{Q,n}$. 
\begin{remark}
For the purposes of this paper, there is some flexibility in the choice of signs $r_\alpha$.  The signs here are defined in such a way that our hypotheisized local Langlands correspondence for covers matches what is known for covers of $\alg{SL}_2$, e.g., metaplectic correspondences of Shimura and Waldspurger.
\end{remark}

Both $\Whit$ and $\sspl(\sheaf{D}_{Q,n}^{\SC})$ are $\sheaf{\hat T}_{\SC}$-torsors, and the following proposition describes how $\omega$ interacts with the torsor structure.
\begin{proposition}
Let $\nu \From \sheaf{\hat T}_{\SC} \To \sheaf{\hat T}_{\SC}$ be the homomorphism corresponding to the unique homomorphism $\sheaf{Y}_{Q,n}^{\SC} \To \sheaf{Y}_{Q,n}^{\SC}$ which sends $\tilde \alpha^\vee$ to $- n_\alpha Q(\alpha^\vee) \tilde \alpha^\vee$ for all simple roots $\alpha$.  Then $\omega$ lies over $\nu$, i.e., the following diagram commutes.
$$\begin{tikzcd}
\sheaf{\hat T_{\SC}} \times \Whit \arrow{r}{\ast} \arrow{d}{\nu \times \omega} & \Whit \arrow{d}{\omega} \\
\sheaf{\hat T_{\SC}} \times \sspl(\sheaf{D}_{Q,n}^{\SC}) \arrow{r}{\ast} & \sspl(\sheaf{D}_{Q,n}^{\SC})
\end{tikzcd}$$
\end{proposition}
\proof
We must trace through the action of $\sheaf{\hat T}_{\SC} = \shom(\sheaf{Y}_{Q,n}^{\SC}, \sheaf{G}_m)$; we work over a finite \'etale cover of $S$ over which $\alg{T}$ splits in what follows.  Then, for any simple root $\alpha \in \Delta$, and any element $h \in \sheaf{G}_m$, there exists a unique element $h_\alpha \in \sheaf{\hat T}_{\SC}$ such that for all $\beta \in \Delta$,
$$h_\alpha(\tilde \beta^\vee) = \begin{cases} 1 & \text{ if } \beta \neq \alpha; \\ h & \text{ if } \beta = \alpha.  \end{cases}$$
If $\eta \in \Whit$ then $[h_\alpha \ast \eta](e_{h_\alpha \ast \eta, \alpha}) = 1$ and so $\eta(e_{h_\alpha \ast \eta, \alpha}) = h$.  Therefore,
$$e_{h_\alpha \ast \eta, \beta} = \begin{cases} e_{\eta, \beta} & \text{ if } \beta \neq \alpha; \\ h \ast  e_{\eta, \alpha} & \text{ if } \beta = \alpha. \end{cases}$$
If $\beta \neq \alpha$, then $\omega(h_\alpha \ast \eta)(\tilde \beta^\vee) = r_\beta f([e_{h_\alpha \ast \eta, \beta}])^{n_\beta} = r_\beta f([e_{\eta, \beta}])^{n_\beta} = \omega(\eta)(\tilde \beta^\vee)$.  On the other hand, in the case $\beta = \alpha$ we compute using \cite[Equation (11.2.1)]{B-D},
\begin{align*}
\omega(h_\alpha \ast \eta)(\tilde \alpha^\vee) = r_\alpha f([e_{h_\alpha \ast \eta, \alpha}])^{n_\alpha} &= r_\alpha f \left( [h \ast e_{\eta, \alpha} ] \right)^{n_\alpha} \\
&= r_\alpha f \left( h^{-Q(\alpha^\vee)} \cdot [e_{\eta, \alpha} ] \right)^{n_\alpha} \\
&= r_\alpha h^{ - n_\alpha Q(\phi_\alpha^\vee) } \cdot f \left( [e_{\eta, \alpha} ] \right)^{n_\alpha} \\ 
&= h^{ - n_\alpha Q(\phi_\alpha^\vee) } \cdot \omega(\eta)(\tilde \alpha^\vee).
\end{align*}
This computation demonstrates that the diagram commutes as desired.
\qed


Now let $\mu \From \sheaf{\hat T}_{\SC} \Onto \sheaf{\hat T}_{\SC}$ be the homomorphism corresponding to the unique homomorphism $\sheaf{Y}_{Q,n}^{\SC} \Into \sheaf{Y}_{Q,n}^{\SC}$ which sends $\tilde \alpha^\vee$ to $- m_\alpha \tilde \alpha^\vee$ for all $\alpha \in \Delta$.  As $Q(\alpha^\vee) n_\alpha = m_\alpha \cdot n$, we find that $\nu = n \circ \mu$, where $n$ denotes the $n^{\th}$-power map.

Let $\mu_\ast \Whit$ denote the pushout of the $\sheaf{\hat T}_{\SC}$-torsor $\Whit$, via $\mu$.  Since $\nu$ factors through $\mu$, we find that $\omega \From \Whit \To \sspl(\sheaf{D}_{Q,n}^{\SC})$ factors uniquely through $\bar \omega \From \mu_\ast \Whit \To \sspl(\sheaf{D}_{Q,n}^{\SC})$, making the following diagram commute.
$$\begin{tikzcd}
\sheaf{\hat T_{\SC}} \times \left( \mu_\ast \Whit \right) \arrow{r}{\ast} \arrow{d}{n \times \bar \omega} & \mu_\ast \Whit \arrow{d}{\bar \omega}  \\
\sheaf{\hat T_{\SC}} \times \sspl(\sheaf{D}_{Q,n}^{\SC}) \arrow{r}{\ast} & \sspl(\sheaf{D}_{Q,n}^{\SC})
\end{tikzcd}$$

The pair $(\mu_\ast \Whit, \bar \omega)$ is therefore an object of the category $\sqrt[n]{\sspl(\sheaf{D}_{Q,n}^{\SC})}[S]$; it {\em neutralizes} the gerbe $\sqrt[n]{\sspl(\sheaf{D}_{Q,n}^{\SC})}$.  In particular,
$$\left[ \sqrt[n]{\sspl(\sheaf{D}_{Q,n}^{\SC})} \right] = 0.$$   
Write $\whit = (\mu_\ast \Whit, \bar \omega)$ for this object.  Pushing out via $\epsilon$, we view $\whit$ as an $S$-object of $\gerb{E}_\epsilon^{\SC}(\alg{\tilde T})$.  This completes the construction of the gerbe $\gerb{E}_\epsilon(\alg{\tilde G}) = \gerb{p}^{-1}(\gerb{w})$ associated to the cover $\alg{\tilde G}$ and character $\epsilon$.

\begin{example}
Suppose that $\sheaf{Y}_{Q,n} / \sheaf{Y}_{Q,n}^{\SC}$ is torsion-free and a constant sheaf (equivalently, the center of $\dgp{\tilde G^\vee}$ is connected and constant as a sheaf on $S_{\et}$).  Then the following short exact sequence splits:
$$\sheaf{Y}_{Q,n}^{\SC} \Into \sheaf{Y}_{Q,n} \Onto \sheaf{Y}_{Q,n} /    \sheaf{Y}_{Q,n}^{\SC}.$$

Given such a splitting, write $\sheaf{Y}_{Q,n}^{\cent} \subset \sheaf{Y}_{Q,n}$ for the image of $\sheaf{Y}_{Q,n} / \sheaf{Y}_{Q,n}^{\SC}$ via the splitting.  The identification $\sheaf{Y}_{Q,n} = \sheaf{Y}_{Q,n}^{\SC} \oplus \sheaf{Y}_{Q,n}^{\cent}$ corresponds to an isomorphism $\sheaf{\hat T} \xrightarrow{\sim} \sheaf{\hat T}_{\SC} \times \sheaf{\hat Z}$.   Let $\sheaf{D}_{Q,n}^{\cent}$ be the pullback of $\sheaf{D}_{Q,n}$ to $\sheaf{Y}_{Q,n}^{\cent}$.  From Example \ref{SplitTorusGerbeNeutral}, the short exact sequence $\sheaf{G}_m \Into \sheaf{D}_{Q,n}^{\cent} \Onto \sheaf{Y}_{Q,n}^{\cent}$ splits, providing an object of $\sqrt[n]{\sspl(\sheaf{D}_{Q,n}^{{\cent}})}$.  

Chasing diagrams gives a map of object sets,
$$\sqrt[n]{\sspl(\sheaf{D}_{Q,n}^{\SC})}[S] \times \sqrt[n]{\sspl(\sheaf{D}_{Q,n}^{{\cent}})}[S] \To \sqrt[n]{\sspl(\sheaf{D}_{Q,n})}[S].$$
A splitting of $\sheaf{D}_{Q,n}^{\cent}$ gives an object of $\sqrt[n]{\sspl(\sheaf{D}_{Q,n}^{{\cent}})}[S]$ and $\whit$ provides an object of $\sqrt[n]{\sspl(\sheaf{D}_{Q,n}^{\SC})}[S]$.  Hence the gerbe $\gerb{E}_\epsilon(\alg{\tilde G})$ is neutral when $\sheaf{Y}_{Q,n} / \sheaf{Y}_{Q,n}^{\SC}$ is torsion-free and a constant sheaf.
\end{example}

\subsection{Well-aligned functoriality}
\label{WAFGerbe}

Consider a well-aligned homomorphism $\tilde \iota \From \alg{\tilde G}_1 \To \alg{\tilde G}_2$ of covers, each endowed with Borel subgroup and maximally split maximal torus, i.e., a morphism in the category $\Cat{WAC}_S$.  Fix $\epsilon$ as before.  We have constructed gerbes $\gerb{E}_\epsilon(\alg{\tilde G}_1)$ and $\gerb{E}_\epsilon(\alg{\tilde G}_2)$ associated to $\alg{\tilde G}_1$ and $\alg{\tilde G}_2$, banded by $\sheaf{\tilde Z}_1^\Vee$ and $\sheaf{\tilde Z}_2^\Vee$, respectively.  We have constructed a homomorphism of dual groups $\iota^\vee \From \dgp{\tilde G}_2^\Vee \To \dgp{\tilde G}_1^\Vee$ in Section \ref{WAFDualGroup}, which (after taking $\CC$-points) restricts to $\iota^\vee \From \sheaf{\tilde Z}_2^\Vee \To \sheaf{\tilde Z}_1^\Vee$.  Here we construct a functor of gerbes $\gerb{i} \From \gerb{E}_\epsilon(\alg{\tilde G}_2) \To \gerb{E}_\epsilon(\alg{\tilde G}_1)$, lying over $\iota^\vee \From \sheaf{\tilde Z}_{2}^\Vee \To \sheaf{\tilde Z}_{1}^\Vee$.

Well-alignedness give a commutative diagram in which the first row is the pullback of the second.
$$\begin{tikzcd}
\alg{K}_2 \inarrow{r} \arrow{d}{=} & \alg{T}_1' \onarrow{r} \arrow{d}{\iota'} & \alg{T}_1 \arrow{d}{\iota} \\
\alg{K}_2 \inarrow{r} & \alg{T}_2' \onarrow{r} & \alg{T}_2
\end{tikzcd}$$
This gives a commutative diagram for the second Brylinski-Deligne invariants.  After pulling back to $\sheaf{Y}_{1,Q_1,n}$ and $\sheaf{Y}_{2,Q_2,n}$, we get a commutative diagram of sheaves of abelian groups on $S_{\et}$, in which the top row is the pullback of the bottom.
$$\begin{tikzcd}
\sheaf{G}_m \inarrow{r} \arrow{d}{=} & \sheaf{D}_{1,Q_1,n} \onarrow{r} \arrow{d}{\iota'} & \sheaf{Y}_{1,Q_1,n} \arrow{d}{\iota} \\
\sheaf{G}_m \inarrow{r} & \sheaf{D}_{2,Q_2,n} \onarrow{r}& \sheaf{Y}_{2,Q_2,n}
\end{tikzcd}$$
(Assumption \ref{OddnEvenQ} is in effect, so $Q_1$ and $Q_2$ are even-valued if $n$ is odd.)

We have homomorphisms of sheaves of abelian groups,
$$\iota \From \sheaf{Y}_{1,Q_1,n} \To \sheaf{Y}_{2,Q_2,n}, \quad \sheaf{Y}_{1,Q_1,n}^{\SC} \To \sheaf{Y}_{2,Q_2,n}^{\SC}, \quad \frac{\sheaf{Y}_{1,Q_1,n}}{\sheaf{Y}_{1,Q_1,n}^{\SC}} \To \frac{\sheaf{Y}_{2,Q_2,n}}{\sheaf{Y}_{2,Q_2,n}^{\SC}}.$$
Applying $\shom(\bullet, \sheaf{G}_m)$ yields homomorphisms of sheaves of abelian groups,
$$\hat \iota \From \sheaf{\hat T}_2 \To \sheaf{\hat T}_1, \quad \sheaf{\hat T}_{\SC,2} \To \sheaf{\hat T}_{\SC,1}, \quad \sheaf{\hat Z}_2 \To \sheaf{\hat Z}_1.$$

A splitting of $\sheaf{D}_{2,Q_2,n}$ pulls back to a splitting of $\sheaf{D}_{1,Q_1,n}$ (call this pullback map $\iota^\ast$), giving a commutative diagram of sheaves on $S_{\et}$.
$$\begin{tikzcd}
\sheaf{\hat T}_2 \times \sspl(\sheaf{D}_{2,Q_2,n}) \arrow{r}{\ast} \arrow{d}{\hat \iota \times \iota^\ast} & \sspl(\sheaf{D}_{2,Q_2,n}) \arrow{d}{\iota^\ast} \\
\sheaf{\hat T}_1 \times \sspl(\sheaf{D}_{1,Q_1,n}) \arrow{r}{\ast} & \sspl(\sheaf{D}_{1,Q_1,n})
\end{tikzcd}$$
This defines a functor of gerbes 
$$\sqrt[n]{\iota^\ast} \From \sqrt[n]{\sspl(\sheaf{D}_{2,Q_2,n})} \To \sqrt[n]{\sspl(\sheaf{D}_{1,Q_1,n})},$$
lying over $\hat \iota \From \sheaf{\hat T}_{2,[n]} \To \sheaf{\hat T}_{1,[n]}$.  Pushing out via $\epsilon$ yields a functor of gerbes,
$$\gerb{i} \From \gerb{E}_\epsilon(\alg{\tilde T}_2) \To \gerb{E}_\epsilon(\alg{\tilde T}_1).$$
The same process applies to $\iota \From \sheaf{Y}_{1,Q_1,n}^{\SC} \To \sheaf{Y}_{2,Q_2,n}^{\SC}$, giving a functor of gerbes, $\gerb{i}^{\SC} \From \gerb{E}_\epsilon^{\SC}(\alg{\tilde T}_2) \To \gerb{E}_\epsilon^{\SC}(\alg{\tilde T}_1)$.  By pulling back in stages, we find a square of gerbes and functors, and a natural isomorphism $\Fun{S} \From \gerb{p}_1 \circ \gerb{i} \xRightarrow{\sim} \gerb{i}^{\SC} \circ \gerb{p}_2$ making the diagram 2-commute.
$$\begin{tikzcd}
 \gerb{E}_\epsilon(\alg{\tilde T}_2) \arrow{r}{\gerb{i}} \arrow{d}{\gerb{p}_2} & \gerb{E}_\epsilon(\alg{\tilde T}_1) \arrow{d}{\gerb{p}_1} \\
\gerb{E}_\epsilon^{\SC}(\alg{\tilde T}_2)  \arrow{r}{\gerb{i}^{\SC}} & \gerb{E}_\epsilon^{\SC}(\alg{\tilde T}_1)
\end{tikzcd}$$

If $\alg{U}_1$ and $\alg{U}_2$ are the unipotent radicals of $\alg{B}_1$ and $\alg{B}_2$, respectively, then pulling back gives a map $\iota^\ast \From \Whit_2 \To \Whit_1$.  Condition (1) of well-alignedness states that $\Ker(\iota)$ is contained in the center of $\alg{G}_1$, and so simple root subgroups in $\alg{U}_1$ map isomorphically to simple root subgroups in $\alg{U}_2$.  Compatibility of quadratic forms $Q_1$ and $Q_2$ implies compatibility in the constants $n_\alpha$, $m_\alpha$, and $r_\alpha$ for simple roots, and so we find a commutative diagram
$$\begin{tikzcd}
(\mu_2)_\ast \Whit_2 \arrow{r}{\bar \omega_2} \arrow{d}{\iota^\ast} & \sspl(\sheaf{D}_{2,Q_2,n}^{\SC}) \arrow{d}{\iota^\ast} \\
(\mu_1)_\ast \Whit_1 \arrow{r}{\bar \omega_1} & \sspl(\sheaf{D}_{1,Q_1,n}^{\SC})  
\end{tikzcd}$$
Write $\whit_1$ for the object $((\mu_1)_\ast \Whit_1, \bar \omega_1)$ of $\gerb{E}_\epsilon^{\SC}(\alg{\tilde T}_1)$, and similarly $\whit_2$ for the object $((\mu_2)_\ast \Whit_2, \bar \omega_2)$ of $\gerb{E}_\epsilon^{\SC}(\alg{\tilde T}_2)$.  The commutative diagram above gives an isomorphism $f \From \gerb{i}^{\SC}( \whit_2) \xrightarrow{\sim} \whit_1$ in $\gerb{E}_\epsilon^{\SC}(\alg{\tilde T}_1)$.

If $(\gerbob{e},j)$ is an object of $\gerb{E}_\epsilon(\alg{\tilde G}_2) = \gerb{p}_2^{-1}(\whit_2)$, i.e., $\gerbob{e}$ is an object of $\gerb{E}_\epsilon(\alg{\tilde T}_2)$ and $j \From \gerb{p}_2(\gerbob{e}) \To \whit_2$ is an isomorphism, then we find a sequence of isomorphisms,
$$\gerb{i}(j) \From \gerb{p}_1(\gerb{i}(\gerbob{e})) \xrightarrow{\Fun{S}} \gerb{i}^{\SC}(\gerb{p}_2(\gerbob{e})) \xrightarrow{j} \gerb{i}^{\SC}(\whit_2) \xrightarrow{f} \whit_1.$$
In this way $(\gerb{i}(\gerbob{e}),\gerb{i}(j))$ becomes an object of $\gerb{E}_\epsilon(\alg{\tilde G}_1) = \gerb{p}_1^{-1}(\whit_1)$.  This extends to a functor of gerbes, $\gerb{i} \From \gerb{E}_\epsilon(\alg{\tilde G}_2) \To \gerb{E}_\epsilon(\alg{\tilde G}_1)$, lying over $\iota^\vee \From \sheaf{\tilde Z}_2^\Vee \To \sheaf{\tilde Z}_1^\Vee$.  

Given a pair of well-aligned homomorphisms,
$$\alg{\tilde G}_1 \xrightarrow{\tilde \iota_1} \alg{\tilde G}_2 \xrightarrow{\tilde \iota_2} \alg{\tilde G}_2,$$
with $\tilde \iota_3 = \tilde \iota_2 \circ \tilde \iota_1$, we find three functors of gerbes.
$$\begin{tikzcd}
\gerb{E}_\epsilon(\alg{\tilde G}_3) \arrow{r}[swap]{\gerb{i}_2} \arrow[bend left=20]{rr}{\gerb{i}_3} & \gerb{E}_\epsilon(\alg{\tilde G}_2) \arrow{r}[swap]{\gerb{i}_1} & \gerb{E}_\epsilon(\alg{\tilde G}_1).
\end{tikzcd}$$
lying over three homomorphisms of sheaves,
$$\begin{tikzcd}
\sheaf{\tilde Z}_3^\Vee \arrow{r}[swap]{\iota_2^\vee} \arrow[bend left=20]{rr}{\iota_3^\vee} & \sheaf{\tilde Z}_2^\Vee \arrow{r}[swap]{\iota_1^\vee} & \sheaf{\tilde Z}_1^\Vee.
\end{tikzcd}$$

The functoriality of pullbacks, pushouts, taking $n^{\th}$ roots of torsors, etc., yields a \textbf{natural isomorphism} of functors:
\begin{equation}
\label{CompositionWellAligned}
N(\iota_1, \iota_2) \From \gerb{i}_3 \xRightarrow{\sim} \gerb{i}_1 \circ \gerb{i}_2.
\end{equation}

We can summarize these results using the language of ``weak functors'' \cite[Tag 003G]{Stacks}.  We have given all the data for such a weak functor, and a reader who wishes to check commutativity of a few diagrams may verify the following. 
\begin{proposition}
The construction of the associated gerbe defines a weak functor from $\Cat{WAC}_S$ (the category of triples $(\alg{\tilde G}, \alg{B}, \alg{T})$ and well-aligned homomorphisms),  to the 2-category of gerbes on $S_{\et}$, functors of gerbes, and natural isomorphisms thereof.
\end{proposition}

\subsection{Well-definedness}
\label{GerbeWellDefined}

The construction of the associated gerbe $\gerb{E}_\epsilon(\alg{\tilde G})$ depended on the choice of Borel subgroup $\alg{B}$ and maximally split maximal torus $\alg{T}$.  Now we demonstrate that $\gerb{E}_\epsilon(\alg{\tilde G})$ is well-defined independently of these choices, in a suitable 2-categorical sense.

Consider another choice $\alg{B}_0 \supset \alg{T}_0$.  Our constructions, with these two choices of tori and Borel subgroups, yield two dual groups $\dgp{\tilde G}_0^\Vee$ and $\dgp{\tilde G}^\Vee$ and two gerbes $\gerb{E}_\epsilon(\alg{\tilde G})$ and $\gerb{E}_{0,\epsilon}(\alg{\tilde G})$.
    
Let $\sheaf{Y}_0$ and $\sheaf{Y}$ be the cocharacter lattices of $\alg{T}_0$ and $\alg{T}$, respectively, and $\Phi_0^\vee, \Phi^\vee$ the coroots therein.  The Borel subgroups provide systems of simple coroots $\Delta_0^\vee$ and $\Delta^\vee$, respectively (sheaves of sets on $S_{\et}$).  Similarly we have character lattices $\sheaf{X}_0, \sheaf{X}$, roots $\Phi_0, \Phi$, and simple roots $\Delta_0, \Delta$.  The cover $\alg{\tilde G}$ yields quadratic forms $Q \From \sheaf{Y} \To \ZZ$ and $Q_0 \From \sheaf{Y}_0 \To \ZZ$.  The second Brylinski-Deligne invariant gives extensions $\sheaf{G}_m \Into \sheaf{D} \Onto \sheaf{Y}$ and $\sheaf{G}_m \Into \sheaf{D}_0 \Onto \sheaf{Y}_0$.  

By Proposition \ref{BTConj}, there exists $g \in \alg{G}(S)$ such that $\Int(g) \alg{T}_0 = \alg{T}$ and $\Int(g) \alg{B}_0 = \alg{B}$.  This automorphism $\Int(g)$ lifts to an automorphism of $\alg{G}'$.  This defines a well-aligned isomorphism of covers.
$$\begin{tikzcd}
\alg{K}_{2} \inarrow{r} \arrow{d}{=} & \alg{G}' \onarrow{r} \arrow{d}{\Int(g)'} & \alg{G} \arrow{d}{\Int(g)} & & \alg{B}_{0} \arrow{d} & \alg{T}_{0} \arrow{d}  \\
\alg{K}_{2} \inarrow{r} & \alg{G}' \onarrow{r} & \alg{G} & & \alg{B} & \alg{T}
\end{tikzcd}$$
This well-aligned isomorphism of covers yields an equivalence of gerbes, 
$$\gerb{Int}(g) \From \gerb{E}_\epsilon(\alg{\tilde G}) \To \gerb{E}_{0,\epsilon}(\alg{\tilde G}),$$ 
lying over $\Int(g)^\vee \From \sheaf{\tilde Z^\vee} \To \sheaf{\tilde Z}_0^\Vee$.

Suppose that $g' \in \alg{G}(S)$ also satisfies $\Int(g') \alg{T}_0 = \alg{T}$ and $\Int(g') \alg{B}_0 = \alg{B}$.  As in the proof of Theorem \ref{WellDefinedDualGroup}, $g' = t g$ for a unique $t \in \alg{T}(S)$ and $\Int(g') = \Int(t) \Int(g)$.  

This gives a natural isomorphism of functors,
$$N(g', g) \From \gerb{Int}(g') \xRightarrow{\sim} \gerb{Int}(g) \circ  \gerb{Int}(t)$$
Our upcoming Proposition \ref{IntTFunctors} will provide a natural isomorphism $\gerb{Int}(t) \xRightarrow{\sim} \gerb{Id}$.  Assuming this for the moment, we find a natural isomorphism $\gerb{Int}(g') \xRightarrow{\sim} \gerb{Int}(g)$.  This demonstrates that the gerbe $\gerb{E}_\epsilon(\alg{\tilde G})$ is well-defined in the following 2-categorical sense:  
\begin{enumerate}
\item
for each pair $i = (\alg{B}, \alg{T})$ consisting of a Borel subgroup of $\alg{G}$ and a maximally split maximal torus therein, we have constructed a gerbe $\gerb{E}_\epsilon^i(\alg{\tilde G})$;
\item
for each pair $i = (\alg{B}, \alg{T})$, $j = (\alg{B}_0, \alg{T}_0)$, we have constructed a family $P(i,j)$ of gerbe equivalences $\gerb{Int}(g) \From \gerb{E}_\epsilon^i(\alg{\tilde G}) \To \gerb{E}_\epsilon^j(\alg{\tilde G})$, indexed by those $g$ which conjugate $i$ to $j$;
\item
for any two elements $g, g'$ which conjugate $i$ to $j$, there is a distinguished natural isomorphism of gerbe equivalences from $\gerb{Int}(g')$ to $\gerb{Int}(g)$.
\end{enumerate}
Once we define the natural isomorphism $\gerb{Int}(t) \xRightarrow{\sim} \gerb{Id}$, we will have defined $\gerb{E}_\epsilon(\alg{\tilde G})$ uniquely up to equivalence, the equivalences being defined uniquely up to unique natural isomorphism (we learned this notion of well-definedness from reading various works of James Milne).  Defining the natural isomorphism requires some computation and is the subject of the section below.

\subsubsection{The isomorphism $\gerb{Int}(t) \xRightarrow{\sim} \gerb{Id}$}
For what follows, define $\delta_Q \From \sheaf{Y} \To n^{-1} \sheaf{X}$ to be the unique homomorphism satisfying
\begin{equation}
\label{DefinitionOfDelta}
\langle \delta_Q(y_1), y_2 \rangle = \beta_Q(y_1, y_2) \text{ for all } y_1, y_2 \in \sheaf{Y}.
\end{equation}
The constants $m_\phi$ and $n_\phi$ arise in the following useful result.
\begin{lemma}
\label{DeltaM}
For all $\phi \in \Phi$, we have $\delta_Q(n_\phi \phi^\vee) = \delta_Q(\tilde \phi^\vee) = m_\phi \phi$.
\end{lemma}
\proof
For all $\phi \in \Phi$, and all $y \in \sheaf{Y}$, we have
\begin{equation}
\label{DeltaBeta}
\langle \delta_Q(\tilde \phi^\vee), y \rangle = \beta_Q(n_\phi \phi^\vee, y).
\end{equation}
Weyl-invariance of the quadratic form (applying the root reflection for $\phi$) implies
$$\beta_Q(\phi^\vee, y) = \beta_Q(-\phi^\vee, y - \langle \phi, y \rangle \phi^\vee) = - \beta_Q(\phi^\vee, y) + \beta_Q(\phi^\vee, \phi^\vee) \langle \phi, y \rangle.$$
Adding $\beta_Q(\phi^\vee, y)$ to both sides and dividing by two, we find that
$$\beta_Q(\phi^\vee, y) = n^{-1} Q(\phi^\vee) \langle \phi, y \rangle.$$
Substituting into \eqref{DeltaBeta} yields
$$\langle \delta_Q(\tilde \phi^\vee), y \rangle = n_\phi  \beta_Q( \phi^\vee, y)= \langle n_\phi n^{-1} Q(\phi^\vee) \cdot \phi, y \rangle.$$
Since this holds for all $y \in \sheaf{Y}$, we have
$$\delta_Q(\tilde \phi^\vee) = n_\phi n^{-1} Q(\phi^\vee) \cdot \phi= m_\phi \cdot \phi.$$
\qed

Consider the homomorphisms of sheaves of abelian groups on $S_{\et}$,
$$\sheaf{T} \xrightarrow{\delta_Q} \sheaf{\hat T} \xtwoheadrightarrow{p} \sheaf{\hat T}_{\SC},$$
obtained by applying $\shom(\bullet, \sheaf{G}_m)$ to the homomorphisms
$$\sheaf{Y}_{Q,n}^{\SC} \Into \sheaf{Y}_{Q,n} \xrightarrow{\delta_Q} \sheaf{X}.$$
An object of $\gerb{E}_\epsilon(\alg{\tilde G})$ is a triple $(\sheaf{H}, h, j)$ where
\begin{itemize}
\item
$\sheaf{H} \xrightarrow{h} \sspl(\sheaf{D}_{Q,n})$ is an $n^{\th}$ root of the $\sheaf{\hat T}$-torsor $\sspl(\sheaf{D}_{Q,n})$.
\item
$j \From p_\ast \sheaf{H} \To \mu_\ast \Whit$ is an isomorphism in the gerbe $\gerb{E}_\epsilon^{\SC}(\alg{\tilde T}) = \epsilon_\ast \sqrt[n]{\sspl(\sheaf{D}_{Q,n}^{\SC})}$.  Thus $j = \tau^\vee \wedge j_0$ where $j_0 \From p_\ast \sheaf{H} \To \mu_\ast \Whit$ is an isomorphism in the gerbe $\sqrt[n]{\sspl(\sheaf{D}_{Q,n}^{\SC})}$ and $\tau^\vee \in \sheaf{\tilde T}_{\SC}^\Vee$.  See \ref{PushoutGerbe} for a general description of morphisms in gerbes obtained by pushing out.
\end{itemize}

Given $\sheaf{H} \xrightarrow{h} \sspl(\sheaf{D}_{Q,n})$ and $\hat t \in \sheaf{\hat T}$, write $\hat t \circ h \From \sheaf{H} \To \sspl(\sheaf{D}_{Q,n})$ for the map obtained by composing $h$ with the automorphism $\hat t$ of the $\sheaf{\hat T}$-torsor $\sspl(\sheaf{D}_{Q,n})$.

Similarly, given $j = \tau^\vee \wedge j_0$, with $p_\ast \sheaf{H} \xrightarrow{j_0} \mu_\ast \Whit$, $\tau^\vee \in \sheaf{\tilde T}_{\SC}^{\Vee}$, and given $\hat t_{\SC} \in \sheaf{\hat T}_{\SC}$, write $\hat t_{\SC} \circ j$ for $\tau^\vee \wedge (\hat t_{\SC} \circ j_0)$.  Here $\hat t_{\SC} \circ j_0 \From p_\ast \sheaf{H} \To \mu_\ast \Whit$ is the map obtained by composing $j$ with the automorphism $\hat t_{\SC}$ of the $\sheaf{\hat T}_{\SC}$-torsor $\mu_\ast \Whit$.  Similarly, write $j \circ \hat t_{\SC}$ for $\tau^\vee \wedge (j_0 \circ \hat t_{\SC})$.  Since $j_0$ is an isomorphism of $\sheaf{\hat T}_{\SC}$-torsors, we have $j \circ \hat t_{\SC} = \hat t_{\SC} \circ j$.

The following result describes the functor $\gerb{Int}(t)$ explicitly.
\begin{lemma}
For all $t \in \alg{T}(S)$, the equivalence of gerbes $\gerb{Int}(t) \From \gerb{E}_\epsilon(\alg{\tilde G}) \To \gerb{E}_\epsilon(\alg{\tilde G})$ sends an object $(\sheaf{H}, h,j)$ to the object $(\sheaf{H}, \delta_Q(t)^{-n} \circ h, p(\delta_Q(t))^{-1} \circ j)$.
\end{lemma}  
\proof
We work locally on $S_{\et}$ throughout the proof.  The action $\Int(t)$ of $t$ on the extension $\sheaf{G}_m \Into \sheaf{D}_{Q,n} \Onto \sheaf{Y}_{Q,n}$ is given by \cite[Equation 11.11.1]{B-D}, which states (in different notation) that
$$\Int(t) d = d \cdot \delta_Q(y)(t)^{-n} = d \cdot y \left( \delta_Q(t)^{-n} \right),$$
for all $d \in \sheaf{D}_{Q,n}$ lying over $y \in \sheaf{Y}_{Q,n}$.  We find that $\Int(t)$ is the automorphism of the extension $\sheaf{D}_{Q,n}$ determined by the element $\delta_Q(t)^{-n} \in \sheaf{\hat T} = \shom(\sheaf{Y}_{Q,n}, \sheaf{G}_m)$.  Hence the functor $\gerb{Int}(t) \From \gerb{E}_\epsilon(\alg{\tilde T}) \To \gerb{E}_\epsilon(\alg{\tilde T})$ sends $(\sheaf{H}, h)$ to $(\sheaf{H}, \delta_Q(t)^{-n} \circ h)$.  It remains to see how $\gerb{Int}(t)$ affects the third term in a triple $(\sheaf{H}, h, j)$.  

Conjugation by $t$ gives a map of torsors $\Whit \To \Whit$.  Let $j_0'$ be the map which makes the following triangle commute:
$$\begin{tikzcd}
p_\ast \sheaf{H} \arrow{r}{j_0} \arrow{dr}[swap]{j_0'} & \mu_\ast \Whit \arrow{d}{\Int(t)} \\
\phantom{a} & \mu_\ast \Whit
\end{tikzcd}$$
We have $\gerb{Int}(t) (\sheaf{H}, h, \tau^\vee \wedge j_0) = (\sheaf{H}, \delta_Q(t)^{-n} \circ h, \tau^\vee \wedge j_0')$ and must describe $j_0'$.

Take an element $\eta= (\eta_\alpha : \alpha \in \Delta) \in \Whit$ (decomposed with respect to the basis of simple roots $\Delta$).  Conjugation by $t$ yields a new element $\Int(t) \eta$ satisfying
$$\Int(t) \eta(u) = \eta(t^{-1} u t), \quad \text{ for all } u \in \sheaf{U}.$$
Decomposing along the simple root subgroups,
$$\Int(t) \eta = (\alpha(t) \ast \eta_\alpha : \alpha \in \Delta) \in \Whit,$$
where $[a \ast \eta_\alpha](e) = \eta_\alpha(a^{-1} e)$ for any $a \in \sheaf{G}_m$.

Comparing to the action of $\sheaf{\hat T}_{\SC}$ on $\Whit$, we find that
$$\Int(t) \eta= \theta(t) \ast \eta,$$
where $\theta \From \sheaf{T} \To \sheaf{\hat T}_{\SC}$ denotes the homomorphism dual to the homomorphism of character lattices sending $\tilde \alpha^\vee \in \sheaf{Y}_{Q,n}^{\SC}$ to $\alpha \in \sheaf{X}$.  

In the quotient $\sheaf{\hat T}_{\SC}$-torsor $\mu_\ast \Whit$, the action is given by $\Int(t) \bar \eta = \mu(\theta(t)) \cdot \bar \eta$, for all $\bar \eta \in \mu_\ast \Whit$.  More explicitly, $\mu \circ \theta \From \sheaf{\hat T}_{\SC} \To \sheaf{T}$ is the homomorphism dual to the map of character lattices sending $\tilde \alpha^\vee$ to $-m_{\alpha} \alpha$.  By Lemma \ref{DeltaM}, we have $\mu(\theta(t)) = p(\delta_Q(t))^{-1}$ (recall that $p$ corresponds to the inclusion $\sheaf{Y}_{Q,n}^{\SC} \Into \sheaf{Y}_{Q,n}$).  It follows that $j_0' = \Int(t) \circ j_0 = p(\delta_Q(t))^{-1} \circ j_0$.  This yields the result:
$$\gerb{Int}(t) (\sheaf{H}, h, j) = (\sheaf{H}, \delta_Q(t)^{-n} \circ h, p(\delta_Q(t))^{-1} \circ j).$$
\qed

\begin{proposition}
\label{IntTFunctors}
Let $(\sheaf{H}, h, j)$ be an object of the category $\gerb{E}_\epsilon(\alg{\tilde G})$.  Then, for all $t \in \alg{T}(S)$, the morphism $\delta_Q(t)^{-1} \From \sheaf{H} \To \sheaf{H}$ defines an isomorphism from $\gerb{Int}(t) (\sheaf{H}, h, j)$ to $(\sheaf{H}, h, j)$.  As objects vary, this defines a natural isomorphism $\Fun{A}(t) \From \gerb{Int}(t) \Rightarrow \gerb{Id}$ of functors from $\gerb{E}_\epsilon(\alg{\tilde G})$ to itself.  For a pair of elements $t_1, t_2 \in \alg{T}(S)$, we have a commutative diagram of functors and natural isomorphisms.
$$\begin{tikzcd}
\gerb{Int}(t_1 t_2) \arrow[Rightarrow]{r}{=} \arrow[Rightarrow]{d}{\Fun{A}(t_1 t_2)} & \gerb{Int}(t_1) \circ \gerb{Int}(t_2) \arrow[Rightarrow]{d}{\Fun{A}(t_1) \circ \Fun{A}(t_2)} \\
\gerb{Id} \arrow[Rightarrow]{r}{=} & \gerb{Id}
\end{tikzcd}$$
\end{proposition}
\proof
As $\sheaf{H} \xrightarrow{h} \sspl(\sheaf{D}_{Q,n})$ is an $n^{\th}$ root of the torsor $\sspl(\sheaf{D}_{Q,n})$, we have $\delta_Q(t)^{-n} \circ h = h \circ \delta_Q(t)^{-1}$.  As $j = \tau^\vee \wedge j_0$, with $j_0$ an isomorphism of $\sheaf{\hat T}_{\SC}$-torsors, we have $p(\delta_Q(t)^{-1}) \circ j = j \circ p(\delta_Q(t)^{-1})$.  

Hence the isomorphism $\delta_Q(t)^{-1} \From \sheaf{H} \To \sheaf{H}$ defines an isomorphism
\begin{align*}
\gerb{Int}(t) (\sheaf{H}, h, j) &= (\sheaf{H}, \delta_Q(t)^{-n} \circ h, p(\delta_Q(t)^{-1}) \circ j) \\
&= (\sheaf{H}, h \circ \delta_Q(t)^{-1}, j \circ p(\delta_Q(t)^{-1}) ) \\
& \xrightarrow{\delta_Q(t)^{-1}} (\sheaf{H}, h, j).
\end{align*}

As their definition depends only on $t$, these isomorphisms $\delta_Q(t)^{-1}$ form a natural isomorphism $\Fun{A}(t) \From \gerb{Int}(t) \Rightarrow \gerb{Id}$.  As $\delta_Q(t_1 t_2)^{-1} = \delta_Q(t_1)^{-1} \delta_Q(t_2)^{-1}$, we find an {\em equality} of functors $\gerb{Int}(t_1 t_2) = \gerb{Int}(t_1) \circ \gerb{Int}(t_2)$, and the commutative diagram of the proposition.  
\qed

\subsection{Change of base scheme}
\label{BaseChangeGerbe}

Let $\alg{\tilde G}$ be a degree $n$ cover of a quasisplit group $\alg{G} \supset \alg{B} \supset \alg{T}$ over $S$, as before.  Let $\gamma \From S_0 \To S$ be a morphism of schemes as in Section \ref{BaseChangeDualGroup}.  Then $\gamma$ gives rise to a pullback functor $\gamma^\ast$ from sheaves on $S_{\et}$ to sheaves on $S_{0,\et}$, and from gerbes on $S_{\et}$ to gerbes on $S_{0, \et}$.  

Assuming that $\amu_n(S)$ and $\amu_n(S_0)$ are cyclic groups of order $n$, we may identify these groups via $\gamma \From \amu_n(S) \To \amu_n(S_0)$, and a character $\epsilon \From \mu_n \Into \CC^\times$ corresponds to a character $\epsilon_0 \From \amu_n(S_0) \Into \CC^\times$.

Pullback via $\gamma$ defines a degree $n$ cover $\alg{\tilde G}_0$ of a quasisplit group $\alg{G}_0 \supset \alg{B}_0 \supset \alg{T}_0$ over $S_0$.  We have constructed gerbes $\gerb{E}_\epsilon(\alg{\tilde G})$ and $\gerb{E}_{\epsilon_0}(\alg{\tilde G}_0)$ associated to this data, banded by $\sheaf{\tilde Z}^\Vee$ and $\sheaf{\tilde Z}_0^\Vee$, respectively.  We also consider the pullback gerbe $\gamma^\ast \gerb{E}_\epsilon(\alg{\tilde G})$, banded by $\gamma^\ast \sheaf{\tilde Z}^\Vee$, and recall from Section \ref{BaseChangeDualGroup} that there is a natural isomorphism $N^\vee \From \gamma^\ast \sheaf{\tilde Z}^\Vee \To \sheaf{\tilde Z}_0^\Vee$.

Similarly, we find isomorphisms of sheaves on $S_{0,\et}$,
\begin{enumerate}
\item $\hat N \From \gamma^\ast \sheaf{\hat T} \xrightarrow{\sim} \sheaf{\hat T}_0$ (from $N \From \gamma^\ast \sheaf{Y}_{Q,n} \xrightarrow{\sim} \sheaf{Y}_{0,Q_0,n}$);
\item $N^\vee \From \gamma^\ast \sheaf{\tilde T}_{\SC}^\Vee \xrightarrow{\sim} \sheaf{\tilde T}_{0,\SC}^\Vee$;
\item $N \From \gamma^\ast \sheaf{D}_{Q,n} \xrightarrow{\sim} \sheaf{D}_{0,Q_0,n}$ (the construction of the second Brylinski-Deligne invariant is compatible with pullbacks); 
\item $N \From \gamma^\ast \Whit \xrightarrow{\sim} \Whit_0$ (since $\alg{B}_0 = \gamma^\ast \alg{B}$).
\end{enumerate}

From these observed isomorphisms, we find that the pullback of \'etale sheaves from $S_{\et}$ to $S_{0,\et}$ defines a functor 
$$\gerb{N}' \From \gerb{E}_\epsilon(\alg{\tilde G})  \To \gamma_\ast \gerb{E}_{\epsilon_0}(\alg{\tilde G}_0)),$$
given on objects by $\gerb{N}'(\sheaf{H}, h, j) = (\gamma^\ast \sheaf{H}, \gamma^\ast h, \gamma^\ast j)$.  For example, if $\sheaf{H}$ is a $\sheaf{\hat T}$-torsor, then $\gamma^\ast \sheaf{H}$, a priori a $\gamma^\ast \sheaf{\hat T}$-torsor, becomes a $\sheaf{\hat T}_0$-torsor via $\hat N$.  

From \cite[Chapitre V, Proposition 3.1.8]{Giraud}, such an equivalence of gerbes $\gerb{N}'$ determines a unique, up to unique natural isomorphism, equivalence of gerbes,
$$\gerb{N} \From \gamma^\ast \gerb{E}_\epsilon(\alg{\tilde G}) \To \gerb{E}_{\epsilon_0}(\alg{\tilde G}_0),$$
lying over the natural isomorphism of bands $N^\vee \From \gamma^\ast \sheaf{\tilde Z}^\Vee \To \sheaf{\tilde Z}_0^\Vee$.  In this way, the construction of the gerbe associated to a cover is compatible with change of base scheme.

 \subsection{Parabolic subgroups}
\label{ParabolicGerbe}

Return to the degree $n$ cover $\alg{\tilde G}$ of a quasisplit group $\alg{G} \supset \alg{B} \supset \alg{T}$ over $S$.  Let $\alg{P} \subset \alg{G}$ be a parabolic subgroup defined over $S$, containing $\alg{B}$.  As in Section \ref{ParabolicDualGroup}, we consider a Levi decomposition $\alg{P} = \alg{M} \alg{N}$ and resulting cover $\alg{\tilde M}$.  Fix $\epsilon \From \mu_n \Into \CC^\times$ as before.

Consider the gerbes $\gerb{E}_\epsilon(\alg{\tilde G})$ and $\gerb{E}_\epsilon(\alg{\tilde M})$.  An object of $\gerb{E}_\epsilon(\alg{\tilde G})$ is a triple $(\sheaf{H}, h, j)$, where $(\sheaf{H}, h)$ is an $n^{\th}$ root of the torsor $\sspl(\sheaf{D}_{Q,n})$, $j = \tau^\vee \wedge j_0$, $\tau^\vee \in \sheaf{\tilde T}_{\SC}^\Vee$, and $j_0 \From p_\ast\sheaf{H} \To \mu_\ast \Whit$ is an isomorphism of $\sheaf{\hat T}_{\SC}$-torsors.

Restriction of characters via $\alg{U} \supset \alg{U}_{\alg{M}} = \alg{U} \cap \alg{M}$ provides a homomorphism of sheaves from $\Whit$ to $\Whit_{\alg{M}}$ (the Whittaker torsor for $\alg{M}$).  The inclusion $\sheaf{Y}_{\alg{M}}^{\SC} \subset \sheaf{Y}^{\SC}$ of coroot lattices provides a homomorphism of sheaves of abelian groups $\sheaf{\hat T}_{\SC} \Onto \sheaf{\hat T}_{\alg{M},\SC}$, where 
$$\sheaf{\hat T}_{\SC} = \shom(\sheaf{Y}^{\SC}, \sheaf{G}_m), \quad \sheaf{\hat T}_{\alg{M},\SC} = \shom(\sheaf{Y}_{\alg{M}}^{\SC}, \sheaf{G}_m).$$
Similarly, it provides a homomorphism $\sheaf{\tilde T}_{\SC}^\Vee \Onto \sheaf{\tilde T}_{\alg{M}, \SC}^\Vee$, where 
$$\sheaf{\tilde T}_{\SC}^\Vee = \shom(\sheaf{Y}^{\SC}, \CC^\times), \quad \sheaf{\tilde T}_{\alg{M},\SC}^\Vee = \shom(\sheaf{Y}_{\alg{M}}^{\SC}, \CC^\times).$$
Define $\tau_{\alg{M}}^\vee$ to be the image of $\tau^\vee$ under this homomorphism.

The constants defining $\mu \From \sheaf{\hat T}_{\SC} \To \sheaf{\hat T}_{\SC}$ are the same as those defining the corresponding map, $\mu_{\alg{M}} \From \sheaf{\hat T}_{\alg{M},\SC} \To \sheaf{\hat T}_{\alg{M},\SC}$.  We find a commutative diagram.
$$\begin{tikzcd}
 \sheaf{\hat T}_{\SC} \times \mu_\ast \Whit  \arrow{r}{\ast} \arrow{d} & \mu_\ast \Whit \arrow{d} \\
  \sheaf{\hat T}_{\alg{M},\SC} \times (\mu_{\alg{M}})_\ast  \Whit_{\alg{M}} \arrow{r}{\ast} &  (\mu_{\alg{M}})_\ast \Whit_{\alg{M}} 
\end{tikzcd}$$

Composing $j_0 \From p_\ast \sheaf{H} \To \mu_\ast \Whit$ with the map to $(\mu_{\alg{M}})_\ast \Whit_{\alg{M}}$ defines a map
$$j_{\alg{M},0} \From \sheaf{H} \To (\mu_{\alg{M}})_\ast \Whit_{\alg{M}}.$$
This defines a functor of gerbes,
$$\gerb{i} \From \gerb{E}_\epsilon(\alg{\tilde G}) \To \gerb{E}_\epsilon(\alg{\tilde M}), \quad (\sheaf{H}, h, \tau^\vee \wedge j_0) \mapsto (\sheaf{H}, h, \tau_{\alg{M}}^\vee \wedge j_{\alg{M},0}),$$
lying over $\iota^\vee \From \sheaf{\tilde Z}^\Vee \Into \sheaf{\tilde Z}_{\alg{M}}^\Vee$.

\subsection{Weyl action on the gerbe associated to the cover of the torus}
\label{WeylGerbe}
We keep the degree $n$ cover $\alg{\tilde G}$ of a quasisplit group $\alg{G} \supset \alg{B} \supset \alg{T}$ over $S$.  Let $\dot w \in \sheaf{N}[S]$ represent an element of the Weyl group $w \in \sheaf{W}[S]$.  As in Section \ref{WeylDualTorus}, $\Int(\dot w)$ defines a well-aligned homomorphism from $\alg{\tilde T}$ to itself.  
$$\begin{tikzcd}
\alg{K}_2 \inarrow{r} \arrow{d}{=} & \alg{T}' \onarrow{r} \arrow{d}{\Int(\dot w)} & \alg{T} \arrow{d}{\Int(w)} \\
\alg{K}_2 \inarrow{r} & \alg{T}' \onarrow{r} & \alg{T}
\end{tikzcd}$$ 
This defines an equivalence of gerbes on $S_{\et}$,
$$\gerb{Int}(\dot w) \From \gerb{E}_\epsilon(\alg{\tilde T}) \To \gerb{E}_\epsilon(\alg{\tilde T}),$$
lying over the homomorphism $\Int(w)^\vee \From \sheaf{\tilde T}^\Vee \To \sheaf{\tilde T}^\Vee$.

Given two such representatives $\dot w, \ddot w \in \sheaf{N}[S]$ for $w$, there exists a unique $t \in \alg{T}(S)$ such that $\ddot w = \dot w \cdot t$.  The natural isomorphism of functors $\Fun{A(t)} \From \gerb{Int}(t) \xRightarrow{\sim} \gerb{Id}$, defined in Proposition \ref{IntTFunctors}, yields a natural isomorphism
$$\gerb{Int}(\ddot w) \xRightarrow{\sim} \gerb{Int}(\dot w).$$
In this way, there is an equivalence of gerbes $\gerb{Int}(w) \From \gerb{E}_\epsilon(\alg{\tilde T}) \To \gerb{E}_\epsilon(\alg{\tilde T})$, defined uniquely up to unique natural isomorphism.  

For what follows later, it will be useful to describe the functor $\gerb{Int}(\dot w)$ explicitly.  So consider an object $\sheaf{H} \xrightarrow{h} \sspl(\sheaf{D}_{Q,n})$ of $\gerb{E}_\epsilon(\alg{\tilde T})$.  Thus $\sheaf{H}$ is a $\sheaf{\hat T}$-torsor on $S_{\et}$, and $h$ is a map of sheaves on $S_{\et}$ satisfying
$$h(\hat \tau \ast a) = \hat \tau^n \ast h(a), \text{ for all } a \in \sheaf{H}, \hat \tau \in \sheaf{\hat T}.$$ 

The well-aligned homomorphism $\Int(\dot w) \From \alg{\tilde T} \To \alg{\tilde T}$ yields a commutative diagram when we take the second Brylinski-Deligne invariant.
$$\begin{tikzcd}
\sheaf{G}_m \inarrow{r} \arrow{d}{=}  & \sheaf{D}_{Q,n} \onarrow{r} \arrow{d}{\Int(\dot w)} & \sheaf{Y}_{Q,n} \arrow{d}{\Int(w)} \\
\sheaf{G}_m \inarrow{r} & \sheaf{D}_{Q,n} \onarrow{r} & \sheaf{Y}_{Q,n}
\end{tikzcd}$$
If $s \in \sspl(\sheaf{D}_{Q,n})$, define $\dot w(s) \in \sspl(\sheaf{D}_{Q,n})$ by
$$\dot w(s) = \Int(\dot w) \circ s \circ \Int(w)^{-1}.$$

Allowing $s$ to vary, this gives a map of sheaves on $S_{\et}$,
$$\dot w \From \sspl(\sheaf{D}_{Q,n}) \To \sspl(\sheaf{D}_{Q,n})$$
which satisfies
$$\dot w( \hat \tau \ast s) = w(\hat \tau) \ast \dot w(s).$$


The functor $\gerb{Int}(\dot w) \From \gerb{E}_\epsilon(\alg{\tilde T}) \To \gerb{E}_\epsilon(\alg{\tilde T})$ sends $(\sheaf{H}, h)$ to $\left( {}^w \sheaf{H}, \dot w \circ h \right)$, where ${}^w \sheaf{H}$ is the $\sheaf{\hat T}$-torsor on $S_{\et}$ which coincides with $\sheaf{H}$ as sheaves of sets on $S_{\et}$, and in which the torsor structure is twisted:
$$\hat \tau \ast_w x \defeq w^{-1} \hat \tau \ast x.$$


An explicit description of $\dot w \From \sspl(\sheaf{D}_{Q,n}) \To \sspl(\sheaf{D}_{Q,n})$ seems unwieldy, in general.  But in a special case, we may describe $\dot w(s)$ explicitly.  For this special case, define $\alg{T}_d$ to be the maximal $S$-split torus in $\alg{T}$; its cocharacter lattice $Y_d$ coincides with $\sheaf{Y}[S]$.  Write $X_d$ for the character lattice of this split torus, and $\Phi_d$ for the subset of $X_d$ consisting of the roots of $\alg{G}$ relative to $\alg{T}_d$.  Write $W_d$ for the Weyl group of the relative root system; it can be identified with $\sheaf{W}[S]$.  If $\beta \in \Phi_d$ is a relative root, then define $w_\beta$ to be the corresponding reflection in $W_d$.  We refer to the seminal paper of Borel and Tits for the results stated here about relative root systems; for example, the following is a result of \cite[Th\'eor\`eme 5.3]{BorelTits}:
$$W_d = \langle w_\beta : \beta \in \Phi_d \text{ and } 2 \beta \not \in \Phi_d \rangle.$$

Suppose that $\beta \in \Phi_d$ and $2 \beta \not \in \Phi_d$ hereafter, i.e., $\beta$ is non-multipliable.  Let $\{ \alpha_1, \ldots, \alpha_\ell \}$ be the absolute roots (in $\Phi[U]$ for some Galois cover $U \To S$) whose restrictions to $\alg{T}_d$ coincide with $\beta$.  These roots are pairwise orthogonal, and form a single Galois orbit.  Let $S' / S$ be the \'etale cover corresponding to the finite Galois module $\{ \alpha_1, \ldots, \alpha_\ell \}$.  As the roots $\alpha_1, \ldots, \alpha_\ell$ lie in the same Galois orbit, the associated constants $n_{\alpha_1}, \ldots, n_{\alpha_\ell}$ are equal, and we write $n_\beta$ for their common value.  Write $q_\beta$ for the common value of $Q(\alpha_1^\vee), \ldots, Q(\alpha_\ell^\vee)$.

Define $\beta^\vee = \sum_{i=1}^\ell \alpha_i^\vee$ for the resulting cocharacter of $\alg{T}_d$, and define $\tilde \beta^\vee = n_\beta \beta^\vee$.  The structure theory of quasisplit groups provides a homomorphism defined over $S$, with finite kernel, $\alg{Res}_{S'/S} \alg{SL}_2 \To \alg{G}$.  Composing with the natural embedding $\alg{SL}_2 \Into \alg{Res}_{S'/S} \alg{SL}_2$, we find a homomorphism with finite kernel,
$$r_\beta \From \alg{SL}_2 \To \alg{G}.$$
The standard coroot $t \mapsto diag(t,t^{-1})$ of $\alg{SL}_2$ maps to the cocharacter $\beta^\vee \in Y_d$.

The cover $\alg{\tilde G} = (\alg{G}', n)$ pulls back via $r_\beta$ to a cover $(\alg{SL}_2', n)$ of $\alg{SL}_2$, defined uniquely up to unique isomorphism by the integer $Q(\beta^\vee) = \ell q_\beta$.  The second Brylinski-Deligne invariant for $\alg{SL}_2'$, an extension $\sheaf{D}_{\ell q_\beta}$, fits into a commutative diagram.
$$\begin{tikzcd}
\sheaf{G}_m \inarrow{r} \arrow{d}{=} & \sheaf{D}_{\ell q_\beta} \onarrow{r} \inarrow{d} & \ZZ \inarrow{d}{a \mapsto a \cdot \beta^\vee} \\
\sheaf{G}_m \inarrow{r} & \sheaf{D} \onarrow{r} & \sheaf{Y}
\end{tikzcd}$$

Let $\alg{U}_\ast$ be the locally closed subscheme of $\alg{SL}_2$ obtained from the closed subgroup of upper-triangular unipotent matrices, by deleting the identity section.  Let $\alg{U}_\ast^-$ be the analogous scheme of lower-triangular unipotent matrices.  If $e \in \alg{U}_\ast(S)$, then there exists $e^- \in \alg{U}_\ast^-(S)$ such that $\dot w_\beta \defeq r_\beta(e e^- e)$ represents the relative root reflection $w_\beta$ (cf. \cite[\S 1.8]{Deodhar}).  This provides an element $[e] \in \sheaf{D}$ lying over $\beta^\vee$ by the construction of \cite[\S 11.2]{B-D} (the construction for $\alg{SL}_2$ is in \cite{B-D}, and apply the inclusion $\sheaf{D}_{\ell q_\beta} \Into \sheaf{D}$).  Raising to the $n_\beta$ power yields an element $[e]^{n_\beta} \in \sheaf{D}_{Q,n}$ lying over $\tilde \beta^\vee \in \sheaf{Y}_{Q,n}$.

The following lemma is a slight adaptation of \cite[Equation 11.11.2]{B-D}.
\begin{lemma}
\label{WDLemma}
Suppose that $\dot w = r_\beta(e e^- e)$ and $d \in \sheaf{D}_{Q,n}$ lies over $y \in \sheaf{Y}_{Q,n}$.  Then
$$\Int(\dot w) d = d \cdot \prod_{i=1}^\ell \left( [e_i]^{-n_\beta \langle \tilde \alpha_i, y \rangle}  \cdot (-1)^{q_\beta \varepsilon( - n_\beta \langle \tilde \alpha_i, y \rangle )  } \right).$$
Here $\varepsilon(N) = N(N+1) / 2$.  The elements $[e_i] \in \sheaf{D}$ are described in the proof. 
\end{lemma}
\proof
Extending scalars to the \'etale cover $S' / S$, we have an isomorphism
$$\left( \alg{Res}_{S' / S} \alg{SL}_{2} \right)_{S'} \isom \prod_{i=1}^\ell \alg{SL}_{2,S'}.$$
The natural embedding $\alg{SL}_2 \Into \alg{Res}_{S' / S} \alg{SL}_2$ is identified with the diagonal embedding, after extending scalars to $S'$.  

In this way, we may write $e = \prod_{i=1}^\ell e_i$ and $e^- = \prod_{i=1}^\ell e_i^-$, where $e_i$ and $e_i^-$ commute with $e_j$ and $e_j^-$ for $i \neq j$.  In this way, $w_\beta = \prod w_{\alpha_i}$ and the latter product consists of commuting reflections in $\sheaf{W}[S']$.  Write $\dot w_i = e_i e_i^- e_i$, and note that $\dot w_\beta = \prod \dot w_i$ (a product of commuting elements of $\alg{G}(S')$).  Tracing through the construction of \cite[\S 11.1]{B-D}, we find that $[e] = \prod_{i=1}^\ell [e_i]$, where the latter product consists of commuting elements of $\sheaf{D}$, with $[e_i]$ lying over $\alpha_i^\vee$.  Define $\tilde \beta = \sum_{i=1}^\ell \tilde \alpha_i$ in what follows.

\cite[Equation 11.11.2]{B-D} gives a formula for each $1 \leq i \leq \ell$,
$$\Int(\dot w_i) d = d \cdot [e_i]^{- \langle \alpha_i, y \rangle } \cdot (-1)^{q_\beta \varepsilon( - \langle \alpha_i, y \rangle )}.$$
Since $\alpha_i = n_\beta \tilde \alpha_i$ for all $i$, we have
\begin{equation}
\label{ConjWd}
\Int(\dot w_i) d = d \cdot [e_i]^{- n_\beta \langle \tilde \alpha_i, y \rangle } \cdot (-1)^{q_\beta \varepsilon( - n_\beta  \langle \tilde \alpha_i, y \rangle ) )}.
\end{equation}
Since $\langle \alpha_i, \alpha_j^\vee \rangle = 0$ for $i \neq j$, repeated application of this equation yields
$$\Int(\dot w) d = d \cdot \prod_{i=1}^\ell \left( [e_i]^{-n_\beta \langle \tilde \alpha_i, y \rangle}  \cdot (-1)^{q_\beta \varepsilon( - n_\beta \langle \tilde \alpha_i, y \rangle )  } \right).$$
\qed

Now, if $s \in \sspl(\sheaf{D}_{Q,n})$, and $\dot w_\beta = r_\beta(e e^- e)$ as above, then we say that $s$ is \defined{aligned with $\dot w_\beta$} if $s(\tilde \alpha_i^\vee) = [e_i]^{n_\beta}$ for all $1 \leq i \leq \ell$.  
  
\begin{thm}
\label{WeylActsOnSplittings}
For $\dot w_\beta = r_\beta(e e^- e)$ as before, and $s \in \sspl(\sheaf{D}_{Q,n})$ aligned with $\dot w_\beta$, we have
$$\dot w_\beta(s) = \begin{cases} \prod_{i=1}^\ell \tilde \alpha_i (-1)^{n/2} \ast s & \text{ if } q_\beta \text{ is odd}; \\ s & \text{ if } q_\beta \text{ is even}. \end{cases}$$
(By Assumption \ref{OddnEvenQ}, $q_\beta$ odd implies $n$ is even.  Each $\tilde \alpha_i \in \sheaf{X}_{Q,n}$ is viewed as a cocharacter of $\sheaf{\hat T}$ here.)
\end{thm}
\proof
From the formula \eqref{ConjWd}, we compute
\begin{align*}
[\dot w_\beta(s)](y) &= \Int(\dot w_\beta) s(\Int(w_\beta)^{-1} y) \\
&= \Int(\dot w_\beta) s \left( y - \sum_j \langle \tilde \alpha_j, y \rangle \tilde \alpha_j^\vee \right) \\
&= \Int(\dot w_\beta) s(y) \cdot \prod_j \Int(\dot w_\beta) s(\tilde \alpha_j^\vee)^{- \langle \tilde \alpha_j, y \rangle } \\
&= s(y) \cdot \prod_{i=1}^\ell [e_i]^{-n_\beta \langle \tilde \alpha_i, y \rangle}  \cdot (-1)^{q_\beta \varepsilon( - n_\beta \langle \tilde \alpha_i, y \rangle ) } \\
& \cdot \prod_{j=1}^\ell \left( s(\tilde \alpha_j^\vee)^{- \langle \tilde \alpha_j, y \rangle} \cdot \prod_{i=1}^\ell \left( [e_i]^{n_\beta \langle \tilde \alpha_j, y \rangle \langle \tilde \alpha_i, \tilde \alpha_j^\vee \rangle}  \cdot (-1)^{- q_\beta \langle \tilde \alpha_j, y \rangle \varepsilon( - n_\beta \langle \tilde \alpha_i, \tilde \alpha_j^\vee \rangle ) } \right) \right) \\
&= s(y) \cdot \prod_{i=1}^\ell (-1)^{q_\beta \epsilon(-n_\beta \langle \tilde \alpha_i, y \rangle) - q_\beta \langle \tilde \alpha_i, y \rangle \varepsilon(-2 n_\beta) }.
\end{align*}
In the last line, we use the fact that $s(\tilde \alpha_i^\vee) = [e_i]^{n_\beta}$ for all $1 \leq i \leq \ell$, and we use the orthogonality relation $\langle \tilde \alpha_j, \tilde \alpha_i^\vee \rangle = 0$ for $i \neq j$, and $\langle \tilde \alpha_j, \tilde \alpha_j^\vee \rangle = 2$.

For all $N \in \ZZ$, $\varepsilon(2N) = N$ modulo $2$, and $\varepsilon(-N) - N = \varepsilon(N)$ modulo $2$.  Hence
$$\varepsilon(- n_\beta \langle \tilde \alpha_i, y \rangle) - \langle \tilde \alpha_i, y \rangle \varepsilon(-2 n_\beta) = \varepsilon(- n_\beta \langle \tilde \alpha_i, y \rangle) - n_\beta \langle \tilde \alpha_i, y \rangle = \varepsilon (n_\beta \langle \tilde \alpha_i, y \rangle), \text{ modulo } 2.$$
We find that
$$[\dot w_\beta(s)](y) = s(y) \cdot \prod_{i=1}^\ell (-1)^{q_\beta \cdot \varepsilon (n_\beta \langle \tilde \alpha_i, y \rangle)}.$$
If $q_\beta$ is even, then the exponent of $(-1)$ is even and $\dot w_\beta(s) = s$ as claimed.  

If $q_\beta$ is odd, then $n$ is even and $n_\beta$ is even.  Moreover, $n$ is a multiple of $4$ if and only if $n_\beta$ is a multiple of $4$.  We find that
$$\prod_{i=1}^\ell (-1)^{q_\beta \epsilon(n_\beta \langle \tilde \alpha_i, y \rangle )} = \prod_{i=1}^\ell (-1)^{n_\beta \langle \tilde \alpha_i, y \rangle / 2} = \prod_{i=1}^\ell (-1)^{n \langle \tilde \alpha_i, y \rangle / 2}.$$
This is precisely the evaluation of $y$ at $\prod_i \tilde \alpha_i(-1)^{n/2} \in \sheaf{\hat T}_{[2]} = \shom(\sheaf{Y}_{Q,n}, \mu_2)$.
\qed

\section{The metaGalois group}
\label{MetaGaloisSection}
We now specialize to three classes of base scheme of arithmetic interest.
\begin{description}
\item[Global]  $S = \Spec(F)$ for a global field $F$;
\item[Local]  $S = \Spec(F)$ for a local field $F$;
\item[Local integral]  $S = \Spec(\OO)$ for the ring of integers $\OO$ in a nonarchimedean local field $F$.
\end{description}

Choose a separable closure $\bar F / F$ in all three cases, and write $\Gal_F = \Gal(\bar F / F)$ for the resulting absolute Galois group.  In the local integral case, the separable closure $\bar F / F$ provides a geometric base point $\bar s$ for $\Spec(\OO)$ as well, and define $\Gal_\OO = \pi_1^{\et}(\Spec(\OO), \bar s)$.  This is a profinite group, topologically generated by a geometric Frobenius element $\Fr$; thus we write $\Gal_\OO = \langle \Fr \rangle_{\prof}$.

When $S = \Spec(F)$ for a local or global field, or $S = \Spec(\OO)$ for the ring of integers in a nonarchimedean local field, write $\Gal_S$ for $\Gal_F$ or $\Gal_\OO$ accordingly.

\subsection{Construction of the metaGalois group}
The metaGalois group will be a profinite group fitting into a central extension,
$$\mu_2 \Into \mGal_S \Onto \Gal_S.$$
When $F$ has characteristic two, we {\em define} the metaGalois group $\mGal_F$ to be the trivial extension $\Gal_F \times \mu_2$.  The metaGalois group $\mGal_\OO$ will {\em not be defined} when $\OO$ has residual characteristic $2$, reflecting the idea that metaGalois representations cannot be ``unramified at $2$'' (though one might propose an alternative notion of ``minimally ramified'').

\subsubsection{Local fields}

When $F$ is a local field (with $2 \neq 0$) the quadratic Hilbert symbol defines a symmetric nondegenerate $\ZZ$-bilinear form
$$\Hilb_2 \From F^\times_{/2} \times F^\times_{/2} \To \mu_2.$$
The abelianized Galois group $\Gal_F^{\ab}$ is defined, up to unique isomorphism, from $F$ alone.  When $F$ is nonarchimedean, we normalize the valuation so that $\val(F^\times) = \ZZ$, and we normalize the reciprocity map of local class field theory to send a geometric Frobenius element to an element of valuation $1$.  The reciprocity map gives a surjective homomorphism
$$\rec_{F/2} \From \Gal_F^{\ab} \Onto F^\times_{/2}.$$
Composing the Hilbert symbol with $\rec_{F/2}$ defines a function
$$h \From \Gal_F^{\ab} \times \Gal_F^{\ab} \To \mu_2,$$
and it is straightforward to verify that $h$ is a (bimultiplicative) continuous symmetric $2$-cocycle.  This incarnates a commutative extension
$$\mu_2 \Into \mGal_F^{\ab} \Onto \Gal_F^{\ab}$$
of profinite groups.  Concretely, $\mGal_F^{\ab} = \Gal_F^{\ab} \times \mu_2$ as sets, and 
$$(\gamma_1, \epsilon_1) \cdot (\gamma_2, \epsilon_2) \defeq \left( \gamma_1 \gamma_2, \epsilon_1 \epsilon_2 \cdot h(\gamma_1, \gamma_2) \right).$$

The pullback of this extension to $\Gal_F$ will be called the \defined{metaGalois group of $F$}, written $\mGal_F$.  It is a central extension of $\Gal_F$ by $\mu_2$,
$$\mu_2 \Into \mGal_F \Onto \Gal_F.$$

\subsubsection{The local integral case}

Suppose that $F$ is a nonarchimedean \textbf{nondyadic} (i.e., $\val(2) = 0$) local field.  Then the quadratic Hilbert symbol satisfies
$$\Hilb_2(u,v) = 1 \text{ for all } u,v \in \OO^\times.$$
Write $\Inertia \subset \Gal_F$ for the inertial subgroup, so that $\rec_F(\Inertia) = \OO^\times$.  The cocycle $h$ is trivial when restricted to $\Inertia \times \Inertia$.  Thus $\gamma \mapsto (\gamma, 1)$ gives a canonical splitting $\sigma^\circ \From \Inertia \Into \mGal_F$.  The natural map $\Gal_F \Onto \Gal_\OO$ identifies $\Gal_\OO$ with $\Gal_F / \Inertia$.  Define $\mGal_\OO = \mGal_F / \sigma^\circ(\Inertia)$ to obtain a commutative diagram with exact rows.
$$\begin{tikzcd}
\mu_2 \arrow{d}{=} \inarrow{r} & \mGal_F \arrow{r} \onarrow{d} & \Gal_F \onarrow{d} \\
\mu_2 \inarrow{r} & \mGal_\OO \arrow{r} & \Gal_\OO
\end{tikzcd}$$
We call $\mGal_\OO$ the \defined{metaGalois group of $\OO$}.  If $\gamma \in \Gal_F$ lifts $\Frob$, then every element of $\mGal_\OO$ is equal to $(\gamma^{\hat n}, \pm 1)$ (mod $\sigma^\circ(\Inertia)$) for some $\hat n \in \hat \ZZ$.  In this way, a Frobenius lift provides an isomorphism from $\mGal_\OO$ to the group with underlying set $\langle \gamma \rangle_{\prof} \times \mu_2$ and multiplication given by
$$(\gamma^{\hat n_1}, \epsilon_1) \cdot (\gamma^{\hat n_2}, \epsilon_2) = \left( \gamma^{\hat n_1 + \hat n_2}, \epsilon_1 \epsilon_2 \cdot (-1)^{ \hat n_1 \hat n_2 (q-1)/2} \right),$$
where $q$ is the cardinality of the residue field of $\OO$.

\subsubsection{Global fields}

When $F$ is a global field (with $2 \neq 0$ as before), the Hilbert symbol defines a symmetric $\ZZ$-bilinear form,
$$\Hilb_2 \From \AA^\times_{/2} \times \AA^\times_{/2} \To \mu_2,$$
obtained as the product of local Hilbert symbols.  This defines a continuous symmetric 2-cocycle, from which we get a commutative extension,
\begin{equation}
\label{ExtensionAdelic}
\mu_2 \Into \widetilde{\AA^\times_{/2}} \Onto \AA^\times_{/2}.
\end{equation}
Global quadratic reciprocity for the Hilbert symbol ($\Hilb_2(u,v) = 1$ for all $u,v \in F^\times$) provides a canonical splitting $\sigma_F \From F^\times_{/2} \Into \widetilde{\AA^\times_{/2}}$.  Taking the quotient yields a commutative extension,
\begin{equation}
\label{QuotientExtension}
\mu_2 \Into \frac{\widetilde{\AA^\times_{/2}}}{\sigma_F(F^\times_{/2})} \Onto \frac{ \AA^\times_{/2} }{F^\times_{/2}}.
\end{equation}
The global reciprocity map of class field theory gives an surjective homomorphism,
$$\rec_{F/2} \From \Gal_F^{\ab} \Onto \AA^\times_{/2} / F^\times_{/2} \ident (\AA^\times / F^\times)_{/2}.$$
Pulling back \eqref{QuotientExtension} via $\rec_{F/2}$ yields a commutative extension,
$$\mu_2 \Into \mGal_F^{\ab} \Onto \Gal_F^{\ab}.$$
Pulling back via $\Gal_F \Onto \Gal_F^{\ab}$ defines the \defined{metaGalois group of $F$},
$$\mu_2 \Into \mGal_F \Onto \Gal_F.$$

\subsubsection{Compatibilities}
\label{CompatibilityMetaGalois}
If $v \in \VV$ is a place of a global field $F$, then an embedding $\bar F \Into \bar F_v$ of separable closures determines an injective homomorphism $\iota_v \From \Gal_{F_v} \Into \Gal_F$.  As the global Hilbert symbol is the product of local ones, we find a homomorphism $\tilde \iota_v$ realizing $\widetilde{\Gal}_{F_v}$ as the pullback of the extension $\widetilde{\Gal}_F$.
$$\begin{tikzcd}
\mu_2 \inarrow{r} \arrow{d}{=} & \mGal_{F_v} \arrow{r} \inarrow{d}{\tilde \iota_v} & \Gal_{F_v} \inarrow{d}{\iota_v} \\
\mu_2 \inarrow{r} & \mGal_F \arrow{r} & \Gal_F
\end{tikzcd}$$

For local and global fields $F$, a choice of separable closure $\bar F / F$ entered the construction of the metaGalois group.  Suppose that $\bar F_0$ is another separable closure of $F$.  Every $F$-algebra isomorphism $\iota \From \bar F_0 \xrightarrow{\sim} \bar F$ yields an isomorphism $\iota \From \Gal(\bar F_0 / F) \xrightarrow{\sim} \Gal(\bar F / F)$.  The resulting isomorphism $\Gal(\bar F_0 / F)^{\ab} \xrightarrow{\sim} \Gal(\bar F / F)^{\ab}$ does not depend on $\iota$.  

The separable closure $\bar F_0$ yields a cocycle $h_0 \From \Gal(\bar F_0 / F)^{\ab} \times \Gal(\bar F_0 / F)^{\ab} \To \mu_2$, and thus a metaGalois group $\mGal(\bar F_0 / F)$.  Since the defining cocycles $h_0$ and $h$ factor through abelianized Galois groups, the isomorphism $\iota \From \Gal(\bar F_0 / F) \To \Gal(\bar F / F)$ lifts canonically to an isomorphism of metaGalois groups.
$$\begin{tikzcd}
\mu_2 \arrow{d}{=} \inarrow{r} & \mGal(\bar F_0 / F)\onarrow{r} \arrow{d}{\tilde \iota} & \Gal(\bar F_0 / F)\arrow{d}{\iota} \\
\mu_2 \inarrow{r} & \mGal(\bar F / F) \onarrow{r} & \Gal(\bar F / F)
\end{tikzcd}$$

\subsection{The Brauer class}
When $F$ is a local or global field, the metaGalois group is an extension $\mu_2 \Into \mGal_F \Onto \Gal_F$.  As such, it has a cohomology class in the Brauer group.
$$\left[ \mGal_F \right] \in H_{\et}^2(F, \mu_2) = \Br(F)_{[2]}.$$

This Brauer class is often trivial -- the metaGalois group often splits, though it rarely has a canonical splitting.
\begin{proposition}
\label{OddRCTrivialBr}
Suppose that $F$ is a \textbf{nondyadic} (i.e., $\val(2) = 0$) nonarchimedean local field.  Then $\left[ \mGal_F \right]$ is the trivial class.
\end{proposition}
\proof
The projection $\mGal_F \To \mGal_\OO$ identifies the metaGalois group of $F$ with the pullback of the metaGalois group of $\OO$.
$$\begin{tikzcd}
\mu_2 \arrow{r} \arrow{d}{=} & \mGal_F \arrow{r} \arrow{d} & \Gal_F \arrow{d} \\
\mu_2 \arrow{r} & \mGal_\OO \arrow{r} & \Gal_\OO = \langle \Fr \rangle_{\prof}
\end{tikzcd}$$
But every extension of $\hat \ZZ$ by $\mu_2$ splits (though not canonically); hence the metaGalois group splits and its Brauer class is trivial.
\qed

\begin{proposition}
\label{RealNontrivBr}
Over $\RR$, the metaGalois group is a nonsplit extension, so $\left[ \mGal_\RR \right]$ is the unique nontrivial class in the Brauer group $\Br(\RR)$.
\end{proposition}
\proof
Let $\sigma$ denote complex conjugation, $\Gal_\RR = \Gal(\CC / \RR) = \{ \Id, \sigma \}$.  The metaGalois group is a cyclic group of order $4$ sitting in an extension
$$\mu_2 \Into \mGal_\RR \Onto \Gal_\RR.$$
Indeed, the cocycle $h$ satisfies $h(\sigma, \sigma) = \Hilb_2(-1,-1) = -1$.  Thus
$$(\sigma,1) \cdot (\sigma,1) = (\Id, -1) \in \mGal_\RR.$$
Hence $(\sigma,1)$ is an element of order $4$ and $\left[ \mGal_\RR \right]$ is nontrivial.
\qed

\begin{proposition}
Let $F_2$ be a \textbf{dyadic} nonarchimedean local field of characteristic zero.  Then $\left[ \mGal_{F_2} \right]$ is trivial if $[F_2 : \QQ_2]$ is even, and is nontrivial if $[F_2 : \QQ_2]$ is odd.  
\end{proposition}
\proof
Let $d = [F_2 : \QQ_2]$.  By approximation, there exists a global field $F$ such that $F \otimes_\QQ \QQ_2$ is isomorphic to $F_2$ as an $F$-algebra.  Indeed, the primitive element theorem allows us to write $F_2 = \QQ_2(z)$ for $z \in F_2$ a root of a monic irreducible polynomial $P \in \QQ_2[X]$.  This gives an isomorphism of $\QQ_2$-algebras from $\QQ_2[X] / (P)$ to $F_2$.  A small change in the coefficients of $P$ will not change the isomorphism class of the field $\QQ_2[X] / (P)$, by Krasner's lemma.  Hence, by density of $\QQ$ in $\QQ_2$, we may assume that $P \in \QQ[X]$.  As $P$ is monic irreducible over $\QQ_2$, it is monic irreducible over $\QQ$.  The field $F = \QQ[X] / (P)$ satisfies the condition that $F \otimes_\QQ \QQ_2$ is isomorphic to $F_2$ as an $F$-algebra.

The global metaGalois group $\mGal_F$ has a Brauer class $\beta_F$ with local components $\beta_{F,v}$ satisfying
\begin{itemize}
\item
$\beta_{F,v}$ is trivial when $F_v$ has odd residual characteristic (by Proposition \ref{OddRCTrivialBr});
\item
$\beta_{F,2} = \left[ \mGal_{F_2} \right]$ at the unique place of even residual characteristic;
\item
$\beta_{F,v}$ is nontrivial at all real places (by Proposition \ref{RealNontrivBr});
\item
(Parity condition) $\beta_{F,v}$ is nontrivial at a set of places of even cardinality.
\end{itemize}
We have $d = [F_2 : \QQ_2] = [F : \QQ] = r_1 + 2 r_2$, where $r_1$ is the number of real places, and $r_2$ the number of complex places.  It follows that $d$ is even if and only if $r_1$ is even.  The parity condition on the global Brauer class implies that $r_1$ is even if and only if $\beta_{F,2}$ is the trivial class.
\qed

\begin{corollary}
Let $F$ be a global field, with $2 \neq 0$ in $F$.  Then the Brauer class of $\mGal_F$ is that of the unique quaternion algebra which is ramified at all real places and all dyadic places of odd degree over $\QQ_2$.  
\end{corollary}
\proof
This follows directly from the previous three propositions, and the local-global compatibility of the metaGalois group.
\qed

In particular, the Brauer class $\left[ \mGal_\QQ \right]$ is that of the quaternion algebra $\frac{(-1,-1)}{\QQ}$ ramified only at $2$ and $\infty$.  If $F$ is a global field of characteristic $p \neq 2$, then $\left[ \mGal_F \right]$ is the trivial class.

\subsection{Splitting by additive characters}

The metaGalois group may be a nonsplit extension of $\Gal_S$ by $\mu_2$, and even when it splits, it rarely splits canonically.  However, an additive character suffices to split the metaGalois group after pushing out via $\mu_2 \Into\mu_4$.  In the three cases of interest, define a $\sheaf{G}_m[S]$-torsor $\Psi_S$ as follows.
\begin{itemize}
\item
When $F$ is local, let $\Psi_F$ be the set of nontrivial continuous homomorphisms from $F$ to $\CC^\times$.  If $u \in F^\times$, $\psi \in \Psi_F$, write $[u \ast \psi](x) = \psi(u^{-1} x)$.  In this way, $\Psi_F$ is a $F^\times$-torsor.
\item
When $F$ is global, let $\Psi_F$  be the set of nontrivial continuous homomorphisms from $\AA / F$ to $\CC^\times$.  If $u \in F^\times$, $\psi \in \Psi_F$, write $[u \ast \psi](x) = \psi(u^{-1} x)$.  In this way, $\Psi_F$ is a $F^\times$-torsor.
\item
When $F$ is local nonarchimedean, with ring of integers $\OO$, let $\Psi_\OO$  be the set of nontrivial continuous homomorphisms from $F / \OO$ to $\CC^\times$.  If $u \in \OO^\times$, $\psi \in \Psi_\OO$, write $[u \ast \psi](x) = \psi(u^{-1} x)$.  In this way, $\Psi_\OO$ is a $\OO^\times$-torsor.
\end{itemize}

Define here $\mu_4 = \amu_4(\CC) = \{ 1,-1, i, -i \}$.  When $F$ is a local field (with $2 \neq 0$), and $\psi \in \Psi_F$, the \defined{Weil index} is a function $\weil_F(\bullet, \psi) \From F^\times_{/2} \To \mu_4$ which satisfies
\begin{equation}
\label{WeilHilbert}
\frac{\weil_F(uv, \psi)}{\weil_F(u, \psi) \weil_F(v,\psi)} = \Hilb_2(u,v).
\end{equation}
Our $\weil_F(u, \psi)$ is defined to be $\gamma(u x^2) / \gamma(x^2)$ in Weil's notation from \cite[\S 29]{WeilActa} and is written $\gamma_F(u, \psi)$ in \cite[\S A.3]{RangaRao} and elsewhere.  From \cite[Proposition A.11]{RangaRao}, the local Weil indices are trivial on $\OO^\times$ at all nondyadic places. 

When $F$ is a global field and $\psi \in \Psi_F$, the \defined{Weil index} is the function
$$\weil_F(\bullet, \psi) \From \AA^\times_{/2} \To \mu_4$$
defined as the product of local Weil indices.  The global Weil index is trivial on $F^\times$ by \cite[\S II.30, Proposition 5]{WeilActa}.  As the global Hilbert symbol $\Hilb_2 \From \AA^\times \times \AA^\times \To \mu_2$ is defined as the product of local Hilbert symbols, the formula \eqref{WeilHilbert} holds in the global setting too.  

Write $\mGal_S^{(4)}$ for the pushout of $\mGal_S$ via the inclusion $\mu_2 \Into \mu_4$ (when $S = \Spec(F)$ or $S = \Spec(\OO)$ as usual).
$$\begin{tikzcd}[column sep = 4em]
\mu_2 \inarrow{r} \inarrow{d}{\iota} & \mGal_S \onarrow{r} \inarrow{d} & \Gal_S \arrow{d}{=} \\
\mu_4 \inarrow{r} & \mGal_S^{(4)} \onarrow{r} & \Gal_S \arrow[bend right=20, dashed]{l}[swap]{s(\psi)?}
\end{tikzcd}$$
The splittings of $\mGal_S^{(4)}$, if they exist, form a $\Hom(\Gal_S, \mu_4)$-torsor.  In what follows, if $u \in \sheaf{G}_m[S]$, define $\chi_u \From \Gal_S \To \mu_2$ to be the quadratic character associated to the \'etale extension $F[\sqrt{u}]$ (in the local or global case) or $\OO[\sqrt{u}]$ (in the nondyadic local integral case).

\begin{proposition}
\label{WeilSplitsMetaGalois}
For each additive character $\psi \in \Psi_S$, the Weil index provides a splitting $s(\psi) \From \Gal_S \To \mGal_S^{(4)}$.  Moreover, this system of splittings satisfies
$$s(u \ast \psi) = \chi_u \ast s(\psi) \text{ for all } u \in \sheaf{G}_m[S].$$
\end{proposition}
The splittings $s(\psi)$ are described in three cases below.

\subsubsection{Local fields}
When $F$ is a local field, the pushout $\mGal_S^{(4)}$ can be identified with the product $\Gal_F \times \mu_4$ as a set, with multiplication given by
$$(\gamma_1, \zeta_1) \cdot (\gamma_2, \zeta_2) = \left( \gamma_1 \gamma_2, \zeta_1 \zeta_2 \cdot \Hilb_2(\rec_{F/2}(\gamma_1), \rec_{F/2}(\gamma_2) ) \right).$$
For $\psi \in \Psi_F$, \eqref{WeilHilbert} provides a splitting $s(\psi) \From \Gal_F \To \mGal_F^{(4)}$,
$$s(\psi)(\gamma) = \left( \gamma, \weil_F( \rec_{F/2}(\gamma), \psi) \right), \text{ for all } \gamma \in \Gal_F.$$

If $u \in F^\times$, then \cite[Corollary A.5]{RangaRao} states that $\weil_F(a, u \ast \psi) = \Hilb_2(a,u) \cdot \weil_F(a, \psi)$.  Since $\Hilb_2(\rec_{F/2} \gamma, u) = \chi_u(\gamma)$, we find $s(u \ast \psi) = \chi_u \ast s(\psi)$.

\subsubsection{The local integral case}

When $F$ is a nonarchimedean, nondyadic local field, the local Weil index is trivial on $\OO^\times$.  Given a character $\psi \in \Psi_\OO \subset \Psi_F$, the splitting $s(\psi) \From \Gal_F \To \mGal_F^{(4)}$ coincides with the canonical splitting $\sigma^\circ$ on inertia,
$$s(\psi) (\gamma) = \sigma^\circ(\gamma) = (\gamma, 1), \text{ for all } \gamma \in \Inertia.$$
It follows that $s(\psi)$ descends to a splitting of $\mGal_\OO^{(4)}$ at nondyadic places.
$$\begin{tikzcd}[column sep = 4em]
\mu_4 \inarrow{r} & \mGal_\OO^{(4)} \onarrow{r} & \Gal_\OO. \arrow[bend right=20]{l}[swap]{s(\psi)}
\end{tikzcd}$$

If $u \in \OO^\times$, write $\bar u$ for its image in the residue field $\FF_q$.  As before, we have $s(u \ast \psi) = \chi_u \cdot s(\psi)$.  But now, the quadratic character $\chi_u$ is restricted to $\Gal_\OO = \langle \Frob \rangle_{\prof}$; we have
$$\chi_u(\Frob) = \bar u^{(q-1) / 2} \in \mu_2.$$
In other words, $\chi_u$ is the character of $\Gal_\OO$ which sends $\Frob$ to the Legendre symbol of the reduction of $u$.  

\subsubsection{Global fields}

In the global setting, pushing out via $\mu_2 \Into \mu_4$ gives a short exact sequence
\begin{equation}
\label{Adelic4}
\mu_4 \Into \widetilde{\AA^\times_{/2}}^{(4)} \Onto \AA^\times_{/2}.
\end{equation}
The middle term is given by $\widetilde{\AA^\times_{/2}}^{(4)} = \AA^\times_{/2} \times \mu_4$ as a set, with multiplication given by
$$(u_1, \zeta_1) \cdot (u_2, \zeta_2) = (u_1 u_2, \zeta_1 \zeta_2 \cdot \Hilb_2(u_1, u_2)).$$
A character $\psi \in \Psi_F$ provides a splitting of the extension \eqref{Adelic4},
$$s_\AA(\psi)(u) = (u, \weil_F(u, \psi) ) \text{ for all } u \in \AA^\times_{/2}.$$
Since $\weil_F(u, \psi) = 1$ for all $u \in F^\times$, this splitting restricts to the canonical splitting $\sigma_F \From F_{/2}^\times \To \widetilde{\AA_{/2}^\times}$.  Thus $s_\AA(\psi)$ descends and pulls back to a splitting $s(\psi) \From \Gal_F \To  \mGal_F^{(4)}$.  If $u \in F^\times$, then our local results and local-global compatibility imply that $s(u \ast \psi) = \chi_u \cdot s(\psi)$.

\subsection{Restriction}

Suppose that $F' / F$ is a finite separable extension with $F' \subset \bar F$.  In the local integral case, suppose that $F' / F$ is unramified and let $\OO'$ be the ring of integers in $F'$.  Write $S' = \Spec(F')$ in the cases of local or global fields, and write $S' = \Spec(\OO')$ in the local integral case.  We have defined metaGalois groups for $S$ and $S'$.
\begin{equation}
\label{MGalInc}
\begin{tikzcd}
\mu_2 \inarrow{r} \arrow{d}{=} & \mGal_{S'} \onarrow{r} \arrow[dashed]{d}{?}& \Gal_{S'} \inarrow{d} \\
\mu_2 \inarrow{r} & \mGal_S \onarrow{r} & \Gal_S
\end{tikzcd}
\end{equation}

The inclusion $F' \subset \bar F$ gives an inclusion of Galois groups $\Gal_{S'} \Into \Gal_S$, but a natural inclusion of metaGalois groups is not obvious.  In particular, the cocycle defining $\mGal_S$ does not restrict to the cocycle defining $\mGal_{S'}$.

Fortunately, a beautiful insight of Wee Teck Gan gives such an inclusion of metaGalois groups, using a ``lifting theorem'' of Edward Bender \cite{Bender}.  We explain this insight here.

In the case of local fields, consider a nonzero element $u \in F'$, and the ``trace form'' (cf. \cite{SerreTrace}) $F' \To F$ given by $x \mapsto \Tr_{F' / F}(u x^2)$.  Viewing this as a quadratic form on a finite-dimensional $F$-vector space $F'$, it has a Hasse-Witt invariant (an element of $\{ \pm 1 \}$).  Define
$$\HW(u) = \frac{ \text{Hasse-Witt invariant of } x \mapsto \Tr_{F'/F}(u x^2) }{ \text{Hasse-Witt invariant of } x \mapsto \Tr_{F'/F}(x^2) }.$$
This function depends only on the square class of $u$.

Bender's theorem \cite[Theorem 1]{Bender} states that
$$\Hilb_{F',2}(u,v) = \frac{ \HW(u) \HW(v) }{ \HW(uv) } \cdot \Hilb_{F,2} ( \Norm_{F'/F} u, \Norm_{F'/F} v).$$
Let $\iota \From \Gal_{F'} \Into \Gal_F$ be the canonical inclusion, so that $\rec_{F}( \iota(\gamma)) = \Norm_{F'/F} \rec_F(\gamma)$ for all $\gamma \in \Gal_{F'}$.
\begin{proposition}
Let $F$ be a local field (with $2 \neq 0$ as usual).  Then the function $\tilde \iota \From \mGal_{F'} \Into \mGal_F$, given by
$$\tilde \iota (\gamma, \pm 1) = \left( \iota(\gamma), \pm \HW( \rec_{F'} \gamma) \right)$$
is a group homomorphism completing the commutative diagram \eqref{MGalInc}.
\end{proposition}
\proof
Consider any $\gamma_1, \gamma_2 \in \Gal_{F'}$ and define $u_1 \defeq \rec_{F'}(\gamma_1)$, $u_2 \defeq \rec_{F'}(\gamma_2)$.  Thus  $\rec_{F}( \iota(\gamma_1)) = \Norm_{F'/F} u_1$ and $\rec_{F}( \iota(\gamma_2)) = \Norm_{F'/F} u_2$.  For all $\epsilon_1, \epsilon_2 \in \{ \pm 1 \}$, we compute
\begin{align*}
\tilde \iota \left( (\gamma_1, \epsilon_1) \cdot (\gamma_2, \epsilon_2) \right) &= \tilde \iota \left( \gamma_1 \gamma_2, \epsilon_1 \epsilon_2 \Hilb_{F',2} \left( \rec_{F'}(\gamma_1), \rec_{F'}(\gamma_2) \right) \right)  \\
& = \left( \iota(\gamma_1 \gamma_2), \epsilon_1 \epsilon_2 \Hilb_{F',2} (u_1, u_2) \cdot \HW(u_1 u_2) \right) \\
& = \left( \iota(\gamma_1) \iota(\gamma_2), \epsilon_1 \epsilon_2 \Hilb_{F,2} (\Norm_{F' / F} u_1, \Norm_{F' / F} u_2) \HW(u_1) \HW(u_2) \right) \\
& = \left( \iota(\gamma_1), \epsilon_1 \HW(u_1) \right) \cdot \left( \iota(\gamma_2), \epsilon_2 \HW(u_2) \right) \\
&= \tilde \iota( \gamma_1, \epsilon_1) \cdot \tilde \iota(\gamma_2, \epsilon_2).
\end{align*}
\qed

In the local integral case, when $\OO$ is the ring of integers in a nondyadic nonarchimedean field,  $\HW(u) = 1$ for all $u \in \OO^\times$.  From this it follows that $\tilde \iota \From \mGal_{F'} \Into \mGal_{F}$ descends to an injective homomorphism.
$$\begin{tikzcd}
\mu_2 \inarrow{r} \arrow{d}{=} & \mGal_{\OO'} \onarrow{r} \arrow{d}{\tilde \iota}& \Gal_{\OO'} \inarrow{d}{\iota} \\
\mu_2 \inarrow{r} & \mGal_\OO \onarrow{r} & \Gal_\OO
\end{tikzcd}$$

In the global case, when $F$ is a number field, we note that $\prod_{v} \HW_v(u) = 1$ for all $u \in F^\times$ (here $\HW_v$ denotes the invariant as the place $v$).  From this it follows that the injective homomorphisms $\tilde \iota_v \From \mGal_{F_v'} \Into \mGal_{F_v}$ yield a injective homomorphism globally.
$$\begin{tikzcd}
\mu_2 \inarrow{r} \arrow{d}{=} & \mGal_{F'} \onarrow{r} \arrow{d}{\tilde \iota}& \Gal_{F'} \inarrow{d}{\iota} \\
\mu_2 \inarrow{r} & \mGal_F \onarrow{r} & \Gal_F
\end{tikzcd}$$

Taken together, these inclusions $\tilde \iota \From \mGal_{S'} \Into \mGal_S$ allow one to canonically ``restrict'' metaGalois representations (representations of $\mGal_S$).  

\section{L-groups, parameters, L-functions}

\subsection{L-groups}

We use the term ``L-group'' to refer to a broad class of extensions of Galois groups by complex reductive groups.  Unlike Langlands, Vogan, and others, we do not assume that our L-groups are endowed with a conjugacy class of splittings.  Our L-groups are more closely related to the ``weak E-groups'' of \cite[Definition 3.24]{VoganLLC}.  But we maintain the letter ``L'' since our L-groups are still connected to L-functions.  

The other difference between our L-groups and those in the literature is that (for reasons which will become clear) we consider our L-groups as objects of a 2-category.  A base scheme $S$ and geometric point $\bar s \To S$ will be fixed as in the previous section.
\begin{definition}
An \defined{L-group} is a pair $(G^\vee, {}^\EL G)$, where $G^\vee$ is a complex linear algebraic group (not necessarily connected) and ${}^\EL G$ is an extension of locally compact groups
$$G^\vee \Into {}^\EL G \Onto \Gal_S,$$
for which the conjugation action of any element of ${}^\EL G$ on $G^\vee$ is complex-algebraic.
\end{definition}

\begin{remark}
For complex linear algebraic groups, we do not distinguish between the underlying variety and its $\CC$-points.  Thus we say $G^\vee$ is a complex linear algebraic group, and also view $G^\vee$ as a locally compact group.
\end{remark}

Of course, Langlands' L-group ${}^\EL G = \Gal_F \ltimes G^\vee$ (associated to a reductive group $\alg{G}$ over a field $F$) is an example.  When $V$ is a finite-dimensional complex vector space, the direct product $\Gal_S \times GL(V)$ is an L-group.  Since we don't assume $G^\vee$ to be connected, our metaGalois group $\mGal_S$ is an L-group. 

\begin{definition}
Given two L-groups,
$$G_1^\vee \Into {}^\EL G_1 \Onto \Gal_S, \quad G_2^\vee \Into {}^\EL G_2 \Onto \Gal_S,$$
an \defined{L-morphism} ${}^\EL \rho \From {}^\EL G_1 \To {}^\EL G_2$ will mean a continuous group homomorphism lying over $\Id \From \Gal_S \To \Gal_S$, which restricts to a complex algebraic homomorphism $\rho^\vee \From G_1^\vee \To G_2^\vee$.  An \defined{L-equivalence} will mean an invertible  L-morphism.  
\end{definition}
In other words, an L-morphism fits into a commutative diagram, with the middle column continuous and the left column complex-algebraic.
$$\begin{tikzcd}
G_1^\vee \inarrow{r} \arrow{d}{\rho^\vee}& {}^\EL G_1 \onarrow{r} \arrow{d}{{}^\EL \rho} & \Gal_S \arrow{d}{=} \\
G_2^\vee \inarrow{r} & {}^\EL G_2 \onarrow{r} & \Gal_S
\end{tikzcd}$$

\begin{definition}
Given two L-morphisms ${}^\EL \rho, {}^\EL \rho' \From {}^\EL G_1 \To {}^\EL G_2$, a natural isomorphism ${}^\EL \rho \xRightarrow{\sim} {}^\EL \rho'$ will mean an element $a \in Z_2^\vee$ (the center of $G_2^\vee$) such that
$${}^\EL \rho'(g) = a \cdot {}^\EL \rho(g) \cdot a^{-1} \text{ for all } g \in {}^\EL G_1.$$
In particular, note that ${}^\EL \rho$ and ${}^\EL \rho'$ coincide on $G_1^\vee$ when they are naturally isomorphic.

While the axioms for a 2-category are not satisfied if one looks at {\em all} L-groups, L-morphisms, and natural isomorphisms, this does define a 2-category of L-groups, L-{\em equivalences}, and natural isomorphisms.
\end{definition}

In many cases of interest (e.g., when ${}^\EL G_2$ arises as the L-group of a split reductive group), the only natural isomorphism is the identity.  However, in some  nonsplit cases, e.g., ${}^\EL G_2 = \Gal_S \ltimes SL_3(\CC)$, the Langlands L-group of a quasisplit $\alg{G} = \alg{PGU}_3$, a nontrivial element $a \in Z_2^\vee$ does not lie in the center of ${}^\EL G_2$.  Such an element $a$ may determine a nonidentity natural isomorphism.

An \defined{L-representation} of an L-group ${}^\EL G$ will mean a pair $(\rho, V)$, where $V$ is a finite-dimensional complex vector space, and $\rho \From {}^\EL G \To GL(V)$ is a continuous homomorphism whose restriction to $G^\vee$ is complex algebraic.  Giving an L-representation of ${}^\EL G$ is the same as giving an L-morphism ${}^\EL \rho \From {}^\EL G \To \Gal_S \times GL(V)$.

\subsection{Parameters}

Write $\Weil_S$ for the \defined{Weil group}.  When $S = \Spec(F)$ for a local or global field, this Weil group $\Weil_S$ is $\Weil_F$ defined as in \cite{ArtinTate}; when $S = \Spec(\OO)$, we define $\Weil_S$ to be the free cyclic group $\langle \Frob \rangle \isom \ZZ$ generated by a geometric Frobenius $\Frob$.  In all cases, the Weil group is endowed with a continuous homomorphism $\Weil_S \To \Gal_S$ with dense image.

Let $G^\vee \Into {}^\EL G \Onto \Gal_S$ be an L-group.  A \defined{Weil parameter} is a continuous homomorphism $\phi \From \Weil_S \To {}^\EL G$ lying over $\Weil_S \To \Gal_S$, such that $\phi(w)$ is semisimple for all $w \in \Weil_S$ (see \cite[\S 8.2]{BorelCorvallis}).  The reader may follow \cite{BorelCorvallis} and \cite{GrossReeder} to define Weil-Deligne parameters in this general context, when working over a local field.

Write $\Par(\Weil_S, {}^\EL G)$ for the set of ${}^\EL G$-valued Weil parameters.  It is endowed with an action of $G^\vee$ by conjugation:  if $g \in G^\vee$ and $\phi$ is a parameter, then define
$${}^ g \phi(w) = \phi(g^{-1} w g).$$
Two parameters are called \defined{equivalent} if they are in the same $G^\vee$-orbit.

Composition with an L-morphism ${}^\EL \rho \From {}^\EL G_1 \To {}^\EL G_2$ defines a map,
$${}^\EL \rho \From \Par(\Weil_S, {}^\EL G_1) \To \Par(\Weil_S, {}^\EL G_2).$$
Moreover, this map is equivariant, in the sense that for all $g_1 \in G_1^\vee$ and all parameters $\phi \in \Par(\Weil_S, {}^\EL G_1)$, we have
$${}^\EL \rho \left( {}^{g_1} \phi \right) = {}^{\rho^\vee(g_1)} \left( {}^\EL \rho(\phi) \right).$$
Thus the L-morphism $\rho$ descends to a well-defined map of equivalence classes
$${}^\EL \rho \From \frac{ \Par(\Weil_S, {}^\EL G_1)}{G_1^\vee-\text{conjugation} } \To \frac{\Par(\Weil_S, {}^\EL G_2)}{G_2^\vee-\text{conjugation} }.$$

Next, consider a natural isomorphism of L-morphisms $\rho \xRightarrow{\sim} \rho'$, with $\rho, \rho' \From {}^\EL G_1 \To {}^\EL G_2$.  We find two maps of parameter spaces,
$${}^\EL \rho, {}^\EL \rho' \From \Par(\Weil_S, {}^\EL G_1) \To \Par(\Weil_S, {}^\EL G_2),$$ 
and an element $a \in Z_2^\vee$ such that ${}^\EL \rho'$ is obtained from ${}^\EL \rho$ by conjugation by $a$.

It follows that ${}^\EL \rho$ and ${}^\EL \rho'$ induce the {\em same} map on equivalence classes,
$${}^\EL \rho = {}^\EL \rho' \From \frac{ \Par(\Weil_S, {}^\EL G_1)}{G_1^\vee-\text{conjugation} } \To \frac{\Par(\Weil_S, {}^\EL G_2)}{G_2^\vee-\text{conjugation} }.$$

Suppose that an L-group ${}^\EL G$ is defined up to L-equivalence, and the L-equivalence defined up to unique natural isomorphism.  Then the set of equivalence classes of parameters
$$\frac{\Par(\Weil_S, {}^\EL G)}{G^\vee-\text{conjugation}}$$    
is uniquely defined up to unique isomorphism.
\begin{remark}
Refinements of the Langlands parameterization for quasisplit groups suggest that one should look not only at equivalence classes of (Weil or Weil-Deligne) parameters, but also irreducible representations of the component group of the centralizer of a parameter.  Or, following Vogan \cite{VoganLLC}, one can look at $G^\vee$-equivariant perverse sheaves on a suitable variety of parameters.  The fact that conjugation by $a \in Z_2^\vee$ commutes with the conjugation action of $G_2^\vee$ implies that conjugation by $a$ preserves not only the equivalence class of a Weil parameter for ${}^\EL G_2$, but also the equivalence class of such a refined parameter.  If an L-group is defined up to L-equivalence, and the L-equivalence defined up to unique natural isomorphism, then the set of equivalence classes of refined parameters is uniquely defined up to unique isomorphism.
\end{remark}

\subsection{L-functions}

Let $G^\vee \Into {}^\EL G \Onto \Gal_S$ be an L-group, and $\phi \From \Weil_S \To {}^\EL G$ a Weil parameter (or we may take $\phi$ to be a Weil-Deligne parameter in the local case).  Let $(\rho, V)$ be an L-representation of ${}^\EL G$.  Then
$$\rho \circ \phi \From \Weil_S \To GL(V)$$
is a Weil representation (or Weil-Deligne representation in the local case).  As such we obtain an L-function (as defined by Weil and discussed in \cite[\S 3.3]{TateCorvallis}),
$$L(\phi, \rho, s) \defeq L(\rho \circ \phi, s).$$
Choosing an additive character $\psi$ as well gives an $\epsilon$-factor (see \cite[\S 3.4]{TateCorvallis}, based on work of Langlands and Deligne),
$$\epsilon(\phi, \rho, \psi, s) \defeq \epsilon(\rho \circ \phi, \psi, s).$$

In the local integral case $S = \Spec(\OO)$, we have $\Weil_S = \langle \Frob \rangle$, and we define the L-functions and $\epsilon$-factors to be those coming from the unramified representation of $\Weil_F$ by pullback.  

In the setting of Langlands L-groups, a zoo of L-representations arises from complex algebraic representations of $G^\vee$, yielding well-known ``standard'' L-functions, symmetric power and exterior power L-functions, etc..  

But in our very broad setting, we limit our discussion to {\em adjoint} L-functions, as these play an important role in representation theory and their definition is ``internal.''  Consider any L-morphism $\rho \From {}^\EL H \To {}^\EL G$ of L-groups.
$$\begin{tikzcd}
H^\vee \inarrow{r} \arrow{d}{\rho^\vee} & {}^\EL H \onarrow{r} \arrow{d}{ {}^\EL \rho} & \Gal_S \arrow{d}{=}  \\
G^\vee \inarrow{r} & {}^\EL G \onarrow{r} & \Gal_S
\end{tikzcd}$$
For example, we might consider the case where $H^\vee$ is a Levi subgroup of $G^\vee$ (as arises in the Langlands-Shahidi method, \cite{Shahidi}).

Let $\Lie{g}^\vee$ be the complex Lie algebra of $G^\vee$.  The homomorphism ${}^\EL \rho$ followed by conjugation gives an adjoint representation:
$$Ad_\rho \From {}^\EL H \To GL\left( \Lie{g}^\vee  \right).$$
Suppose we have a decomposition of $\Lie{g}^\vee$ as a representation of ${}^\EL H$,
\begin{equation}
\label{AdSummands}
\Lie{g}^\vee = \bigoplus_{i = 0}^h \Lie{g}_{i}^\vee.
\end{equation}
For example, when $H^\vee$ is a Levi subgroup of a parabolic $P^\vee \subset G^\vee$, we may decompose $\Lie{g}^\vee$ into $\Lie{h}^\vee$ and the steps in the nilradical of the Lie algebra of $P^\vee$ and its opposite.   

A decomposition \eqref{AdSummands} gives representations $Ad_{i} \From {}^\EL H \To GL\left( \Lie{g}_{i}^\vee  \right)$.  When $\phi \From \Weil_S \To {}^\EL H$ is a Weil parameter, we obtain L-functions
$$L(\phi, Ad_i, s) \defeq L(Ad_i \circ \phi, s).$$
In particular, when ${}^\EL H = {}^\EL G$, and $\rho = \Id$, we write $Ad$ for the adjoint representation of ${}^\EL G$ on $\Lie{g}^\vee$.  This yields \defined{the adjoint L-function} $L(\phi, Ad, s)$ for any Weil parameter $\phi \From \Weil_S \To {}^\EL G$.  When $H^\vee$ is a Levi subgroup of a parabolic in $G^\vee$, and $Ad_i$ arises from a step in the nilradical of the parabolic, we call $L(\phi, Ad_i, s)$ a \defined{Langlands-Shahidi L-function}.

\begin{remark}
The importance of such L-functions for covering groups is suggested by recent work of D. Szpruch \cite{Szp}, who demonstrates that the Langlands-Shahidi construction of L-functions carries over to the metaplectic group.  But it is not clear how to extend the Langlands-Shahidi method to other covering groups, where uniqueness of Whittaker models often fails.  The thesis work of Gao Fan \cite{GaoThesis} takes some promising steps in this direction.  The general machinery of adjoint L-functions also suggests an analogue, for covering groups, of the Hiraga-Ichino-Ikeda conjecture \cite[Conjecture 1.4]{HII} on formal degrees (see Ichino-Lapid-Mao \cite{IchinoLapidMao}).  It is also supported by the simpler observation that theta correspondence for the metaplectic group $\widetilde{Sp}_{2n}$ provides a definition of adjoint L-functions independently of choices of additive characters.
\end{remark}

\subsection{The L-group of a cover}
\label{LGroupCover}
Now we define the L-group of a cover.  Let $\alg{\tilde G}$ be a degree $n$ cover of a quasisplit group $\alg{G}$ over $S$.  Fix an injective character $\epsilon \From \mu_n \Into \CC^\times$.  Choose a separable closure $\bar F / F$, yielding a geometric base point $\bar s \To S$ and the absolute Galois group $\Gal_S = \pi_1^{\et}(S, \bar s)$.

Recall the constructions of the previous three sections.
\begin{itemize}
\item
$\dgp{\tilde G^\vee}$ denotes the dual group of $\alg{\tilde G}$, a local system on $S_{\et}$ of pinned reductive groups over $\ZZ$,  with center $\dgp{\tilde Z^\vee}$.  It is endowed with a homomorphism $\tau_Q \From \amu_2 \To \dgp{\tilde Z^\vee}$.
\item
$\gerb{E}_\epsilon(\alg{\tilde G})$ is the gerbe associated to $\alg{\tilde G}$, a gerbe on $S_{\et}$ banded by $\sheaf{\tilde Z^\vee} = \dgp{\tilde Z^\vee}(\CC)$.  
\item
$\mu_2 \Into \mGal_S \Onto \Gal_S$ is the metaGalois group.
\end{itemize}

Define $\tilde Z^\vee = \sheaf{\tilde Z}_{\bar s}^\Vee = \dgp{\tilde Z}_{\bar s}^\Vee(\CC)$.  This is the center of the complex pinned reductive group $\tilde G^\vee = \sheaf{\tilde G}_{\bar s}^\Vee = \dgp{\tilde G}_{\bar s}^\Vee(\CC)$.  Note that $\Gal_S$ acts by pinned automorphisms on $\tilde G^\vee$.  

Pushing out $\mGal_S$ via $\tau_Q \From \mu_2 \To \tilde Z^\vee$ defines an L-group,
\begin{equation}
\label{Twist1}
\tilde Z^\vee \Into (\tau_Q)_\ast \mGal_S \Onto \Gal_S.
\end{equation}
From Theorem \ref{AppxFundGpDefined}, the fundamental group of the gerbe $\gerb{E}_\epsilon(\alg{\tilde G})$ at the base point $\bar s$ is an L-group, well-defined up to L-equivalence, and the L-equivalence well-defined up to unique natural isomorphism,
\begin{equation}
\label{Twist2}
\tilde Z^\vee \Into \pi_1^{\et}(\gerb{E}_\epsilon(\alg{\tilde G}), \bar s) \Onto \Gal_S.
\end{equation}

\begin{remark}
The extensions \eqref{Twist1} and \eqref{Twist2} play the role of the first and second twist in \cite{MWCrelle}.  In fact \eqref{Twist1} is canonically isomorphic to the first twist in the split case; the extension \eqref{Twist2} may not coincide with the second twist under some circumstances, and the construction here is more general than \cite{MWCrelle} in both cases.
\end{remark}

The Baer sum of \eqref{Twist1} and \eqref{Twist2} is an L-group which will be called ${}^\EL \tilde Z$,
$$\tilde Z^\vee \Into {}^\EL \tilde Z \Onto \Gal_S.$$
\defined{The L-group} of $\alg{\tilde G}$ is defined to be the pushout of ${}^\EL \tilde Z$ via the inclusion $\tilde Z^\vee \Into \tilde G^\vee$,
$$\tilde G^\vee \Into {}^\EL \tilde G \Onto \Gal_S.$$
More explicitly, this pushout is in the $\Gal_S$-equivariant sense.  In other words,
$${}^\EL \tilde G = \frac{\tilde G^\vee \rtimes {}^\EL \tilde Z}{ \langle (z, z^{-1}) : z \in \tilde Z^\vee \rangle },$$
where the semidirect product action ${}^\EL \tilde Z \To \Aut(\tilde G^\vee)$ is given by projection ${}^\EL \tilde Z \To \Gal_S$ followed by the action $\Gal_S \To \Aut(\tilde G^\vee)$. 

By construction, ${}^\EL \tilde G$ is well-defined by $\alg{\tilde G}$ and $\epsilon$ up to L-equivalence, and the equivalence defined uniquely up to unique natural isomorphism.  To describe ${}^\EL \tilde G$ on the nose (not ``up to L-equivalence''), one must choose a geometric base point $\bar z$ for the gerbe $\gerb{E}_\epsilon(\alg{\tilde G})$ over $\bar s$. 

\subsection{Well-aligned functoriality}
\label{WAFLGroup}

Our notation here follows that of Section \ref{WAFGerbe}:  consider a well-aligned homomorphism $\tilde \iota \From \alg{\tilde G}_1 \To \alg{\tilde G}_2$ of covers, each endowed with Borel subgroup and maximally split maximal torus, i.e., a morphism in the category $\Cat{WAC}_S$.  Fix $\epsilon$ as before.  We have constructed gerbes $\gerb{E}_\epsilon(\alg{\tilde G}_1)$ and $\gerb{E}_\epsilon(\alg{\tilde G}_2)$ associated to $\alg{\tilde G}_1$ and $\alg{\tilde G}_2$, banded by $\sheaf{\tilde Z}_1^\Vee$ and $\sheaf{\tilde Z}_2^\Vee$, respectively.  We have constructed a homomorphism of dual groups $\iota^\vee \From \dgp{\tilde G}_2^\Vee \To \dgp{\tilde G}_1^\Vee$ in Section \ref{WAFDualGroup}, which (after taking $\CC$-points) restricts to $\iota^\vee \From \sheaf{\tilde Z}_2^\Vee \To \sheaf{\tilde Z}_1^\Vee$.  This homomorphism $\iota^\vee$ is compatible with the 2-torsion elements, i.e. $\iota^\vee \circ \tau_{Q_2} = \tau_{Q_1}$.  In Section \ref{WAFGerbe}, we constructed a functor of gerbes $\gerb{i} \From \gerb{E}_\epsilon(\alg{\tilde G}_2) \To \gerb{E}_\epsilon(\alg{\tilde G}_1)$, lying over $\iota^\vee \From \sheaf{\tilde Z}_{2}^\Vee \To \sheaf{\tilde Z}_{1}^\Vee$.

Define $\tilde Z_1^\vee = \dgp{\tilde Z}_{1,\bar s}^\Vee(\CC)$ and $\tilde Z_2^\vee = \dgp{\tilde Z}_{2,\bar s}^\Vee(\CC)$.  These are the centers of the complex pinned reductive groups $\tilde G_1^\vee = \dgp{\tilde G}_{1,\bar s}^\Vee(\CC)$ and $\tilde G_2^\vee = \dgp{\tilde G}_{2,\bar s}^\Vee(\CC)$.

The compatibility $\iota^\vee \circ \tau_{Q_2} = \tau_{Q_1}$ defines an L-morphism,
\begin{equation}
\label{Twist1Diag}
\begin{tikzcd}
\tilde Z_2^\vee \inarrow{r} \arrow{d} & (\tau_{Q_2})_\ast \mGal_S \onarrow{r} \arrow{d} & \Gal_S \arrow{d}{=} \\
\tilde Z_1^\vee \inarrow{r} & (\tau_{Q_1})_\ast \mGal_S \onarrow{r} & \Gal_S 
\end{tikzcd}
\end{equation}
The functor of gerbes $\gerb{i} \From \gerb{E}_\epsilon(\alg{\tilde G}_2) \To \gerb{E}_\epsilon(\alg{\tilde G}_1)$ defines an L-morphism, well-defined up to unique natural isomorphism.
\begin{equation}
\label{Twist2Diag}
\begin{tikzcd}
\tilde Z_2^\vee \inarrow{r} \arrow{d} & \pi_1^{\et}(\gerb{E}_\epsilon(\alg{\tilde G}_2), \bar s) \onarrow{r} \arrow{d} & \Gal_S \arrow{d}{=} \\
\tilde Z_1^\vee \inarrow{r} & \pi_1^{\et}(\gerb{E}_\epsilon(\alg{\tilde G}_1), \bar s) \onarrow{r} &  \Gal_S 
\end{tikzcd}
\end{equation}
Applying Baer sums to \eqref{Twist1Diag} and \eqref{Twist2Diag} yields an L-morphism,
$$\begin{tikzcd}
\tilde Z_2^\vee \inarrow{r} \arrow{d} & {}^\EL \tilde Z_2 \onarrow{r} \arrow{d} & \Gal_S \arrow{d}{=} \\
\tilde Z_1^\vee \inarrow{r} & {}^\EL \tilde Z_1 \onarrow{r} & \Gal_S 
\end{tikzcd}$$
Pushing out yields an L-morphism, well-defined up to natural isomorphism,
$$\begin{tikzcd}
\tilde G_2^\vee \inarrow{r} \arrow{d}{\iota^\vee} & {}^\EL \tilde G_2 \onarrow{r} \arrow{d}{ {}^\EL \iota} & \Gal_S \arrow{d}{=} \\
\tilde G_1^\vee \inarrow{r} & {}^\EL \tilde G_1 \onarrow{r} & \Gal_S 
\end{tikzcd}$$

In this way, the construction of the L-group is contravariantly functorial for well-aligned homomorphisms.

\subsection{Local-global compatibility}

Suppose that $\gamma \From F \Into F_v$ is the inclusion of a global field $F$ into its localization at a place.  Let $\alg{\tilde G}$ be a degree $n$ cover of a quasisplit group $\alg{G}$ over $F$, and let $\epsilon \From \mu_n \Into \CC^\times$ be an injective character.  Let $\bar F \Into \bar F_v$ be an inclusion of separable closures, inducing an inclusion $\Gal_{F_v} \Into \Gal_F$ of absolute Galois groups.  Write $S = \Spec(F)$ and $S_v = \Spec(F_v)$, and $\bar s \To S$ and $\bar s_v \To S_v$ for the geometric base points arising from $\bar F \Into \bar F_v$.

Write $\alg{\tilde G}_v$ for the pullback of $\alg{\tilde G}$ via $\Spec(F_v) \To \Spec(F)$.  Similarly, write $Q_v$ for its Brylinski-Deligne invariant.  The results of Section \ref{BaseChangeDualGroup} identify $\tilde G^\vee = \dgp{\tilde G}_{\bar s}^\Vee(\CC)$ with the corresponding dual group for $\alg{\tilde G}_v$ (relative to the separable closure $\bar F_v$).  Thus we simply write $\tilde G^\vee$ for their dual groups and $\tilde Z^\vee$ for the centers thereof.  

The results of Sections \ref{BaseChangeDualGroup} and \ref{CompatibilityMetaGalois} together provide an L-morphism, unique up to unique natural isomorphism.
\begin{equation}
\label{LGTwist1}
\begin{tikzcd}
\tilde Z^\vee \inarrow{r} \arrow{d}{=} & (\tau_{Q_v})_\ast {}^\EL \mGal_{F_v} \onarrow{r} \inarrow{d} & \Gal_{F_v} \inarrow{d} \\
\tilde Z^\vee \inarrow{r} & (\tau_Q)_\ast {}^\EL \mGal_F \onarrow{r} & \Gal_F
\end{tikzcd}
\end{equation}

The results of Section \ref{BaseChangeGerbe} and following Theorem \ref{AppxFundGpDefined} give an L-morphism, unique up to unique natural isomorphism.
\begin{equation}
\label{LGTwist2}
\begin{tikzcd}
\tilde Z^\vee \inarrow{r} \arrow{d}{=} & \pi_1^{\et}(\gerb{E}_\epsilon(\alg{\tilde G}_v), \bar s_v) \onarrow{r} \inarrow{d} & \Gal_{F_v} \inarrow{d} \\
\tilde Z^\vee \inarrow{r} & \pi_1^{\et}(\gerb{E}_\epsilon(\alg{\tilde G}), \bar s) \onarrow{r} & \Gal_F
\end{tikzcd}
\end{equation}

Applying Baer sums to \eqref{LGTwist1} and \eqref{LGTwist2}, and pushing out via $\tilde Z^\vee \Into \tilde G^\vee$, yields an L-morphism, unique up to unique natural isomorphism.
$$\begin{tikzcd}
\tilde G^\vee \inarrow{r} \arrow{d}{=} & {}^\EL \tilde G_v \onarrow{r} \inarrow{d} & \Gal_{F_v} \inarrow{d} \\
\tilde G^\vee \inarrow{r} & {}^\EL \tilde G \onarrow{r} & \Gal_F
\end{tikzcd}$$
In this way, we identify the L-group of $\alg{\tilde G}_v$ with the pullback of the L-group of $\alg{\tilde G}$, via the inclusion of Galois groups $\Gal_{F_v} \Into \Gal_F$.

\subsection{Parabolic subgroups}
\label{ParabolicLGroup}
Return to a degree $n$ cover $\alg{\tilde G}$ of a quasisplit group $\alg{G}$ over $S$, and let $\alg{P} \subset \alg{G}$ be a parabolic subgroup defined over $S$.  As before, consider a Levi decomposition $\alg{P} = \alg{M} \alg{N}$ and the resulting cover $\alg{\tilde M}$.  Fix $\epsilon \From \mu_n \Into \CC^\times$.

Compatibility of the dual groups from Section \ref{ParabolicDualGroup} gives inclusions
$$ \tilde Z^\vee \Into \tilde Z_M^\vee \Into \tilde M^\vee \Into \tilde G^\vee,$$
where $\tilde Z^\vee$ denotes the center of $\tilde G^\vee$, and $\tilde Z_M^\vee$ denotes the center of $\tilde M^\vee$.  As these inclusions are compatible with the 2-torsion elements in $\tilde Z^\vee$ and $\tilde Z_M^\vee$, we find an L-morphism.
\begin{equation}
\label{PSTwist1}
\begin{tikzcd}
\tilde Z^\vee \inarrow{r}  \inarrow{d} &  (\tau_Q)_\ast {}^\EL \mGal_F \onarrow{r} \inarrow{d}  & \Gal_S \arrow{d}{=} \\
\tilde Z_M^\vee \inarrow{r} &  (\tau_{Q_{\alg{M}}})_\ast {}^\EL \mGal_F \onarrow{r}  & \Gal_S \\
\end{tikzcd}
\end{equation}
Section \ref{ParabolicGerbe} provided a functor of gerbes $\gerb{i} \From \gerb{E}_\epsilon(\alg{\tilde G}) \To \gerb{E}_\epsilon(\alg{\tilde M})$, lying over $\sheaf{\tilde Z}^\Vee \Into \sheaf{\tilde Z}_{\alg{M}}^\Vee$.  This defines an L-morphism of \'etale fundamental groups, up to unique natural isomorphism.
\begin{equation}
\label{PSTwist2}
\begin{tikzcd}
\tilde Z^\vee \inarrow{r} \inarrow{d} & \pi_1^{\et}(\gerb{E}_\epsilon(\alg{\tilde G}), \bar s) \onarrow{r} \inarrow{d} & \Gal_S \arrow{d}{=} \\
\tilde Z_M^\vee \inarrow{r} & \pi_1^{\et}(\gerb{E}_\epsilon(\alg{\tilde M}), \bar s) \onarrow{r}  & \Gal_S \\
\end{tikzcd}
\end{equation}

Applying Baer sums to \eqref{PSTwist1} and \eqref{PSTwist2} yields an L-morphism, uniquely defined up to unique natural isomorphism.
$$\begin{tikzcd}
\tilde Z^\vee \inarrow{r} \inarrow{d} & {}^\EL \tilde Z \onarrow{r} \inarrow{d} & \Gal_S \arrow{d}{=} \\
\tilde Z_M^\vee \inarrow{r} & {}^\EL \tilde Z_M \onarrow{r} & \Gal_S \\
\end{tikzcd}$$
The universal property of pushouts yields an L-morphism.
$$\begin{tikzcd}
\tilde M^\vee \inarrow{r} \inarrow{d}& {}^\EL \tilde M \onarrow{r} \inarrow{d} & \Gal_S  \arrow{d}{=} \\
\tilde G^\vee \inarrow{r}  & {}^\EL \tilde G \onarrow{r} & \Gal_S
\end{tikzcd}$$
This L-morphism is well-defined up to conjugation by $\tilde Z_M^\vee$.

As a special case, when $\alg{P} = \alg{B} = \alg{T} \alg{U}$, we find an L-morphism
$$\begin{tikzcd}
\tilde T^\vee \inarrow{r} \inarrow{d}& {}^\EL \tilde T \onarrow{r} \inarrow{d} & \Gal_S  \arrow{d}{=} \\
\tilde G^\vee \inarrow{r}  & {}^\EL \tilde G \onarrow{r} & \Gal_S
\end{tikzcd}$$
A \defined{principal series} parameter for $\alg{\tilde G}$ is a Weil parameter $\phi \From \Weil_S \To {}^\EL \tilde G$, which factors through a Weil parameter $\Weil_S \To {}^\EL \tilde T$ via the L-morphism above.  By \cite[Theorem 2]{Rajan}, principal series parameters exist, i.e., there exist Weil parameters $\Weil_S \To {}^\EL \tilde T$.  

\subsection{The Weyl-group action on the L-group of a cover of torus}
\label{WeylLTorus}
We keep the degree $n$ cover $\alg{\tilde G}$ of a quasisplit group $\alg{G} \supset \alg{B} \supset \alg{T}$ over $S$.  As before, $\tilde Z^\vee = \dgp{\tilde Z}_{\bar s}(\CC)$ and similarly $\tilde T^\vee = \dgp{\tilde T}_{\bar s}(\CC)$ and $\tilde G^\vee = \dgp{\tilde G}_{\bar s}(\CC)$.  Write $\tilde N^\vee$ for the normalizer of $\tilde T^\vee$ in $\tilde G^\vee$.  Define $W = \sheaf{W}_{\bar s}$, a finite group endowed with $\Gal_S$-action.  As a group endowed with $\Gal_S$-action, $W$ is identified with $\tilde W^\vee = \tilde N^\vee / \tilde T^\vee$.  Assume in this section (as in Section \ref{WeylDualTorus}) that $\alg{T}$ splits over a {\em cyclic} Galois cover of $S$.

For convenience, write $\OO_{\bar s}^\times = \sheaf{G}_{m, \bar s}$; thus $\OO_{\bar s}^\times = \bar F^\times$ if $F$ is a local or global field, and $\OO_{\bar s}^\times$ is endowed with an action of $\Gal_S$.  

Choose a geometric base point $\bar z = (\sheaf{H}, h, j) \in \gerb{E}_\epsilon(\alg{\tilde G})_{\bar s}$.  Such a geometric base point defines a geometric base point $\bar x = (\sheaf{H}, h) \in \gerb{E}_\epsilon(\alg{\tilde T})_{\bar s}$, and a commutative diagram
$$\begin{tikzcd}
\tilde Z^\vee \inarrow{r} \inarrow{d} & \pi_1( \gerb{E}_\epsilon(\alg{\tilde G}), \bar z) \onarrow{r} \inarrow{d} & \Gal_S \arrow{d}{=} \\
\tilde T^\vee \inarrow{r} &  \pi_1( \gerb{E}_\epsilon(\alg{\tilde T}), \bar x ) \onarrow{r} & \Gal_S
\end{tikzcd}$$
Taking the Baer sum with $(\tau_Q)_\ast \mGal_S$, we find an L-morphism
$$\begin{tikzcd}
\tilde Z^\vee \inarrow{r} \inarrow{d} & {}^\EL \tilde Z_{\bar z} \onarrow{r} \inarrow{d} & \Gal_S \arrow{d}{=} \\
\tilde T^\vee \inarrow{r} & {}^\EL \tilde T_{\bar x} \onarrow{r} & \Gal_S
\end{tikzcd}$$
Here we place $\bar z$ and $\bar x$ in subscripts to emphasize the dependence on geometric base point.  We must be careful here to discuss L-morphisms ``on the nose'' and not just up to natural isomorphism.

This L-morphism identifies ${}^\EL \tilde T_{\bar x}$ with $({}^\EL \tilde Z_{\bar z} \ltimes \tilde T^\vee) /  \langle (u,u^{-1}) : u \in \tilde Z^\vee \rangle$.  As described in the previous section, the universal property of pushouts yields an L-morphism,
$$\begin{tikzcd}
\tilde T^\vee \inarrow{r} \inarrow{d} & {}^\EL \tilde T_{\bar x} \onarrow{r} \inarrow{d} & \Gal_S \arrow{d}{=} \\
\tilde G^\vee \inarrow{r} & {}^\EL \tilde G_{\bar z} \onarrow{r} & \Gal_S
\end{tikzcd}$$
In this way, we find a chain of subgroups, ${}^\EL \tilde Z_{\bar z} \Into {}^\EL \tilde T_{\bar x} \Into {}^\EL \tilde G_{\bar z}$ lying over $\Gal_S$.  The inclusion ${}^\EL \tilde T_{\bar x} \Into {}^\EL \tilde G_{\bar z}$ is an inclusion of semidirect products,  
$${}^\EL \tilde T_{\bar x} = \frac{ {}^\EL \tilde Z_{\bar z} \ltimes \tilde T^\vee}{ \langle (u,u^{-1}) : u \in \tilde Z^\vee \rangle} \Into \frac{ {}^\EL \tilde Z_{\bar z} \ltimes \tilde G^\vee}{ \langle (u,u^{-1}) : u \in \tilde Z^\vee \rangle} = {}^\EL \tilde G_{\bar z}.$$
For all $\zeta \in {}^\EL \tilde Z_{\bar z}$, $g^\vee \in \tilde G^\vee$, we write $[\zeta,g^\vee]$ for the corresponding element of ${}^\EL \tilde G_{\bar z}$. 

Let $\dot w \in \sheaf{N}[S]$ represent an element of the Weyl group $w \in \sheaf{W}[S] = W^{\Gal_S}$.  Let $w^\vee$ denote the element of $\tilde W^\vee$ corresponding to the element $w$.  There exists a $\Gal_S$-invariant representative $n^\vee \in \tilde N^\vee$ for $w^\vee$ (by \cite[Lemma 6.2]{BorelCorvallis}).  Conjugation by $n^\vee$ gives an L-morphism $\Int(n^\vee) \From {}^\EL \tilde T_{\bar x} \To {}^\EL \tilde T_{\bar x}$.  

The $\Gal_S$-invariant representative $n^\vee$ of $w^\vee$ is unique up to multiplication by $(\tilde T^\vee)^{\Gal_S}$.  Hence the L-morphism $\Int(n^\vee)$ depends only on $w^\vee$, which in turn depends only on $w$.  Therefore we write $\Int(w^\vee)$ instead of $\Int(n^\vee)$, as we describe the L-morphism below.
\begin{equation}
\label{LMorphW1}
\begin{tikzcd}
\tilde T^\vee \inarrow{r} \arrow{d}{\Int(w^\vee)} & {}^\EL \tilde T_{\bar x} \onarrow{r} \arrow{d}{\Int(w^\vee)} & \Gal_S \arrow{d}{=} \\
\tilde T^\vee \inarrow{r} & {}^\EL \tilde T_{\bar x} \onarrow{r} & \Gal_S
\end{tikzcd}
\end{equation}

\begin{proposition}
\label{WFix1}
For all $w \in \sheaf{W}[S]$, and all $\zeta \in {}^\EL \tilde Z_{\bar z} \subset {}^\EL \tilde T_{\bar x}$, we have $\Int(w^\vee) \zeta = \zeta$.
\end{proposition}
\proof
We compute directly, writing $\gamma$ for the image of $\zeta$ in $\Gal_S$, and $n^\vee$ for a $\Gal_S$-invariant representative of $w^\vee$ in $\tilde N^\vee$.
\begin{align*}
\Int(w^\vee) \zeta &= \Int(n^\vee) [\zeta, 1] \quad \text{(in the group ${}^\EL \tilde G_{\bar z}$)}, \\
&= [1, n^\vee] \cdot [\zeta,1] \cdot [1,(n^\vee)^{-1}], \\
&= [1, n^\vee] \cdot [\zeta,1] \cdot [1,(n^\vee)^{-1}] \cdot [\zeta,1]^{-1} [\zeta,1], \\
&= [1, n^\vee] \cdot [1, {}^\gamma (n^\vee)^{-1} ] \cdot [\zeta, 1], \\
&= [\zeta, 1] = \zeta \quad \text{(since $n^\vee$ is $\Gal_S$-invariant)}.
\end{align*}
\qed

The previous proposition describes the automorphism $\Int(w^\vee)$ of ${}^\EL \tilde T_{\bar x}$.  For any element $[\zeta, t^\vee] \in {}^\EL \tilde T_{\bar x}$, we find that
$$\Int(w^\vee) [\zeta, t^\vee] = [\zeta, \Int(w^\vee) t^\vee].$$

There is another action of the Weyl group on ${}^\EL \tilde T_{\bar x}$, arising as in Sections \ref{WeylDualTorus} and \ref{WeylGerbe} from the well-aligned homomorphism $\Int(\dot w)$.
$$\begin{tikzcd}
\alg{K}_2 \inarrow{r} \arrow{d}{=} & \alg{T}' \onarrow{r} \arrow{d}{\Int(\dot w)} & \alg{T} \arrow{d}{\Int(w)} \\
\alg{K}_2 \inarrow{r} & \alg{T}' \onarrow{r} & \alg{T}
\end{tikzcd}$$ 
From Section \ref{WAFLGroup}, the well-aligned homomorphism defines an L-equivalence lying over $\Int(w)^\vee$, up to natural isomorphism.
\begin{equation}
\label{LMorphW2}
\begin{tikzcd}
\tilde T^\vee \inarrow{r} \arrow{d}{\Int(w)^\vee} & {}^\EL \tilde T_{\bar x} \onarrow{r} \arrow{d}{{}^\EL \Int(\dot w)} & \Gal_S \arrow{d}{=} \\
\tilde T^\vee \inarrow{r} & {}^\EL \tilde T_{\bar x} \onarrow{r} & \Gal_S
\end{tikzcd}
\end{equation}

From Section \ref{WeylDualTorus}, we know that $\Int(w^\vee) = \Int(w)^\vee$, as automorphisms of the dual torus $\tilde T^\vee$.  The rest of this section will be devoted to demonstrating that the L-morphism $\Int(w^\vee)$ in \eqref{LMorphW1} is naturally isomorphic to the L-morphism ${}^\EL \Int(\dot w)$ in \eqref{LMorphW2}.  For this, we must describe ${}^\EL \Int(\dot w)$ in much more detail.

The well-aligned homomorphism $\Int(\dot w)$ gives an equivalence of gerbes,
$$\gerb{Int}(\dot w) \From \gerb{E}_\epsilon(\alg{\tilde T}) \To \gerb{E}_\epsilon(\alg{\tilde T}),$$
lying over the homomorphism $\Int(w)^\vee \From \sheaf{\tilde T}^\Vee \To \sheaf{\tilde T}^\Vee$.  This equivalence sends the object $\bar x = (\sheaf{H}, h)$ to the object $\gerb{Int}(\dot w) \bar x = ({}^w \sheaf{H}, \dot w \circ h)$.  The latter object is described in Section \ref{WeylGerbe}.

Suppose that $\dot \rho \From \bar x \to \gerb{Int}(\dot w) \bar x$ is an isomorphism in the gerbe $\gerb{E}_\epsilon(\alg{\tilde T})$.  Then $\gerb{Int}(\dot w)$ and the isomorphism $\dot  \rho$ define an L-equivalence we call $I(\dot w, \dot  \rho)$,
$$\begin{tikzcd}
\tilde T^\vee \inarrow{r} \arrow{d}{\Int(w)^\vee} & \pi_1( \gerb{E}_\epsilon(\alg{\tilde T}), \bar x ) \onarrow{r} \arrow{d}{I(\dot w, \dot  \rho)} & \Gal_S \arrow{d}{=} \\
\tilde T^\vee \inarrow{r} &  \pi_1( \gerb{E}_\epsilon(\alg{\tilde T}), \bar x ) \onarrow{r} & \Gal_S
\end{tikzcd}$$
Explicitly, if $\phi \From \bar x \To {}^\gamma \bar x$ is an element of $\pi_1( \gerb{E}_\epsilon(\alg{\tilde T}), \bar x )$ lying over $\gamma \in \Gal_S$,
$$I(\dot w, \dot  \rho) \phi = {}^\gamma \dot  \rho^{-1} \circ \gerb{Int}(\dot w) \phi \circ \dot  \rho.$$
$I(\dot w, \dot \rho) \phi$ fits into the commutative diagram below.
$$\begin{tikzcd}[column sep = 5.5em]
\bar x \arrow{d}{\dot  \rho} \arrow{r}{I(\dot w, \dot \rho) \phi} & {}^\gamma \bar x \arrow{d}{ {}^\gamma \dot  \rho} \\
\gerb{Int}(\dot w) \bar x \arrow{r}{\gerb{Int}(\dot w) \phi} & \gerb{Int}(\dot w)  {}^{\gamma} \bar x
\end{tikzcd}$$
On the right side of the diagram, we use the fact that $\gamma(\dot w) = \dot w$, and therefore
$$\gerb{Int}(\dot w) {}^\gamma \bar x = {}^\gamma \left( \gerb{Int}(\dot w)  \bar x \right).$$

\begin{proposition}
Allowing $\dot w$ to vary over representatives for $w$ in $\sheaf{N}[S]$, and allowing $\dot \rho$ to vary over isomorphisms from $\bar x$ to $\gerb{Int}(\dot w) \bar x$, the family of L-equivalences $I(\dot w, \dot \rho)$ defines an L-equivalence $I(w)$ up to unique natural isomorphism.
$$\begin{tikzcd}
\tilde T^\vee \inarrow{r} \arrow{d}{\Int(w)^\vee} & \pi_1( \gerb{E}_\epsilon(\alg{\tilde T}), \bar x ) \onarrow{r} \arrow{d}{I(w)} & \Gal_S \arrow{d}{=} \\
\tilde T^\vee \inarrow{r} &  \pi_1( \gerb{E}_\epsilon(\alg{\tilde T}), \bar x ) \onarrow{r} & \Gal_S
\end{tikzcd}$$
\end{proposition}
\proof
We first vary representatives of $w$.  If $\dot w, \ddot w \in \sheaf{N}[S]$ are two representatives of $w \in \sheaf{W}[S]$, then the unique element $t \in \sheaf{T}[S]$ satisfying $\ddot w = \dot w t$ determines a natural isomorphism of functors $\Fun{A}(t) \From \gerb{Int}(\ddot w) \Rightarrow \gerb{Int}(\dot w)$, according to Proposition \ref{IntTFunctors}.  

If $\dot \rho \From \bar x \To \gerb{Int}(\dot w) \bar x$ is an isomorphism, define $\ddot \rho \From \bar x \To \gerb{Int}(\ddot w) \bar x$ by $\ddot \rho = \Fun{A}(t)_{\bar x}^{-1} \circ \dot \rho$, as in the commutative triangle below.
$$\begin{tikzcd}
\phantom{\bar x} & \gerb{Int}(\ddot w) \bar x \arrow{d}{\Fun{A}(t)_{\bar x}} \\
\bar x \arrow{r}{\dot \rho} \arrow{ur}{\ddot \rho} & \gerb{Int}(\dot w) \bar x 
\end{tikzcd}$$
Naturality of $\Fun{A}(t)$ implies that the following diagram commutes.
$$\begin{tikzcd}[column sep = 4em]
\gerb{Int}(\ddot w) \bar x \arrow{r}{ \gerb{Int}(\ddot w) \phi} \arrow{d}{\Fun{A}(t)_{\bar x}} & \gerb{Int}(\ddot w) {}^\gamma \bar x \arrow{d}{\Fun{A}(t)_{{}^\gamma \bar x} = {}^\gamma \Fun{A}(t)_{\bar x}} \\
\gerb{Int}(\dot w) \bar x \arrow{r}{\gerb{Int}(\dot w) \phi} & \gerb{Int}(\dot w) {}^\gamma \bar x
\end{tikzcd}$$
It follows that
\begin{align*}
I(\ddot w, \ddot \rho) \phi &= {}^\gamma \ddot \rho^{-1} \circ \gerb{Int}(\ddot w) \phi \circ \ddot \rho \\
&= {}^\gamma \dot \rho^{-1} \circ {}^\gamma \Fun{A}(t)_{\bar x} \circ \gerb{Int}(\ddot w) \phi \circ \Fun{A}(t)_{\bar x}^{-1} \circ \dot \rho \\
&= {}^\gamma \dot \rho^{-1} \circ \gerb{Int}(\dot w) \phi \circ \dot \rho = I(\dot w, \dot \rho) \phi.
\end{align*}
Hence $I(\ddot w, \ddot \rho) = I(\dot w, \dot \rho)$.  

Next, we consider varying the isomorphism $\dot \rho$ while keeping $\dot w$ fixed.  Consider another isomorphism $\breve \rho \From \bar x \To \gerb{Int}(\dot w) \bar x$.  There exists a unique $\tau^\vee \in \tilde T^\vee$ such that $\breve \rho = \dot \rho \circ \tau^\vee$.  It follows that ${}^\gamma \breve \rho = \dot \rho \circ \gamma(\tau^\vee)$.  Therefore,
\begin{align*}
I(\dot w, \breve \rho) \phi &= {}^\gamma \breve \rho^{-1} \circ \gerb{Int}(\dot w) \phi \circ \breve \rho \\
&= \frac{1}{\gamma(\tau^\vee)} \circ {}^\gamma \dot \rho^{-1} \circ \gerb{Int}(\dot w) \phi \circ \dot \rho \circ \tau^\vee \\
&= \frac{\tau^\vee}{\gamma(\tau^\vee)} \circ \left( {}^\gamma \dot \rho^{-1} \circ \gerb{Int}(\dot w) \phi \circ \dot \rho \right) \\
&= \frac{\tau^\vee}{\gamma(\tau^\vee)} \circ I(\dot w, \dot \rho) \phi.
\end{align*}
Therefore, $\tau^\vee$ defines a natural isomorphism from the L-equivalence $I(\dot w, \dot \rho)$ to the L-equivalence $I(\dot w, \breve \rho)$.

Hence, as $\dot w$ and $\dot \rho$ vary, the L-equivalence $I(\dot w, \dot \rho)$ varies by uniquely-determined natural isomorphism.  Thus we find an L-equivalence $I(w)$ defined uniquely up to unique natural isomorphism.  
\qed

The following lemma will allow us to study $I(\dot w, \dot \rho)$, in an important special case.  As we work in the geometric fibre throughout, recall that $\tilde T^\vee = \sheaf{\tilde T}_{\bar s}^\vee$, and define $H = \sheaf{H}_{\bar s}$, $\Spl(\sheaf{D}_{Q,n}) = \sspl(\sheaf{D}_{Q,n})_{\bar s}$, and $\hat T = \sheaf{\hat T}_{\bar s}$.  Define $\Phi = \Phi_{\bar s}$ for the set of absolute roots.  Write $Y_{Q,n} = \sheaf{Y}_{Q,n, \bar s}$, so that 
$$\hat T = \Hom(Y_{Q,n}, \OO_{\bar s}^\times), \quad \tilde T^\vee = \Hom(Y_{Q,n}, \CC^\times).$$
\begin{lemma}
\label{RhoAA}
Suppose that $\dot w = r_\beta(e e^- e)$ arises from a relative root $\beta$, for which $2 \beta$ is not a root, as in Theorem \ref{WeylActsOnSplittings}.  Then there exists a pair $(\dot \rho, \dot a)$, with $\dot \rho$ an isomorphism $\bar x \xrightarrow{\sim} \gerb{Int}(\dot w) \bar x$ in the gerbe $\sqrt[n]{\sspl(\sheaf{D}_{Q,n})}$, and $\dot a \in H = \sheaf{H}_{\bar s}$, such that $\dot \rho(\dot a) = \dot a$.
\end{lemma}
\proof
Since we work in the geometric fibre, there exists a splitting $\sigma \in \Spl(\sheaf{D}_{Q,n})$ such that $\sigma$ is aligned with $\dot w$ (aligned as in Theorem \ref{WeylActsOnSplittings}).  Moreover, $h \From H \To \Spl(\sheaf{D}_{Q,n})$ is surjective (Definition \ref{DefCover}(3) implies that $\hat T \xrightarrow{n} \hat T$ is surjective), so there exists $a \in H$ such that $h(a) = \sigma$.

Let $\rho \From \bar x \To \gerb{Int}(\dot w) \bar x$ be any isomorphism.  We find that $\rho \From H \To {}^w H$ satisfies
$$\sigma = h(a) = \dot w(h(\rho(a)).$$
There exists a unique $r \in \hat T$ such that $\rho(a) = r \ast a$.  From this it follows that
$$\sigma = \dot w(h(r \ast a)) = \dot w(r^n \ast \sigma) = w^{-1}(r)^n \ast \dot w(\sigma).$$

Let $\{ \alpha_1, \ldots, \alpha_\ell \}$ denote the roots in $\Phi$ which restrict to the relative root $\beta$.  Define $\tilde \beta \From \OO_{\bar s}^\times \To \hat T = \Hom(Y_{Q,n}, \OO_{\bar s}^\times)$ by
$$\tilde \beta(u) (y) = \prod_{i=1}^\ell u^{\langle \tilde \alpha_i, y \rangle}.$$
Since $w = \prod w_{\alpha_i}$ (orthogonal root reflections), $w(\tilde \beta(u)) = \tilde \beta(u)^{-1}$

Theorem \ref{WeylActsOnSplittings} gives the formula $\dot w(\sigma) = \tilde \beta(\pm 1) \ast \sigma$, where the sign is $-1$ if and only if $q_\beta$ and $n/2$ are both odd.  Hence
$$\sigma = w^{-1}(r^n) \ast \tilde \beta(\pm 1) \ast \sigma.$$
Since $w(\tilde \beta(\pm 1)) = \tilde \beta(\pm 1)$, we have
$$r^n = \tilde \beta(\pm 1).$$

Let $\xi \in \OO_{\bar s}^\times$ be an element satisfying $\xi^n = \pm 1$ (the same sign as above).  Here we use the fact that $n$ is coprime to the characteristics of all residue fields of prime ideals in $S$.  Thus $\tilde \beta(\xi)/r \in \hat T_{[n]}$.  Choose a square root $\sqrt{\xi}$ of $\xi$ in $\OO_{\bar s}^\times$.  (In characteristic $2$, we have $-1 = 1$, and we may take $\xi = 1$ and $\sqrt{\xi} = 1$.)

View $\tilde \beta(\xi)/r$ as an automorphism of $\gerb{Int}(\dot w) \bar x$ in $\sqrt[n]{\sspl(\sheaf{D}_{Q,n})}$ and define
$$\dot \rho = (\tilde \beta(\xi) / r) \circ \rho , \quad \dot a = \tilde \beta(\sqrt{\xi}) \ast a.$$
Then we find
\begin{align*}
\dot \rho(\dot a) &= \dot \rho \left( \tilde \beta(\sqrt{\xi}) \ast a \right) \\
&= \tilde \beta(\sqrt{\xi})^{-1} \ast \dot \rho(a) \\
&= \tilde \beta(\sqrt{\xi})^{-1} \ast \frac{\tilde \beta(\xi)}{r} \ast \rho(a) \\
&= \tilde \beta(\sqrt{\xi})^{-1} \cdot \tilde \beta(\xi) \ast a \\
&=  \tilde \beta(\sqrt{\xi})^{-1} \cdot \tilde \beta(\xi) \cdot \tilde \beta(\sqrt{\xi})^{-1} \ast \dot a = \dot a.
\end{align*}

\qed

The functor of gerbes $\gerb{E}_\epsilon(\alg{\tilde G}) \To \gerb{E}_\epsilon(\alg{\tilde T})$ sends $\bar z$ to $\bar x$ and induces an L-morphism,
$$\begin{tikzcd}
\tilde Z^\vee \inarrow{r} \inarrow{d} & \pi_1( \gerb{E}_\epsilon(\alg{\tilde G}), \bar z ) \onarrow{r} \inarrow{d} & \Gal_S \arrow{d}{=} \\
\tilde T^\vee \inarrow{r} &  \pi_1( \gerb{E}_\epsilon(\alg{\tilde T}), \bar x ) \onarrow{r} & \Gal_S
\end{tikzcd}$$
In this way we view $\pi_1( \gerb{E}_\epsilon(\alg{\tilde G}), \bar z )$ as a subset of $\pi_1( \gerb{E}_\epsilon(\alg{\tilde T}), \bar x )$.

\begin{proposition}
Suppose that $\phi \From \bar z \To {}^\gamma \bar z$ is an element of $\pi_1(\gerb{E}_\epsilon(\alg{\tilde G}), \bar z)$ lying over $\gamma \in \Gal_S$.  Suppose that $\dot w = r_\beta(e e^- e)$ arises from a relative root $\beta$ as in the previous lemma, and choose $\dot \rho \From \bar x \To \gerb{Int}(\dot w) \bar x$ and $\dot a \in \sheaf{H}_{\bar s}$ as in that lemma.  Then
$$I(\dot w, \dot \rho) \phi = \phi.$$  
\end{proposition}
\proof
We continue to work in the geometric fibre throughout this proof.  Note that $\dot a \in H$ and every point of $H$ can be expressed as $\hat \tau \ast \dot a$ for some $\hat \tau \in \hat T$.  Since $\dot \rho(\dot a) = \dot a$, we have for all $\hat \tau \in \hat T$,
$$\dot \rho(\hat \tau \ast \dot a) = w^{-1}(\hat \tau) \ast \dot a.$$
Since $\dot \rho$ intertwines $h$ and $\dot w \circ h$, we find that
$$h(\dot a) = \dot w(h(\dot \rho(\dot a))) = \dot w(h(\dot a)).$$

Now we consider an isomorphism $\phi \From \bar z \To {}^\gamma \bar z$ in the gerbe $\gerb{E}_{\epsilon}(\alg{\tilde G})$.  Such an isomorphism can be written as a contraction $\phi = t^\vee \wedge \phi_0$, where $\phi_0 \From \bar x \To {}^\gamma \bar x$ is an isomorphism in $\sqrt[n]{\sspl(\sheaf{D}_{Q,n})}$ and $t^\vee \in \tilde T^\vee = \Hom(Y_{Q,n}, \CC^\times)$.

The object ${}^\gamma \bar z$ is equal to $({}^\gamma \sheaf{H}, {}^\gamma h, {}^\gamma j)$, where ${}^\gamma \sheaf{H}$ is the torsor with the same underlying sheaf of sets $\sheaf{H}$, and torsor structure given by
$$\hat \tau \ast_\gamma a = \gamma^{-1}(\hat \tau) \ast a \text{ for all } a \in H.$$

The isomorphism $\phi_0$, viewed as a morphism of $\hat T$-torsors from $H$ to ${}^\gamma H$ lying over $\gamma \From \hat T \To \hat T$, satisfies $\phi_0(\dot a) = \hat f \ast \dot a$ for some $\hat f \in \hat T$.  For all $\hat \tau \in \hat T$,
$$\phi_0(\hat \tau \ast \dot a) = \gamma^{-1}(\hat \tau) \ast \hat f \ast \dot a.$$
The isomorphism $\phi = t^\vee \wedge \phi_0$ is therefore determined by the element
$$t^\vee \wedge \hat f  \in \tilde T^\vee \wedge_{\hat T_{[n]}} \hat T = \frac{\tilde T^\vee \times \hat T}{ \langle (\epsilon(\zeta), \zeta^{-1}) : \zeta \in \hat T_{[n]} \rangle}.$$

Recall that $j \From p_\ast \sheaf{H} \To \mu_\ast \Whit$ is an isomorphism in the gerbe $\gerb{E}_{\epsilon}^{\SC}(\alg{\tilde T})$.  Since $\phi$ is an isomorphism in $\gerb{E}_\epsilon(\alg{\tilde G})$, the following diagram commutes in $\gerb{E}_\epsilon^{\SC}(\alg{\tilde T})$.
$$\begin{tikzcd}
p_\ast \sheaf{H} \arrow{r}{j} \arrow{d}{p_\ast \phi} & \mu_\ast \Whit \arrow{d}{\gamma} \\
p_\ast {}^\gamma \sheaf{H} \arrow{r}{{}^\gamma j} & \mu_\ast \Whit
\end{tikzcd}$$
The commutativity of this diagram is equivalent to the fact that
$$t^\vee \wedge \hat f  \in \Ker \left( p \From \tilde T^\vee \wedge_{\hat T_{[n]}} \hat T  \To \tilde T_{\SC}^\vee \wedge_{\hat T_{\SC, [n]}} \hat T_{\SC} \right).$$
Note that $\hat T = \Hom(Y_{Q,n}, \bar F^\times)$ and $\tilde T^\vee = \Hom(Y_{Q,n}, \CC^\times)$, and they are contracted over $\hat T_{[n]} = \Hom(Y_{Q,n}, \mu_n)$ via $\epsilon$.  The map $p$ is given by restricting homomorphisms from $Y_{Q,n}$ to $Y_{Q,n}^{\SC}$.  Recall that $\{ \alpha_1, \ldots, \alpha_\ell \}$ are the roots restricting to the given relative root $\beta$.

Since $t^\vee \wedge \hat f \in \Ker(p)$, there exist $\zeta_i \in \mu_n$ for all $1 \leq i \leq \ell$, such that
$$t^\vee(\tilde \alpha_i^\vee) = \epsilon(\zeta_i)^{-1}, \quad \hat f(\tilde \alpha_i^\vee) = \zeta_i .$$
For each $y \in Y_{Q,n}$ the simple root reflection acts by $w_{\alpha_i}(y) = y - \langle \tilde \alpha_i, y \rangle \tilde \alpha_i^\vee$.  Therefore, for all $1 \leq i \leq \ell$, we have
$$w_{\alpha_i}(t^\vee) = t^\vee \cdot \tilde \alpha_i(\epsilon(\zeta_i))^{-1}, \quad w_{\alpha_i}(\hat f) = \hat f \cdot \tilde \alpha_i(\zeta_i).$$
Here $\tilde \alpha_i(\zeta_i)$ denotes the element of $\hat T_{[n]} = \Hom(Y_{Q,n}, \mu_n)$ given by
$$\tilde \alpha_i(\zeta_i)(y) = \zeta_i^{\langle \tilde \alpha_i, y \rangle},$$
and similarly $\tilde \alpha_i(\epsilon(\zeta_i))$ is an element of $\tilde T_{[n]}^\vee$.

Hence, for all $1 \leq i \leq \ell$, we have
$$w_{\alpha_i}(t^\vee \wedge \hat f) =  t^\vee \cdot \tilde \alpha_i(\epsilon(\zeta_i))^{-1} \wedge \hat f \cdot \tilde \alpha_i(\zeta_i) = t^\vee \wedge \hat f.$$
As $w = \prod_{i=1}^\ell w_{\alpha_i}$, we have $w(t^\vee \wedge \hat f) = t^\vee \wedge \hat f$.  

Now we compute
\begin{align*}
I(\dot w, \dot \rho) \phi_0 (\dot a)&= \left[ {}^\gamma \dot \rho^{-1} \circ \gerb{Int}(\dot w) \phi_0 \circ \dot \rho \right] (\dot a) \\
&= \dot \rho^{-1} \left( \phi_0 \left( \dot \rho(\dot a) \right) \right) \\
&=  \dot \rho^{-1} \left( \phi_0 (\dot a) \right) \\
&= \dot \rho^{-1} \left( \hat f \ast \dot a \right) \\
&= w(\hat f) \ast \dot a \\
&= w(\hat f) \cdot \hat f^{-1} \ast \phi_0(\dot a).
\end{align*}
Therefore,
$$I(\dot w, \dot \rho) \phi_0 = \frac{w(\hat f)}{\hat f} \ast \phi_0.$$
Since $\phi = t^\vee \wedge \phi_0$, we compute
\begin{align*}
I(\dot w, \dot \rho) \phi &= I(\dot w, \dot \rho) (t^\vee \wedge \phi_0) \\
&= w(t^\vee) \wedge I(\dot w, \dot \rho) \phi_0 \\
&=  w(t^\vee) \wedge \left( \frac{w(\hat f)}{\hat f} \ast \phi_0 \right) \\
&=  t^\vee \frac{w(t^\vee)}{t^\vee} \wedge \left( \frac{w(\hat f)}{\hat f} \ast \phi_0 \right) \\
&= t^\vee \wedge \phi_0 = \phi.
\end{align*} 
The last step follows from the identity $w(t^\vee) \wedge w(\hat f) = t^\vee \wedge \hat f$.  
\qed
\begin{remark}
The notation may be difficult to follow in the above argument.  In the contraction of torsors $\tilde T^\vee \wedge_{\hat T_{[n]}} \Hom(\bar x, {}^\gamma \bar x)$, we have
$$\tau_1^\vee \wedge (\hat u_1 \ast \phi_0) = \tau_2^\vee \wedge (\hat u_2 \ast \phi_0)$$
whenever $\tau_1^\vee \wedge u_1 = \tau_2^\vee \wedge u_2$ in $\tilde T^\vee \wedge_{\hat T_{[n]}} \hat T$.
\end{remark}

\begin{thm}
\label{BigWTheorem}
For any $w \in W^{\Gal_S}$, the L-morphism ${}^\EL \Int(\dot w) \From {}^\EL \tilde T \To {}^\EL \tilde T$ is naturally isomorphic to the L-morphism $\Int(w^\vee)$.
\end{thm}
\proof
Suppose first that $\dot w = r_\beta(e e^- e)$, and $\dot \rho$ is chosen as in the previous proposition.  Then the L-morphism
$$I(\dot w, \dot \rho) \From \pi_1( \gerb{E}_\epsilon(\alg{\tilde T}), \bar x) \To \pi_1( \gerb{E}_\epsilon(\alg{\tilde T}), \bar x)$$
restricts to the identity on $\pi_1( \gerb{E}_\epsilon(\alg{\tilde G}), \bar z)$.

The Weyl action on $(\tau_Q)_\ast \mGal_F$ is also trivial, and so we find that the L-morphism
$${}^\EL \Int(\dot w) \From {}^\EL \tilde T_{\bar x} \To {}^\EL \tilde T_{\bar x}$$
restricts to the identity on ${}^\EL \tilde Z_{\bar z}$.  (Strictly speaking, ${}^\EL \Int(\dot w)$ is only defined up to natural isomorphism; but fixing $\dot \rho$ defines ${}^\EL \Int(\dot w)$ on the nose.)  

The semidirect product decomposition, 
$${}^\EL \tilde T_{\bar x} = \frac{{}^\EL \tilde Z_{\bar z} \ltimes \tilde T^\vee}{ \langle (u, u^{-1}) : u \in \tilde z^\vee \rangle},$$
implies that there is a unique L-morphism from ${}^\EL \tilde T_{\bar x}$ to itself which acts as the identity on ${}^\EL \tilde Z_{\bar z}$ and acts via $\Int(w^\vee) = \Int(w)^\vee$ on $\tilde T^\vee$.  Proposition \ref{WFix1} implies that
$${}^\EL \Int(\dot w) \text{ is naturally isomorphic to } \Int(w^\vee).$$

Since $W^{\Gal_S}$ is generated by relative root reflections (for relative roots $\beta$ such that $2 \beta$ is not a relative root), we find that ${}^\EL \Int(\dot w)$ is naturally isomorphic to $\Int(w^\vee)$ for all $w \in W^{\Gal_S}$. 
\qed

\part{Genuine representations}

In this part, we review basic facts about the $\epsilon$-genuine representations (admissible and automorphic) of covering groups $\tilde G$.  Most of these facts are corollaries of previous works by many authors, consolidated and organized for convenience.  However, a new feature of representation theory is the organization of representations into pouches based on their character.  As we review features of the $\epsilon$-genuine representation theory, we introduce corresponding features of L-groups. 

\section{Local fields}

Let $F$ be a local field, and let $\alg{G}$ be a quasisplit reductive group over $F$.  Let $n$ be a positive integer and let $\alg{\tilde G} = (\alg{G}', n)$ be a degree $n$ cover of $\alg{G}$ over $F$ (Definition \ref{DefCover}).  Fix an injective character $\epsilon \From \mu_n \Into \CC^\times$.

The cover $\alg{\tilde G}$ yields a short exact sequence
\begin{equation}
\label{algcover}
\alg{K}_2(F) \Into \alg{G}'(F) \Onto \alg{G}(F).
\end{equation}
Define $G = \alg{G}(F)$.  The Hilbert symbol provides a homomorphism $\Hilb_n \From \alg{K}_2(F) \To \mu_n$.  (When $F \isom \CC$, the Hilbert symbol is trivial; otherwise the Hilbert symbol is surjective).  Pushing out \eqref{algcover} via $\Hilb_n$ yields a central extension of groups
\begin{equation}
\label{loccover}
\mu_n \Into \tilde G \Onto G.
\end{equation}

As a topological group, $\tilde G$ is described in \cite[Construction 10.3]{B-D}.  When $F$ is nonarchimedean, $\tilde G$ is a totally disconnected locally compact group.  When $F$ is archimedean, $\tilde G$ is a real-reductive group in the sense of Harish-Chandra.  When $F \isom \CC$, the Hilbert symbol is trivial, and the extension $\tilde G \To G$ splits canonically.

Let $K$ be a maximal compact subgroup of $G$, and let $\tilde K$ denote its preimage in $\tilde G$.  Then $\tilde K$ is a maximal compact subgroup of $\tilde G$.  Define $\Irr(\tilde K)$ to be the set of equivalence classes of continuous irreducible finite-dimensional representations of $\tilde K$ on complex vector spaces (the \defined{irreps} of $\tilde K$).

\begin{definition}
Suppose that $F$ is nonarchimedean.  An  \defined{admissible representation} of $\tilde G$ is a pair $(\pi, V)$ where $V$ is a complex vector space and $\pi \From \tilde G \To GL(V)$ is a group homomorphism, such that 
\begin{enumerate}
\item
For all $[\chi] \in \Irr(\tilde K)$, the $(\tilde K, [\chi])$-isotypic subspace $V_{[\chi]}$ is finite-dimensional;
\item
As a representation of $\tilde K$, $V = \bigoplus_{[\chi] \in \Irr(\tilde K)} V_{[\chi]}$.
\end{enumerate}
\end{definition}

When $F$ is archimedean, write $\Lie{g}$ for the complexified Lie algebra of $\tilde G$; this is naturally identified with the complexified Lie algebra of $G$.  
\begin{definition}
Suppose that $F$ is archimedean.  An \defined{admissible representation} of $\tilde G$ will mean an admissible $(\Lie{g}, \tilde K)$-module $V$.  In particular, $V$ decomposes as a direct sum of finite-dimensional $\tilde K$-isotypic representations, $V = \bigoplus_{[\chi] \in \Irr(\tilde K)} V_{[\chi]}$ as in the nonarchimedean case.
\end{definition}

In both nonarchimedean and archimedean cases, an admissible representation $V$ of $\tilde G$ is called \defined{$\epsilon$-genuine} if for all $\zeta \in \mu_n$ and all $v \in V$, we have $\pi(\zeta) v = \epsilon(\zeta) \cdot v$.  Note that $\mu_n \subset \tilde K$, so an admissible representation is $\epsilon$-genuine if and only if all of its $\tilde K$-isotypic components are $\epsilon$-genuine.  Define $\Irr_\epsilon(\alg{\tilde G})$ (or $\Irr_\epsilon(\alg{\tilde G} / F)$ in case of confusion) to be the set of equivalence classes of irreducible admissible $\epsilon$-genuine representations of $\tilde G$.

In the nonarchimedean case, write $\hecke_\epsilon(\tilde G)$ for the convolution algebra
$$\hecke_\epsilon(\tilde G) = \{ f \in C_c^\infty(\tilde G) : f(\zeta \tilde g) = \epsilon(\zeta) \cdot f(\tilde g) \text{ for all } \zeta \in \mu_n, \tilde g \in \tilde G \}.$$
Here we fix a Haar measure on $\tilde G$, and convolution is given by
$$[f_1 \ast f_2](g) = \int_G f_1(x) \cdot f_2(x^{-1} g) dx.$$
As $f_1$ and $f_2$ are ``$\epsilon$-genuine functions,'' the integrand is a well-defined function on $G = \tilde G / \mu_n$.  

In the archimedean case, write $\hecke_\epsilon(\tilde G)$ for the convolution algebra of $\epsilon$-genuine, left and right $\tilde K$-finite distributions on $\tilde G$ supported on $\tilde K$.  In both cases, we have a faithful functor from the category of admissible $\epsilon$-genuine representations of $\tilde G$ to the category of (nondegenerate) $\hecke_\epsilon(\tilde G)$-modules.  

\subsection{Unitary, discrete series, and tempered representations}

\begin{definition}
A \defined{unitary representation} of $\tilde G$ is a pair $(\pi, \hat V)$, where $\hat V$ is a separable Hilbert space, and $\pi \From \tilde G \To U(\hat V)$ is an action of $\tilde G$ on $\hat V$ by unitary transformations, such that the resulting map $\tilde G \times \hat V \To \hat V$ is continuous.

A unitary representation $(\pi, \hat V)$ of $\tilde G$ is called \defined{irreducible} if the only closed $\tilde G$-invariant subspaces of $\hat V$ are $0$ and $\hat V$.  Write $\Irr_\epsilon^{\uni}(\alg{\tilde G})$ for the set of equivalence classes of irreducible $\epsilon$-genuine unitary representations of $\tilde G$.
\end{definition}

The following theorem combines a few fundamental results, essentially due to Harish-Chandra.  See \cite{WallachRRGI} for proofs in the (more difficult) archimedean case.  
\begin{thm}
If $(\pi, \hat V)$ is an irreducible unitary representation of $\tilde G$, then the subspace $V \subset \hat V$ of $\tilde K$-finite vectors inherits the structure of an irreducible admissible representation of $\tilde G$.    This defines an injective function,
$$\Irr_\epsilon^{\uni}(\alg{\tilde G}) \Into \Irr_\epsilon(\alg{\tilde G}).$$
\end{thm}

\begin{definition}
Suppose that $(\pi, \hat V)$ is an irreducible unitary $\epsilon$-genuine representation of $\tilde G$.  For all $v_1, v_2 \in \hat V$, the \defined{matrix coefficient} $m_{v_1, v_2} \From \tilde G \To \CC$ is the function
$$m_{v_1, v_2}(g) = \langle v_1, \pi(g) v_2 \rangle.$$
We say that $(\pi, \hat V)$ is a \defined{discrete series} (respectively, \defined{tempered}) representation if for all $\tilde K$-finite vectors $v_1, v_2 \in \hat V$,
$$m_{v_1, v_2} \in L^2(\tilde G / Z(\tilde G)), \quad (\text{resp., } L^{2 + \epsilon}(\tilde G / Z(\tilde G)) \text{ for all } \epsilon > 0).$$
Write $\Irr_\epsilon^{\dis}(\alg{\tilde G})$ (respectively $\Irr_\epsilon^{\temp}(\alg{\tilde G})$) for the set of equivalence classes of discrete series (resp. tempered) $\epsilon$-genuine unitary representations of $\tilde G$.
\end{definition}

In this way, we organize the irreducible admissible $\epsilon$-genuine representations of $\tilde G$ in a nested fashion.
$$\Irr_\epsilon^{\dis}(\alg{\tilde G}) \subset \Irr_\epsilon^{\temp}(\alg{\tilde G}) \subset \Irr_\epsilon^{\uni}(\alg{\tilde G}) \subset \Irr_\epsilon(\alg{\tilde G}).$$
This mirrors the familiar ``uncovered'' case, where there are inclusions,
$$\Irr^{\dis}(\alg{G}) \subset \Irr^{\temp}(\alg{G}) \subset \Irr^{\uni}(\alg{G}) \subset \Irr(\alg{G}).$$

\subsection{Tori}

In previous papers \cite{MWToriOld} \cite{MWToriNew}, we have studied the $\epsilon$-genuine representations of covers of tori over local and global fields.  Here we review the main results over local fields.  Let $\alg{T}$ be a torus over a local field $F$, and let $\alg{\tilde T}$ be a degree $n$ cover of $\alg{T}$.  The resulting central extension
$$\mu_n \Into \tilde T \Onto T$$
is a ``Heisenberg'' type group, and its irreducible representations are therefore determined by their central character.  Let $Z(\tilde T)$ denote the center of $\tilde T$.  For any $[\pi] \in \Irr_\epsilon(\alg{\tilde T})$, let $\chi_{[\pi]}$ be its central character.  Then we have
$$\chi_{[\pi]} \in \Hom_\epsilon( Z(\tilde T), \CC^\times),$$
the set of $\epsilon$-genuine continuous homomorphisms.
\begin{thm}
The map $[\pi] \To \chi_{[\pi]}$ gives a bijection $\Irr_\epsilon(\alg{\tilde T}) \xrightarrow{\sim} \Hom_\epsilon( Z(\tilde T), \CC^\times)$.
\end{thm}
\proof
A proof of this analogue of the Stone von-Neumann theorem can be found in \cite[Theorem 3.1]{MWToriOld} and \cite[Proposition 2.2]{ABPTV} (for the archimedean case).  
\qed

Let $Y$ be the cocharacter lattice of $\alg{T}$.  Let $Q$ be the first Brylinski-Deligne invariant of the cover $\alg{\tilde T}$, and define
$$Y_{Q,n} = \{ y \in Y : n^{-1} B_Q(y, y') \in \ZZ \text{ for all } y' \in Y \},$$
as in \cite{MWToriOld}.  Write $\alg{T}_{Q,n}$ for the $F$-torus with cocharacter lattice $Y_{Q,n}$.  The inclusion $Y_{Q,n} \Into Y$ corresponds to an isogeny $\alg{T}_{Q,n} \To \alg{T}$, and we define $C^\dag(T) = \Im( \alg{T}_{Q,n}(F) \To \alg{T}(F) )$.  Define $C(\tilde T)$ for the preimage of $C^\dag(T)$ in $\tilde T$.  The group $C(\tilde T)$ is called the \defined{central core} of $\tilde T$.

The following result is contained in \cite[Theorem 1.3]{MWToriOld}.
\begin{thm}
\label{CenterLocal}
$C(\tilde T)$ is a finite-index subgroup of $Z(\tilde T)$.  If $\alg{T}$ is split then $C(\tilde T) = Z(\tilde T)$.
\end{thm}

Let $\alg{\tilde T}_{Q,n}$ be the pullback of the cover $\alg{\tilde T}$ to $\alg{T}_{Q,n}$.  Then $\tilde T_{Q,n}$ is abelian and $C(\tilde T) = \Im(\tilde T_{Q,n} \To \tilde T)$.  If $[\pi] \in \Irr_\epsilon(\alg{\tilde T})$, then define $\omega_\pi \in\Irr_\epsilon(\alg{\tilde T}_{Q,n})$ to be the pullback of the central character of $\pi$ to $\tilde T_{Q,n}$.  The character $\omega_\pi$ is called the \defined{central core character} of $[\pi]$.  

\begin{definition}
\label{DefPouchTori}
If $[\pi_1], [\pi_2] \in \Irr_\epsilon(\alg{\tilde T})$, we say that $[\pi_1]$ and $[\pi_2]$ belong to the same \defined{pouch} if they share the same central core character.  In other words, pouches for covers of tori are the fibres of the central core character map
$$\omega \From \Irr_\epsilon(\alg{\tilde T}) \To \Irr_\epsilon(\alg{\tilde T}_{Q,n}).$$
\end{definition}
If $\alg{T}$ is split, then the pouches are singletons.  The image of $\omega$ is a subject for another paper.

On the L-group side, note that the map $\alg{\tilde T}_{Q,n} \To \alg{\tilde T}$ is well-aligned.  Hence it corresponds to an L-morphism,
$$\begin{tikzcd}
\tilde T^\vee \inarrow{r} \arrow{d} & {}^\EL \tilde T \onarrow{r} \arrow{d} & \Gal_F \arrow{d}{=} \\
\tilde T_{Q,n}^\vee \inarrow{r} & {}^\EL \tilde T_{Q,n} \onarrow{r} & \Gal_F
\end{tikzcd}$$

Here the dual groups $\tilde T^\vee$ and $\tilde T_{Q,n}^\vee$ are both equal to $\Hom(Y_{Q,n}, \CC^\times)$, the homomorphism $\tilde T^\vee \To \tilde T_{Q,n}^\vee$ is the identity, and the above diagram is an L-isomorphism.  Hence there is a natural bijection of Weil parameters,
$${}^\EL \omega \From  \Irr_\epsilon(\alg{\tilde T}) \To \Irr_\epsilon(\alg{\tilde T}_{Q,n}).$$

\subsection{Central core character}

An element $g \in G$ is called \defined{regular semisimple} if it is semisimple and the neutral component of its centralizer $\alg{Z}_{\alg{G}}^\circ(g)$ is a maximal torus in $\alg{G}$.  In this case, $g$ is certainly an element of this maximal torus since it centralizes the maximal torus.  Write $G^{\reg}$ for the locus of regular semisimple elements of $G = \alg{G}(F)$; this is a dense open subset of $G$.  Define $\tilde G^{\reg}$ to be the preimage of $G^{\reg}$ in $\tilde G$.

Let $\alg{S} = \alg{Z}^\circ(\alg{G})$ be the neutral component of the center of $\alg{G}$, i.e., the maximal torus contained in the center of $\alg{G}$.  Let $\alg{\tilde S}$ the degree $n$ cover of $\alg{S}$ obtained by pulling back $\alg{\tilde G}$.  If $\alg{T}$ is any maximal torus of $\alg{G}$, then $\alg{S} \subset \alg{T}$.   The cocharacter lattice of $\alg{S}$ is naturally identified with the sublattice $Y^W$ of Weyl-invariants in $Y$.  Define $Y_{Q,n}^W = Y_{Q,n} \cap Y^W$, and let $\alg{S}_{Q,n}$ be the $F$-torus with cocharacter lattice $Y_{Q,n}^W$.  The inclusion $Y_{Q,n}^W \Into Y^W$ corresponds to an $F$-isogeny $\alg{S}_{Q,n} \To \alg{S}$, and we define 
$$C^\dag(G) = \Im(\alg{S}_{Q,n}(F) \To \alg{S}(F)).$$
\begin{definition}
The \defined{central core} of $\tilde G$ is the preimage $C(\tilde G)$ of $C^\dag(G)$ in $\tilde G$.
\end{definition}
Alternatively, pulling back $\alg{\tilde G}$ yields a cover $\alg{\tilde S}_{Q,n}$ and a continuous map $\tilde S_{Q,n} \To \tilde S$.  The central core of $\tilde G$ is the image:
$$C(\tilde G) = \Im \left( \tilde S_{Q,n} \To \tilde S \right).$$
Note that
$$Y_{Q,n}^W \subset (Y^W)_{Q,n} = \{y_1 \in Y^W : B_Q(y_1, y_2) \in n \ZZ \text{ for all } y_2 \in Y^W \}.$$
It follows that $C^\dag(G) \subset C^\dag(S)$, and so $C(\tilde G) \subset C(\tilde S)$.

\begin{proposition}
The central core $C(\tilde G)$ is a finite-index subgroup of $Z(\tilde G)$.
\end{proposition}
\proof
Suppose that $g \in G^{\reg}$, and let $\alg{T} = \alg{Z}_{\alg{G}}^\circ(g)$ denote the neutral component of its centralizer (a maximal $F$-torus in $\alg{G}$).  Since $Y_{Q,n}^W \subset Y_{Q,n} \subset Y$, we find that
$$C(\tilde G) \subset Z(\tilde T),$$
by Theorem \ref{CenterLocal}.  Hence if $\tilde g \in \tilde G^{\reg}$ is a lift of $g$, then $\tilde g \in \tilde T$, and so $C(\tilde G)$ commutes with $\tilde g$.  Since $C(\tilde G)$ commutes with every regular semisimple element $\tilde g \in \tilde G^{\reg}$, the density of $\tilde G^{\reg}$ in $\tilde G$ implies that $C(\tilde G) \subset Z(\tilde G)$.

Since the isogeny $\alg{S}_{Q,n} \To \alg{S}$ has degree coprime to the characteristic of $F$, the image $C^\dag(G)$ has finite index in $S$.  We find a chain of subgroups of $G$ and finite-index inclusions $C^\dag(G) \subset S \subset Z(G)$.  Writing $Z^\dag(G)$ for the image of $Z(\tilde G)$ in $G$, we have $C^\dag(G) \subset Z^\dag(G) \subset Z(G)$ and therefore $C^\dag(G)$ has finite index in $Z^\dag(G)$.  Thus $C(\tilde G)$ has finite index in $Z(\tilde G)$.
\qed

\begin{definition}
For any $[\pi] \in \Irr_\epsilon(\alg{\tilde G})$, let $\omega_\pi$ be the pullback of the central character of $[\pi]$ to $\tilde S_{Q,n}$.  The character $\omega_\pi \in \Irr_\epsilon(\alg{\tilde S}_{Q,n})$ is called the \defined{central core character} of $[\pi]$.
\end{definition}

In this way, the central core character provides a map,
$$\omega \From  \Irr_\epsilon(\alg{\tilde G}) \To  \Irr_\epsilon(\alg{\tilde S}_{Q,n}).$$
Observe that the map $\alg{\tilde S}_{Q,n} \To \alg{\tilde G}$ is well-aligned.  Therefore, it induces an L-morphism,
$$\begin{tikzcd}
\tilde G^\vee \inarrow{r} \arrow{d} & {}^\EL \tilde G \onarrow{r} \arrow{d} & \Gal_F \arrow{d}{=} \\
\tilde S_{Q,n}^\vee \inarrow{r} & {}^\EL \tilde S_{Q,n} \onarrow{r} & \Gal_F
\end{tikzcd}$$
Hence there is a natural map of Weil parameters, generalizing the case of tori,
$${}^\EL \omega \From  \Irr_\epsilon(\alg{\tilde G}) \To \Irr_\epsilon(\alg{\tilde S}_{Q,n}).$$

\subsection{Characters and pouches}

Just as we partitioned $\Irr_\epsilon(\alg{\tilde T})$ into pouches, we can partition $\Irr_\epsilon(\alg{\tilde G})$ into pouches using the character distribution (occasionally assuming $\Char(F) = 0$).  When $[\pi] \in \Irr_\epsilon(\alg{\tilde G})$, the character of $[\pi]$ is a conjugation-invariant distribution $\Tr[\pi]$ on $\tilde G$.  Let $\Theta_\pi$ denote the restriction of $[\pi]$ to $\tilde G^{\reg}$, i.e., to test functions supported in the regular locus.  The following theorem in harmonic analysis is a result of many extensions to a deep result of Harish-Chandra.  
\begin{thm}
\label{CharacterFunction}
$\Theta_\pi$ coincides, after choice of Haar measure, with a smooth (real analytic in the archimedean case, locally constant in the nonarchimedean case) function on $\tilde G^{\reg}$.  In other words, there exists a smooth function $\Theta_\pi \From \tilde G^{\reg} \To \CC$ such that for every test function $f \in C_c^\infty(\tilde G^{\reg})$ we have
$$\Tr[\pi](f) = \int_{\tilde G^{\reg}} \Theta_\pi(x) f(x) dx.$$
\end{thm}
\proof
The key ideas are in the work of Harish-Chandra for semisimple groups at first; see \cite{HCBull}.  The idea extends in a straightforward way for real-reductive groups in the Harish-Chandra class (including our covering groups); see \cite[\S 8.4.1]{WallachRRGI} for a treatment. For a large class of $p$-adic groups, see Clozel \cite{Clo}; the extension to covering groups requires no new ideas.  In characteristic $p$, we refer to Gopal Prasad's appendix to \cite{A-D} and the treatment of Bushnell and Henniart \cite[Corollary A.11]{B-H}.
\qed

\begin{thm}
\label{CharacterInvariant}
Suppose that $F$ has characteristic zero.  If $[\pi_1], [\pi_2] \in \Irr_\epsilon(\alg{\tilde G})$ and $\Theta_{\pi_1} = \Theta_{\pi_2}$, then $[\pi_1] = [\pi_2]$.  In other words, the character distribution, restricted to the regular locus, determines the isomorphism class of a genuine irreducible representation.
\end{thm}
\proof
The character $\Tr[\pi]$ determines the equivalence class $[\pi] \in \Irr_\epsilon(\alg{\tilde G})$.  This result follows from \cite[\S X.1, Theorem 10.6]{Kna} when $F$ is archimedean, and from \cite[Corollary 2.20]{B-Z} when $F$ is nonarchimedean.  In characteristic zero, it is known \cite[Corollaire 4.3.3]{WWL2} that the character distribution $\Tr[\pi]$ is locally integrable, and determined by the function $\Theta_{\pi}$:  for all test functions $f \in C_c^\infty(\tilde G)$ (not necessarily supported on the regular locus!),
$$\Tr[\pi](f) = \int_{\tilde G^{\reg}} \Theta_\pi(x) f(x) dx.$$
\qed

\begin{remark}
If $F$ has characteristic $p \neq 0$, the local integrability of characters is an open question (see Rodier \cite{Rod} and Lemaire \cite{Lem} for two cases where the characters are proven locally integrable).  For this reason, in characteristic $p$, we do not know whether each equivalence class $[\pi] \in \Irr_\epsilon(\alg{\tilde G})$ is determined by $\Theta_\pi$.
\end{remark}

\begin{definition}
Suppose that $\tilde x \in \tilde G^{\reg}$, mapping to $x \in G^{\reg}$.  Let $\alg{T} = \alg{Z}_{\alg{G}}^\circ(x)$ be the maximal torus centralizing $x$.  We say that $\tilde x$ is \defined{genuinely supportive} if $\tilde x \in Z(\tilde T)$.  Let $\tilde G^{\greg}$ denote the set of genuinely supportive elements of $\tilde G^{\reg}$.
\end{definition}

\begin{proposition}
If $[\pi] \in \Irr_\epsilon(\alg{\tilde G})$ then $\Theta_\pi$ is supported on $\tilde G^{\greg}$.
\end{proposition}
\proof
$\Theta_\pi$ is a conjugation-invariant function on $\tilde G^{\reg}$.  If $\tilde x \in \tilde G^{\reg}$ and $\tilde x$ is not genuinely supportive, then there exists $t \in T \subset Z_G(x)$ and $\zeta \in \mu_n$ such that
$$\Int(t) \tilde x = \zeta \cdot \tilde x, \text{ and } 1 \neq \zeta.$$
It follows that $\Theta_\pi(\tilde x) = \Theta_\pi(\Int(t) \tilde x) = \epsilon(\zeta) \cdot \Theta_\pi(\tilde x)$.  Hence $\Theta_\pi(x) = 0$.
\qed

For purposes of parameterization, it seems appropriate to place irreducible representations in the same ``pouch'' if their characters agree on a locus called the ``regular core.''
\begin{definition}
Suppose that $\tilde x \in G^{\reg}$, mapping to $x \in G^{\reg}$.  Let $\alg{T} = \alg{Z}_{\alg{G}}^\circ(x)$.  We say that $\tilde x$ is in the \defined{regular core} of $G$ if $\tilde x \in C(\tilde T)$, the central core of $\tilde T$.  Define $\tilde G^{\creg}$ to be the regular core of $\tilde G$.
\end{definition}
Write $G^{\creg}$ for the image of $\tilde G^{\creg}$ in $G$.  Since $C(\tilde T) \subset Z(\tilde T)$, we find that 
$$C(\tilde G) \cap \tilde G^{\reg} \subset \tilde G^{\creg} \subset \tilde G^{\greg}.$$
The following definition generalizes Definition \ref{DefPouchTori} from tori to reductive groups.
\begin{definition}
Suppose that $[\pi_1], [\pi_2] \in \Irr_\epsilon(\alg{\tilde G})$.  We say that $[\pi_1]$ and $[\pi_2]$ belong to the same \defined{pouch} if $\Theta_{\pi_1}$ and $\Theta_{\pi_2}$ (functions on $\tilde G^{\reg}$ by Theorem \ref{CharacterFunction}) coincide on $\tilde G^{\creg}$.
\end{definition}

\begin{lemma}
Assume that $F$ has characteristic zero.  Suppose that $[\pi] \in \Irr_\epsilon(\alg{\tilde G})$.  Then there exists $\tilde g \in \tilde G^{\creg}$ such that $\Theta_\pi(\tilde g) \neq 0$.
\end{lemma}
\proof
From \cite[Th\'eor\`eme 4.3.2]{WWL2}, in characteristic zero, there exists an open neighborhood $0 \in U \subset \Lie{g}_F$ in the Lie algebra of $\alg{G}$ on which the local character expansion of Harish-Chandra and Howe is valid:
$$\Theta_\pi( \exp(X) ) = \sum_{\OO} c_\OO(\pi) \cdot \hat \mu_\OO(X), \text{ for all } X \in U \cap \Lie{g}_F^{\reg}.$$
Here $\Lie{g}^{\reg}$ denotes the regular semisimple locus.  We use the fact that the cover $\tilde G \To G$ splits over an open subgroup, and any two such splittings coincide on a (possibly smaller) open subgroup.  This allows us to interpret $\exp(X)$ as an element of $\tilde G$, when $X$ is close to $0 \in \Lie{g}_F$.

On the other hand, for any maximal torus $\alg{T} \subset \alg{G}$, the subgroup $C(\tilde T)$ is open of finite index in $\tilde T$.  There are finitely many $\alg{G}(F)$-conjugacy classes of $F$-rational maximal tori in $\alg{G}$.  It follows that for a sufficiently small open neighborhood $0 \in U \subset \Lie{g}_F$, $\exp(U \cap \Lie{g}^{\reg}) \subset \tilde G^{\creg}$.

Hence, if $\Theta_\pi(\tilde g) = 0$ for all $\tilde G^{\creg}$, then we find that $c_\OO(\pi) = 0$ for all nilpotent orbits $\OO$.  But this cannot occur; the character $\Theta_\pi$ of an irreducible admissible representation cannot vanish in a neighborhood of the identity.  
\qed

\begin{proposition}
Suppose that $F$ has characteristic zero.  If $[\pi_1]$ and $[\pi_2]$ belong to the same pouch in $\Irr_\epsilon(\alg{\tilde G})$, then the central core character of $[\pi_1]$ equals the central core character of $[\pi_2]$.
\end{proposition}
\proof
Fix an element $\tilde g \in \tilde G^{\creg}$ for which $\Theta_{\pi_1}(\tilde g) \neq 0$.  If $\tilde c \in C(\tilde G)$, we claim that $\tilde c \cdot \tilde g \in \tilde G^{\creg}$ as well.  Indeed, since $\tilde c \in C(\tilde G)$ and $C(\tilde G) \subset Z(\tilde G)$, we find that $\alg{T} = \alg{Z}_{\alg{G}}^\circ(\tilde g) = \alg{Z}_{\alg{G}}^\circ(\tilde c \cdot \tilde g)$.  In particular, $\tilde c \cdot \tilde g \in \tilde G^{\reg}$.  Since $\tilde g \in C(\tilde T)$ and $\tilde c \in C(\tilde G) \subset C(\tilde T)$, we find that $\tilde c \cdot \tilde g \in C(\tilde T)$.  Thus $\tilde c \cdot \tilde g \in \tilde G^{\creg}$ as claimed.

We now find that
$$\omega_{\pi_1}(\tilde c) \cdot \Theta_{\pi_1}(\tilde g) = \Theta_{\pi_1}(\tilde c \cdot \tilde g) = \Theta_{\pi_2}(\tilde c \cdot \tilde g) = \omega_{\pi_2}(\tilde c) \cdot \Theta_{\pi_2}(\tilde g) = \omega_{\pi_2}(\tilde c) \cdot \Theta_{\pi_1}(\tilde g).$$
Hence $\omega_{\pi_1}(\tilde c) = \omega_{\pi_2}(\tilde c)$.  
\qed

\subsection{Twisting}
Let $\xi \From G \To \CC^\times$ be a smooth homomorphism.  Then $\xi$ pulls back to a smooth character of $\tilde G$ as well.  When $[\pi] \in \Irr_\epsilon(\alg{\tilde G})$, the twist $[\xi \cdot \pi] \in \Irr_\epsilon(\alg{\tilde G})$ as well.  This gives an action
$$\Tw \From \Hom(G, \CC^\times) \times \Irr_\epsilon(\alg{\tilde G}) \To \Irr_\epsilon(\alg{\tilde G}).$$
It is possible that $[\pi] = [\xi \cdot \pi]$, even if $\xi$ is nontrivial.  

\begin{definition}
Let $\xi \From G \To \CC^\times$ be a smooth homomorphism.  We say that $\xi$ is \defined{genuinely trivial} if $\xi(g) = 1$ for all $g \in G^{\greg}$.  We say that $\xi$ is core-trivial if $\xi(g) = 1$ for all $g \in G^{\creg}$.
\end{definition}

\begin{proposition}
Assume that $F$ has characteristic zero.  Suppose that $[\pi] \in \Irr_\epsilon(\alg{\tilde G})$ and $\xi \in \Hom(G, \CC^\times)$ is genuinely trivial.  Then $[\pi] = [\xi \cdot \pi]$.
\end{proposition}
\proof
If $\xi(g) = 1$ for all $g \in G^{\greg}$, then
$$\Theta_\pi(\tilde g) = \Theta_\pi(\tilde g) \cdot \xi(g) = \Theta_{\xi \cdot \pi}(\tilde g)$$
for all $\tilde g \in \tilde G^{\greg}$ mapping to $g \in G^{\greg}$.  But $\Theta_\pi$ and $\Theta_{\xi \cdot \pi}$ are supported on $\tilde G^{\greg}$, and so we find that $\Theta_\pi = \Theta_{\xi \cdot \pi}$.  Hence $[\pi] = [\xi \cdot \pi]$ by Theorem \ref{CharacterInvariant}. 
\qed

\begin{proposition}
\label{TwistingMatters}
Suppose that $[\pi] \in \Irr_\epsilon(\alg{\tilde G})$ and $\xi \in \Hom(G, \CC^\times)$ is core-trivial.  Then $[\pi]$ and $[\xi \cdot \pi]$ belong to the same pouch.
\end{proposition}
\proof
This follows almost from the definition; if $\xi(g) = 1$ for all $g \in G^{\creg}$ then
$$\Theta_\pi(\tilde g) = \Theta_\pi(\tilde g) \cdot \xi(g) = \Theta_{\xi \cdot \pi}(\tilde g)$$
for all $\tilde g \in \tilde G^{\creg}$ mapping to $g \in G^{\creg}$.  Hence $[\pi]$ and $[\xi \cdot \pi]$ belong to the same pouch.
\qed

On the dual group side, let $Z^\vee = \Hom(Y / Y^{\SC}, \CC^\times)$ denote the center of the dual group $G^\vee$.  If $\eta \in H^1(\Weil_F, Z^\vee)$, then one may associate a character $\xi \in \Hom(G, \CC^\times)$.  Let $\tilde Z^\vee = \Hom(Y_{Q,n} / Y_{Q,n}^{\SC}, \CC^\times)$ denote the center of the dual group $\tilde G^\vee$.  Then restriction gives a Galois-equivariant homomorphism $Z^\vee \To \tilde Z^\vee$, and thus a homomorphism
$$H^1(\Weil_F, Z^\vee) \To H^1(\Weil_F, \tilde Z^\vee).$$
Twisting a ${}^\EL \tilde G$-valued Weil parameter by a $\tilde Z^\vee$-valued cocycle, we find a map
$${}^\EL \Tw \From H^1(\Weil_F, Z^\vee) \times \WP_\epsilon(\alg{\tilde G}) \To \WP_\epsilon(\alg{\tilde G}).$$
We might expect $\eta \in H^1(\Weil_F, Z^\vee)$ to correspond to a core-trivial character $\xi \in \Hom(G, \CC^\times)$, if its image in $H^1(\Weil_F, \tilde Z^\vee)$ is trivial.

\subsection{Langlands classification}

Let $\alg{B}$ be a Borel subgroup of $\alg{G}$, defined over $F$.  Let $\alg{A}$ be a maximal $F$-split torus in $\alg{G}$, whose centralizer $\alg{T}$ is a maximal torus in $\alg{B}$.  Let $\Phi_F = \Phi_F(\alg{G}, \alg{A})$ be the resulting set of relative roots, and $\Phi_F^+$ the positive roots determined by $\alg{B}$, and $\Delta_F$ the simple positive roots therein.  Our treatment of the Langlands classification closely follows Ban and Jantzen \cite{BanJantzen}, who extend the Langlands classification to coverings of $p$-adic groups.  

A \defined{standard parabolic} will mean an $F$-parabolic subgroup of $\alg{G}$ containing $\alg{B}$.  If $\alg{P} = \alg{M} \alg{N}$ is a standard parabolic subgroup of $\alg{G}$, write $\Delta_{F,M} \subset \Delta_F$ for the corresponding set of simple roots.  Let $\alg{A}_M$ be the maximal split torus contained in the center of $\alg{M}$, a subgroup of $\alg{A}$.  Define $\Lie{a}_M^\ast = \Hom_F(\alg{M}, \alg{G}_m) \otimes_\ZZ \RR$.  Restriction of a character from $\alg{M}$ to $\alg{A}_M$ extends to an identification of real vector spaces,
$$\Lie{a}_M^\ast = \Hom_F(\alg{M}, \alg{G}_m) \otimes_\ZZ \RR \ident \Hom(\alg{A}_M, \alg{G}_m) \otimes_\ZZ \RR.$$

An inclusion of standard Levi subgroups $\alg{L} \subset \alg{M}$ gives an inclusion $\alg{A}_{M} \subset \alg{A}_L$.  The identifications above yield $\RR$-linear maps,
$$r_M^L \From \Lie{a}_L^\ast \Onto \Lie{a}_{M}^\ast, \quad i_M^L \From \Lie{a}_{M}^\ast \Into \Lie{a}_L^\ast.$$
Write $\Lie{a}^\ast = \Hom(\alg{A}, \alg{G}_m) \otimes \RR$.  We find an injective linear map $i_M \From \Lie{a}_M^\ast \Into \Lie{a}^\ast$ and a surjective linear map $r_M \From \Lie{a}^\ast \Onto \Lie{a}_M^\ast$.  The relative Weyl group $W_F = W_F(\alg{G}, \alg{A})$ acts on $\Lie{a}^\ast$, and there is a unique-up-to-scaling $W_F$-invariant inner product $\langle \cdot, \cdot \rangle$ on $\Lie{a}^\ast$.  

If $\nu = \xi \otimes r \in \Lie{a}_M^\ast$ is a basic tensor, and $\tilde m \in \tilde M$ lies over $m \in M$, define
$$\exp(\nu) \tilde m = \vert \xi(m) \vert^r.$$
This extends to give a homomorphism
$$\exp \From \Lie{a}_M^\ast \To \Hom(\tilde M, \RR_{>0}^\times).$$
Elements $\exp(\nu) \in \Hom(\tilde M, \RR_{>0}^\times)$ are \defined{unramified characters} of $\tilde M$.  

Within the real vector space $\Lie{a}_M^\ast$, define an open cone,
$$(\Lie{a}_M)_+^\ast = \{ \nu \in \Lie{a}_M^\ast : \langle \nu, r_M(\alpha) \rangle > 0 \text{ for all } \alpha \in \Delta_F - \Delta_{F,M} \}.$$

Now, fix a standard parabolic subgroup $\alg{P} = \alg{M} \alg{N}$.  The central extension $\alg{K}_2 \Into \alg{G}' \Onto \alg{G}$ splits uniquely over $\alg{N}$, and thus we have a canonical splitting $N \Into \tilde N$.  Pulling back $\alg{\tilde G}$ to covers $\alg{\tilde P}$ and $\alg{\tilde M}$, we find a semidirect product structure
$$\tilde P = \tilde M \ltimes N.$$
The adjoint action of $\tilde M$ on $N$ has the same modular character $\delta_P$ as the adjoint action of $M$ on $N$. 

Suppose that $(\sigma, \hat H)$ is an $\epsilon$-genuine irreducible unitary representation of $\tilde M$.  Choose $\nu \in \Lie{a}_M^\ast$, so $\exp(\nu) \From \tilde M \To \RR_{>0}^\times$ is an unramified character.  Let $H$ be the subset of $\hat H$ consisting of {\em smooth vectors}.  Define
\begin{equation*}
\begin{split}
I^\infty(\alg{\tilde P}, \sigma, \nu) = \{ & f \From \tilde G \To H : f \text{ is smooth, and } \\
& f(n \tilde m \tilde g) = \delta_P(\tilde m)^{1/2} \sigma(\tilde m) \exp(\nu) (\tilde m) 
\text{ for all } n \in N, \tilde m \in \tilde M, \tilde g \in \tilde G \}.
\end{split}
\end{equation*}
Let $I(\alg{\tilde P}, \sigma, \nu)$ denote the space of $K$-finite vectors therein.  We find that $I(\alg{\tilde P}, \sigma, \nu)$ is an $\epsilon$-genuine admissible representation of $\tilde G$.  This is a result of Ban and Jantzen \cite[\S 3]{BanJantzen} in the nonarchimedean case, and follows from Borel and Wallach \cite[\S III.3.2]{BorelWallach} in the archimedean case.

\begin{thm}[Langlands classification]
If $\sigma$ is an $\epsilon$-genuine irreducible {\em tempered} representation of $\tilde M$ and $\nu \in (\Lie{a}_M)_+^\ast$, then the admissible representation $I(\alg{\tilde P}, \sigma, \nu)$ has a unique irreducible quotient, which we call $J(\alg{\tilde P}, \sigma, \nu)$.  This $J(\alg{\tilde P}, \sigma, \nu)$ is an $\epsilon$-genuine irreducible admissible representation of $\tilde G$.

If $\pi$ is an $\epsilon$-genuine irreducible admissible representation of $\tilde G$, then there exists a unique triple $(\alg{P}, [\sigma], \nu)$, where $\alg{P} = \alg{M} \alg{N}$ is a standard parabolic subgroup of $\alg{G}$, $[\sigma]$ is an equivalence class of $\epsilon$-genuine irreducible tempered representations of $\tilde M$, and $\nu \in (\Lie{a}_M)_+^\ast$ such that $\pi$ is equivalent to $J(\alg{\tilde P}, \sigma, \nu)$.
\end{thm}
\proof
For covering groups in the nonarchimedean case, this is the main result of Ban and Jantzen \cite[Theorem 1.1]{BanJantzen}.  Their restrictions on the characteristic of $F$ are unnecessary here as we work with covers $\alg{\tilde G}$ coming from central extensions of $\alg{G}$ by $\alg{K}_2$ rather than general topological extensions.  Ban and Jantzen rely on the uniqueness of a splitting of $\tilde N \To N$, whereas we can exploit the uniqueness of a splitting of $\alg{N}' \To \alg{N}$ in the algebraic category.  

Indeed, for the canonical splitting $\can \From N \Into \tilde N$, and every $m \in M$, we have
$$\Int(m)^{-1} \circ \can \circ \Int(m) = \can,$$
since both arise from {\em algebraic} splittings of $\alg{K}_2 \Into \alg{N}' \Onto \alg{N}$.  This suffices to demonstrate that \cite[Lemma 2.7, Proposition 2.11]{BanJantzen} hold in arbitrary characteristic, for our class of covers.

Similarly, if $\alg{U} \subset \alg{G}$ is {\em any} unipotent subgroup, and $a \in G$ centralizes $\alg{U}$, then the canonical splitting $\can \From U \Into \tilde U$ satisfies
$$\Int(a) \circ \can = \can.$$
This suffices to demonstrate that \cite[Lemma 2.13]{BanJantzen} holds in arbitrary characteristic, for our class of covers.  Thus the main results of \cite{BanJantzen} hold in arbitrary characteristic, for the nonarchimedean covering groups considered in this paper.

In the archimedean case, we refer to Borel and Wallach \cite[IV, Corollary 4.6, Theorem 4.11]{BorelWallach} (based on earlier work of Harish-Chandra, Casselman, Mili\v{c}i\'c and Langlands).
\qed

On the L-group side, a Weil parameter $\phi \in \WP_\epsilon(\alg{\tilde G})$ is called tempered if its image is bounded, i.e., if the closure of its image in ${}^\EL \tilde G$ is compact.  Recall from Section \ref{ParabolicLGroup} that a parabolic subgroup $\alg{P} = \alg{M} \alg{N}$ yields an embedding of L-groups,
$$\begin{tikzcd}
\tilde M^\vee \inarrow{r} \inarrow{d} & {}^\EL \tilde M \onarrow{r} \inarrow{d} & \Gal_F \arrow{d}{=} \\
\tilde G^\vee \inarrow{r} &  {}^\EL \tilde G \onarrow{r} & \Gal_F
\end{tikzcd}$$
This arises by pushing out an embedding of L-groups,
$$\begin{tikzcd}
\tilde Z^\vee \inarrow{r} \inarrow{d} & {}^\EL \tilde Z \onarrow{r} \inarrow{d} & \Gal_F \arrow{d}{=} \\
\tilde Z_M^\vee \inarrow{r} &  {}^\EL \tilde Z_M \onarrow{r} & \Gal_F
\end{tikzcd}$$

Within $\tilde Z_M^\vee$, let $\tilde Z_{M, \hyp}^\vee$ be the subgroup of \defined{hyperbolic} elements -- those whose eigenvalues are positive real numbers.  The following result gives a Langlands classification for parameters.
\begin{proposition}
For all $\phi \in \WP_\epsilon(\alg{\tilde G})$, there exists a unique triple $(\alg{P}, \phi_M, \eta)$, where $\alg{P} = \alg{M} \alg{N}$ is a standard parabolic subgroup of $\alg{G}$, $\phi_M \in \WP_\epsilon^{\temp}(\alg{\tilde M})$ is a tempered Weil parameter, and $\eta \in H^1(\Weil_F, \tilde Z_{M,\hyp}^\vee)$, such that $\phi$ is equivalent to $\eta \ast \phi_M$ (the parameter obtained by twisting $\phi_M$ by $\eta$, then including ${}^\EL \tilde M \Into {}^\EL \tilde G$).
\end{proposition}
\proof
The proof in the uncovered setting, sketched in various sources, and treated thoroughly by Silberger and Zink in \cite{SilbergerZink} carries over without significant changes.
\qed

\section{Spherical representations}

Now suppose that $\OO$ the ring of integers in a nonarchimedean local field $F$, with residue field $\FF_q$.  Let $\alg{G}$ be a quasisplit reductive group over $\OO$.  Let $\alg{A}$ be a maximal $\OO$-split torus in $\alg{G}$, and let $\alg{T}$ be the centralizer of $\alg{A}$; then $\alg{T}$ is a maximal $\OO$-torus in $\alg{G}$.  Let $\alg{B} = \alg{T} \alg{U}$ be an $\OO$-Borel subgroup of $\alg{G}$ containing $\alg{T}$.  For convenience, let $E/F$ be an unramified extension over which $\alg{T}$ splits, and let ${\Fr}$ be a generator of $\Gal(E/F)$ which corresponds to the geometric Frobenius automorphism over $\FF_q$.  Let $Y$ be the cocharacter lattice of $\alg{T}$, viewed as a $\ZZ[\Fr]$-module.  Let $X$ be the character lattice of $\alg{T}$.  Let $W$ be the (absolute) Weyl group of $\alg{G}$ with respect to $\alg{T}$.  Note that $Y^{\Fr}$ is the group of cocharacters of the split torus $\alg{A}$.  

Let $\alg{\tilde G} = (\alg{G}', n)$ be a degree $n$ cover of $\alg{G}$ over $\OO$; such covers are studied in \cite{MWIntegral}.  The requirement that $\alg{\tilde G}$ is defined over $\OO$ implies that $n$ divides $q-1$.  This ``tameness'' gives an exact sequence,
$$\alg{K}_2(\OO) \Into \alg{K}_2(F) \xtwoheadrightarrow{\Hilb_n} \mu_n.$$

Define $G^\circ = \alg{G}(\OO)$ and $G = \alg{G}(F)$.  Pushing out $\alg{G}'(F)$ via $\Hilb_n$, the central extension $\mu_n \Into \tilde G \Onto G$ is endowed with a splitting over the hyperspecial maximal compact subgroup $G^\circ \subset G$.  We view $G^\circ$ as a subgroup of $\tilde G$ throughout.  Similarly, we write $T^\circ = \alg{T}(\OO)$ and $T = \alg{T}(F)$, and we find a central extension $\mu_n \Into \tilde T \Onto T$ endowed with a splitting over $T^\circ$.  Define $W^\circ = N_{G^\circ}(T^\circ) / T^\circ$.  The inclusion $N_{G^\circ}(T^\circ) \Into \alg{N}_{\alg{G}}(\alg{T})(\OO)$ yields an identification of $W^\circ$ with $W^{\Fr}$.

Fix an injective character $\epsilon \From \mu_n \To \CC^\times$ as usual.  Suppose that $(\pi, V)$ is an irreducible admissible $\epsilon$-genuine representation of $\tilde G$.  We say that $(\pi, V)$ is \defined{spherical} if the space of $G^\circ$-fixed vectors $V^{G^\circ}$ is nonzero.  Define $\Irr_\epsilon(\alg{\tilde G} / \OO)$ to be the set of equivalence classes of spherical irreducible admissible $\epsilon$-genuine representations of $\tilde G$.

\subsection{Hecke algebras}

The $\epsilon$-genuine \defined{spherical Hecke algebra}, denoted $\hecke_\epsilon(\tilde G, G^\circ)$ is the space of compactly supported functions $f \From \tilde G \To \CC$ which satisfy
$$f(\zeta k_1 \tilde g k_2) = \epsilon(\zeta) \cdot f(\tilde g) \text{ for all } \zeta \in \mu_n, k_1,k_2 \in G^\circ, \tilde g \in \tilde G.$$
Here, we fix the Haar measure on $\tilde G$ for which $G^\circ$ has measure $1$.  Then $\hecke_\epsilon(\tilde G, G^\circ)$ is a subalgebra of $\hecke_\epsilon(\tilde G)$,
$$\hecke_\epsilon(\tilde G, G^\circ) = \Char(G^\circ) \ast \hecke_\epsilon(\tilde G) \ast \Char(G^\circ).$$

When $(\pi, V)$ is a smooth $\epsilon$-genuine spherical representation of $\tilde G$, the set $V^{G^\circ}$ of $G^\circ$-fixed vectors forms a module over the Hecke algebra $\hecke_\epsilon(\tilde G, G^\circ)$.  In this way, $\Irr_\epsilon(\alg{\tilde G} / \OO)$ is identified with the set of equivalence classes of irreducible $\hecke_\epsilon(\tilde G, G^\circ)$-modules.

The \defined{support} of $\hecke_\epsilon(\tilde G, G^\circ)$ is
$$\{ \tilde g \in \tilde G : f(\tilde g) \neq 0 \text{ for some } f \in \hecke_\epsilon(\tilde G, G^\circ) \}.$$
Similarly, we may define the spherical Hecke algebra $\hecke_\epsilon(\tilde T, T^\circ)$, its support, and the set of its irreducible modules $\Irr_\epsilon(\alg{\tilde T} / \OO)$.

The support of $\hecke_\epsilon(\tilde G, G^\circ)$ is a set of the form $\mu_n G^\circ \Lambda G^\circ$, for some subset $\Lambda \subset \tilde G$.  The Cartan decomposition gives
$$G = G^\circ Y^{\Fr} G^\circ;$$
in the sense that for any choice of uniformizing parameter $\varpi$, every element of $G$ can be expressed as $k_1 y(\varpi) k_2$ for some $y \in Y^{\Fr}$.  In this way, the support of $\hecke_\epsilon(\tilde G, G^\circ)$ corresponds to a $W^\circ$-stable subset of $Y^{\Fr}$.  We call this subset of $Y^{\Fr}$ the \defined{combinatorial support} of $\hecke_\epsilon(\tilde G, G^\circ)$.  Similarly, the support of $\hecke_\epsilon(\tilde T, T^\circ)$ is determined by its combinatorial support, a subset of $Y^{\Fr}$.  

The supports of $\hecke_\epsilon(\tilde G, G^\circ)$ and $\hecke_\epsilon(\tilde T, T^\circ)$ are given by Mcnamara \cite[Theorem 10.1]{McNamara}, and by Wen Wei Li in \cite[\S 3.2]{WWL1}.
\begin{proposition}
\label{HeckeSupport}
The combinatorial supports of $\hecke_\epsilon(\tilde G, G^\circ)$ and $\hecke_\epsilon(\tilde T, T^\circ)$ coincide.  They equal $Y_{Q,n}^{\Fr}$.  
\end{proposition}

From \cite[Lemma 3.2.3]{WWL1}, it is known that the Hecke algebra $\hecke_\epsilon(\tilde T, T^\circ)$ is commutative.  It is (noncanonically) isomorphic to a polynomial ring over $\CC$ in $r$ variables, where $r$ is the rank of $\alg{A}$.  The support of $\hecke_\epsilon(\tilde T, T^\circ)$ is contained in the centralizer of $T^\circ$ in $\tilde T$ (see \cite[\S 3.2]{MWToriNew}).  It follows that $W^\circ = N_{G^\circ}(T^\circ) / T^\circ$ acts on $\hecke_\epsilon(\tilde T, T^\circ)$ by the rule
$${}^w f(\tilde t) \defeq f( \dot w^{-1} \tilde t \dot w ) \text{ for all } w \in W^\circ \text{ represented by } \dot w \in N_{G^\circ}(T^\circ).$$

The Satake isomorphism is a ring homomorphism from the spherical Hecke algebra of an unramified $p$-adic reductive group to the corresponding Hecke algebra of a maximal torus.  For covering groups, this has been considered by McNamara \cite{McNamara} and Wen Wei Li \cite{WWL1}.

Write $\delta: \tilde T \To \RR_{>0}^\times$ for the modular character for the adjoint action of $\tilde T$ on $U$.  Suppose that $f \in \hecke_\epsilon(\tilde G, G^\circ)$.  Define the \defined{Satake transform} of $f$, a function $\Sat f \From \tilde T \To \CC$, by
$$[\Sat f](\tilde t) = \delta(\tilde t)^{-1/2} \int_U f(u \tilde t ) du.$$

From \cite[Proposition 3.2.5]{WWL1}, we have the following result.
\begin{thm}
The Satake transform is a ring isomorphism from $\hecke_\epsilon(\tilde G, G^\circ)$ to $\hecke_\epsilon(\tilde T, T^\circ)^{W^\circ}$.
\end{thm}

\begin{corollary}
\label{UniqueSph}
The Hecke algebra $\hecke_\epsilon(\tilde G, G^\circ)$ is a finitely-generated commutative $\CC$-algebra.  If $(\pi, V)$ is an irreducible $\epsilon$-genuine spherical representation of $\tilde G$, then $V^{G^\circ}$ is one-dimensional.
\end{corollary}

\begin{corollary}
\label{SatBij}
The Satake isomorphism defines a bijection,
$$\Irr_\epsilon(\alg{\tilde G} / \OO) \ident \Hom(\hecke_\epsilon(\tilde G, G^\circ), \CC) \xrightarrow{\Sat^\ast} \Hom(\hecke_\epsilon(\tilde T, T^\circ)^{W^\circ}, \CC) \ident \Irr_\epsilon(\alg{\tilde T} / \OO) / W^{\Fr}.$$
\end{corollary}

In other words, the Satake isomorphism gives a bijection from the set of equivalence classes of $\epsilon$-genuine spherical irreps of $\tilde G$ to the $W^{\Fr}$-orbits on the set of equivalence classes of $\epsilon$-genuine spherical irreps of $\tilde T$.  In applying the Satake isomorphism, the $W^{\Fr}$ action arises from the $W^\circ$ action on the Hecke algebra.  This coincides with another action described below.

For $w \in W^{\Fr}$, conjugation by $w$ defines a homomorphism $\Int(w) \From \alg{G} \To \alg{G}$, which restricts to a homomorphism $\Int(w) \From \alg{T} \To \alg{T}$.  For any representative $\dot w \in \alg{N}_{\alg{G}}(\alg{T})(\OO)$, we find a well-aligned homomorphism of covers from $\alg{\tilde T}$ to itself, as in Section \ref{WeylDualTorus}.
$$\begin{tikzcd}
\alg{K}_2 \inarrow{r} \arrow{d}{=} & \alg{T}' \onarrow{r} \arrow{d}{\Int(\dot w)} & \alg{T} \arrow{d}{\Int(w)}  \\
\alg{K}_2 \inarrow{r} & \alg{T}' \onarrow{r} & \alg{T} 
\end{tikzcd}$$
Pulling back via $\Int(\dot w)$, we find a map, $\Int(\dot w)^\ast \From \Irr_\epsilon(\tilde T / \OO) \To  \Irr_\epsilon(\tilde T / \OO)$.  Explicitly, if $(\pi, V)$ is an $\epsilon$-genuine spherical irrep of $\tilde T$, then $\Int(\dot w)^\ast (\pi, V) = (\pi \circ \dot w^{-1}, V)$, where
$$[\pi \circ \dot w^{-1}](\tilde t) \defeq \pi(\dot w^{-1} \tilde t \dot w).$$
If $\ddot w = \dot w \tau$ for some $\tau \in T^\circ$, then 
$$\pi \circ \ddot w^{-1}(\tilde t) = \pi(\ddot w^{-1} \tilde t \ddot w) = \pi(\tau^{-1} \dot w^{-1} \tilde t \dot w \tau) = \pi(\tau)^{-1} [\pi \circ \ddot w^{-1}](\tilde t) \pi(\tau).$$
Hence $\pi(\tau)$ gives an equivalence of representations from $(\pi \circ \dot w^{-1}, V)$ to $(\pi \circ \ddot w^{-1}, V)$.  Therefore, we find a well-defined map (independent of representative),
$$\Int(w)^\ast \From \Irr_\epsilon(\tilde T / \OO) \To  \Irr_\epsilon(\tilde T / \OO).$$
This gives an action of $W^{\Fr}$ on $\Irr_\epsilon(\tilde T / \OO)$ which coincides with the action obtained from the $W^\circ$ action on the Hecke algebra $\hecke_\epsilon(\tilde T, T^\circ)$.

\section{Automorphic representations}

Now let $F$ be a global field, $\VV$ the set of places of $F$, and $\alg{G}$ a quasisplit reductive group over $F$.  Write $G_F = \alg{G}(F)$, $G_v = \alg{G}(F_v)$ for all $v \in \VV$, and $G_\AA = \alg{G}(\AA)$, where $\AA$ denotes the ring of adeles of $F$.  

Let $\alg{\tilde G} = (\alg{G}', n)$ be a degree $n$ cover of $\alg{G}$ over $F$.  We may assume (see \cite[\S 10.4]{B-D}) that $\alg{G}$ and $\alg{G}'$ are defined over the ring $\OO_\SS$ of $\SS$-integers in $F$, for some finite subset $\SS \subset \VV$.  

In this way, for every place $v \in \VV$, we find a central extension of locally compact groups, $\mu_n \Into \tilde G_v \Onto G_v$, endowed with splittings $G_v^\circ = \alg{G}_v(\OO_v) \Into \tilde G_v$ for all $v \in \VV - \SS$.  Define $\dtilde{G}_\AA$ to be the restricted direct product of the extensions $\tilde G_v$, with respect to the compact open subgroups $G_v^\circ$ (defined at almost all places).  
$$\bigoplus_{v \in \VV} \mu_n \Into \dtilde{G}_{\AA} \Onto G_\AA.$$
Pushing out via the product map $\bigoplus_{v \in \VV} \mu_n \To \mu_n$ yields a central extension $\tilde G_\AA$ as in the diagram below.
$$\begin{tikzcd}
\bigoplus_{v \in \VV} \mu_n \inarrow{r} \onarrow{d}{\Pi} & \dtilde{G}_\AA \onarrow{r} \onarrow{d} & G_\AA \arrow{d}{=} \\
\mu_n \inarrow{r} & \tilde G_\AA \onarrow{r} & G_\AA
\end{tikzcd}$$
As explained in \cite[\S 10.4]{B-D}, $\mu_n \Into \tilde G_\AA \Onto G_\AA$ is a central extension of locally compact groups, split canonically over $G_F$.  Write $\tilde Z_\AA$ for the center of $\tilde G_\AA$.  Then $G_F$ is discrete in $\tilde G_\AA$ and $\tilde Z_\AA G_F \backslash \tilde G_\AA$ has finite volume with respect to a $\tilde G_\AA$-invariant measure.

\subsection{Admissible and unitary representations}

Our approach to the admissible and automorphic representations of $\tilde G_\AA$ coincides with the toral case in \cite[\S 4]{MWToriNew}, based on Flath \cite{Flath} and Borel-Jacquet \cite{BorelJacquet}.  Let $\AA_{\fin}$ be the finite adeles.  Write $\dtilde G_{\AA_{\fin}}$ for the restricted direct product of $\tilde G_v$ with respect to the compact open subgroups $G_v^\circ$, where the product is indexed by only the nonarchimedean places $v \in \VV_{\fin}$.  Then $\dtilde G_{\AA_{\fin}}$ is a totally disconnected locally compact group, a central extension as below.
$$\bigoplus_{v \in \VV_{\fin}} \mu_n \Into \dtilde G_{\AA_{\fin}} \Onto G_{\AA_{\fin}}.$$
Let $\dtilde K_{\fin}$ be an open compact subgroup of $\dtilde G_{\AA_{\fin}}$.  

Define the $\epsilon$-genuine Hecke algebra,
\begin{equation*}
\begin{split}
\hecke_{\fin, \epsilon} = \{ & f \in C_c^\infty(\dtilde G_{\AA_{\fin}}) : f( (\zeta_v) \cdot \dtilde g) = \prod_{v \in \VV_{\fin}} \epsilon(\zeta_v) \cdot f(\dtilde g) \\
&  \text{ for all } (\zeta_v) \in \bigoplus_{v \in \VV_{\fin}} \mu_n, \quad \dtilde g \in \dtilde G_{\AA_{\fin}} \}.
\end{split}
\end{equation*}
As in \cite[Example 2]{Flath}, the Hecke algebra $\hecke_{\fin, \epsilon}$ is isomorphic to the tensor product of the $\epsilon$-genuine Hecke algebras $\hecke_{v, \epsilon} = \hecke_\epsilon(\tilde G_v)$ with respect to the system of idempotents $\Char(G_v^\circ)$ at almost all places.  An \defined{admissible} $\epsilon$-genuine representation of $\dtilde G_{\AA_{\fin}}$ is a representation $(\pi, V)$ of $\dtilde G_{\AA_{\fin}}$ such that
$$\pi( (\zeta_v) \dtilde g) = \prod_{v \in \VV_{\fin}} \epsilon(\zeta_v) \cdot \pi(\dtilde g), \text{ for all } (\zeta_v) \in \bigoplus_{v \in \VV_{\fin}} \mu_n, \dtilde g \in \dtilde G_{\AA_{\fin}},$$
and which decomposes as a direct sum of irreducible representations of $\dtilde K_{\fin}$, each occurring with finite multiplicity.

The commutativity of the spherical Hecke algebras $\hecke(\tilde G_v, G_v^\circ)$ at almost all places implies the factorization of irreducible admissible representations.  The following is adapted from \cite[Proposition 4.1]{MWToriNew}, and is a direct result of \cite[Theorem 2]{Flath}.
\begin{proposition}
For every irreducible admissible $\epsilon$-genuine representation $\pi_{\fin}$ of $\dtilde G_{\AA_{\fin}}$, there exists a unique family of equivalence classes $([\pi_v])_{v \in \VV_{\fin}}$ of irreducible admissible $\epsilon$-genuine representations of each $\tilde G_v$, such that $\pi_v$ is spherical for almost all $v \in \VV_{\fin}$ and $\pi_{\fin}$ is isomorphic to the restricted tensor product of representations $\bigotimes_{v \in \VV_{\fin}} \pi_v$ with respect to some choice of nonzero spherical vectors at almost all places.
\end{proposition}

At the archimedean places $v \in \VV_\infty$, we have Lie groups $\mu_n \Into \tilde G_v \Onto G_v$; for each such place, choose a maximal compact subgroup $K_v$ with pullback $\tilde K_v \subset \tilde G_v$. Let $\Lie{g}_v$ denote the complexified Lie algebra of $G_v$.  The admissible $\epsilon$-genuine representations of $\tilde G_v$ (i.e. admissible $\epsilon$-genuine $(\Lie{g}_v, \tilde K_v)$-modules) may be identified with the admissible modules over $\hecke_{v,\epsilon} = \hecke_\epsilon(\tilde G_v)$ -- the algebra of $\epsilon$-genuine, left and right $\tilde K_v$-finite distributions on $\tilde G_v$ with support in $\tilde K_v$.  Define $\hecke_{\infty, \epsilon} = \bigotimes_{v \in \VV_\infty} \hecke_{v, \epsilon}$.

An $\epsilon$-genuine \defined{admissible representation} of $\tilde G_\AA$ is an admissible module over the Hecke algebra $\hecke_{\AA, \epsilon} \defeq \hecke_{\fin, \epsilon} \otimes \hecke_{\infty, \epsilon}$.

A \defined{unitary representation} of $\tilde G_\AA$ will mean a continuous representation of $\tilde G_\AA$ by unitary operators on a Hilbert space.  As in \cite[\S 4]{MWToriNew}, following \cite[Theorem 4]{Flath} and Moore \cite[Lemma 6.3]{Moore}, the unitary representations of $\tilde G_\AA$ can be factored as Hilbert space tensor products of unitary representations at all places (spherical almost everywhere).  We summarize the results in the following theorem.
\begin{thm}
The restricted tensor product gives a bijection between two sets:
\begin{itemize}
\item
the set of families $([\pi_v])_{v \in \VV}$, where $[\pi_v] \in \Irr_\epsilon(\alg{\tilde G} / F_v)$ for all $v \in \VV$ and $[\pi_v] \in \Irr_\epsilon(\alg{\tilde G} / \OO_v)$ for almost all $v \in \VV_{\fin}$;
\item
the set of irreducible admissible $\epsilon$-genuine representations of $\tilde G_\AA$.
\end{itemize}
Similarly, in the unitary setting, we find a bijection between two sets:
\begin{itemize}
\item
the set of families $([\pi_v])_{v \in \VV}$, where $[\pi_v] \in \Irr_\epsilon^{\uni}(\alg{\tilde G} / F_v)$ for all $v \in \VV$ and $[\pi_v] \in \Irr_\epsilon^{\uni}(\alg{\tilde G} / \OO_v)$ for almost all $v \in \VV_{\fin}$;
\item
the set of irreducible admissible $\epsilon$-genuine unitary representations of $\tilde G_\AA$.
\end{itemize}
\end{thm}

\subsection{Automorphic representations}

\begin{definition}
An \defined{$\epsilon$-genuine automorphic central character} ($\epsilon$-genuine \defined{ACC}, for short) for $\tilde G_\AA$ is a continuous homomorphism $\chi \From \tilde Z_\AA \To \CC^\times$ such that $\chi$ is trivial on $G_F \cap \tilde Z_\AA$ and the restriction of $\chi$ to $\mu_n$ coincides with $\epsilon$.
\end{definition}

\begin{proposition}
Every $\epsilon$-genuine ACC $\chi \From \tilde Z_\AA \To \CC^\times$ extends uniquely to a continuous character $G_F \tilde Z_\AA \To \CC^\times$ which is trivial on $G_F$.
\end{proposition}
\proof
Multiplication gives a surjective continuous homomorphism,
$$G_F \times \tilde Z_\AA \Onto G_F \tilde Z_\AA.$$
Since $G_F \times \tilde Z_\AA$ is a $\sigma$-compact locally compact topological group, the open mapping theorem implies that the multiplication map above is an open map.  Hence $\tilde Z_\AA$ (the image of the open subgroup $\{ 1 \} \times \tilde Z_\AA$ under multiplication) is an open subgroup of $G_F \tilde Z_\AA$.  Thus the inclusion $Z_\AA \Into G_F Z_\AA$ induces a topological isomorphism of locally compact abelian groups,
$$\frac{\tilde Z_\AA}{G_F \cap \tilde Z_\AA} \xrightarrow{\sim} \frac{G_F \cdot \tilde Z_\AA}{G_F}.$$
The result follows immediately.
\qed

The following definition is lifted with slight adaptation from \cite[\S 4.2]{BorelJacquet} (and generalizes \cite[Definition 4.9]{MWToriNew}
\begin{definition}
A function $f \From \tilde G_\AA \To \CC$ is an $\epsilon$-genuine \defined{automorphic form} if it satisfies the following conditions.
\begin{enumerate}
\item
$f(\zeta \gamma \tilde g) = \epsilon(\zeta) \cdot f(\tilde g)$ for all $\gamma \in G_F$, $\zeta \in \mu_n$, and $\tilde g \in \tilde G_\AA$.
\item
Pulling back, $f$ is locally constant on $\tilde G_{\AA_{\fin}}$ and smooth on $\tilde G_v$ for all $v \in \VV_\infty$.
\item
There is a simple element $\xi = \xi_{\fin} \otimes \xi_\infty \in \hecke_{\fin, \epsilon} \otimes \hecke_{\infty, \epsilon}$ such that $f \ast \xi = f$.
\item
For every $y \in \tilde G_{\AA}$ and every place $v \in \VV_\infty$, the function $\tilde G_v \To \CC$, $x \mapsto f(x \cdot y)$ has moderate growth.
\item
The span of $\{f \ast \xi : \xi \in \hecke_{\fin} \otimes \hecke_\infty \}$ is an admissible representation of $\tilde G_\AA$.
\end{enumerate}
The space of $\epsilon$-genuine automorphic forms will be written $\AF_\epsilon(\tilde G_\AA)$.  
\end{definition}

\begin{definition}
An $\epsilon$-genuine \defined{automorphic representation} of $\tilde G_\AA$ is an irreducible admissible subquotient of $\AF_\epsilon(\tilde G_\AA)$.  Define $\Irr_\epsilon(\alg{\tilde G} / F)$ to be the set of equivalence classes (as admissible representations of $\tilde G_\AA$) of $\epsilon$-genuine automorphic representations of $\tilde G_\AA$.
\end{definition}

By Schur's Lemma, every $\epsilon$-genuine automorphic representation of $\tilde G_\AA$ has an $\epsilon$-genuine automorphic central character.  If $\chi$ is an $\epsilon$-genuine ACC, then we also write $\chi$ for its extension to $G_F \cdot \tilde Z_\AA$.  If $\chi$ is {\em unitary}, define $L_\chi^2(G_F \backslash \tilde G_\AA)$ to be the unitary representation of $\tilde G_\AA$, induced (in the sense of Blattner \cite{Bla}) from the unitary character $\chi$ of $G_F \cdot \tilde Z_\AA$.  Let $L_\chi^{2, \cusp}(G_F \backslash \tilde G_\AA)$ be the cuspidal subspace, discussed for example in the work of Moeglin and Waldspurger \cite{MoeglinWaldspurger}.  

An \defined{$L^2$ automorphic representation} of $\alg{\tilde G}$ with central character $\chi$ is an irreducible unitary subrepresentation of $L_\chi^2(G_F \backslash \tilde G_\AA)$.  An automorphic representation is called \defined{cuspidal} if it lies within $L_\chi^{2, \cusp}(G_F \backslash \tilde G_\AA)$.  Define $\Irr_\chi^{\uni}(\alg{\tilde G} / F)$ to be the set of equivalence classes (as unitary representations of $\tilde G_\AA$) of $L^2$ automorphic representations of $\alg{\tilde G}$ with central character $\chi$.  Let $\Irr_\chi^{\cusp}(\alg{\tilde G} / F)$ be the subset consisting of cuspidal automorphic representations.

From the work of Harish-Chandra, Gelfand, and Piatetski-Shapiro, see also \cite[\S I.2.18]{MoeglinWaldspurger}, the cuspidal automorphic spectrum decomposes as nicely as possible.
\begin{proposition}
Let $\chi$ be a unitary $\epsilon$-genuine ACC.  Then $L_\chi^{2, \cusp}(G_F \backslash \tilde G_\AA)$ decomposes discretely into Hilbert space direct sum of $\epsilon$-genuine unitary irreps, each with finite multiplicity.
$$L_\chi^{2,\cusp}(G_F \backslash \tilde G_\AA) \isom \widehat \bigoplus_{\pi \in \Irr_\epsilon^{\cusp}(\alg{\tilde G} / F) } m_\pi \cdot [\pi].$$
\end{proposition}

\part{Parameters for split tori}

\section{A tale of two functors}

Here we consider the simplest case in which we can give a one-to-one parameterization of irreducible genuine representations by Weil parameters:  the case of {\em sharp} covers of split tori over local fields.  This case underlies more difficult cases to follow.  

Here $F$ will denote a local field, and $\alg{T}$ a split torus over $F$.  Fix a separable closure $\bar F / F$.  Let $Y$ be the cocharacter lattice of $\alg{T}$, and $X$ the character lattice.  We identify $T = \alg{T}(F)$ with $Y \otimes F^\times$; for $y \in Y$ and $u \in F^\times$, we write $u^y$ for $y \otimes u \in T$.  Similarly, we write $T^\vee = X \otimes \CC^\times$ for the complex dual torus; for $x \in X$ and $z \in \CC^\times$, we write $z^x$ instead of $x \otimes z$.

Let $n$ be a positive integer for which $\mu_n = \{ \zeta \in F^\times : \zeta^n = 1 \}$ has $n$ elements.  Fix an injective character $\epsilon \From \mu_n \Into \CC^\times$.

Let $\Cat{Cov}_n(\alg{T})$ be the category of covers $\alg{\tilde T} = (\alg{T}', n)$ of $\alg{T}$ of degree $n$ over $F$, as in Definition \ref{DefCover}.  Each such cover yields a quadratic form $Q \From Y \To \ZZ$ (its first Brylinski-Deligne invariant) and hence a lattice $Y_{Q,n} \subset Y$.  Let $\Cat{Cov}_n^\sharp(\alg{T})$ be the full subcategory whose objects are those extensions $(\alg{T}', n)$ for which $Y = Y_{Q,n}$, and for which $Q$ is even-valued if $n$ is odd.  These are the \defined{sharp covers} of degree $n$, and they form a Picard groupoid.

If $\alpha \From \alg{S} \xrightarrow{\sim} \alg{T}$ is an isomorphism of split tori, then pulling back gives an equivalence of Picard groupoids $\alpha^\ast \From \Cat{Cov}_n^\sharp(\alg{T}) \To \Cat{Cov}_n^\sharp(\alg{S})$.

\subsection{The functor of genuine characters}

Suppose that $\alg{\tilde T} \in \Cat{Cov}_n^\sharp(\alg{T})$ is a sharp degree $n$ cover of $\alg{T}$ over $F$.  Let $Q \From Y \To \ZZ$ be its first Brylinski-Deligne invariant.  The sharp cover $\alg{\tilde T}$ defines a short exact sequences of locally compact \textbf{abelian} groups (the construction of \cite[\S 10.3]{B-D}, abelian by \cite[Proposition 4.1]{MWToriOld}),
$$\mu_n \Into \tilde T \Onto T.$$
Define $\Irr(\alg{T}) = \Hom(T, \CC^\times)$ (continuous homomorphisms).  Recall that 
$$\Irr_\epsilon(\alg{\tilde T}) = \Hom_\epsilon(\tilde T, \CC^\times),$$ 
the set of continuous $\epsilon$-genuine homomorphisms (those that restrict to $\epsilon$ on $\mu_n$).  Then $\Irr_\epsilon(\alg{\tilde T})$ is a $\Irr(\alg{T})$-torsor by twisting.

A morphism $f \From \alg{\tilde T}_1 \To \alg{\tilde T}_2$ in $\Cat{Cov}_n^\sharp(\alg{T})$ yields an isomorphism $f \From \tilde T_1 \To \tilde T_2$ (an isomorphism of extensions of $T$ by $\mu_n$).  The map $\chi \mapsto \chi \circ f$ gives a morphism of $\Irr(\alg{T})$-torsors from $\Irr_\epsilon(\alg{\tilde T}_2)$ to $\Irr_\epsilon(\alg{\tilde T}_1)$.

\begin{proposition}
The construction above defines a contravariant additive functor of Picard groupoids,
$$\Irr_\epsilon \From \Cat{Cov}_n^\sharp(\alg{T}) \To \Cat{Tors}(\Irr(\alg{T})).$$
\end{proposition}
\proof
The only thing to check is compatibility with Baer sums.  If $\alg{\tilde T} = \alg{\tilde T}_1 \Baer \alg{\tilde T}_2$ is the Baer sum of two objects of $\Cat{Cov}_n^\sharp(\alg{T})$, then the resulting extension $\mu_n \Into \tilde T \Onto T$ is the Baer sum of $\tilde T_{1}$ and $\tilde T_{2}$.  If $t \in T$, then the elements of $\tilde T$ lying over $t$ are given by pairs $(\tilde t_1, \tilde t_2) \in \tilde T_{1} \times \tilde T_{2}$ lying over $t$, modulo the relation $(\tilde t_1 \zeta, \tilde t_2) = (\tilde t_1, \zeta \tilde t_2)$ for all $\zeta \in \mu_n$.  

If $\chi_1 \in \Irr_\epsilon(\alg{\tilde T}_1)$ and $\chi_2 \in \Irr_\epsilon(\alg{\tilde T}_2)$, then define
$$\chi(\tilde t_1, \tilde t_2) = \chi_1(\tilde t_1) \cdot \chi_2(\tilde t_2).$$
Since both $\chi_1$ and $\chi_2$ are $\epsilon$-genuine, $\chi(\tilde t_1 \zeta, \tilde t_2) = \chi(\tilde t_1, \zeta \tilde t_2)$, and so $\chi$ descends to an $\epsilon$-genuine character of $\tilde T$.  The map $(\chi_1, \chi_2) \mapsto \chi$ defines a function $\Irr_\epsilon(\alg{\tilde T}_1) \times \Irr_\epsilon(\alg{\tilde T}_2) \To \Irr_\epsilon(\alg{\tilde T})$.  Checking compatibility with twisting by $\Irr(\alg{T})$, this defines an isomorphism of $\Irr(\alg{T})$-torsors,
$$\Irr_\epsilon(\alg{\tilde T}_1) \Baer \Irr_\epsilon(\alg{\tilde T}_2) \xrightarrow{\sim} \Irr_\epsilon(\alg{\tilde T}_1 \Baer \alg{\tilde T}_2).$$
This isomorphism defines an additive structure on the functor $\Irr_\epsilon$.
\qed

If $\alpha \From \alg{S} \xrightarrow{\sim} \alg{T}$ is an isomorphism of split tori, and $\alg{\tilde T} \in \Cat{Cov}_n^\sharp(\alg{T})$, and $\alg{\tilde S} = \alpha^\ast \alg{\tilde T}$, then we find a commutative diagram of locally compact abelian groups with exact rows.
$$\begin{tikzcd}
\mu_n \inarrow{r} \arrow{d}{=} & \tilde S \onarrow{r} \arrow{d}{\tilde \alpha} & S \arrow{d}{\alpha} \\
\mu_n \inarrow{r} & \tilde T \onarrow{r} & T
\end{tikzcd}
$$
Pulling back characters gives a group isomorphism $\alpha^\ast \From \Irr(\alg{T}) \To \Irr(\alg{S})$.  Pushing out torsors via $\alpha^\ast$ gives an equivalence of Picard groupoids, which we also write $\alpha^\ast \From \Cat{Tors}(\Irr(\alg{T})) \To \Cat{Tors}(\Irr(\alg{S}))$.  Pulling back genuine characters gives a bijection $\tilde \alpha^\ast \From \Irr_\epsilon(\alg{\tilde T}) \To \Irr_\epsilon(\alg{\tilde S})$.  Allowing covers to vary, these bijections define a natural isomorphism $\tilde \alpha^\ast \From \alpha^\ast \circ \Irr_\epsilon \Rightarrow \Irr_\epsilon \circ \alpha^\ast$ making the following diagram 2-commute.
\begin{equation}
\label{PullbackIrr}
\begin{tikzcd}
\Cat{Cov}_n^\sharp(\alg{T}) \arrow{r}{\Irr_\epsilon} \arrow{d}{\alpha^\ast} & \Cat{Tors}(\Irr(\alg{T})) \arrow{d}{\alpha^\ast} \\
\Cat{Cov}_n^\sharp(\alg{S}) \arrow{r}{\Irr_\epsilon} & \Cat{Tors}(\Irr(\alg{S}))
\end{tikzcd}
\end{equation}

\subsection{The functor of Weil parameters}

Since $\alg{T}$ is a split torus over $F$, and $Y = Y_{Q,n}$, the L-group of $\alg{T}$ (see Section \ref{LGroupCover}) is a central extension of locally compact groups,
$$T^\vee \Into {}^\EL \tilde T \Onto \Gal_F,$$
with $T^\vee = X \otimes \CC^\times$.  Since $\alg{T}$ is split, the L-group is well-defined up to unique isomorphism, without specifying a base point in the associated gerbe.

Define $\WP(\alg{T}) = \Hom(\Weil_F, T^\vee)$, the {\em abelian group} of Weil parameters for $\alg{T}$.  Define $\WP_\epsilon(\alg{\tilde T})$ to be the set of Weil parameters from $\Weil_F$ to ${}^\EL \tilde T$.    Then $\WP_\epsilon(\alg{\tilde T})$ is naturally a torsor for $\WP(\alg{T})$.

A morphism $f \From \alg{\tilde T}_1 \To \alg{\tilde T}_2$ in $\Cat{Cov}_n^\sharp(\alg{T})$ is automatically a well-aligned isomorphism.  It induces an equivalence of L-groups, by Section \ref{WAFLGroup}.
$$\begin{tikzcd}
T^\vee \inarrow{r} \arrow{d}{=} & {}^\EL \tilde T_2 \arrow{d}{{}^\EL f} \onarrow{r} & \Gal_F \arrow{d}{=} \\
T^\vee \inarrow{r} & {}^\EL \tilde T_1 \onarrow{r} & \Gal_F
\end{tikzcd}$$
Composition with ${}^\EL f$ defines an isomorphism of $\WP(\alg{T})$-torsors, from ${}^\EL \tilde T_2$-valued parameters to ${}^\EL \tilde T_1$-valued parameters.

\begin{proposition}
This defines a contravariant additive functor of Picard groupoids,
$$\WP_\epsilon \From \Cat{Cov}_n^\sharp(\alg{T}) \To \Cat{Tors}(\WP(\alg{T})).$$
\end{proposition}
\proof
We have described such a functor $\WP_\epsilon$, at the level of objects and morphisms.  To check that $\WP_\epsilon$ respects the additive structure, we must see how the construction of the L-group behaves for the Baer sum of two covers.  So consider two objects $\alg{\tilde T}_1$ and $\alg{\tilde T}_2$ of $\Cat{Cov}_n^\sharp(\alg{T})$, and write $\alg{\tilde T}$ for their Baer sum.  Write $\sheaf{D}_1, \sheaf{D}_2, \sheaf{D}$ for the second Brylinski-Deligne invariants of $\alg{\tilde T}_1$, $\alg{\tilde T}_2$ and $\alg{\tilde T}$, respectively.  The construction of the second Brylinski-Deligne invariant respects Baer sums, and so we find extensions of sheaves on $F_{\et}$,
$$\sheaf{G}_m \Into \sheaf{D}_{i} \Onto Y, \text{ for } i = 1,2,$$
and $\sheaf{D} = \sheaf{D}_1 \Baer \sheaf{D}_2$.  The $\shom(Y, \sheaf{G}_m)$-torsor of splittings $\sspl(\sheaf{D})$ is naturally identified with the contraction of the torsors $\sspl(\sheaf{D}_1)$ and $\sspl(\sheaf{D}_2)$.  This gives an equivalence from the gerbe of $n^{\th}$ roots of $\sspl(\sheaf{D})$ to the contraction of gerbes,
$$\sqrt[n]{\sspl(\sheaf{D})} \xrightarrow{\sim} \sqrt[n]{\sspl(\sheaf{D}_1)} \Baer \sqrt[n]{\sspl(\sheaf{D}_2)}.$$
If $z_1, z_2$ are base points (over some finite Galois $F' / F$), objects of $\sqrt[n]{\sspl(\sheaf{D}_1)}$ and $\sqrt[n]{\sspl(\sheaf{D}_2)}$ respectively, then we find a base point $z = z_1 \wedge z_2$ of $\sqrt[n]{\sspl(\sheaf{D})}$.  

This defines a map, for all $\gamma \in \Gal_F$
$$\Hom(z_1, {}^\gamma z_1) \times \Hom(z_2, {}^\gamma z_2) \To \Hom(z, {}^\gamma z).$$
Assembling these, we find an isomorphism of extensions of $\Gal_F$ by $T^\vee$,
\begin{equation}
\label{BSp1}
\pi_1(\gerb{E}_\epsilon(\alg{\tilde T}), \bar z) \xrightarrow{\sim} \pi_1(\gerb{E}_\epsilon(\alg{\tilde T}_1), \bar z_1) \Baer \pi_1(\gerb{E}_\epsilon(\alg{\tilde T}_2), \bar z_2).
\end{equation}

Let $Q, Q_1, Q_2$ be the first Brylinski-Deligne invariants of $\alg{\tilde T}$, $\alg{\tilde T}_1$, $\alg{\tilde T}_2$, respectively.  Then $Q = Q_1 + Q_2$.  This gives an isomorphism of extensions of $\Gal_F$ by $T^\vee$,
\begin{equation}
\label{BSp2}
(\tau_Q)_\ast \mGal_F \xrightarrow{\sim} (\tau_{Q_1})_\ast \mGal_F \Baer (\tau_{Q_2})_\ast \mGal_F.
\end{equation}

Taking the Baer sums of \eqref{BSp1} and \eqref{BSp2}, we find an L-equivalence of L-groups,
$${}^\EL \tilde T \xrightarrow{\sim} {}^\EL \tilde T_1 \Baer {}^\EL \tilde T_2,$$
This L-equivalence defines an isomorphism of $\WP(\alg{T})$-torsors, from the ${}^\EL \tilde T$-valued Weil parameters $\WP_\epsilon(\alg{\tilde T}) $ to the contraction of $\WP(\alg{T})$-torsors $\WP_\epsilon(\alg{\tilde T}_1) \Baer \WP_\epsilon(\alg{\tilde T}_2)$.  This defines the additive structure on the functor $\WP_\epsilon$ as required.
\qed

If $\alpha \From \alg{S} \xrightarrow{\sim} \alg{T}$ is an isomorphism of split tori, and $\alg{\tilde T} \in \Cat{Cov}_n^\sharp(\alg{T})$, and $\alg{\tilde S} = \alpha^\ast \alg{\tilde T}$, then we find a well-aligned isomorphism,
$$\begin{tikzcd}
\alg{K}_2 \inarrow{r} \arrow{d}{=} & \alg{S}' \onarrow{r} \arrow{d}{\alpha'} & \alg{S} \arrow{d}{\alpha} \\
\alg{K}_2 \inarrow{r} & \alg{T}' \onarrow{r} & \alg{T}
\end{tikzcd}
$$
This defines an equivalence of L-groups,
$$\begin{tikzcd}
T^\vee \inarrow{r} \arrow{d}{\alpha^\vee} & {}^\EL \tilde T \onarrow{r} \arrow{d}{{}^\EL \alpha} & \Gal_F \arrow{d}{=} \\
S^\vee \inarrow{r} & {}^\EL \tilde S \onarrow{r} & \Gal_F
\end{tikzcd}$$

The isomorphism $\alpha^\vee \From T^\vee \To S^\vee$ gives a group isomorphism $\alpha^\vee \From \WP(\alg{T}) \To \WP(\alg{S})$.  Composition with ${}^\EL \alpha$ gives a bijection ${}^\EL \alpha \From \WP_\epsilon(\alg{\tilde T}) \To \WP_\epsilon(\alg{\tilde S})$.  Allowing covers to vary, these bijections define a natural isomorphism ${}^\EL \alpha \From \alpha^\vee \circ \WP_\epsilon \Rightarrow \WP_\epsilon \circ \alpha^\vee$ making the following diagram 2-commute.
\begin{equation}
\label{PullbackWP}
\begin{tikzcd}
\Cat{Cov}_n^\sharp(\alg{T}) \arrow{r}{\WP_\epsilon} \arrow{d}{\alpha^\ast} & \Cat{Tors}(\WP(\alg{T})) \arrow{d}{\alpha^\vee} \\
\Cat{Cov}_n^\sharp(\alg{S}) \arrow{r}{\WP_\epsilon} & \Cat{Tors}(\WP(\alg{S}))
\end{tikzcd}
\end{equation}

\subsection{The goal}

Class field theory gives a group isomorphism (the simplest case of \cite{LanglandsAbelian}),
$$\Lang \From \Irr(\alg{T}) \xrightarrow{\sim} \WP(\alg{T}).$$
Indeed, if $\rec \From \Weil_F^{\ab} \xrightarrow{\sim} F^\times$ is the isomorphism of class field theory (normalized so that a geometric Frobenius corresponds to a uniformizing parameter, in the nonarchimedean case), then $\Lang$ is the isomorphism given by composing the sequence of isomorphisms,
\begin{align*}
\Lang \From \Irr(\alg{T}) &= \Hom(T, \CC^\times) = \Hom(Y \otimes F^\times, \CC^\times) \\
& \xrightarrow{\can_Y}  \Hom(F^\times, X \otimes \CC^\times) \\
& \xrightarrow{\rec^\ast}  \Hom(\Weil_F, X \otimes \CC^\times) = \Hom(\Weil_F, T^\vee) = \WP(\alg{T}). 
\end{align*}
Here and later, $\can_Y \From \Hom(Y \otimes F^\times, \CC^\times) \xrightarrow{\sim} \Hom(F^\times, X \otimes \CC^\times)$ is the natural isomorphism arising from the duality of $X$ and $Y$.

The goal for the next few sections is the construction of a natural isomorphism $\Lang_\epsilon \From \Lang \circ \Irr_\epsilon \Rightarrow \WP_\epsilon$ of additive functors, which makes the following diagram of Picard groupoids and additive functors 2-commute.
\begin{equation}
\label{GoalParam}
\begin{tikzcd}
\phantom{a} &  \phantom{a} & \Cat{Tors}(\Irr(\alg{T}))  \arrow{dd}{\Lang} \\
\Cat{Cov}_n^\sharp(\alg{T}) \arrow{urr}{\Irr_\epsilon} \arrow{drr}[swap]{\WP_\epsilon} & \phantom{1000} \Downarrow \Lang_\epsilon?  &  \phantom{a} \\
\phantom{a} & \phantom{a} & \Cat{Tors}(\WP(\alg{T}))
\end{tikzcd}
\end{equation}
Here we abuse notation slightly, and write $\Lang$ not only for the isomorphism $\Irr(\alg{T}) \xrightarrow{\sim} \WP(\alg{T})$, but also for the equivalence of Picard groupoids given by pushing out via $\Lang$, $\Cat{Tors}(\Irr(\alg{T})) \xrightarrow{\sim} \Cat{Tors}(\Irr(\WP(\alg{T})))$.

This pushout isomorphism $\Lang$ is naturally isomorphic to the functor which sends an $\Irr(\alg{T})$-torsor $V$ to the $\WP(\alg{T})$-torsor with underlying set $V$ and action given by
$$\phi \ast v = \Lang^{-1}(\phi) \ast v \text{ for all } v \in V, \phi \in \WP(\alg{T}).$$
In this way, viewing $\Lang$ as the identity map on underlying sets, the natural isomorphism $\Lang_\epsilon$ will give a bijection from $\Irr_\epsilon(\alg{\tilde T})$ to $\WP_\epsilon(\alg{\tilde T})$, for any cover $\alg{\tilde T} \in \Cat{Cov}_n^\sharp(\alg{T})$.

We will also prove that $\Lang_\epsilon$ is compatible with pullbacks, in the following sense.  If $\alpha \From \alg{S} \To \alg{T}$ is an isomorphism of split tori, and and $\alg{\tilde T} \in \Cat{Cov}_n^\sharp(\alg{T})$, and $\alg{\tilde S} = \alpha^\ast \alg{\tilde T}$, then $\Lang$ is compatible with pullbacks in the sense that the following diagram of groups and isomorphisms commutes.
$$\begin{tikzcd}
\Irr(\alg{T}) \arrow{r}{\Lang} \arrow{d}{\alpha^\ast} & \WP(\alg{T}) \arrow{d}{\alpha^\ast} \\
\Irr(\alg{S}) \arrow{r}{\Lang} & \WP(\alg{S})
\end{tikzcd}$$

We will find a commutative diagram of additive functors and natural isomorphisms.
\begin{equation}
\label{PullbackLang}
\begin{tikzcd}[column sep = 6em]
\alpha^\ast \circ \Lang \circ \Irr_\epsilon \arrow[Rightarrow]{r}{\Lang_\epsilon} \arrow[Rightarrow]{d} & \alpha^\ast \circ \WP_\epsilon \arrow[Rightarrow]{dd}{{}^\EL \alpha} \\
\Lang \circ \alpha^\ast \circ \Irr_\epsilon \arrow[Rightarrow]{d}{\tilde \alpha^\ast} & \\
\Lang \circ \Irr_\epsilon \circ \alpha^\ast \arrow[Rightarrow]{r}{\Lang_\epsilon} & \WP_\epsilon \circ \alpha^\ast
\end{tikzcd}
\end{equation}
In other words, given an $\epsilon$-genuine character of a cover $\tilde T$, we may proceed in two directions:  first, we may parameterize it by a ${}^\EL \tilde T$-valued Weil parameter, and then push that parameter to find a ${}^\EL \tilde S$-valued Weil parameter.  Second, we may pull back the character to find an $\epsilon$-genuine character of the cover $\tilde S$, and then parameterize by an ${}^\EL \tilde S$-valued Weil parameter.  The resulting ${}^\EL \tilde S$-valued Weil parameters will be the same, according to the diagram above.

Thus $\Lang_\epsilon \From \Lang \circ \Irr_\epsilon \xRightarrow{\sim} \WP_\epsilon$ will parameterize $\epsilon$-genuine characters by Weil parameters, for sharp covers of split tori.  Formulating $\Lang_\epsilon$ as a natural isomorphism of additive functors means that our parameterization of $\epsilon$-genuine representations will be compatible with morphisms of covers, Baer sums, twisting, and pullbacks (via isomorphisms of tori).  The proof requires some intermediate steps, and we will proceed through two more functors and three natural isomorphisms, according to the roadmap below.
\begin{equation}
\label{FunctorRoadmap}
\begin{tikzcd}
\Lang \circ \Irr_\epsilon \arrow[Rightarrow]{r}{\Del}&  \rec^\ast \circ \Irr_\epsilon^\DEE  \arrow[Leftrightarrow]{r}{\P} & \rec^\ast \circ \WP_\epsilon^{\A} & \WP_\epsilon \arrow[Rightarrow]{l}[swap]{\AB}.
\end{tikzcd}
\end{equation}

\section{Deligne's construction}

Here we utilize an idea of Deligne (personal communication) which relates the genuine characters of one extension to splittings of a ``dual'' extension.  This provides a natural isomorphism from the functor $\can_Y \circ \Irr_\epsilon$ to a functor $\Irr_\epsilon^\DEE$ which will be easier to work with later.  We fully credit Deligne with the elegant idea of this construction, and we take responsibility for any sloppiness in exposition.

\subsection{Dual extensions}

We give a very general construction here; suppose that $J$ and $C$ are locally compact abelian groups (with composition written multiplicatively), and $X$ and $Y$ are finite-rank free abelian groups in duality as before.  There is a group isomorphism,
$$\Ext(Y \otimes J, C) \isom \Ext(J, X \otimes C),$$
where $\Ext$ denotes the group of isomorphism classes of extensions (of locally compact abelian groups).  This follows from the fact that $\Ext$ distributes over finite direct sums.  We can beef this up to an equivalence of Picard groupoids as follows:  Let $\Cat{Ext}(Y \otimes J, C)$ be the Picard groupoid of extensions of locally compact abelian groups, and $E$ an object therein,
$$C \Into E \Onto Y \otimes J.$$
Tensoring with $X$ yields an extension of locally compact abelian groups,
$$X \otimes C \Into X \otimes E \Onto X \otimes Y \otimes J.$$
Pulling back via the canonical inclusion $\iota \From \ZZ \Into X \otimes Y$ yields an extension
$$X \otimes C \Into \iota^\ast(X \otimes E) \Onto J.$$
\begin{proposition}
This defines an equivalence of Picard groupoids,
$$\DEE_Y = \iota^\ast(X \otimes \bullet) \From \Cat{Ext}(Y \otimes J, C) \To \Cat{Ext}(J, X \otimes C).$$
\end{proposition}
\proof
We begin by checking that the construction $\DEE_Y$ gives a bijection on isomorphism classes.  Choosing a $\ZZ$-basis $\{ y_1, \ldots, y_r \}$ of $Y$, and dual basis $\{ x_1, \ldots, x_r \}$ of $X$ yields  group isomorphisms $\Ext(Y \otimes J, C) \isom \bigoplus_{i=1}^r \Ext(J, C)$ and $\Ext(J, X \otimes C) \isom \bigoplus_{i=1}^r \Ext(J, C)$.  

The canonical map $\iota \From \ZZ \Into X \otimes Y$ satisfies $\iota(1) = \sum_{i=1}^r x_i \otimes y_i$.  It follows that the construction gives a commutative diagram.
$$\begin{tikzcd}
\Ext(Y \otimes J, C) \arrow{r}{\DEE_Y} \arrow{d}{\sim} & \Ext(J, X \otimes C) \arrow{d}{\sim} \\
\bigoplus_{i=1}^r \Ext(J, C) \arrow{r}{\bigoplus \DEE_\ZZ} & \bigoplus_{i=1}^r \Ext(J, C).
\end{tikzcd}$$
But the functor $\DEE_\ZZ$ is (naturally isomorphic to) the identity functor on the category $\Cat{Ext}(J, C)$.  Thus the bottom row of the commutative diagram is an isomorphism, and so the top row is an isomorphism.

To demonstrate that $\DEE_Y$ is an equivalence of categories, we trace through an automorphism of an extension $C \Into E \Onto Y \otimes J$.  Such an automorphism is given by an element $a \in \Hom(Y \otimes J, C)$ (a continuous homomorphism).  With respect to the chosen basis, $a = \prod \xi_i^{x_i}$ for some family of homomorphisms $\xi_i \From J \To C$.  In other words, $a(j^y) = \prod_i \xi_i(j)^{\langle x_i, y \rangle}$.  

If $e \in E$ lies over $j^y \in Y \otimes J$, then the automorphism $\eta_a \in \Aut(E)$ is given by
$$\eta_a(e) = e \cdot \prod_i \xi_i(j)^{\langle x_i, y \rangle}.$$
The resulting automorphism of $X \otimes C \Into X \otimes E \Onto X \otimes Y \otimes J$ is characterized by
$$\eta_{a,X}(x \otimes e) = (x \otimes e) \cdot \prod_i \xi_i(j)^{\langle x_i, y \rangle \cdot x},$$
if $x \in X$, and $e \in E$ lies over $j^y \in Y \otimes J$.  Thus in the pullback via $\iota \From \ZZ \Into X \otimes Y$, we find that the automorphism of $\DEE_Y E \defeq \iota^\ast(X \otimes E)$ is given by
$$\DEE_Y \eta_a(\tilde j) = \tilde j \cdot \prod_i \xi_i(j)^{\langle x_i, y_i \rangle \cdot x_i} = \tilde j \cdot \prod_i \xi_i(j)^{x_i},$$
for all $\tilde j \in \DEE_Y E$ lying over $j \in J$.

The automorphism group of $E$ is canonically identified with $\Hom(Y \otimes J, C)$.  The automorphism group of $\DEE_Y E$ is canonically identified with $\Hom(J, X \otimes C)$.  The computations above demonstrate that the functor $\DEE_Y$ coincides on automorphism groups with the canonical isomorphism,
$$\DEE_Y = \can_Y \From \Hom(Y \otimes J, C) \xrightarrow{\sim} \Hom(J, X \otimes C).$$
Therefore $\DEE_Y$ is an equivalence of Picard groupoids.  Compatibility with Baer sums is left to the reader.
\qed

The equivalence $\DEE_Y$ depends functorially on $C$ and $J$ as well.  For example, if $j \From J_1 \To J_2$ is a homomorphism of locally compact abelian groups, then a natural isomorphism makes the square below 2-commute.
\begin{equation}
\label{PullbackDEE}
\begin{tikzcd}
\Cat{Ext}(Y \otimes J_2, C) \arrow{r}{\DEE_Y} \arrow{d}{j^\ast} & \Cat{Ext}(J_2, X \otimes C) \arrow{d}{j^\ast} \\
\Cat{Ext}(Y \otimes J_1, C) \arrow{r}{\DEE_Y} & \Cat{Ext}(J_1, X \otimes C)
\end{tikzcd}
\end{equation}
Here the vertical arrows are additive functors given by pullback via $j$.

The equivalence $\DEE_Y$ also depends functorially on $Y$, at least for isomorphisms.  Let $\alpha \From Y_1 \To Y_2$ be an isomorphism of finite-rank free abelian groups, and write $\alpha^\vee \From X_2 \To X_1$ for the dual isomorphism.  Then a natural isomorphism makes the square below 2-commute.
\begin{equation}
\label{PullbackDEEY}
\begin{tikzcd}
\Cat{Ext}(Y_2 \otimes J, C) \arrow{r}{\DEE_Y} \arrow{d}{\alpha^\ast} & \Cat{Ext}(J, X_2 \otimes C) \arrow{d}{\alpha_\ast^\vee} \\
\Cat{Ext}(Y_1 \otimes J, C) \arrow{r}{\DEE_Y} & \Cat{Ext}(J, X_1 \otimes C)
\end{tikzcd}
\end{equation}

\subsection{Genuine characters and splittings}

Now consider an extension of locally compact abelian groups,
$$\mu_n \Into E \Onto Y \otimes J,$$
Pushing out via $\epsilon \From \mu_n \Into \CC^\times$ yields a commutative diagram
$$\begin{tikzcd}
\mu_n \inarrow{r} \inarrow{d}{\epsilon} &E\onarrow{r} \inarrow{d} & Y \otimes J \arrow{d}{=} \\
\CC^\times \inarrow{r} & \epsilon_\ast E \onarrow{r} & Y \otimes J
\end{tikzcd}$$
Giving an $\epsilon$-genuine character of $E$ is the same, by the universal property of pushouts, as giving a splitting of the bottom row of this diagram.  In this way, we find an isomorphism of $\Hom(Y \otimes J, \CC^\times)$-torsors,
\begin{equation}
\label{GenuineCharacterAsSplitting}
\Hom_\epsilon(E, \CC^\times) \xrightarrow{\sim} \Spl(\epsilon_\ast E).
\end{equation}

Applying the functor $\DEE_Y = \iota^\ast (X \otimes \bullet)$ to $\epsilon_\ast E$ yields an extension
$$X \otimes \CC^\times \Into \DEE_Y \epsilon_\ast E \Onto J.$$

A splitting of $\epsilon_\ast E$ is the same as an isomorphism from $\epsilon_\ast E$ to the trivial extension (i.e., zero object) in the Picard groupoid $\Cat{Ext}(Y \otimes J, \CC^\times)$.  Since $\DEE_Y$ is an equivalence of Picard groupoids, it defines a bijection,
\begin{equation}
\label{SplittingsToSplittings}
\DEE_Y \From \Spl(\epsilon_\ast E) \xrightarrow{\sim} \Spl(\DEE_Y \epsilon_\ast E).
\end{equation}
Assembling \eqref{GenuineCharacterAsSplitting} and \eqref{SplittingsToSplittings}, we have a bijection,
\begin{equation}
\Hom_\epsilon(E, \CC^\times) \xrightarrow{\sim} \Spl( \DEE_Y \epsilon_\ast E),
\end{equation}
which intertwines the $\Hom(Y \otimes J, \CC^\times)$-action on $\Hom_\epsilon(E, \CC^\times)$ with the $\Hom(J, X \otimes \CC^\times)$ action on $\Spl( \DEE_Y \epsilon_\ast E)$, the intertwining passing through the isomorphism $\can_Y$.

\subsection{A natural isomorphism}

If $\mu_n \Into \tilde T \Onto T$ arises from a sharp cover $\alg{\tilde T}$ as before, then define ${}^\DEE \tilde T = \DEE_Y \epsilon_\ast \tilde T$.  This is an extension of locally compact abelian groups,
$$X \otimes \CC^\times \Into {}^\DEE \tilde T \Onto F^\times.$$
Define $\Irr_\epsilon^\DEE(\alg{\tilde T}) = \Spl( {}^\DEE \tilde T)$.  We have described a bijection,
$$\Del \From \Irr_\epsilon(\alg{\tilde T}) = \Hom_\epsilon(\tilde T, \CC^\times) \xrightarrow{\sim} \Spl( {}^\DEE \tilde T) = \Irr_\epsilon^\DEE(\alg{\tilde T}).$$

Define $\Irr^\DEE(\alg{T}) = \Hom(F^\times, X \otimes \CC^\times)$, so $\can_Y \From \Irr(\alg{T}) \To \Irr^\DEE(\alg{T})$ is an isomorphism.  We find that $\Irr_\epsilon^\DEE(\alg{\tilde T})$ is naturally a $\Irr^\DEE(\alg{T})$-torsor.  Moreover, $\Del$ intertwines the $\Irr(\alg{T})$-torsor structure on $\Irr_\epsilon(\alg{\tilde T})$ with the $\Irr^\DEE(\alg{T})$-torsor structure on $\Irr_\epsilon^\DEE(\alg{\tilde T})$, via $\can_Y$.  

Putting everything together, we find the following diagram of Picard groupoids and additive functors, made 2-commutative by $\Del \From \can_Y \circ  \Irr_\epsilon \Rightarrow \Irr_\epsilon^\DEE$.
$$\begin{tikzcd}[column sep = 4em]
\Cat{Cov}_n^\sharp(\alg{T}) \arrow{r} \arrow[bend left=20]{rr}{\Irr_\epsilon} \arrow[bend right=35]{rrd}[swap]{\Irr_\epsilon^\DEE} & \Cat{Ext}(Y \otimes F^\times, \mu_n) \arrow{d}{\DEE_Y \epsilon_\ast} \arrow{r}{\Hom_\epsilon(\bullet, \CC^\times)} & \Cat{Tors}( \Irr(\alg{T}) ) \arrow{d}{\can_Y} \\
\phantom{a} & \Cat{Ext}(F^\times, X \otimes \CC^\times) \arrow{r}{\Spl} & \Cat{Tors}(\Irr^\DEE(\alg{T}))
\end{tikzcd}$$
For reference, we trace an object $\alg{\tilde T} \in \Cat{Cov}_n^\sharp(\alg{T})$ through this diagram.
 $$\begin{tikzcd}[column sep = 7em]
\alg{\tilde T} \arrow[mapsto]{r} \arrow[bend left=15, mapsto]{rr}{\Irr_\epsilon} \arrow[bend right=35, mapsto]{rrd}[swap]{\Irr_\epsilon^\DEE} & \tilde T \arrow[mapsto]{d} \arrow[mapsto]{r} & \Hom_\epsilon(\tilde T, \CC^\times) \\
\phantom{a} & {}^\DEE \tilde T \arrow[mapsto]{r} & \Spl \left( {}^\DEE \tilde T \right)
\end{tikzcd}$$

Compressing the diagram down a bit, we find that $\Del \From \can_Y \circ \Irr_\epsilon \Rightarrow \Irr_\epsilon^\DEE$ is a natural isomorphism of additive functors, making the following diagram of Picard groupoids and additive functors 2-commute.
\begin{equation}
\label{TautologicalParam}
\begin{tikzcd}
\phantom{a} &  \phantom{a} & \Cat{Tors}(\Irr(\alg{T})) \arrow{dd}{\can_Y} \\
\Cat{Cov}_n^\sharp(\alg{T}) \arrow{urr}{\Irr_\epsilon} \arrow{drr}[swap]{\Irr_\epsilon^\DEE} & \phantom{1000} \Downarrow \Del  &  \phantom{a} \\
\phantom{a} & & \Cat{Tors}(\Irr^\DEE(\alg{T}))
\end{tikzcd}
\end{equation}

Composing with $\rec^\ast$, and recalling that $\Lang = \rec^\ast \circ \can_Y$, we have constructed a natural isomorphism
\begin{equation}
\label{LeftStep}
\Del \From \Lang \circ \Irr_\epsilon \Longrightarrow \rec^\ast \circ \Irr_\epsilon^\DEE.
\end{equation}

The natural isomorphism $\Del$ is compatible with pullbacks in the following way.  Suppose that $\alpha \From \alg{S} \To \alg{T}$ is an isomorphism of split tori, $\alg{\tilde T} \in \Cat{Cov}_n^\sharp(\alg{T})$ and $\alg{\tilde S}$ is the pullback of $\alg{\tilde T}$.  The following diagram of additive functors and natural isomorphisms commutes.
\begin{equation}
\label{PullbackDel}
\begin{tikzcd}
\alpha^\ast \circ \Lang \circ \Irr_\epsilon \arrow[Rightarrow]{r}{\Del} \arrow[Rightarrow]{d} & \alpha^\ast \circ \rec^\ast \circ \Irr_\epsilon^\DEE \arrow[Rightarrow]{dd}{{}^\DEE \tilde \alpha} \\
\Lang \circ \alpha^\ast \circ \Irr_\epsilon \arrow[Rightarrow]{d}{\tilde \alpha^\ast} & \\
\Lang \circ \Irr_\epsilon \circ \alpha^\ast \arrow[Rightarrow]{r}{\Del} & \rec^\ast  \circ \Irr_\epsilon^\DEE \circ \alpha^\ast
\end{tikzcd}
\end{equation}
Indeed, $\rec^\ast$ and $\Lang$ are compatible with pullback.  From \eqref{PullbackDEEY}, the construction of the ``dual extension'' is compatible with pullback.  In particular, $\alpha$ induces an isomorphism ${}^\DEE \tilde \alpha \From {}^\DEE \tilde T \To {}^\DEE \tilde S$ lying over $\alpha^\vee \From T^\vee \To S^\vee$.

Given an $\epsilon$-genuine character of $\tilde T$, we may proceed in two directions:  first, one can construct a splitting of ${}^\DEE \tilde T$, and then push via ${}^\DEE \tilde \alpha$ to find a splitting of ${}^\DEE \tilde S$.  Second, we can pull back to find an $\epsilon$-genuine character of $\tilde S$, and use that to construct a splitting of ${}^\DEE \tilde S$.  The above diagram states that the resulting splitting of ${}^\DEE \tilde S$ will be the same.

\section{Weil parameters as splittings}

The natural isomorphism $\Del \From \can_Y \circ \Irr_\epsilon \Rightarrow \Irr_\epsilon^\DEE$ relates genuine characters to splittings of a central extension.  Here we define a natural isomorphism which relates Weil parameters to splittings of another central extension.  

Begin with the L-group $T^\vee \Into {}^\EL \tilde T \Onto \Gal_F$.  Since $H^2(\Gal_F, \CC^\times) = 0$, this extension splits.  While there is no canonical splitting of ${}^\EL \tilde T$, the splitting over the commutator subgroup $[\Gal_F, \Gal_F]$ is canonical; indeed, two splittings of ${}^\EL \tilde T$ differ by an element of $\Hom(\Gal_F, T^\vee)$, and such a homomorphism is trivial on $[\Gal_F, \Gal_F]$ since $T^\vee$ is abelian.  From this observation, we get an extension of {\em abelian} groups,
\begin{equation}
\label{AlmostA}
T^\vee \Into {}^\EL \tilde T / [\Gal_F, \Gal_F] \Onto \Gal_F^{\ab}.
\end{equation}
Let ${}^{\A} \tilde T$ be the pullback of \eqref{AlmostA} via $F^\times \xrightarrow{\rec^{-1}} \Weil_F^{\ab} \Into \Gal_F^{\ab}$,
$$T^\vee \Into {}^{\A} \tilde T \Onto F^\times.$$
($\A$ stands for abelian.)  Define $\WP_\epsilon^{\A}(\alg{\tilde T}) = \Spl({}^{\A} \tilde T)$, viewed as a $\Hom(F^\times, T^\vee)$-torsor.

\begin{proposition}
Every Weil parameter $\phi \From \Weil_F \To {}^\EL \tilde T$ descends uniquely to a homomorphism $\Weil_F^{\ab} \To {}^{\A} \tilde T$.  This defines a bijection,
$$\AB \From \WP_\epsilon(\alg{\tilde T}) \rightarrow \WP_\epsilon^\A(\alg{\tilde T}),$$
which intertwines the $\WP(\alg{T}) = \Hom(\Weil_F, T^\vee)$-torsor structure on the left with the $\Irr^\DEE(\alg{T}) = \Hom(F^\times, T^\vee)$-torsor structure on the right, via the reciprocity isomorphism $\rec^\ast \From \Irr^\DEE(\alg{T}) \xrightarrow{\sim} \WP(\alg{T})$.
\end{proposition}
\proof
Let $\sigma \From \Gal_F \To {}^\EL \tilde T$ be any splitting and write $\sigma' \From \Weil_F \To {}^\EL \tilde T$ for its pullback via the canonical map $\Weil_F \To \Gal_F$.  If $\phi \in \WP_\epsilon(\alg{\tilde T})$ is a Weil parameter, then both $\sigma'$ and $\phi$ are homomorphisms from $\Weil_F$ to ${}^\EL \tilde T$, and there exists a unique function $\tau \From \Weil_F \To T^\vee$ such that $\phi(w) = \sigma(w) \cdot \tau(w), \text{ for all } w \in \Weil_F$.

Since $T^\vee$ is contained in the center of ${}^\EL \tilde T$, we find that $\tau \From \Weil_F \To T^\vee$ is a homomorphism.  Hence $\phi = \sigma \cdot \tau$ sends the commutator subgroup of $\Weil_F$ to the image of $[\Gal_F, \Gal_F]$ in ${}^\EL \tilde T$.  Thus $\phi$ descends uniquely to a homomorphism $\Weil_F^{\ab} \To {}^{\A} \tilde T$.  Pulling back via the reciprocity isomorphism $\rec^{-1} \From F^\times \xrightarrow{\sim} \Weil_F^{\ab}$ yields a homomorphism $\AB(\phi) \From F^\times \To {}^\A \tilde T$ as claimed.  This homomorphism $\AB(\phi)$ is a splitting of ${}^{\A} \tilde T$, and the construction of $\AB(\phi)$ from $\phi$ is compatible with twisting by $\WP(\alg{T})$ throughout.
\qed

The construction of the L-group gives an additive functor of Picard categories, $\Cat{Cov}_n^\sharp(\alg{T}) \To \Cat{CExt}(\Gal_F, T^\vee)$ (the latter category is a Picard groupoid with one isomorphism class, since $H^2(\Gal_F, T^\vee)$ is trivial).  The construction of ${}^{\A} \tilde T$ from the (split) extension ${}^\EL \tilde T$ comes from another additive functor of Picard categories, $\Cat{CExt}(\Gal_F, T^\vee) \To \Cat{Ext}(F^\times, T^\vee)$.  Here $\Cat{Ext}(F^\times, T^\vee)$ is the Picard groupoid of extensions of locally compact abelian groups.  It can be proven that every such extension splits.  Finally we apply the additive functor of splittings $\Spl$, 
$$\Spl \From \Cat{Ext}(F^\times, T^\vee) \To \Cat{Tors}(\Hom(F^\times, T^\vee)).$$  
Recall that $\Irr^\DEE(\alg{T})$ denotes $\Hom(F^\times, T^\vee)$.  As the composition of three additive functors, $\WP_\epsilon^{\A}$ is an additive functor,
$$\WP_\epsilon^{\A} \From \Cat{Cov}_n^\sharp(\alg{T}) \To \Cat{Tors}(\Irr^\DEE(\alg{T})).$$
The bijection $\phi \mapsto \AB(\phi)$, from the set of Weil parameters to $\Spl({}^{\A} \tilde T)$, gives a natural isomorphism of additive functors, $\AB \From \WP_\epsilon \xRightarrow{\sim} \rec^\ast \circ \WP_\epsilon^\A$.
\begin{equation}
\label{ParametersAsSplittings}
\begin{tikzcd}
\phantom{a} &  \phantom{a} & \Cat{Tors}(\WP(\alg{T}))  \\
\Cat{Cov}_n^\sharp(\alg{T}) \arrow{urr}{\WP_\epsilon} \arrow{drr}[swap]{\WP_\epsilon^{\A}} & \phantom{1000} \Downarrow \AB  &  \phantom{a} \\
\phantom{a} & & \Cat{Tors}(\Irr^\DEE(\alg{T})) \arrow{uu}{\rec^\ast}
\end{tikzcd}
\end{equation}

The natural isomorphism $\AB$ is compatible with pullbacks in the following way.  Suppose that $\alpha \From \alg{S} \To \alg{T}$ is an isomorphism of split tori, $\alg{\tilde T} \in \Cat{Cov}_n^\sharp(\alg{T})$ and $\alg{\tilde S}$ is the pullback of $\alg{\tilde T}$.  Then we find the L-morphism ${}^\EL \tilde \alpha \From {}^\EL \tilde T \To {}^\EL \tilde S$.  This descends to a morphism of extensions ${}^\A \tilde \alpha$ as below.
$$\begin{tikzcd}
T^\vee \arrow{d}{\alpha^\vee} \inarrow{r} & {}^\A \tilde T \arrow{d}{ {}^\A \tilde \alpha} \onarrow{r} & F^\times \arrow{d}{=} \\
S^\vee \inarrow{r} & {}^\A \tilde S \onarrow{r} & F^\times
\end{tikzcd}$$
In this way, $\alpha$ yields a bijection ${}^\A \alpha \From \Spl({}^\A \tilde T) \To \Spl({}^\A \tilde S)$.  

The following diagram of additive functors and natural isomorphisms commutes.
\begin{equation}
\label{PullbackAB}
\begin{tikzcd}[row sep = 5em, column sep = 6em]
\alpha^\ast \circ \WP_\epsilon \arrow[Rightarrow]{r}{\AB} \arrow[Rightarrow]{d}{{}^\EL \tilde \alpha} & \alpha^\ast \circ \rec^\ast \circ \WP_\epsilon^\A \arrow[Rightarrow]{d}{{}^\A \tilde \alpha} \\
\WP_\epsilon \circ \alpha^\ast \arrow[Rightarrow]{r}{\AB} & \rec^\ast  \circ \WP_\epsilon^\A \circ \alpha^\ast
\end{tikzcd}
\end{equation}
Given a Weil parameter valued in ${}^\EL \tilde T$, we may proceed in two directions:  first, one can push via ${}^\EL \tilde \alpha$ to find a Weil parameter valued in ${}^\EL \tilde S$, which descends to a splitting of ${}^\A \tilde S$.  Second, one can descend the original parameter to find a splitting of ${}^\A \tilde T$, and then push via ${}^\A \tilde \alpha$ to find a splitting of ${}^\A \tilde S$.  The above diagram states that the resulting splitting of ${}^\A \tilde S$ will be the same.

\section{Parameterization}

With \eqref{LeftStep} and \eqref{ParametersAsSplittings}, we have now constructed additive functors and natural isomorphisms as below:
$$\begin{tikzcd}
\Lang \circ \Irr_\epsilon \arrow[Rightarrow]{r}{\Del} & \rec^\ast \circ \Irr_\epsilon^\DEE & \rec^\ast \circ \WP_\epsilon^{\A}  & \WP_\epsilon \arrow[Rightarrow]{l}[swap]{\AB}
\end{tikzcd}$$
The real work comes here, as we define a natural isomorphism of functors $\Irr_\epsilon^\DEE \Leftrightarrow \WP_\epsilon^{\A}$ linking the middle terms.  For this, we will define isomorphisms $\P(\alg{\tilde T}) \From {}^\DEE \tilde T \xrightarrow{\sim} {}^{\A} \tilde T$, in the category of extensions of $F^\times$ by $T^\vee$.  

A morphism $j \From \alg{\tilde T}_1 \To \alg{\tilde T}_2$ in $\Cat{Cov}_n^\sharp(\alg{T})$ yields two morphisms in $\Cat{Ext}(F^\times, T^\vee)$.
$${}^\DEE j \From {}^\DEE \tilde T_2 \To {}^\DEE \tilde T_1, \quad {}^\A j \From {}^\A \tilde T_2 \To {}^\A \tilde T_1.$$
We will demonstrate that our isomorphisms $\P$ are compatible with morphisms in $\Cat{Cov}_n^\sharp(\alg{T})$, i.e., the following diagram commutes for all $j \From \alg{\tilde T}_1 \To \alg{\tilde T}_2$.
$$\begin{tikzcd}[column sep = 6em]
{}^\DEE \tilde T_2 \arrow{r}{{}^\DEE j} \arrow{d}{\P(\alg{\tilde T}_2)} & {}^\DEE \tilde T_1 \arrow{d}{\P(\alg{\tilde T}_1)} \\
{}^\A \tilde T_2 \arrow{r}{{}^\A j} & {}^\A \tilde T_1
\end{tikzcd}$$

This will provide a natural isomorphism of torsors 
$$\Irr_\epsilon^\DEE(\alg{\tilde T}) = \Spl({}^\DEE \tilde T) \xrightarrow{\P(\alg{\tilde T})} \Spl({}^{\A} \tilde T) = \WP_\epsilon^{\A}(\alg{\tilde T}).$$

In the work that follows, we frequently refer to extensions which are ``incarnated'' by bimultiplicative cocycles, so we fix some terminology here.  Suppose that $J$ and $C$ are locally compact abelian groups.  A bimultiplicative cocycle
$$\theta \From J \times J \To C$$
is a continuous function which factors through $J \otimes_\ZZ J \To C$.  (The usual 2-cocycle identity follows from this).  Such a cocycle allows one to define a central extension 
$$C \Into J \times_\theta C \Onto J,$$
where $J \times_\theta C = J \times C$ as sets, and multiplication is given by
$$(j_1, c_1) \cdot (j_2, c_2) = (j_1 j_2, c_1 c_2 \cdot \theta(j_1, j_2) ).$$
We often write $c$ instead of $(1,c)$, viewing $C$ as a subgroup of $J \times_\theta C$.  But we typically write $s(j) = (j,1)$, using the letter $s$ for the section (not often a homomorphism) from $J$ to $J \times_\theta C$.  

For a bimultiplicative cocycle $\theta$, induction gives a formula for powers,
\begin{equation}
\label{gpower}
(j,1)^g = (j^g, 1) \cdot \theta(j,j)^{\frac{g(g-1)}{2}} \text{ for all } j \in J, g \in \ZZ_{>0}.
\end{equation}

\subsection{Incarnated covers}

Suppose that $\alg{\tilde T} = (\alg{T}', n)$ is a sharp cover \defined{incarnated} by $C \in X \otimes X$ as in \cite[\S 3]{B-D}.  We often utilize an ordered basis $\Basis = ( y_1, \ldots, y_r )$ of $Y$ and dual basis $( x_1, \ldots, x_r )$ of $X$ in what follows.  Write $C = \sum_{i,j} c_\Basis^{ij} x_i \otimes x_j$ with respect to this basis.  Thus $\alg{T}' = \alg{T} \times_{\theta_C} \alg{K}_2$ (as sheaves of sets on $F_{\Zar}$), where $\theta_C \From \alg{T} \times \alg{T} \To \alg{K}_2$ is the bimultiplicative cocycle (a map of sheaves of sets),
$$\theta_C(t_1, t_2) = \prod_{i,j} \{ x_i(t_1), x_j(t_2) \}^{c_\Basis^{ij}}.$$
Note that the right side depends on $C$, but not on the choice of basis.  Hence $\alg{T}'$ depends only on $C$ and not on the choice of basis.

The isomorphism class of $\alg{\tilde T} \in \Cat{Cov}_n^\sharp(\alg{T})$ is determined by the associated quadratic form $Q(y) = C(y,y)$.  If $C_0 \in X \otimes X$, and $C_0(y,y) = C(y,y)$ for all $y \in Y$, define $A = C - C_0$, the alternating bilinear form represented by the matrix $(a_\Basis^{ij})$, $a_\Basis^{ij} = c_{\Basis}^{ij} - c_{\Basis,0}^{ij}$.  Write $\alg{\tilde T}_0$ for the degree $n$ cover incarnated by $C_0$.  Associated to $A$ is an isomorphism $\iota_{\Basis,A} \From \alg{\tilde T}_0 \To \alg{\tilde T}$ given by
$$\iota_{\Basis,A}(t, \kappa) = \left( t, \kappa \cdot \prod_{i < j} \{ x_i(t), x_j(t) \}^{a_\Basis^{ij}} \right).$$
For any element $\hat u = \prod_i u_i^{x_i} \in X \otimes F^\times$, define an automorphism $j_{\hat u} \in \Aut(\alg{\tilde T})$ by
$$j_{\hat u}(t, \kappa) = \left( t, \kappa \cdot \prod_i \{ x_i(t), u_i \} \right).$$

From \cite[\S 3]{B-D}, it suffices to study incarnated covers and morphisms as above.
\begin{proposition}
\label{IncarnatedOK}
The category $\Cat{Cov}_n^\sharp(\alg{T})$ is equivalent to its full subcategory whose objects are the covers incarnated by elements of $X \otimes X$.  Every morphism in this full subcategory can be expressed as the  composition of a morphism $\iota_{\Basis,A}$ (for some alternating form $A$) and an automorphism of type $j_{\hat u}$ (for some $\hat u \in X \otimes F^\times$).
\end{proposition}
\proof
The isomorphism classes in $\Cat{Cov}_n^\sharp(\alg{T})$ are in bijection with quadratic forms $Q \From Y \To \ZZ$ such that $Y = Y_{Q,n}$ and such that $Q$ is even-valued if $n$ is odd.  Since every quadratic form $Q \From Y \To \ZZ$ can be written as $Q(y) = C(y,y)$ for some $C \in X \otimes X$, we find that every isomorphism class contains an incarnated cover.

The proposition now follows, since the category $\Cat{Cov}_n^\sharp(\alg{T})$ is a groupoid, the isomorphisms of type $\iota_{\Basis, A}$ link all objects within the same isomorphism class, and the automorphism group of each object is identified with $X \otimes F^\times$ by \cite[\S 3]{B-D}.
\qed

\subsection{The extension ${}^\DEE \tilde T$}

For $\alg{\tilde T}$ a sharp cover incarnated by $C \in X \otimes X$, every ordered basis $\Basis$ of $Y$ yields a section  of ${}^\DEE \tilde T$.  We describe this here, and track the dependence on basis and the effect of morphisms of covers.

First, $\alg{\tilde T}$ yields an extension $\CC^\times \Into \epsilon_\ast \tilde T \Onto T$, with $\epsilon_\ast \tilde T = T \times_{\theta_{C, \epsilon}} \CC^\times$ and incarnating cocycle
$$\theta_{C, \epsilon}(t_1, t_2) = \prod_{i,j} \Hilb_n^\epsilon( x_i(t_1), x_j(t_2) )^{c_\Basis^{ij}}.$$
Hereafter we write $\Hilb_n^\epsilon = \epsilon \circ \Hilb_n \From F^\times \times F^\times \To \alg{\mu}_n(\CC)$.

An element of $T = Y \otimes F^\times$ can be written uniquely as $\prod_i u_i^{y_i}$ with $u_i \in F^\times$ for all $1 \leq i \leq r$.  Similarly, an element of $T^\vee = X \otimes \CC^\times$ can be written uniquely as $\prod_i z_i^{x_i}$ with $z_i \in \CC^\times$ for all $1 \leq i \leq r$.

To construct ${}^\DEE \tilde T$, we first tensor with $X$ to obtain
$$X \otimes \CC^\times \Into X \otimes \epsilon_\ast \tilde T \Onto X \otimes Y \otimes F^\times.$$
The basis $\Basis$ determines a section $s_\Basis^\otimes \From X \otimes Y \otimes F^\times \To X \otimes \epsilon_\ast \tilde T = X \otimes (T \times_{\theta_{C,\epsilon}} \CC^\times)$,
$$s_\Basis^\otimes \left( \prod_{i=1}^r \prod_{j=1}^r u_{ij}^{x_i \otimes y_j} \right) = \prod_{i=1}^r  {\underbrace{\left( \prod_{j=1}^r u_{ij}^{y_j},1 \right)}_{\text{an element of } \epsilon_\ast \tilde T} }^{x_i}  \in X \otimes \epsilon_\ast \tilde T.$$  
Define a bimultiplicative cocycle $\theta_\Basis^\otimes \From (X \otimes Y \otimes F^\times) \times (X \otimes Y \otimes F^\times) \To X \otimes \CC^\times,$
$$\theta_{\Basis}^\otimes \left( u^{x_i \otimes y_j}, v^{x_k \otimes y_\ell} \right) = \begin{cases}  \Hilb_n^\epsilon( u, v)^{c_{\Basis}^{j \ell} x_i}, & \text{ if } i = k; \\ 1 & \text{ otherwise.} \end{cases}$$

The section $s_\Basis^\otimes$ defines a group isomorphism,
$$\sigma_\Basis^\otimes \From (X \otimes Y \otimes F^\times) \times_{\theta_\Basis^\otimes} (X \otimes \CC^\times) \To X \otimes \epsilon_\ast \tilde T,$$
characterized by $\sigma_\Basis^\otimes(u^{x_i \otimes y_j}, t^\vee) = t^\vee \cdot s_\Basis^\otimes(u^{x_i \otimes y_j})$ for all $1 \leq i,j \leq r$, $u \in F^\times$.

Recall that ${}^\DEE \tilde T = \iota^\ast(X \otimes \epsilon_\ast \tilde T)$, where $\iota \From \ZZ \Into X \otimes Y$ sends $1$ to $\sum_i x_i \otimes y_i$.  Define $\theta_{\Basis} \From F^\times \times F^\times \To X \otimes \CC^\times$ to be the pullback of $\theta_{\Basis}^\otimes$ via $\iota$.  Thus
\begin{equation}
\label{thetaB}
\begin{split}
\theta_{\Basis}(u,v) &= \theta_\Basis^\otimes \left( \prod_i u^{x_i \otimes y_i} , \prod_j v^{x_j \otimes y_j} \right) \\
&= \prod_{i,j} \theta_\Basis^\otimes \left( u^{x_i \otimes y_i} , v^{x_j \otimes y_j}\right) \\
&= \prod_i \theta_\Basis^\otimes  \left( u^{x_i \otimes y_i} , v^{x_i \otimes y_i}\right) \\
&= \prod_i \Hilb_n^\epsilon(u,v)^{c_{\Basis}^{ii} x_i} = \prod_i \Hilb_n^\epsilon(u,v)^{Q(y_i) x_i}.
\end{split}
\end{equation}
Pulling back $\sigma_\Basis^\otimes$ via $\iota$ gives a group isomorphism,
$$\sigma_{\Basis} \From F^\times \times_{\theta_{\Basis}} T^\vee \xrightarrow{\sim} {}^\DEE \tilde T,$$
characterized by $\sigma_{\Basis}(u, t^\vee) = t^\vee \cdot s_{\Basis}^\otimes(\iota(u))$ for all $u \in F^\times$.

\begin{lemma}
\label{cocycleDEE}
For all $u,v \in F^\times$, we have
$$\theta_{\Basis}(u,v) = \prod_{i=1}^r \Hilb_n^\epsilon(u,v)^{Q(y_i) x_i} = \tau_Q( \Hilb_2(u,v) ).$$
In particular, $\theta_\Basis$ is independent of basis.
\end{lemma}
\begin{remark}
If $\Char(F) = 2$, then $n$ is odd, so $\tau_Q(\pm 1) = 1$, and one should interpret the lemma as $\theta_{\Basis}(u,v) \ident 1$ since the quadratic Hilbert symbol is not defined.
\end{remark}
\proof
Recall that $\tau_Q \From \mu_2 \To T^\vee = X \otimes \CC^\times$ is dual to the homomorphism $Y \To \frac{1}{2} \ZZ / \ZZ$ given by $y \mapsto n^{-1} Q(y)$.  Thus $\tau_Q$ can be expressed using the basis $\Basis$,
$$\tau_Q(-1) = \prod_i  e^{2 \pi i n^{-1} Q(y_i) \cdot x_i}.$$
Since $\alg{\tilde T}$ is a sharp cover, $Y = Y_{Q,n}$, and so $2 Q(y) \in n \ZZ$ for all $y \in Y$.  If $n$ is odd (in particular if $\Char(F) = 2$), this implies that $Q(y_i) \in n \ZZ$ for all $1 \leq i \leq r$, and so $\theta_{\Basis}(u,v) = 1 = \tau_Q(\Hilb_2(u,v))$ as claimed.  

If $n$ is even, then $\Hilb_2(u,v) = \Hilb_n(u,v)^{n/2}$ and so
$$\Hilb_n^\epsilon(u,v)^{Q(y_i)} =  \Hilb_2(u,v)^{2 Q(y_i) / n} = \begin{cases} \Hilb_2(u,v) , & \text{ if } n^{-1} Q(y_i) \not \in \ZZ; \\ 1, & \text{ if } n^{-1} Q(y_i) \in \ZZ. \end{cases}$$
The result follows from our previous computation of $\theta_\Basis(u,v)$ in \eqref{thetaB}.
\qed

Hereafter, we write simply $\theta$ rather than $\theta_{\Basis}$.  Thus $\sigma_{\Basis}$ gives an isomorphism in $\Cat{Ext}(F^\times, T^\vee)$,
$$\sigma_{\Basis} \From F^\times \times_\theta T^\vee \To {}^\DEE \tilde T.$$
The cocycle $\theta(u,v) = \tau_Q(\Hilb_2(u,v))$ is independent of basis, but the isomorphism $\sigma_{\Basis}$ depends on the basis.

We trace through morphisms of extensions here.  First, consider an element $\hat u = \prod_\ell u_\ell^{x_\ell} \in X \otimes F^\times = \Aut(\alg{\tilde T})$.  This defines an automorphism of $\epsilon_\ast \tilde T$, given by
$$j_{\hat u}(t, z) = \left( t, z \cdot \prod_\ell \Hilb_n^\epsilon(x_\ell(t), u_\ell) \right).$$

On the other hand, for any extension $E \in \Cat{Ext}(F^\times, T^\vee)$, $\hat u$ defines an automorphism $j_{\hat u}^\vee \in \Aut(E)$ by 
\begin{equation}
\label{jCheck}
j_{\hat u}^\vee(e) = e \cdot \prod_\ell \Hilb_n^\epsilon(v,u_\ell)^{x_\ell}, \text{ for all } e \in E \text{ lying over } v \in F^\times.
\end{equation}
In particular, the automorphism $j_{\hat u}^\vee$ coincides with ${}^\DEE j_{\hat u}$ on the extension ${}^\DEE \tilde T$.  If $E_1, E_2 \in \Cat{Ext}(F^\times, T^\vee)$, and $\eta \From E_1 \To E_2$ is a morphism of extensions, then
$$j_{\hat u}^\vee \circ \eta = \eta \circ j_{\hat u}^\vee.$$

From this, we directly obtain compatibility of $\sigma_\Basis$ with automorphisms.
\begin{proposition}
\label{AutomorphDParam}
The isomorphisms $\sigma_\Basis$ are compatible with automorphisms $j_{\hat u}$, i.e., the following diagram commutes.
$$\begin{tikzcd}
F^\times \times_\theta T^\vee \arrow{r}{\sigma_\Basis} \arrow{d}{j_{\hat u}^\vee} & {}^\DEE \tilde T \arrow{d}{j_{\hat u}^\vee = {}^\DEE j_{\hat u}} \\
F^\times \times_\theta T^\vee \arrow{r}{\sigma_\Basis} & {}^\DEE \tilde T
\end{tikzcd}$$
\end{proposition}

Next we consider a morphism $\iota_{\Basis,A} \From \alg{\tilde T}_0 \To \alg{\tilde T}$, given by two elements $C_0, C \in X \otimes X$ with $C(y,y) = C_0(y,y)$ for all $y \in Y$.  The matrix $(a_{\Basis}^{ij})$ is defined by $a_{\Basis}^{ij} = c_{\Basis}^{ij} - c_{\Basis,0}^{ij}$, and the morphism of covers $\iota_{\Basis,A}$ induces an isomorphism $\iota_{\Basis,A} \From \epsilon_\ast \tilde T_0 \To \epsilon_\ast \tilde T$,
$$\iota_{\Basis,A}(t,z) = \left( t,z \cdot \prod_{i < j} \Hilb_n^\epsilon(x_i(t), x_j(t))^{a_{\Basis}^{ij}} \right).$$
Write $\sigma_{0, \Basis} \From F^\times \times_\theta T^\vee \To {}^\DEE \tilde T_0$ for the isomorphism defined analogously to $\sigma_\Basis$.
\begin{proposition}
\label{MorphDParam}
The isomorphisms $\sigma_\Basis$ and $\sigma_{0, \Basis}$ fit into a commutative diagram.
$$\begin{tikzcd}
F^\times \times_\theta T^\vee \arrow{r}{\sigma_\Basis} \arrow{d}{=} & {}^\DEE \tilde T \arrow{d}{{}^\DEE \iota_{\Basis,A}} \\
F^\times \times_\theta T^\vee \arrow{r}{\sigma_{0,\Basis}} & {}^\DEE \tilde T_0
\end{tikzcd}$$
\end{proposition}
\proof
It suffices to check that the diagram commutes for elements of the form $(u,1) \in F^\times \times_\theta T^\vee$.  Observe that the vertical map ${}^\DEE \iota_{\Basis,A}$ is the restriction of a homomorphism $[X \otimes \iota_{\Basis,A}] \From X \otimes \epsilon_\ast \tilde T_0 \To X \otimes \epsilon_\ast \tilde T$,
$$[X \otimes \iota_{\Basis,A}] \left( \prod_\ell (t_\ell, 1)^{x_\ell} \right) = \prod_\ell \left( t_\ell, \prod_{i < j} \Hilb_n^\epsilon(x_i(t_\ell), x_j(t_\ell))^{a_{\Basis}^{ij}} \right)^{x_\ell}.$$
If $(u,1) \in F^\times \times_\theta T^\vee$, then
$$\sigma_{\Basis}(u,1) = s_{\Basis}^\otimes \left( \prod_\ell u^{x_\ell \otimes y_\ell} \right) = \prod_\ell (u^{y_\ell}, 1)^{x_\ell} \in X \otimes \epsilon_\ast \tilde T.$$
Similarly,
$$\sigma_{0,\Basis}(u,1) = \prod_\ell (u^{y_\ell}, 1)^{x_\ell} \in X \otimes \epsilon_\ast \tilde T_0.$$

Thus we find
\begin{align*}
{}^\DEE \iota_{\Basis,A} (\sigma_{0,\Basis} (u,1)) &= [X \otimes \iota_{\Basis,A}] \left( \prod_\ell (u^{y_\ell}, 1)^{x_\ell} \right) \\
&= \prod_\ell \left( u^{y_\ell}, \prod_{i < j} \Hilb_n^\epsilon(x_i(u^{y_\ell}), x_j(u^{y_\ell}))^{a_{\Basis}^{ij}} \right)^{x_\ell} \\
&= \prod_\ell \left( u^{y_\ell}, 1 \right)^{x_\ell} = \sigma_{\Basis}(u,1).
\end{align*}
The simplification of the inner product is based on the following observation:  $x_i(u^{y_\ell}) = 1$ unless $i = \ell$ and $x_j(u^{y_\ell}) = 1$ unless $j = \ell$, but the product is indexed by only those $i,j$ with $i < j$.
\qed

\subsection{Change of basis}

If $\Basis' = ( y_1', \ldots, y_r')$ is another ordered basis of $Y$, then we find another isomorphism, $\sigma_{\Basis'} \From F^\times \times_{\theta} T^\vee \To {}^\DEE \tilde T$.  Consider the change of basis matrix $(g_{ij}) \in GL_r(\ZZ)$, so that $x_i' = \sum_j g_{ij} x_j$ and $\sum_i g_{ij} y_i' = y_j$.  For what follows, define
\begin{equation}
\label{DGChi}
\begin{split}
\Delta_{ij} = c_{\Basis'}^{ii} \frac{g_{ij} (g_{ij} - 1)}{2}, \quad \Gamma_{j}^{k \ell} = c_{\Basis'}^{k \ell} g_{kj} g_{\ell j}, \\
\chi(u) = \Hilb_n^\epsilon(u,u) \text{ for all } u \in F^\times.
\end{split}
\end{equation}

\begin{proposition}
\label{BasisChangeSigma}
Define the automorphism $f_{\Basis, \Basis'} \in \Aut(F^\times \times_{\theta} T^\vee)$, in the category $\Cat{Ext}(F^\times, T^\vee)$,
$$f_{\Basis, \Basis'}(u, t^\vee) =  \left( u, t^\vee \cdot \prod_j \chi(u)^{ \left( \sum_i \Delta_{ij} + \sum_{k < \ell} \Gamma_j^{k \ell} \right) x_j } \right).$$
This fits into a commutative diagram in $\Cat{Ext}(F^\times, T^\vee)$.
$$\begin{tikzcd}
F^\times \times_{\theta} T^\vee \arrow{drr}{\sigma_{\Basis}} \arrow{dd}{f_{\Basis, \Basis'}} & & \\
 \phantom{a} & & {}^\DEE \tilde T \\
 F^\times \times_{\theta} T^\vee \arrow{urr}[swap]{\sigma_{\Basis'}} & &
\end{tikzcd}$$
\end{proposition}
\proof
There exists a unique isomorphism $f_{\Basis, \Basis'}$ making the diagram above commute, given by $f_{\Basis, \Basis'}(u, t^\vee) = (u, \delta(u) t^\vee)$, where $\delta(u) = s_{\Basis}^\otimes( \iota(u)) \cdot  s_{\Basis'}^\otimes(\iota(u))^{-1}$.  Now we compute
\begin{align*}
s_{\Basis'}^\otimes(\iota(u)) &=  s_{\Basis'}^\otimes  \left( \prod_i u^{x_i' \otimes y_i'} \right), \\
&=  \prod_{i=1}^r (u^{y_i'}, 1)^{x_i'}, \quad \text{(noting that $(u^{y_i'}, 1) \in \epsilon_\ast \tilde T = T \times_{\theta_{C,\epsilon}} \CC^\times$)} \\
&=  \prod_{i=1}^r (u^{y_i'}, 1)^{\sum_{j=1}^r g_{ij} x_j}, \\
&= \prod_{i,j = 1}^r  (u^{y_i'}, 1)^{g_{ij} x_j},  \\
&=  \prod_j \prod_i \left( u^{g_{ij} y_i'}, 1 \right)^{x_j}  \cdot \prod_{i,j} \chi(u)^{\Delta_{ij} \cdot x_j}, \quad \text{ by \eqref{gpower}} \\
&= \prod_j \left( \prod_i u^{g_{ij} y_i'}, 1 \right)^{x_j} \cdot \prod_{i,j} \chi(u)^{\Delta_{ij} \cdot x_j} \cdot \prod_j \prod_{k < \ell} \chi(u)^{\Gamma_{j}^{k \ell} x_j}, \\
&= \prod_j \left( u^{y_i}, 1 \right)^{x_j} \cdot \prod_j \chi(u)^{ \left( \sum_i \Delta_{ij} + \sum_{k < \ell} \Gamma_j^{k \ell} \right) x_j }, \\
&= s_{\Basis}^\otimes \left( \prod_{j} u^{x_j \otimes y_j} \right)  \cdot \prod_j \chi(u)^{ \left( \sum_i \Delta_{ij} + \sum_{k < \ell} \Gamma_j^{k \ell} \right) x_j }, \\
&= s_{\Basis}^\otimes(\iota(u)) \cdot \prod_j \chi(u)^{ \left( \sum_i \Delta_{ij} + \sum_{k < \ell} \Gamma_j^{k \ell} \right) x_j }.
\end{align*}
The proposition follows immediately.
\qed

\subsection{Parameters}

We begin with the same data, a degree $n$ sharp cover $\alg{\tilde T}$ of a split torus $\alg{T}$ over $F$, incarnated by $C \in X \otimes X$, and an ordered basis $\Basis = (y_1, \ldots, y_r)$ of the cocharacter lattice $Y$.  The second Brylinski-Deligne invariant of $\alg{\tilde T}$ is an extension $\sheaf{G}_m \Into \sheaf{D} \Onto Y$ of sheaves on $F_{\et}$, and we can describe $\sheaf{D}$ in terms of $C \in X \otimes X$.  

Tracing through the construction of the second Brylinski-Deligne invariant \cite[\S 3.12]{B-D}, consider $\alg{T}'(F(\!(\form)\!))$ for a formal parameter $\form$.  This fits into a short exact sequence
$$\alg{K}_2(F(\!(\form)\!)) \Into \alg{T}'(F(\!(\form)\!)) \Onto \alg{T}(F(\!(\form)\!)).$$
This extension is incarnated by the bimultiplicative cocycle $\theta_\form$ satisfying
$$\theta_{\Basis,\form}(f^{y_i}, g^{y_j}) = \{ f,g \}^{c_{\Basis}^{ij}}, \text{ for all } f,g \in F(\!(\form)\!), 1 \leq i,j \leq r.$$
Pushing out via the tame symbol $\partial \From \alg{K}_2(F(\!(\form)\!)) \To \alg{K}_1(F) = F^\times$, and pulling back via $Y \To \alg{T}(F(\!(\form)\!))$, $y \mapsto \form^y$, we find the extension $F^\times \Into \sheaf{D}[F] \Onto Y$.  This extension is incarnated by the bimultiplicative cocycle $\theta_{\Basis,D}$ satisfying
$$\theta_{\Basis,D}(y_i, y_j) = \partial \{ \form, \form \}^{c_{\Basis}^{ij}} = (-1)^{c_{\Basis}^{ij}}.$$
This construction applies equally well over a finite Galois extension of $F$.  Thus the basis $\Basis$ determines an isomorphism 
$$\sigma_{\Basis,D} \From Y \times_{\theta_{\Basis,D}} \sheaf{G}_m \To \sheaf{D},$$
of extensions of sheaves of groups on $F_{\et}$.  Sharpness of $\alg{\tilde T}$ implies that $\sheaf{D}$ is a sheaf of \textbf{abelian} groups on $F_{\et}$.

The choice of basis $\Basis$ gives a splitting of $\sheaf{D}$:  since $Y$ is a constant sheaf of free abelian groups, there exists a unique splitting $d_\Basis \From Y \To \sheaf{D}$ which satisfies $d_\Basis(y_i) = \sigma_{\Basis,D}(y_i, 1)$, for all $1 \leq i \leq r$.  Thus the choice of basis trivializes the $\shom(Y, \sheaf{G}_m) = \sheaf{\hat T}$-torsor $\sspl(\sheaf{D})$.  This trivialization, in turn, neutralizes the gerbe of $n^{\th}$ roots $\sqrt[n]{\sspl(\sheaf{D})}$ over $F$.  Explicitly, define $z_\Basis \in \sqrt[n]{\sspl(\sheaf{D})}[F]$ by $z_\Basis = (\sheaf{\hat T}, h_\Basis)$ with $h_\Basis \From \sheaf{\hat T} \To \sspl(\sheaf{D})$ given by
\begin{equation}
\label{hb}
h_\Basis(\hat a) = \hat a^n \ast d_\Basis.
\end{equation}

This neutralizes the pushout $\gerb{E}_{\epsilon}(\alg{\tilde T})$, and we find an isomorphism in $\Cat{CExt}(\Gal_F, T^\vee)$,
\begin{equation}
\label{LambdaGerbe}
\lambda_\Basis \From  \Gal_F \times T^\vee \To \pi_1^{\et}(\gerb{E}_\epsilon(\alg{\tilde T}), z_\Basis),
\end{equation}
The Baer sum with $(\tau_Q)_\ast \mGal_F$ yields an isomorphism in $\Cat{CExt}(\Gal_F, T^\vee)$,
\begin{equation}
\label{LambdaGal}
\lambda_\Basis \From (\tau_Q)_\ast \mGal_F \To {}^\EL \tilde T_{z_\Basis}.
\end{equation}
Here ${}^\EL \tilde T_{z_{\Basis}}$ denotes the L-group with respect to the base point $z_{\Basis}$.  Note that $(\tau_Q)_\ast \mGal_F = \Gal_F \times T^\vee$ as sets, with multiplication given by the cocycle,
$$(\gamma_1, \gamma_2) \mapsto \tau_Q \left(\Hilb_2 (\rec \gamma_1, \rec \gamma_2 ) \right).$$
From \eqref{LambdaGal}, taking the quotient by $[\Gal_F, \Gal_F]$, and pulling back via $F^\times \To \Weil_F^{\ab} \To \Gal_F^{\ab}$, we find an isomorphism in $\Cat{Ext}(F^\times, T^\vee)$,
\begin{equation}
\label{LambdaA}
\lambda_\Basis \From F^\times \times_\theta T^\vee \To {}^{\A} \tilde T_{z_{\Basis}},
\end{equation}
where we recall that $\theta(u,v) = \tau_Q(\Hilb_2(u,v))$.  For what follows, it helps to describe $\lambda_\Basis \From F^\times \times_\theta T^\vee \To {}^{\A} \tilde T_{z_{\Basis}}$ in more detail.  For an element $(u, 1) \in F^\times \times_{\theta} T^\vee$, take the following steps:
\begin{itemize}
\item
Choose $\gamma \in \Weil_F \subset \Gal_F$ such that $\rec(\gamma) = u$.
\item
Since $z_\Basis$ is defined over $F$, $\Hom(z_\Basis, {}^\gamma z_\Basis) = \Hom(z_\Basis, z_\Basis)$.  Thus we have an element $\Id_{\Basis,\gamma} \in \Hom(z_\Basis, {}^\gamma z_\Basis) \subset \pi_1^{\et}(\gerb{E}_\epsilon(\alg{\tilde T}), z_\Basis)$.
\item
We also have an element $(\gamma, 1) \in (\tau_Q)_\ast \mGal_F$ lying over $\gamma \in \Gal_F$.  Hence 
$$((\gamma,1), \Id_{\Basis, \gamma}) \in {}^\EL \tilde T_{z_{\Basis}} = (\tau_Q)_\ast \mGal_F \Baer \pi_1^{\et}(\gerb{E}_\epsilon(\alg{\tilde T}), z_\Basis).$$
\item
$\lambda_\Basis(u, 1)$ is the image of $((\gamma,1), \Id_{\Basis, \gamma})$ in the quotient ${}^\A \tilde T_{z_\Basis}$.
\end{itemize}

As in Proposition \ref{AutomorphDParam}, the isomorphisms $\lambda_\Basis$ are compatible with automorphisms of covers.  If $\hat u \in X \otimes F^\times$, then we have the resulting automorphism $j_{\hat u} \in \Aut(\alg{\tilde T})$.  As before, let $j_{\hat u}^\vee$ denote the automorphism of any extension in $\Cat{Ext}(F^\times, T^\vee)$ associated to $\hat u$.  Since $j_{\hat u}$ is well-aligned, we also get an automorphism of the L-group ${}^\EL \tilde T_{z_{\Basis}}$, which descends to an automorphism ${}^\A j_{\hat u}$ of ${}^\A \tilde T_{z_{\Basis}}$.
\begin{proposition}
\label{AutomorphAParam}
For all $\hat u \in X \otimes F^\times$, $j_{\hat u}^\vee = {}^\A j_{\hat u}$, making the following diagram commute.
$$\begin{tikzcd}
F^\times \times_\theta T^\vee \arrow{r}{\lambda_{\Basis}} \arrow{d}{j_{\hat u}^\vee} & {}^\A \tilde T_{\Basis} \arrow{d}{j_{\hat u}^\vee = {}^{\A} j_{\hat u}}  \\
F^\times \times_\theta T^\vee \arrow{r}{\lambda_\Basis} &  {}^\A \tilde T_{\Basis}
\end{tikzcd}$$
\end{proposition}
\proof
We trace the automorphism $j_{\hat u}$ through the construction of ${}^\EL \tilde T_{z_{\Basis}}$.  It has no effect on $(\tau_Q)_{\ast} \mGal_F$.  On the other hand, the automorphism $j_{\hat u}$ corresponds to the automorphism of $\sheaf{G}_m \Into \sheaf{D} \Onto Y$,
$$d \mapsto d \cdot \hat u(y), \text{ for all } d \in \sheaf{D} \text{ lying over } y \in Y.$$
Here and later, $\hat u \in X \otimes F^\times = \Hom(Y, F^\times)$, so $\hat u(y) \in F^\times$.

If $s \in \sspl(\sheaf{D})$, then $j_{\hat u}^\ast(s) = \hat u^{-1} \ast s$ is the pullback of $s$ via $j_{\hat u}$.  (See Section \ref{WAFGerbe} for conventions on well-aligned functoriality).  From this, $j_{\hat u}$ gives an equivalence of gerbes $\sqrt[n]{j_{\hat u}^\ast} \From \sqrt[n]{\sspl(\sheaf{D})} \To \sqrt[n]{\sspl(\sheaf{D})}$, which sends the base point $z_{\Basis} = (\sheaf{\hat T}, h_\Basis)$ to $z_{\Basis}' = (\sheaf{\hat T}, \hat u \circ h_\Basis)$.  Here $\hat u \circ h_\Basis \From \sheaf{\hat T} \To \sspl(\sheaf{D})$ satisfies
\begin{equation}
\label{uhb}
[\hat u \circ h_\Basis](\hat a) = \hat a^n \hat u^{-1} \ast d_\Basis.
\end{equation}
Choose any isomorphism $\rho \From z_\Basis \To z_\Basis'$ in $\sqrt[n]{\sspl(\sheaf{D})}[\bar F]$.  Concretely, $\rho \From \sheaf{\hat T} \To \sheaf{\hat T}$ must be a map of $\sheaf{\hat T}$-torsors, so there exists $\hat \tau \in \sheaf{\hat T}$ such that $\rho(\hat a) = \hat \tau \ast \hat a$, for all $\hat a \in \sheaf{\hat T}$.  Moreover, since $\rho$ intertwines $h_{\Basis}$ and $\hat u \circ h_{\Basis}$, \eqref{hb} and \eqref{uhb} imply that $\hat \tau^n = \hat u$.  

Now, for any $f \in \pi_1^{\et}(\sqrt[n]{\sspl(\sheaf{D})}, z_{\Basis}) = \Hom(z_{\Basis}, z_{\Basis})$, define an isomorphism $f' \in \pi_1^{\et}(\sqrt[n]{\sspl(\sheaf{D})}, z_{\Basis}')$ by the rule $f'(\hat a) = f(\hat a)$ for all $\hat a \in \sheaf{\hat T}$.  (Here we note that $z_{\Basis} = (\sheaf{\hat T}, h_{\Basis})$ and $z_{\Basis}' = (\sheaf{\hat T}, \hat u \circ h_{\Basis})$, so $f$ and $f'$ are both given by functions from $\sheaf{\hat T}$ to itself).  Then $f'$ is the image of $f$, under the equivalence of gerbes $\sqrt[n]{j_{\hat u}}$.

Define 
$${}^\rho f \in \pi_1^{\et}(\sqrt[n]{\sspl(\sheaf{D})}, z_{\Basis}) = {}^\gamma \rho^{-1} \circ f' \circ \rho \in \Hom(z_\Basis, {}^\gamma z_\Basis).$$
Then we compute
$${}^\rho f(\hat a) = {}^\gamma \rho^{-1}( f' ( \rho(\hat a) ) ) = {}^\gamma \rho^{-1}(f'(\hat \tau \ast \hat a)) = \frac{\hat \tau}{ \gamma(\hat \tau)} \ast f(\hat a).$$
Pushing out via $\epsilon$ and taking the Baer sum with $(\tau_Q)_\ast \mGal_F$, we find that $f \mapsto {}^\rho f$ induces an automorphism ${}^\EL j_{\hat u}$ of the extension
$$T^\vee \Into {}^\EL \tilde T_{z_{\Basis}} \Onto \Gal_F,$$
given by
$${}^\EL j_{\hat u}(f) = \epsilon \left( \frac{\sqrt[n]{\hat u}}{\gamma (\sqrt[n]{\hat u})} \right) \cdot f.$$
for all $f$ lying over $\gamma \in \Gal_F$.
For $u_\ell \in F^\times$ and $\gamma \in \Weil_F \subset \Gal_F$, the relationship between the Artin symbol and Hilbert symbol gives
$$\frac{ \gamma \sqrt[n]{u_\ell} }{ \sqrt[n]{u_\ell} } = \Hilb_n(u_\ell, \rec \gamma) = \Hilb_n( \rec \gamma, u_\ell)^{-1}.$$

Hence, ${}^\EL j_{\hat u}$ descends to an automorphism ${}^\A j_{\hat u}$ of the extension
$$T^\vee \Into {}^\A \tilde T_{z_{\Basis}} \Onto F^\times,$$
given by
$${}^\A j_{\hat u}(e) = e \cdot \prod_\ell \Hilb_n^\epsilon(v, u_\ell)^{x_\ell},$$
for all $e \in {}^\A \tilde T_{z_{\Basis}}$ lying over $v \in F^\times$.
Hence ${}^\A j_{\hat u} = j_{\hat u}^\vee$ as defined in \eqref{jCheck}.
\qed

Next, as in Proposition \ref{MorphDParam}, consider two elements $C_0, C \in X \otimes X$ such that $C_0(y,y) = C(y,y)$ for all $y \in Y$.  Let $A = C - C_0$, giving an isomorphism $\iota_{\Basis,A} \From \alg{\tilde T}_0 \To \alg{\tilde T}$ of the incarnated covers.  Write $\sheaf{D}_0$, $d_{0, \Basis}$, $z_{0, \Basis}$, ${}^\A \tilde T_{0, z_{0,\Basis}}$, etc., for the analogues of $\sheaf{D}$, $d_\Basis$, $z_\Basis$, ${}^\A \tilde T_{z_{\Basis}}$, using $C_0$ instead of $C$.  

\begin{proposition}
\label{MorphAParam}
The homomorphism ${}^\A \iota_{\Basis,A} \From {}^\A \tilde T_{z_\Basis} \To {}^\A \tilde T_{0,z_{0,\Basis}}$, induced by the well-aligned isomorphism $\iota_{\Basis,A}$, fits into a commutative square.
$$\begin{tikzcd}
F^\times \times_\theta T^\vee \arrow{r}{\lambda_\Basis} \arrow{d}{=} &{}^\A \tilde T_{z_\Basis} \arrow{d}{{}^\A \iota_{\Basis,A}}  \\
F^\times \times_\theta T^\vee \arrow{r}{\lambda_{0,\Basis}} & {}^\A \tilde T_{0,z_{0,\Basis}}
\end{tikzcd}$$
\end{proposition}
\proof
The isomorphism $\iota_{\Basis,A} \From \alg{\tilde T}_0 \To \alg{\tilde T}$ gives an isomorphism of extensions below.
$$\begin{tikzcd}
\sheaf{G}_m \arrow{d}{=} \inarrow{r} & \sheaf{D}_0 \onarrow{r} \arrow{d}{\iota_{\Basis,A}} & Y \arrow{d}{=}  \\
\sheaf{G}_m \inarrow{r} & \sheaf{D} \onarrow{r} & Y
\end{tikzcd}$$
Using the presentation of $\sheaf{D}_0$ and $\sheaf{D}$ by cocycles from $C_0$ and $C$, respectively, we find that
$$\iota_{\Basis,A}(y,1) = \left( y, \prod_{i < j} (-1)^{\langle x_i, y \rangle \cdot \langle x_j, y \rangle \cdot a_{\Basis}^{ij} } \right),$$
for all $y \in Y$.  Recall that the basis $\Basis$ gave splittings of $\sheaf{D}_0$ and $\sheaf{D}$ respectively, satisfying $d_0(y_i) = (y_i, 1)$ and $d(y_i) = (y_i, 1)$.  Since $\iota_{\Basis,A}(y_i, 1) = (y_i, 1)$, we find that $\iota_{\Basis,A} \circ d_0 = d$.

We may use $\iota_{\Basis,A}$ to pull back splittings of $\sheaf{D}$ to splittings of $\sheaf{D}_0$.  We find that this map of $\sheaf{\hat T}$-torsors
$$\iota_{\Basis,A}^\ast \From \sspl(\sheaf{D}) \To \sspl(\sheaf{D}_0)$$
sends the trivialization $d$ to the trivialization $d_0$.  It follows that $\iota_{\Basis,A}^\ast$ gives a functor of gerbes,
$$\iota_{\Basis,A}^\ast \From \sqrt[n]{\sspl(\sheaf{D})} \To \sqrt[n]{\sspl(\sheaf{D}_0)}$$
which sends the base point $z_\Basis$ to the base point $z_{0, \Basis}$.  As both are base points defined over $F$, we find that $\iota_{\Basis,A}^\ast$ sends the identity morphism $\Id_{\Basis, \gamma} \in \Hom(z_\Basis, {}^\gamma z_\Basis) = \Hom(z_\Basis, z_\Basis)$ to the corresponding identity morphism $\Id_{0,\Basis, \gamma} \in \Hom(z_{0, \Basis}, {}^\gamma z_{0, \Basis})$.  

It follows that ${}^\A \iota_{\Basis,A}$ sends $\lambda_{\Basis}(u,1)$ to $\lambda_{0,\Basis}(u,1)$ for all $u \in F^\times$.  The proposition follows.
\qed

\subsection{Change of basis}

If $\Basis' = (y_1', \ldots, y_r')$ is another basis, then we find from it another splitting $d_{\Basis'} \From Y \To \sheaf{D}$ characterized by $d_{\Basis}'(y_i') = (y_i', 1)$.  This provides another trivialization of the torsor $\sspl(\sheaf{D})$, and thus another object $z_{\Basis'} = (\sheaf{\hat T}, h_{\Basis'})$ of $\sqrt[n]{\sspl(\sheaf{D})}[F]$.    

An isomorphism $\iota \From z_{\Basis} \To z_{\Basis'}$ in $\sqrt[n]{\sspl(\sheaf{D})}[\bar F]$ produces an isomorphism ${}^\EL \iota \From {}^\EL T_{z_\Basis} \To {}^\EL T_{z_{\Basis'}}$ in $\Cat{CExt}(\Gal_F, T^\vee)$, and any two such isomorphisms $\iota_1, \iota_2$ produce the {\em same} isomorphism ${}^\EL \iota_1 = {}^\EL \iota_2$.  It follows that there is a unique automorphism $r_{\Basis, \Basis'}$ of $F^\times \times_{\theta} T^\vee$ in the category $\Cat{Ext}(F^\times, T^\vee)$, which makes the following diagram commute.
$$\begin{tikzcd}
F^\times \times_{\theta} T^\vee \arrow{r}{\lambda_\Basis} \arrow{d}{r_{\Basis, \Basis'}} & {}^\A T \arrow{d}{{}^\A \iota}_{z_\Basis} \\
F^\times \times_{\theta} T^\vee \arrow{r}{\lambda_{\Basis'}} & {}^\A T_{z_{\Basis'}} 
\end{tikzcd}$$

Consider the change of basis matrix $(g_{ij}) \in GL_r(\ZZ)$, so that $x_i' = \sum_j g_{ij} x_j$ and $\sum_i g_{ij} y_i' = y_j$.  Define $\Delta_{ij}$, $\Gamma_j^{k \ell}$, and $\chi$ as in \eqref{DGChi}.
\begin{proposition}
\label{BasisChangeLambda}
The automorphism $r_{\Basis, \Basis'} \in \Aut(F^\times \times_{\theta} T^\vee)$ satisfies
$$r_{\Basis, \Basis'}(u, t^\vee ) = \left( u, t^\vee \cdot  \prod_j \chi(u)^{\left( \sum_i \Delta_{ij} + \sum_{k < \ell} \Gamma_j^{k \ell} \right) \cdot x_j} \right).$$
\end{proposition}
\proof
Giving an isomorphism $\iota \From z_{\Basis} \To z_{\Basis'}$ is the same as giving a morphism $f \From \sheaf{\hat T} \To \sheaf{\hat T}$ of $\sheaf{\hat T}$-torsors, $f(\hat a) = \hat e \cdot \hat a$, satisfying $d_{\Basis'} =  \hat e^n \ast d_{\Basis}$.  This condition on $\hat e$ is equivalent to the following, for all $1 \leq j \leq r$.
\begin{align*}
\hat e(y_j)^n &= d_\Basis(y_j)^{-1} \cdot d_{\Basis'}(y_j),  \\
&= (y_j, 1)^{-1} \cdot d_{\Basis'}\left( \sum_i g_{ij} y_i' \right), \\
&= (y_j, 1)^{-1}  \cdot  \prod_i (y_i',1)^{g_{ij}}, \\
&= (y_j, 1)^{-1} \cdot \prod_i \left( g_{ij} y_i', 1 \right) \cdot (-1)^{\Delta_{ij}}, \\
&= (y_j,1)^{-1} \cdot \left( \sum_i g_{ij} y_i', 1 \right) \cdot \prod_i (-1)^{\Delta_{ij}} \cdot \prod_{k < \ell} (-1)^{\Gamma_j^{k \ell} }, \\
&= \prod_i (-1)^{\Delta_{ij}} \cdot \prod_{k < \ell} (-1)^{\Gamma_j^{k \ell} }.
\end{align*}

In other words, giving a morphism $\iota \From z_{\Basis} \To z_{\Basis'}$ (defined over $\bar F$) is the same as giving an element $\hat e \in \sheaf{\hat T}[\bar F] = X \otimes \bar F^\times$ satisfying
\begin{equation}
\label{eFormula}
\hat e^n = \prod_j (-1)^{\left( \sum_i \Delta_{ij} + \sum_{k < \ell} \Gamma_j^{k \ell} \right) \cdot x_j}.
\end{equation}

Now consider an element $(u,1) \in F^\times \times_{\theta} T^\vee$.  Let $\gamma \in \Weil_F$ be an element such that $\rec(\gamma) = u$, and $\Id_{\Basis,\gamma} \in \Hom(z_\Basis, {}^\gamma z_\Basis)$ the resulting element of $\pi_1^{\et}(\gerb{E}_\epsilon(\alg{\tilde T}), z_{\Basis})$.  Similarly, let $\Id_{\Basis', \gamma} \in  \Hom(z_{\Basis'}, {}^\gamma z_{\Basis'})$ be the resulting element of $\pi_1^{\et}(\gerb{E}_\epsilon(\alg{\tilde T}), z_{\Basis'})$

The following diagram of objects and morphisms in $\sqrt[n]{\sspl(\sheaf{D})}$ commutes.
$$\begin{tikzcd}[ column sep = 6em]
z_\Basis \arrow{r}{\Id_{\Basis, \gamma}} \arrow{d}[swap]{\hat e} & {}^\gamma z_\Basis = z_\Basis \arrow{d}{ {}^\gamma \hat e} \\
z_{\Basis'} \arrow{r}{ ({}^\gamma \hat e / \hat e) \cdot \Id_{\Basis', \gamma} } & {}^\gamma z_{\Basis'} = z_{\Basis'}
\end{tikzcd}$$
Here we label morphisms among $z_\Basis$ and $z_\Basis'$ by elements of $\sheaf{\hat T}$, since all morphisms are morphisms of $\sheaf{\hat T}$-torsors from $\sheaf{\hat T}$ to itself.  Projecting from ${}^\EL \tilde T_{z_\Basis}$ to ${}^{\A} \tilde T_{z_\Basis}$, we find that
$${}^\A \iota \left( \lambda_{\Basis}^\A(u,1) \right) = \frac{ {}^\gamma \hat e}{\hat e} \cdot \lambda_{\Basis'}^\A(u,1).$$
The relationship between the Artin symbol and Hilbert symbol gives
$$\epsilon \left( \frac{\gamma( \sqrt[n]{-1})}{\sqrt[n]{-1}} \right) = \Hilb_n^\epsilon(-1,u) = \Hilb_n^\epsilon(u,u) = \chi(u),$$
since $\rec(\gamma) = u$.  Hence equation \eqref{eFormula} yields
$${}^\A \iota( \lambda_{\Basis}^\A(u,1) ) =  \lambda_{\Basis'}^\A(u,1)  \cdot  \prod_j \chi(u)^{\left( \sum_i \Delta_{ij} + \sum_{k < \ell} \Gamma_j^{k \ell} \right) \cdot x_j}.$$
The proposition follows.
\qed

\subsection{Isomorphism of functors}

For each basis $\Basis$ of $Y$, and each cover $\alg{\tilde T}$ of $\alg{T}$ incarnated by $C$, we have constructed isomorphisms $\sigma_\Basis \From F^\times \times_\theta T^\vee \To {}^\DEE \tilde T$ and $\lambda_\Basis \From F^\times \times_\theta T^\vee \To {}^\A \tilde T$.  Here we write ${}^\A \tilde T$ rather than ${}^\A \tilde T_{z_\Basis}$, since ${}^\A \tilde T$ is defined up to unique isomorphism by $\alg{\tilde T}$ without reference to base point.  

There exists a unique isomorphism $\P_\Basis$ in $\Cat{Ext}(F^\times, T^\vee)$ which makes the following diagram commute.
$$\begin{tikzcd}
\phantom{a} & {}^\DEE \tilde T \arrow{dd}{\P_\Basis} \\
F^\times \times_\theta T^\vee \arrow{ur}{\sigma_\Basis} \arrow{dr}[swap]{\lambda_\Basis} &  \phantom{a} \\
\phantom{a} & {}^\A \tilde T
\end{tikzcd}$$

\begin{thm}
If $\Basis$ and $\Basis'$ are two bases of $Y$, then $\P_\Basis = \P_{\Basis'}$.
\end{thm}
\proof
Propositions \ref{BasisChangeSigma} and \ref{BasisChangeLambda} imply that the change of basis functions $f_{\Basis, \Basis'}$ and $r_{\Basis, \Basis'}$ coincide, yielding
\begin{align*}
\P_{\Basis'} &= \lambda_{\Basis'} \circ \sigma_{\Basis'}^{-1} \\
&= \lambda_\Basis \circ r_{\Basis, \Basis'} \circ (\sigma_\Basis \circ f_{\Basis, \Basis'}^{-1})^{-1} \\
&= \lambda_\Basis \circ r_{\Basis, \Basis'} f_{\Basis, \Basis'} \circ \sigma_\Basis^{-1} \\
&= \lambda_\Basis \circ \sigma_\Basis^{-1} = \P_\Basis.
\end{align*}
Here we use the fact that $r_{\Basis, \Basis'}$ and $f_{\Basis, \Basis'}$ not only coincide, but they have order $1$ or $2$ (the character $\chi$ has values in $\pm 1$).  Thus $r_{\Basis, \Basis'} f_{\Basis, \Basis'} = \Id$ in $\Aut(F^\times \times_\theta T^\vee)$.
 \qed
 
From this theorem, we write simply $\P \From {}^\DEE \tilde T \xrightarrow{\sim} {}^\A \tilde T$ without reference to a basis.  Recall that $\Irr_\epsilon^\DEE(\alg{\tilde T}) = \Spl({}^\DEE \tilde T)$, and $\WP_\epsilon^{\A}(\alg{\tilde T}) = \Spl({}^\A \tilde T)$.  Thus $\P$ defines an isomorphism of $\Irr^\DEE(\alg{T})$-torsors,
$$\P(\alg{\tilde T}) \From \Irr_\epsilon^\DEE(\alg{\tilde T}) \xrightarrow{\sim} \WP_\epsilon^{\A}(\alg{\tilde T}).$$
\begin{thm}
Allowing the cover to vary, the system of isomorphisms $\P$ defines a natural isomorphism of additive functors,
$$\P \From \Cat{Cov}_n^\sharp(\alg{T}) \xRightarrow{\sim} \Cat{Tors}(\Irr^\DEE(\alg{T})).$$
\end{thm}
\proof
It suffices to work with the full subcategory of $\Cat{Cov}_n^\sharp$ given by the incarnated covers -- those that arise from elements $C \in X \otimes X$.  We have constructed isomorphisms $\P \From \Irr_\epsilon^\DEE(\alg{\tilde T}) \xrightarrow{\sim} \WP_\epsilon^{\A}(\alg{\tilde T})$ for each such incarnated cover $\alg{\tilde T}$.  For naturality, recall that all morphisms among incarnated covers are composites of automorphisms $j_{\hat u}$ for $\hat u \in X \otimes F^\times$ and isomorphisms $\iota_{\Basis, A}$ for alternating forms $A \in X \otimes X$.

Consider first an automorphism $j_{\hat u}$ and a fixed cover $\alg{\tilde T}$.  In the diagram below, Proposition \ref{AutomorphDParam} gives the commutativity of the top face and Proposition \ref{AutomorphAParam} the bottom face.  Front and back faces commute by the definition of $\P$.  Hence the right face commutes.  
$$\begin{tikzcd}
\phantom{a} & \phantom{a} & \phantom{a} & \phantom{a} & {}^\DEE \tilde T \arrow{dd}{\P} \\
\phantom{a} & F^\times \times_\theta T^\vee \arrow{urrr}[swap]{\sigma_\Basis} \arrow[dashed]{drrr}{\lambda_\Basis} & \phantom{a} & {}^\DEE \tilde T \arrow{ur}[swap]{{}^\DEE j_{\hat u}} \arrow{dd}{\P} & \phantom{a} \\
F^\times \times_\theta T^\vee \arrow{ur}{j_{\hat u}^\vee} \arrow{urrr}[swap]{\sigma_\Basis} \arrow{drrr}{\lambda_\Basis} & \phantom{a} & \phantom{a} & \phantom{a} & {}^\A \tilde T \\
\phantom{a} & \phantom{a} & \phantom{a} & {}^\A \tilde T \arrow{ur}[swap]{{}^\A j_{\hat u}} & \phantom{a}
\end{tikzcd}$$

Next consider $C_0, C \in X \otimes X$ such that $C_0(y,y) = C(y,y)$, and let $A = C - C_0$.  This gives an isomorphism $\iota_{\Basis, A} \From \alg{\tilde T}_0 \To \alg{\tilde T}$.  In the diagram below, Proposition \ref{MorphDParam} gives commutativity of the top triangle and Proposition \ref{MorphAParam} the bottom.  Left and right triangles commute by the definition of $\P$.  Hence the outer square commutes.
$$\begin{tikzcd}[row sep = 4em, column sep = 4em]
{}^\DEE \tilde T \arrow{rr}{{}^\DEE \iota_{\Basis, A}} \arrow{dd}[swap]{\P(\alg{T})} & & {}^\DEE \tilde T_0 \arrow{dd}{\P(\alg{T}_0)} \\
\phantom{a} & F^\times \times_\theta T^\vee 
\arrow{ul}[swap]{\sigma_\Basis}
\arrow{ur}{\sigma_{0,\Basis}}
\arrow{dl}{\lambda_\Basis}
\arrow{dr}[swap]{\lambda_{0,\Basis}}
& \phantom{a} \\
{}^\A \tilde T \arrow{rr}[swap]{{}^\A \iota_{\Basis, A}} & & {}^\A \tilde T_0 
\end{tikzcd}$$
By Proposition \ref{IncarnatedOK}, commutativity of the previous two diagrams implies that $\P$ is a natural isomorphism of functors.  

To check that $\P$ is a natural isomorphism of {\em additive} functors, we may use the monoidal structure on incarnated covers arising from addition of elements of $X \otimes X$.  In other words, if $\alg{\tilde T}_1$ is incarnated by $C_1 \in X \otimes X$, and $\alg{\tilde T}_2$ is incarnated by $C_2 \in X \otimes X$, then $\alg{\tilde T} = \alg{\tilde T}_1 \Baer \alg{\tilde T}_2$ may be identified with the cover incarnated by $C_1 + C_2$.  Such identifications allow us to identify the following, for incarnated covers:
$${}^\DEE \tilde T_1 \Baer {}^\DEE \tilde T_2 \ident {}^\DEE \tilde T, \quad {}^\A \tilde T_1 \Baer {}^\A \tilde T_2 \ident {}^\A \tilde T.$$
As $\P$ is defined by identifying extensions given by the same cocycle $\theta$, we find that $\P$ is compatible with the additive structure of the functors $\Irr_\epsilon^\DEE$ and $\WP_\epsilon^\A$.
\qed

The natural isomorphism $\P$ is compatible with pullbacks in the following way.  If $\alpha \From \alg{S} \To \alg{T}$ is an isomorphism of split tori, and $\alg{\tilde T} \in \Cat{Cov}_n^\sharp(\alg{T})$ is incarnated by $C \in X \otimes X$, and $\alg{\tilde S}$ is the pullback of $\alg{\tilde T}$, then $\alg{\tilde S}$ is incarnated by the pullback of $C$.  If $\Basis$ is a basis of $Y$, then we may pull back $\Basis$ to form a basis of the cocharacter lattice of $\alg{S}$.  We find a commutative diagram of groups and isomorphisms.
$$\begin{tikzcd}
\phantom{a} & \phantom{a} & \phantom{a} & \phantom{a} & {}^\DEE \tilde S \arrow{dd}{\P_S} \\
\phantom{a} & F^\times \times_\theta S^\vee \arrow{urrr}[swap]{\sigma_{\Basis,S}} \arrow[dashed]{drrr}{\lambda_{\Basis,S}} & \phantom{a} & {}^\DEE \tilde T \arrow{ur}[swap]{{}^\DEE \tilde \alpha} \arrow{dd}{\P_T} & \phantom{a} \\
F^\times \times_\theta T^\vee \arrow{ur}{\alpha^\vee} \arrow{urrr}[swap]{\sigma_{\Basis,T}} \arrow{drrr}{\lambda_{\Basis,T}} & \phantom{a} & \phantom{a} & \phantom{a} & {}^\A \tilde S \\
\phantom{a} & \phantom{a} & \phantom{a} & {}^\A \tilde T \arrow{ur}[swap]{{}^\A \tilde \alpha} & \phantom{a}
\end{tikzcd}$$
Indeed, the definition of $\sigma_\Basis$ and $\lambda_\Basis$ (for $S$ and $T$) yields the commutativity of the top and bottom squares.  The front and back faces commute by definitions of $\P_T$ and $\P_S$, respectively.  This makes the right face commute as well.  

We find a commutative diagram of additive functors and natural isomorphisms.
\begin{equation}
\label{PullbackP}
\begin{tikzcd}
\alpha^\ast \circ \Irr_\epsilon^\DEE \arrow[Rightarrow]{r}{\P} \arrow[Rightarrow]{d}{{}^\DEE \tilde \alpha} & \alpha^\ast \circ \WP_\epsilon^\A \arrow[Rightarrow]{d}{{}^\A \tilde \alpha} \\
\Irr_\epsilon^\DEE \circ \alpha^\ast \arrow[Rightarrow]{r}{\P} & \WP_\epsilon^\A \circ \alpha^\ast
\end{tikzcd}
\end{equation}

This completes the chain of additive functors and natural isomorphisms, compatible with pullbacks throughout by \eqref{PullbackDel}, \eqref{PullbackP}, \eqref{PullbackAB}.
$$\begin{tikzcd}
\rec_F^\ast \circ \can_Y \circ \Irr_\epsilon \arrow[Rightarrow]{r}{\Del} & \rec_F^\ast \circ  \Irr_\epsilon^\DEE \arrow[Rightarrow]{r}{\P} & \rec_F^\ast \circ \WP_\epsilon^A & \WP_\epsilon \arrow[Rightarrow]{l}[swap]{\AB}
\end{tikzcd}$$

Assembling commutative diagrams, we have proven the following result.
\begin{thm}
\label{ParameterizationNatural}
Define $\Lang = \rec_F^\ast \circ \can_Y$ and $\Lang_\epsilon = \AB^{-1} \circ \P \circ \Del$.  Then $\Lang_\epsilon \From \Lang \circ \Irr_\epsilon \Rightarrow \WP_\epsilon$ is a natural isomorphism of additive functors.
$$\begin{tikzcd}
\phantom{a} &  \phantom{a} & \Cat{Tors}(\Irr(\alg{T})) \arrow{dd}{\Lang} \\
\Cat{Cov}_n^\sharp(\alg{T}) \arrow{urr}{\Irr_\epsilon} \arrow{drr}[swap]{\WP_\epsilon} & \phantom{1000} \Downarrow \Lang_\epsilon  &  \phantom{a} \\
\phantom{a} & & \Cat{Tors}(\WP(\alg{T}))
\end{tikzcd}$$
When $\alpha \From \alg{S} \To \alg{T}$ is an isomorphism of split tori, the natural isomorphisms $\Lang_\epsilon$ for $\alg{S}$ and $\alg{T}$ commute with pullback as in the diagram \eqref{PullbackLang}.
\end{thm}

This theorem gives a parameterization of genuine characters by Weil parameters, functorial in the (sharp) cover, and compatibly with Langlands' parameterization of characters of algebraic tori over local fields.

\section{The integral case}

Now consider $\alg{T}$, a split torus over $\OO$, the ring of integers in a nonarchimedean local field $F$.  We may consider the category $\Cat{Cov}_n^\sharp(\alg{T}/\OO)$ of sharp covers defined over $\OO$.  By \cite{MWIntegral}, the objects of $\Cat{Cov}_n(\alg{T}/\OO)$ are classified by pairs $(Q, \sheaf{D})$ just as in the case of tori over fields, and each object is isomorphic to one incarnated by an element $C \in X \otimes X$.  The automorphism group of a cover $\alg{\tilde T} = (\alg{T}', n)$ is naturally isomorphic to $\Hom(Y, \OO^\times) = X \otimes \OO^\times$.  

Given a sharp cover $\alg{\tilde T} = (\alg{T}', n)$ over $\OO$, write $T^\circ = \alg{T}(\OO)$, $T = \alg{T}(F)$, and $\mu_n \Into \tilde T \Onto T$ for the resulting extension of locally compact abelian groups.  As $\alg{\tilde T}$ is defined over $\OO$, this comes with a splitting $\sigma^\circ \From T^\circ \Into \tilde T$.  Thus we may consider the set $\Hom_\epsilon(\tilde T / \sigma^\circ(T^\circ), \CC^\times) \subset \Hom_\epsilon(\tilde T, \CC^\times)$ of $T^\circ$-spherical $\epsilon$-genuine characters.  

\subsection{Parameterization by splittings}

We abbreviate
\begin{align*}
\Irr^\circ(\alg{T}) &= \Irr(\alg{T} / \OO) = \Hom(T / T^\circ, \CC^\times) \ident T^\vee, \\
\Irr_\epsilon^\circ(\alg{\tilde T}) &= \Irr_\epsilon(\alg{\tilde T} / \OO) = \Hom_\epsilon(\tilde T / \sigma^\circ(T^\circ), \CC^\times).
\end{align*}
Then $\Irr_\epsilon^\circ(\alg{\tilde T})$ is a $\Irr^\circ(\alg{T})$-torsor.  Write 
$$\iota^\circ \From \Irr^\circ(\alg{T}) \Into \Irr(\alg{T}) = \Irr(\alg{T} / F),$$ 
for the inclusion of the $T^\circ$-spherical characters of $T$ into the group of all characters of $T$.  There is a natural isomorphism of $\Irr(\alg{T})$-torsors,
$$(\iota^\circ)_\ast \Irr_\epsilon^\circ(\alg{\tilde T}) \xrightarrow{\sim} \Irr_\epsilon(\alg{\tilde T}).$$

We identified $\Irr_\epsilon(\alg{\tilde T})$ with the set of splittings of an extension ${}^\DEE \tilde T$.  This works as well for $T^\circ$-spherical characters as follows.  The cover $\alg{\tilde T}$ over $\OO$ yields a commutative diagram with exact rows, with top row split by $\sigma^\circ$.
$$\begin{tikzcd}[column sep = 4em]
\CC^\times \inarrow{r} \arrow{d}{=} & \epsilon_\ast \tilde T^\circ \onarrow{r} \inarrow{d} & T^\circ \arrow[bend right=20]{l}[swap]{\sigma^\circ} \inarrow{d} \\
\CC^\times \inarrow{r} & \epsilon_\ast \tilde T \onarrow{r} & T 
\end{tikzcd}$$

Recalling that $\iota \From \ZZ \Into X \otimes Y$ was the canonical inclusion, and $\DEE_Y = \iota^\ast(X \otimes \bullet)$, we may tensor the above diagram by $X$ and pull back via $\iota$ to form a commutative diagram with exact rows, again with top row split by a homomorphism we call ${}^\DEE \sigma^\circ$,
$$\begin{tikzcd}[column sep = 4.5em]
T^\vee \inarrow{r} \arrow{d}{=} & \DEE_Y \epsilon^\ast \tilde T^\circ \onarrow{r} \inarrow{d} & \OO^\times   \arrow[bend right=20]{l}[swap]{{}^\DEE \sigma^\circ} \inarrow{d} \\
T^\vee \inarrow{r} &{}^\DEE \tilde T \onarrow{r} & F^\times
\end{tikzcd}$$

Just as there is a natural bijection $\Del \From \Irr_\epsilon(\alg{\tilde T}) \To \Irr_\epsilon^\DEE(\alg{\tilde T}) = \Spl({}^\DEE \tilde T)$, the $T^\circ$-spherical characters correspond to those splittings from $\Spl({}^\DEE \tilde T)$ which pull back to the splitting ${}^\DEE \sigma^\circ$ of the top row above.  Thus define $\Irr_\epsilon^{\DEE,\circ}(\alg{\tilde T})$ to be the subset of splittings pulling back to ${}^\DEE \sigma^\circ$.  The commutative diagrams above extend to a diagram of categories and functors below.
\begin{proposition}
\label{UnramifiedDGroup}
There is a natural isomorphism of additive functors, 
$$\Del^\circ \From \can_Y^\circ \circ \Irr_\epsilon^\circ \Rightarrow \Irr_\epsilon^{\DEE,\circ}.$$
\end{proposition}
\proof
Consider the following diagram of Picard categories and additive functors.
$$\begin{tikzcd}
\phantom{a} & \phantom{a} & \Cat{Tors}(T^\vee) \arrow{dd}{\iota_\ast^\circ} \\
\Cat{Cov}_n^\sharp(\alg{T} / \OO) \arrow{dd}{\iota_\ast^\circ} \arrow{r}[swap]{\Irr_\epsilon^\circ} \arrow{rru}[near start]{\Irr_\epsilon^{\DEE, \circ}}  & \Cat{Tors}(\Irr^\circ(\alg{T})) \arrow{ur}[swap]{\can_Y^\circ} \arrow{dd}{\iota_\ast^\circ} & \\
\phantom{a} & \phantom{a} & \Cat{Tors}(\Irr^\DEE(\alg{T})) \\
\Cat{Cov}_n^\sharp(\alg{T} / F) \arrow{r}[swap]{\Irr_\epsilon} \arrow{rru}[near start]{\Irr_\epsilon^{\DEE}}  & \Cat{Tors}(\Irr(\alg{T})) \arrow{ur}[swap]{\can_Y} & \\
\end{tikzcd}$$

The bottom face is the diagram \eqref{TautologicalParam}, in which we found a natural isomorphism
$$\Del \From \can_Y \circ \Irr_\epsilon \Rightarrow \Irr_\epsilon^\DEE.$$
The lateral faces of the prism 2-commute by the commutative diagrams discussed just above.  There are natural isomorphisms
$$\alpha \From \iota_\ast^\circ \circ \can_Y^\circ \circ \Irr_\epsilon^\circ \Rightarrow \can_Y \circ \Irr_\epsilon \circ \iota_\ast^\circ, \quad \beta \From \iota_\ast^\circ \circ \Irr_\epsilon^{\DEE, \circ} \Rightarrow \Irr_\epsilon^\DEE \circ \iota_\ast^\circ.$$
There is a unique natural isomorphism $\Del^\circ \From \can_Y^\circ \circ \Irr_\epsilon^\circ \Rightarrow \Irr_\epsilon^{\DEE,\circ}$ making the following diagram of functors and natural isomorphisms commute.
$$\begin{tikzcd}
\iota_\ast^\circ \circ \can_Y^\circ \circ \Irr_\epsilon^\circ \arrow[Rightarrow]{r}{\alpha} \arrow[Rightarrow]{d}{\Del^\circ} &  \can_Y \circ \Irr_\epsilon \circ \iota_\ast^\circ  \arrow[Rightarrow]{d}{\Del}  \\
\iota_\ast^\circ \circ \Irr_\epsilon^{\DEE, \circ} \arrow[Rightarrow]{r}{\beta}  & \Irr_\epsilon^\DEE \circ \iota_\ast^\circ 
\end{tikzcd}$$ 
\qed

\subsection{Unramified Weil parameters}

We continue to write $\iota^\circ \From \OO \Into F$ for the inclusion, and for all maps and functors resulting from this inclusion.  Given a sharp cover $\alg{\tilde T}$ over $\OO$, and fixing $\epsilon \From \mu_n \Into \CC^\times$ as always, we have constructed the L-group
$$T^\vee \Into {}^\EL T^\circ \Onto \Gal_\OO = \langle \Frob \rangle_{\prof}.$$
This is an extension of {\em abelian} groups, defined uniquely up to unique isomorphism from $\alg{\tilde T}$ and $\epsilon$.  

The construction of the L-group is compatible with base extension from $\OO$ to $F$ (see, e.g., Section \ref{BaseChangeGerbe}), and we find a commutative diagram,
\begin{equation}
\label{LgroupUnramifiedCompatible}
\begin{tikzcd}[column sep = 5.5em]
T^\vee \inarrow{r} \arrow{d}{=} & {}^\EL \tilde T \onarrow{r} \onarrow{d} & \Gal_F \onarrow{d} \arrow[bend right=20, dashed]{l}[swap]{\tau^\circ} \\
T^\vee \inarrow{r} & {}^\EL \tilde T^\circ \onarrow{r} & \Gal_\OO
\end{tikzcd}
\end{equation}
This identifies the complex L-group ${}^\EL \tilde T$ with the pullback of ${}^\EL \tilde T^\circ$.  As such, the top row is endowed with a splitting $\tau^\circ \From \Inertia_F \To {}^\EL \tilde T$ over the inertial subgroup of $\Gal_F$.

The set of Weil parameters $\WP_\epsilon(\alg{\tilde T}) = \WP_\epsilon(\alg{\tilde T} / F)$ consists of continuous homomorphisms $\rho \From \Weil_F \To {}^\EL \tilde T$, lying over the canonical map $\Weil_F \To \Gal_F$.  A subset is given by the unramified Weil parameters $\WP_\epsilon^\circ(\alg{\tilde T}) = \WP_\epsilon(\alg{\tilde T} / \OO)$, consisting of those $\rho \From \Weil_F \To {}^\EL \tilde T$ whose restriction to inertia coincides with $\tau^\circ$.  The unramified parameters $\WP_\epsilon^\circ(\alg{\tilde T})$ are identified, in turn, with the homomorphisms $\Weil_\OO = \langle \Frob \rangle \To {}^\EL \tilde T^\circ$ lying over the canonical map $\Weil_\OO \To \Gal_\OO$.  The set $\WP_\epsilon^\circ(\alg{\tilde T})$ is naturally a torsor for $\WP^\circ(\alg{T}) = \Hom(\Weil_\OO, T^\vee)$, the set of unramified parameters for $\alg{T}$.

These sets of Weil parameters, $\WP_\epsilon^\circ(\alg{\tilde T}) \subset \WP_\epsilon(\alg{\tilde T})$ can be identified with splittings sequences of exact sequences as follows:  the canonical splitting of ${}^\EL \tilde T$ over $[\Gal_F, \Gal_F]$ coincides with the splitting $\tau^\circ$ on $[\Gal_F, \Gal_F] \cap \Inertia$.  In this way, the previous commutative diagram \eqref{LgroupUnramifiedCompatible} gives a commutative diagram,
$$\begin{tikzcd}
T^\vee \inarrow{r} \arrow{d}{=} & {}^\EL \tilde T / [\Gal_F, \Gal_F] \onarrow{r} \onarrow{d} & \Gal_F^{\ab} \onarrow{d} \\
T^\vee \inarrow{r} & {}^\EL \tilde T^\circ \onarrow{r} & \Gal_\OO
\end{tikzcd}$$

Pulling back via $\Weil_F^{\ab} \Into \Gal_F^{\ab}$, and $\Weil_\OO \Into \Gal_\OO$, and via the reciprocity isomorphisms $\rec_F \From \Weil_F^{\ab} \xrightarrow{\sim} F^\times$, $\rec_\OO \From \Weil_\OO \xrightarrow{\sim} \ZZ$, we obtain a commutative diagram of locally compact abelian groups,
$$\begin{tikzcd}
T^\vee \inarrow{r} \arrow{d}{=} & {}^\A \tilde T \onarrow{r} \onarrow{d} & F^\times \onarrow{d}{\val} \\
T^\vee \inarrow{r} & {}^\A \tilde T^\circ \onarrow{r} & \ZZ
\end{tikzcd}$$

This identifies ${}^\A \tilde T$ wth the pullback of ${}^\A \tilde T^\circ$ via $\val \From F^\times \To \ZZ$.  In particular, ${}^\A \tilde T$ is endowed with a canonical splitting ${}^\A \lambda^\circ \From \OO^\times \To {}^\A \tilde T$.  Recall that $\WP_\epsilon^\A(\alg{\tilde T})$ denotes the set of splittings $\Spl({}^\A \tilde T)$, and define $\WP_\epsilon^{\A, \circ}(\alg{\tilde T})$ to be the set of splittings of ${}^\A \tilde T^\circ$.  Equivalently, $\WP_\epsilon^{\A, \circ}(\alg{\tilde T})$ is the set of splittings in $\Spl({}^\A \tilde T)$ which restrict to ${}^\A \lambda^\circ$ on $\OO^\times$.  This gives an inclusion
$\WP_\epsilon^{\A,\circ}(\alg{\tilde T}) \Into \WP_\epsilon^{\A}(\alg{\tilde T})$.
\begin{proposition}
There is a natural isomorphism of additive functors,
$$\AB^\circ \From \WP_\epsilon^\circ \xRightarrow{\sim} \rec_\OO^\ast \circ \WP_\epsilon^{\A, \circ}.$$
\end{proposition}
\proof
Consider the following diagram of Picard categories and additive functors.
$$
\begin{tikzcd}
\phantom{a} & \phantom{a} & \Cat{Tors}(T^\vee) \arrow{dd}{\iota_\ast^\circ} \arrow{dl}{\rec_\OO^\ast} \\
\Cat{Cov}_n^\sharp(\alg{T} / \OO) \arrow{dd}{\iota_\ast^\circ} \arrow{r}[swap]{\WP_\epsilon^\circ} \arrow{rru}[near start]{\WP_\epsilon^{\A, \circ}}  & \Cat{Tors}(\WP^\circ(\alg{T}))  \arrow{dd}{\iota_\ast^\circ} & \\
\phantom{a} & \phantom{a} & \Cat{Tors}(\Irr^\DEE(\alg{\tilde T})) \arrow{dl}{\rec_F^\ast} \\
\Cat{Cov}_n^\sharp(\alg{T} / F) \arrow{r}[swap]{\WP_\epsilon} \arrow{rru}[near start]{\WP_\epsilon^{\A}}  & \Cat{Tors}(\WP(\alg{T}))  & \\
\end{tikzcd}
$$
The proof mirrors Proposition \ref{UnramifiedDGroup}.  The bottom face 2-commutes via the natural isomorphism $\AB$ from \eqref{ParametersAsSplittings}.  The lateral faces 2-commute by examining the commutative diagrams just above.  This family of natural isomorphisms determines a natural isomorphism $\AB^{\circ}$ making the top face 2-commute as well.
\qed

\subsection{Spherical/Unramified comparison}

The connection between genuine characters and Weil parameters was made through a sequence of three natural isomorphisms,
$$\begin{tikzcd}
\rec_F^\ast \circ \can_Y \circ \Irr_\epsilon \arrow[Rightarrow]{r}{\Del} & \rec_F^\ast \circ  \Irr_\epsilon^\DEE \arrow[Rightarrow]{r}{\P} & \rec_F^\ast \circ \WP_\epsilon^A & \WP_\epsilon  \arrow[Rightarrow]{l}[swap]{\AB}.
\end{tikzcd}$$

We have found a spherical analogue of $\Del$ and an unramified analogue of $\AB$, giving the following sequence of functors and natural isomorphisms,
$$\begin{tikzcd}
\rec_\OO^\ast \circ \can_Y \circ \Irr_\epsilon^\circ \arrow[Rightarrow]{r}{\Del^\circ} & \rec_\OO^\ast \circ  \Irr_\epsilon^{\DEE,\circ}  & \rec_\OO^\ast \circ \WP_\epsilon^{\A,\circ}  & \WP_\epsilon^\circ \arrow[Rightarrow]{l}[swap]{\AB^\circ}.
\end{tikzcd}$$

The natural isomorphism $\P \From \Irr_\epsilon^\DEE \Rightarrow \WP_\epsilon^\A$ arose from the full subcategory of incarnated covers, where any basis of $Y$ yields a pair of isomorphisms $\sigma_\Basis$ and $\lambda_\Basis$,
$$\begin{tikzcd}
\phantom{a} & {}^\DEE \tilde T \arrow{dd}{\P} \\
F^\times \times_\theta T^\vee \arrow{ur}{\sigma_\Basis} \arrow{dr}[swap]{\lambda_\Basis} & \phantom{a} \\
\phantom{a} & {}^\A \tilde T
\end{tikzcd}$$
A unique isomorphism $\P$ (independent of the choice of basis) makes this diagram commute; it yields a natural isomorphism from $\Spl({}^\DEE \tilde T)$ to $\Spl({}^\A \tilde T)$.

For sharp covers $\alg{\tilde T}$ over $\OO$ incarnated by $C \in X \otimes X$, the same steps yield isomorphisms $\sigma_\Basis^\circ$ and $\lambda_\Basis^\circ$,
$$\begin{tikzcd}
\phantom{a} & {}^\DEE \tilde T \arrow{dd}{\P} \\
\OO^\times \times_\theta T^\vee \arrow{ur}{\sigma_\Basis^\circ} \arrow{dr}[swap]{\lambda_\Basis^\circ} & \phantom{a} \\
\phantom{a} & {}^\A \tilde T
\end{tikzcd}$$
As the Hilbert symbol of order $n$ is trivial on $\OO^\times \times \OO^\times$ (as $\alg{\tilde T}$ is defined over $\OO$, $n$ is coprime to the residue characteristic), we find that $\OO^\times \times_\theta T^\vee = \OO^\times  \times T^\vee$.  In this way, $\sigma_\Basis^\circ$ is a splitting of ${}^\DEE \tilde T^\circ$ over $\OO^\times$, and $\lambda_\Basis^\circ$ is a splitting of ${}^\A \tilde T^\circ$ over $\OO^\times$.  These coincide with the canonical splittings ${}^\DEE \sigma^\circ$ and ${}^\A \lambda^\circ$ described in the previous two sections.  

In this way, the isomorphism $\P$ sets up a bijection between those splittings of ${}^\DEE \tilde T$ which restrict to ${}^\DEE \sigma^\circ$, and those splittings of ${}^\A \tilde T$ which restrict to ${}^\A \lambda^\circ$.  In other words, $\P$ restricts to a bijection
$$\P^\circ \From \Irr_\epsilon^{\DEE, \circ}(\alg{\tilde T}) \xrightarrow{\sim} \WP_\epsilon^{\A, \circ}(\alg{\tilde T}).$$

\begin{proposition}
The following diagram of Picard categories and additive functors 2-commutes.
$$
\begin{tikzcd}[row sep = 5em]
\Cat{Cov}_n^\sharp(\alg{T} / \OO) \arrow{d}{\iota_\ast^\circ} \arrow[bend left=20]{rr}[swap]{\Irr_\epsilon^{\DEE,\circ}} \arrow[bend right=20]{rr}{\WP_\epsilon^{\A, \circ}}  & & \Cat{Tors}(\Hom(F^\times / \OO^\times, T^\vee)) \arrow{d}{\iota_\ast^\circ} \\
\Cat{Cov}_n^\sharp(\alg{T} / F) \arrow[bend left=20]{rr}[swap]{\Irr_\epsilon^\DEE} \arrow[bend right=20]{rr}{\WP_\epsilon^{\A}}  & & \Cat{Tors}(\Hom(F^\times, T^\vee)) 
\end{tikzcd}
$$
\end{proposition}
\proof
The natural isomorphism $\P$ makes the bottom face of the cylinder 2-commute.  As $\P$ pulls back to a bijection $\P^\circ$ from $\Irr_\epsilon^{\DEE, \circ}(\alg{\tilde T})$ to $\WP_\epsilon^{\A, \circ}(\alg{\tilde T})$, for any sharp cover $\alg{\tilde T}$, we find a unique natural isomorphism $\P^\circ \From \Irr_\epsilon^{\DEE, \circ} \Rightarrow \WP_\epsilon^{\A, \circ}$ making the diagram 2-commute.
\qed

Assembling these propositions, we find that our parameterization makes spherical characters correspond to unramified parameters.  We have two sequences of functors and natural isomorphisms,
$$\begin{tikzcd}
\rec_F^\ast \circ \can_Y \circ \Irr_\epsilon \arrow[Rightarrow]{r}{\Del} & \rec_F^\ast \circ  \Irr_\epsilon^\DEE \arrow[Rightarrow]{r}{\P} & \rec_F^\ast \circ \WP_\epsilon^A  \arrow[Rightarrow]{r}{\AB^{-1}} & \WP_\epsilon \\
\rec_\OO^\ast \circ \can_Y^\circ \circ \Irr_\epsilon^\circ \arrow[Rightarrow]{r}{\Del^\circ} & \rec_\OO^\ast \circ  \Irr_\epsilon^{\DEE,\circ}  \arrow[Rightarrow]{r}{\P^\circ} & \rec_\OO^\ast \circ \WP_\epsilon^{\A,\circ}  \arrow[Rightarrow]{r}{(\AB^\circ)^{-1}} & \WP_\epsilon^\circ.
\end{tikzcd}$$
Recall that $\Lang = \rec_F^\ast \circ \can_Y \From \Hom(T, \CC^\times) \To \Hom(\Weil_F, T^\vee)$.  Define analogously $\Lang^\circ = \rec_\OO^\ast \circ \can_Y^\circ \From \Hom(T / T^\circ, \CC^\times) \To \Hom(\Weil_\OO, T^\vee)$.  Recall that $\Lang_\epsilon = \AB^{-1} \circ \P \circ \Del$, and define analogously $\Lang_\epsilon^\circ = (\AB^\circ)^{-1} \circ \P^\circ \circ \Del^\circ$.

\begin{thm}
\label{UnramToriParam}
The following diagram of Picard categories and additive functors 2-commutes, via the natural isomorphisms $\Lang_\epsilon^\circ$ (on the top face) and $\Lang_\epsilon$ (on the bottom face).
$$
\begin{tikzcd}[column sep = 3.5em]
\phantom{a} & \phantom{a} & \Cat{Tors}(\WP^\circ(\alg{T})) \arrow{dd}{\iota_\ast^\circ} \\
\Cat{Cov}_n^\sharp(\alg{T} / \OO) \arrow{dd}{\iota_\ast^\circ} \arrow{r}[swap]{\WP_\epsilon^\circ} \arrow{rru}[near start]{\Irr_\epsilon^\circ}  & \Cat{Tors}(\Irr^\circ(\alg{T})) \arrow{ur}[swap]{\Lang^\circ} \arrow{dd}{\iota_\ast^\circ} & \\
\phantom{a} & \phantom{a} & \Cat{Tors}(\WP(\alg{T})) \\
\Cat{Cov}_n^\sharp(\alg{T} / F) \arrow{r}[swap]{\WP_\epsilon} \arrow{rru}[near start]{\Irr_\epsilon}  & \Cat{Tors}(\Irr(\alg{T})) \arrow{ur}[swap]{\Lang} & \\
\end{tikzcd}
$$
\end{thm}

This theorem gives a parameterization of spherical genuine characters of $\tilde T$ by unramified Weil parameters, functorial in the choice of sharp cover, and compatible with the parameterization of all genuine characters by all Weil parameters.  As in the case of local fields, one may verify that the natural isomorphisms $\Lang_\epsilon^\circ$ are compatible with pullbacks, for isomorphisms of split tori $\alg{S} \To \alg{T}$ over $\OO$.  

\section{Global case}

Now we consider sharp covers $\alg{\tilde T}$ of a split torus $\alg{T}$ over a {\em global} field $F$.  The methods are much the same as the previous section, with the inclusion $\OO \Into F$ replaced by the inclusion $F \Into \AA$.  Thus we leave a few details to the reader, to adapt proofs from the previous section as needed.

Given a sharp cover $\alg{\tilde T} = (\alg{T}', n)$ over $F$, write $T_F = \alg{T}(F)$ and $T_\AA = \alg{T}(\AA)$, and $\mu_n \Into \tilde T_\AA \Onto T_\AA$ for the resulting extension of locally compact abelian groups.  As $\alg{\tilde T}$ is defined over $F$, this comes with a splitting $\sigma_F \From T_F \Into \tilde T_\AA$.

\subsection{Parameterization by splittings}

Define $\Irr_\AA(\alg{T}) = \Hom(T_\AA, \CC^\times)$, the group of continuous characters, and $\Irr_{\AA,\epsilon}(\alg{\tilde T}) = \Hom_\epsilon(\tilde T_\AA, \CC^\times)$ for the $\Irr_\AA(\alg{T})$-torsor of $\epsilon$-genuine continuous characters.  To give such a character, it is equivalent to give genuine characters of $\tilde T_v$ for all places $v$ of $F$, almost all of which are $T_v^\circ$-spherical; see \cite{MWToriOld} and \cite[\S 4]{MWToriNew} for details.  

We abbreviate,
\begin{align*}
\Irr_F(\alg{T}) &= \Irr(\alg{T} / F) = \Hom(T_\AA / T_F, \CC^\times), \\
\Irr_{F,\epsilon}(\alg{\tilde T}) &= \Irr_\epsilon(\alg{\tilde T} / F) =  \Hom_\epsilon(\tilde T / \sigma_F(T), \CC^\times).
\end{align*}
Then $\Irr_{F, \epsilon}(\alg{\tilde T})$ is a $\Irr_F(\alg{T})$-torsor.  Write $\iota_F \From \Irr_F(\alg{T}) \Into \Irr_\AA(\alg{T})$ for the inclusion of the automorphic characters of $T_\AA$ into the group of all continuous characters of $T_\AA$.  There is a natural isomorphism of $\Irr(\alg{T})$-torsors,
$$(\iota_F)_\ast \left(\Irr_{F,\epsilon}(\alg{\tilde T}) \right) \xrightarrow{\sim} \Irr_{\AA,\epsilon}(\alg{\tilde T}).$$

Define
$$\Irr_\AA^\DEE(\alg{T}) = \Hom(\AA^\times, T^\vee), \quad \Irr_F^\DEE(\alg{T}) = \Hom(\AA^\times / F^\times, T^\vee),$$
so there are natural isomorphisms
$$\can_{\AA,Y} \From \Irr_\AA(\alg{T}) \xrightarrow{\sim} \Irr_\AA^\DEE(\alg{T}), \quad \can_{F,Y} \From \Irr_F(\alg{T}) \xrightarrow{\sim} \Irr_F^\DEE(\alg{T}).$$

We may identify $\Irr_{\AA,\epsilon}(\alg{\tilde T})$ with the set of splittings of an extension ${}^\DEE \tilde T_\AA$, just as in the local case.  Indeed, begin with the extension $\mu_n \Into \tilde T_\AA \Onto T_\AA$.  Push out via $\epsilon$ to get an extension $\CC^\times \Into \epsilon_\ast \tilde T_\AA \Onto T_\AA$.  As $T_\AA = Y \otimes \AA^\times$, we may tensor with $X$ and pull back via $\ZZ \Into X \otimes Y$ to obtain an extension,
$$T^\vee \Into {}^\DEE \tilde T_\AA \Onto \AA^\times.$$
Define $\Irr_{\AA,\epsilon}^\DEE(\alg{\tilde T}) = \Spl({}^\DEE \tilde T_\AA)$, a $\Irr_\AA^\DEE(\alg{T})$-torsor.  As in the local case, there is a natural bijection $\Del_\AA \From \Irr_{\AA,\epsilon}(\alg{\tilde T}) \To \Irr_{\AA,\epsilon}^\DEE(\alg{\tilde T})$, and this is compatible with pullbacks to local fields.  The splitting $\sigma_F$ defines a splitting ${}^\DEE \sigma_F$ in the following commutative diagram.  
$$\begin{tikzcd}[column sep = 4.5em]
T^\vee \inarrow{r} \arrow{d}{=} & {}^\DEE \tilde T_F \onarrow{r} \inarrow{d} & F^\times   \arrow[bend right=20]{l}[swap]{{}^\DEE \sigma_F} \inarrow{d} \\
T^\vee \inarrow{r} & {}^\DEE \tilde T_\AA \onarrow{r} & \AA^\times
\end{tikzcd}$$

The automorphic genuine characters $\Irr_{F,\epsilon}(\alg{\tilde T}) \subset \Irr_\epsilon(\alg{\tilde T})$ correspond to the subset $\Irr_{F,\epsilon}^\DEE(\alg{\tilde T}) \subset \Irr_{\AA,\epsilon}^\DEE(\alg{\tilde T})$ of splittings which restrict to ${}^\DEE \sigma_F$ on $F^\times$.  Then $\Irr_{F,\epsilon}^\DEE(\alg{\tilde T})$ is a $\Irr_F^\DEE(\alg{T})$-torsor.  The commutative diagrams above extend to a diagram of categories and functors below.
\begin{proposition}
There is a natural isomorphism of additive functors,
$$\Del_F \From \can_{F,Y} \circ \Irr_{F, \epsilon} \xRightarrow{\sim} \Irr_{F, \epsilon}^\DEE.$$
\end{proposition}
\proof
Consider the following diagram of Picard categories and additive functors.
$$\begin{tikzcd}
\phantom{a} & \phantom{a} & \Cat{Tors}( \Irr_F^\DEE(\alg{T})) \arrow{dd}{\iota_{F,\ast}} \\
\Cat{Cov}_n^\sharp(\alg{T}) \arrow{dd}{=} \arrow{r}[swap]{\Irr_{F,\epsilon}} \arrow{rru}[near start]{\Irr_{F,\epsilon}^{\DEE}}  & \Cat{Tors}(\Irr_F(\alg{T})) \arrow{ur}[swap]{\can_{F,Y}} \arrow{dd}{\iota_{F,\ast}} & \\
\phantom{a} & \phantom{a} & \Cat{Tors}(\Irr_\AA^\DEE(\alg{T})) \\
\Cat{Cov}_n^\sharp(\alg{T}) \arrow{r}[swap]{\Irr_{\AA,\epsilon}} \arrow{rru}[near start]{\Irr_{\AA,\epsilon}^{\DEE}}  & \Cat{Tors}(\Irr_\AA(\alg{T})) \arrow{ur}[swap]{\can_{\AA,Y}} & \\
\end{tikzcd}$$

The natural isomorphisms $\Del$ and $\Del^\circ$ define a natural isomorphism $\Del_\AA$ making the bottom face 2-commute.  Indeed, $\Irr_{\AA,\epsilon}(\alg{\tilde T})$ gives the set of continuous genuine characters of $\tilde T_\AA$ -- these are described in turn by families of genuine characters of $\tilde T_v$ for all places, almost all of which are $T_v^\circ$-spherical.  But these are described, by splittings of ${}^\DEE \tilde T_v$ for all places, coinciding with the splitting $\sigma_v^\circ$ at all $T_v^\circ$-spherical places.

As in the local integral case, the lateral faces 2-commute, from which we find a unique natural isomorphism $\Del_F$ making the top face 2-commute by pulling back $\Del_\AA$.
\qed

\subsection{Global Weil parameters}

Given a sharp cover $\alg{\tilde T}$ over $F$, we have constructed the L-group,
$$T^\vee \Into {}^\EL T \Onto \Gal_F.$$
The construction of the L-group is compatible with base extension from $F$ to $F_v$, and we find a commutative diagram for every place $v$ of $F$,
\begin{equation}
\label{LgroupUnramifiedCompatible}
\begin{tikzcd}
T^\vee \inarrow{r} \arrow{d}{=} & {}^\EL \tilde T_v \inarrow{r} \onarrow{d} & \Gal_{F_v} \inarrow{d} \\
T^\vee \inarrow{r} & {}^\EL \tilde T \onarrow{r} & \Gal_F
\end{tikzcd}
\end{equation}
This identifies the L-group ${}^\EL \tilde T_v$ with the pullback of ${}^\EL \tilde T$.  

The set of Weil parameters $\WP_{F,\epsilon}(\alg{\tilde T}) = \WP_\epsilon(\alg{\tilde T} / F)$ consists of continuous homomorphisms $\rho \From \Weil_F \To {}^\EL \tilde T$, lying over the canonical map $\Weil_F \To \Gal_F$.  As in the local case, the cohomology $H^2(\Gal_F, T^\vee)$ is trivial, yielding the (noncanonical) splitting of ${}^\EL \tilde T$.  From this it follows that the L-groups split canonically over $[\Gal_F, \Gal_F]$.  Taking the quotient ${}^\EL \tilde T / [\Gal_F, \Gal_F]$ and pulling back via $\AA^\times / F^\times \xrightarrow{\rec_F^{-1}} \Weil_F^{\ab} \To \Gal_F^{\ab}$ yields an extension
$$T^\vee \Into {}^\A \tilde T_F \Onto \AA^\times / F^\times.$$
Pulling back via $F_v^\times \Into \AA^\times$ gives a commutative diagram,
$$\begin{tikzcd}
T^\vee \inarrow{r} \arrow{d}{=} & {}^\A \tilde T_v \onarrow{r} \arrow{d} & F_v^\times \arrow{d} \\
T^\vee \inarrow{r} & {}^\A \tilde T_F \onarrow{r} & \AA^\times / F^\times
\end{tikzcd}$$

The local extensions $T^\vee \Into {}^\A \tilde T_v \Onto F_v^\times$ are endowed with splittings $\sigma_v^\circ$ over $\OO_v^\times$ for almost all $v$.  Define the restricted product,
$${}^\A \tilde T_\AA = \frac{ \{ (\tilde t_v)_{v \in \VV} \in \prod_v {}^\A \tilde T_v : \tilde t_v \in \sigma_v^\circ(\OO_v^\times) \text{ for almost all } v \} }{ \Ker \left( \bigoplus_v T^\vee \xrightarrow{\Pi} T^\vee \right) }.$$
This fits into a commutative diagram with exact rows,
$$\begin{tikzcd}
T^\vee \inarrow{r} \arrow{d}{=} & {}^\A \tilde T_\AA \onarrow{r} \onarrow{d} & \AA^\times \onarrow{d} \\
T^\vee \inarrow{r} & {}^\A \tilde T_F \onarrow{r} & \AA^\times / F^\times
\end{tikzcd}$$
The construction of ${}^\A \tilde T_\AA$ gives a splitting ${}^\A \lambda_F \From F^\times \To {}^\A \tilde T_\AA$ of the top row over $F^\times$

Define $\WP_{F,\epsilon}^\A(\alg{\tilde T}) = \Spl({}^\A \tilde T_F)$; this is in natural bijection with the set $\WP_{F,\epsilon}(\alg{\tilde T})$ of Weil parameters for ${}^L \tilde T_F$.  Analogously, define $\WP_{\AA,\epsilon}^\A(\alg{\tilde T}) = \Spl( {}^\A \tilde T_\AA)$.  Pullback yields an inclusion $\Spl({}^\A \tilde T_F) \Into \Spl( {}^\A \tilde T_\AA)$, whose image consists of splittings of ${}^\A \tilde T_\AA$ which restrict to the splitting ${}^\A \lambda_F$ on $F^\times$. 

Tracing through the definitions, $\Spl( {}^\A \tilde T_\AA)$ is in natural bijection with the set of families $(\phi_v)_{v \in \VV}$ of Weil parameters in $\WP_{v, \epsilon}(\alg{\tilde T})$ at every place, for which $\phi_v$ is unramified almost everywhere.  

Define $\WP_F(\alg{T})$ to be the set of Weil parameters for $\alg{T}$, i.e. $\WP_F(\alg{T}) = \Hom(\Weil_F, T^\vee)$.  Define $\WP_\AA(\alg{T})$ to be the set of families $(\phi_v)_{v \in \VV}$ of Weil parameters in $\WP_v(\alg{T}) = \Hom(\Weil_{F_v}, T^\vee)$, which are unramified almost everywhere.

\begin{proposition}
There is a natural isomorphism of additive functors,
$$\AB_F  \From \WP_{F, \epsilon} \xRightarrow{\sim} \rec_F^\ast \circ \WP_{F,\epsilon}^\A.$$
\end{proposition}
\proof
Consider the following diagram of Picard categories and additive functors.
$$\begin{tikzcd}
\phantom{a} & \phantom{a} & \Cat{Tors}(\Irr_F^\DEE(\alg{T})) \arrow{dd}{\iota_{F,\ast}} \arrow{dl}{\rec_F^\ast}\\
\Cat{Cov}_n^\sharp(\alg{T}) \arrow{dd}{=} \arrow{r}[swap]{\WP_{F,\epsilon}} \arrow{rru}[near start]{\WP_{F,\epsilon}^{\A}}  & \Cat{Tors}(\WP_F(\alg{T}))  \arrow{dd}{\iota_{F,\ast}} & \\
\phantom{a} & \phantom{a} & \Cat{Tors}(\Irr_\AA^\DEE(\alg{T})) \arrow{dl}{\prod_v \rec_v^\ast} \\
\Cat{Cov}_n^\sharp(\alg{T}) \arrow{r}[swap]{\WP_{\AA,\epsilon}} \arrow{rru}[near start]{\WP_{\AA,\epsilon}^{\A}}  & \Cat{Tors}(\WP_\AA(\alg{T}))  & \\
\end{tikzcd}$$
The natural isomorphisms $\AB$ and $\AB^\circ$ in the local and unramified cases yield a natural isomorphism $\AB_\AA$ which makes the bottom face 2-commute.  The lateral faces 2-commute by the local-global compatibility of our constructions.  A natural isomorphism $\AB_F$ making the top face 2-commute follows.
\qed

\subsection{Global comparison}

Over local fields, we have found a connection between genuine characters and Weil parameters through a sequence of three natural isomorphisms,
$$\begin{tikzcd}
\rec_{F_v}^\ast \circ \can_Y \circ \Irr_\epsilon \arrow[Rightarrow]{r}{\Del_v} & \rec_{F_v}^\ast \circ  \Irr_\epsilon^\DEE \arrow[Rightarrow]{r}{\P_v} & \rec_{F_v}^\ast \circ \WP_{v,\epsilon}^\A  \arrow[Rightarrow]{r}{\AB_v^{-1}} & \WP_{v,\epsilon} 
\end{tikzcd}$$
Each term has an unramified counterpart.
$$\begin{tikzcd}
\rec_{\OO_v}^\ast \circ \can_Y \circ \Irr_\epsilon^\circ \arrow[Rightarrow]{r}{\Del_v^\circ} & \rec_{\OO_v}^\ast \circ  \Irr_{v,\epsilon}^{\DEE,\circ}   \arrow[Rightarrow]{r}{\P_v^\circ} & \rec_{\OO_v}^\ast \circ \WP_{v,\epsilon}^{\A,\circ}  \arrow[Rightarrow]{r}{(\AB_v^\circ)^{-1}} & \WP_{v,\epsilon}^\circ 
\end{tikzcd}$$
Local-global compatibility gives an adelic version.
$$\begin{tikzcd}
\rec_\AA^\ast \circ \can_Y \circ \Irr_\epsilon \arrow[Rightarrow]{r}{\Del_\AA} & \rec_\AA^\ast \circ  \Irr_{\AA,\epsilon}^\DEE \arrow[Rightarrow]{r}{\P_\AA} & \rec_\AA^\ast \circ \WP_{\AA,\epsilon}^\A  \arrow[Rightarrow]{r}{\AB_\AA^{-1}} & \WP_{\AA,\epsilon} 
\end{tikzcd}$$

We have found an automorphic version, as below.
$$\begin{tikzcd}
\rec_F^\ast \circ \can_Y \circ \Irr_\epsilon \arrow[Rightarrow]{r}{\Del_F} & \rec_F^\ast \circ  \Irr_{F,\epsilon}^\DEE  & \rec_\AA^\ast \circ \WP_{F,\epsilon}^\A  \arrow[Rightarrow]{r}{\AB_F^{-1}} & \WP_{F,\epsilon} 
\end{tikzcd}$$

To link the middle terms, it suffices as before to consider the incarnated covers, on which any basis of $Y$ yields a pair of isomorphisms $\sigma_{\AA,\Basis}$ and $\lambda_{\AA,\Basis}$,
$$\begin{tikzcd}
\phantom{a} & {}^\DEE \tilde T_\AA \arrow{dd}{\P_\AA} \\
\AA^\times \times_\theta T^\vee \arrow{ur}{\sigma_{\AA,\Basis}} \arrow{dr}[swap]{\lambda_{\AA,\Basis}} & \phantom{a} \\
\phantom{a} & {}^\A \tilde T_\AA
\end{tikzcd}$$
A unique isomorphism $\P_{\AA}$ makes this diagram commute, and is independent of basis; it yields an isomorphism from $\Spl({}^\DEE \tilde T_\AA)$ to $\Spl({}^\A \tilde T_\AA)$, from which we find the natural isomorphism $\P_\AA \From \Irr_{\AA,\epsilon}^\DEE \Rightarrow \WP_{\AA,\epsilon}^\A$.

As the global Hilbert symbol of order $n$ is trivial on $F^\times \times F^\times$ (Hilbert reciprocity), we find that $F^\times \times_\theta T^\vee = F^\times  \times T^\vee$.  In this way, $\sigma_{\AA,\Basis}$ gives a splitting $\sigma_{F,\Basis}$ of ${}^\DEE \tilde T_\AA$ over $F^\times$, and $\lambda_{\AA,\Basis}$ gives a splitting $\lambda_{F,\Basis}$ of ${}^\A \tilde T_\AA$ over $F^\times$.  These coincide with the canonical splittings ${}^\DEE \sigma_F$ and ${}^\A \lambda_F$ described before. 

In this way, the isomorphism $\P_\AA$ sets up a bijection $\P_F$ between those splittings of ${}^\DEE \tilde T_\AA$ which restrict to ${}^\DEE \sigma_F$ on $F^\times$, and those splittings of ${}^\A \tilde T_\AA$ which restrict to ${}^\A \lambda_F$ on $F^\times$.  In other words, $\P_F$ gives a bijection,
$$\P_F \From  \Irr_{F,\epsilon}^{\DEE}(\alg{\tilde T}) \xrightarrow{\sim} \WP_{F,\epsilon}^{\A}(\alg{\tilde T}).$$
Naturality of this bijection is the following.
\begin{proposition}
There is a natural isomorphism of additive functors, $\P_F \From  \Irr_{F,\epsilon}^{\DEE} \Rightarrow \WP_{F,\epsilon}^{\A}$.
\end{proposition}
\proof
Consider the following diagram of Picard categories and additive functors.
$$
\begin{tikzcd}[row sep = 5em]
\Cat{Cov}_n^\sharp(\alg{T}) \arrow{d}{=} \arrow[bend left=20]{rr}[swap]{\Irr_{F,\epsilon}^\DEE} \arrow[bend right=20]{rr}{\WP_{F,\epsilon}^{\A}}  & & \Cat{Tors}(\Hom(\AA^\times / F^\times, T^\vee)) \arrow{d}{\iota_{F,\ast}} \\
\Cat{Cov}_n^\sharp(\alg{T}) \arrow[bend left=20]{rr}[swap]{\Irr_{\AA,\epsilon}^\DEE} \arrow[bend right=20]{rr}{\WP_{\AA,\epsilon}^{\A}}  & & \Cat{Tors}(\Hom(\AA^\times, T^\vee)) 
\end{tikzcd}
$$

The natural isomorphism $\P_\AA$ makes the bottom face of the cylinder 2-commute.  As $\P_\AA$ pulls back to a bijection $\P_F$ from $\Irr_{F,\epsilon}^{\DEE}(\alg{\tilde T})$ to $\WP_{F,\epsilon}^{\A}(\alg{\tilde T})$, we find a unique natural isomorphism $\P_F \From \Irr_{F,\epsilon}^{\DEE} \Rightarrow \WP_{F,\epsilon}^{\A}$ making the diagram 2-commute.
\qed

\begin{thm}
\label{GlobalParameterizationNatural}
Define $\Lang_F = \rec_F \circ \can_{F,Y}$ and $\Lang_\epsilon = \AB_F^{-1} \circ \P_F \circ \Del_F$.  This gives a 2-commutative diagram of Picard groupoids and additive functors.
$$\begin{tikzcd}
\phantom{a} &  \phantom{a} & \Cat{Tors}(\Irr_F(\alg{T})) \arrow{dd}{\Lang_F} \\
\Cat{Cov}_n^\sharp(\alg{T}) \arrow{urr}{\Irr_{F,\epsilon}} \arrow{drr}[swap]{\WP_{F,\epsilon}} & \phantom{1000} \Downarrow \Lang_\epsilon  &  \phantom{a} \\
\phantom{a} & & \Cat{Tors}(\WP_F(\alg{T}))
\end{tikzcd}$$
\end{thm}

This theorem gives a parameterization of genuine automorphic characters of $\tilde T_\AA$ by global Weil parameters with values with ${}^\EL \tilde T$, functorial in the choice of sharp cover, and compatible with the previous local parameterization.  As in the case of local fields, one may verify that the natural isomorphisms $\Lang_\epsilon$ are compatible with pullbacks, for isomorphisms of split tori $\alg{S} \To \alg{T}$ over $F$.  

\section{Split tori}

Let $S$ be the spectrum of a local field, a global field, or the ring of integers in a nonarchimedean local field.  Let $\alg{T}$ be a split torus over $S$, with character lattice $X$ and cocharacter lattice $Y$ (constant sheaves on $S_{\et}$).  Let $\alg{\tilde T}$ be a degree $n$ cover of $\alg{T}$ over $S$.  Let $Q \From Y \To \ZZ$ be the first Brylinski-Deligne invariant.  Following Assumption \ref{OddnEvenQ}, we assume that $Q$ is even-valued if $n$ is odd.  But we do \textbf{not} assume $\alg{\tilde T}$ is a sharp cover here, i.e., we do not assume $Y = Y_{Q,n}$.

We can parameterize the following sets of irreducible genuine representations.
\begin{description}
\item[Local case]  The cover $\alg{\tilde T}$ over a local field $F$ yields a central extension $\mu_n \Into \tilde T \Onto T$.  Recall that $\Irr_\epsilon(\alg{\tilde T} / F)$ is the set of equivalence classes of irreducible $\epsilon$-genuine admissible representations of $\tilde T$.
\item[Local integral case]  The cover $\alg{\tilde T}$ over $\OO$ (the ring of integers in a nonarchimedean local field) yields a central extension $\mu_n \Into \tilde T \Onto T$ and a splitting $T^\circ = \alg{T}(\OO) \Into \tilde T$.  Recall that $\Irr_\epsilon(\alg{\tilde T} / \OO)$ is the set of equivalence classes of irreducible $T^\circ$-spherical $\epsilon$-genuine representations of $\tilde T$.
\item[Global case]  The cover $\alg{\tilde T}$ over a global field $F$ yields a central extension $\mu_n \Into \tilde T_\AA \Onto T_\AA$ and a splitting $T_F = \alg{T}(F) \Into \tilde T_\AA$.  Recall that $\Irr_\epsilon(\alg{\tilde T} / F)$ is the set of equivalence classes of automorphic $\epsilon$-genuine representations of $\tilde T_\AA$.  
\end{description}

Let $\alg{T}_{Q,n}$ denote the split torus over $S$ with cocharacter lattice $Y_{Q,n}$, and define 
$$\iota_{Q,n} \From \alg{T}_{Q,n} \To \alg{T}$$
to be the isogeny (of tori over $S$) corresponding to the inclusion $Y_{Q,n} \Into Y$.  Pulling back the cover $\alg{\tilde T}$ via $\iota_{Q,n}$ yields a degree $n$ cover $\alg{\tilde T}_{Q,n}$ of $\alg{T}_{Q,n}$.  Note that $\alg{\tilde T}_{Q,n}$ is a sharp cover of $\alg{T}_{Q,n}$.  Viewing $\iota_{Q,n}$ as a well-aligned homomorphism of covers $\alg{\tilde T}_{Q,n} \To \alg{\tilde T}$, we find a canonical identification of L-groups, ${}^\EL \tilde T \ident {}^\EL \tilde T_{Q,n}$.  

The results of this section and of \cite{MWToriNew} imply the following.
\begin{thm}
In all three cases above, pulling back the central character via the isogeny $\iota_{Q,n}$ gives a one-to-one function (bijective in the local integral case)
$$\Irr_\epsilon(\alg{\tilde T} / S) \To \Irr_\epsilon(\alg{\tilde T}_{Q,n} / S)$$
Composing with the parameterization $\Lang_\epsilon$, this gives a one-to-one parameterization (bijective in the local integral case)
$$\Irr_\epsilon(\alg{\tilde T} / S) \To \WP_\epsilon(\alg{\tilde T} / S).$$
\end{thm}

In the local integral case, the parameterization is bijective.  In the case of local fields or global fields, we would like to characterize the images of $\Irr_\epsilon(\alg{\tilde T} / S) \To \WP_\epsilon(\alg{\tilde T} / S)$ -- to find the ``relevant'' parameters for covers of tori.  But we leave such a characterization for a future paper.

\part{Other parameterizations}

\section{Spherical/Unramified parameterization}

Let $\alg{G}$ be a quasisplit reductive group over $\OO$, the ring of integers in a nonarchimedean local field $F$.  Let $\alg{\tilde G} = (\alg{G}', n)$ be a degree $n$ cover of $\alg{G}$ defined over $\OO$.  Write $G = \alg{G}(F)$ and $G^\circ = \alg{G}(\OO)$.  Then we have a central extension of locally compact groups,
$$\mu_n \Into \tilde G \Onto G,$$
and a splitting $G^\circ \hookrightarrow \tilde G$.  

We have constructed an L-group,
$$\tilde G^\vee \Into {}^\EL \tilde G \Onto \Gal_\OO,$$
where $\tilde G^\vee$ is a complex reductive group and $\Gal_\OO = \langle \Fr \rangle_{\prof}$.  This L-group is well-defined up to L-equivalence, and the L-equivalence is uniquely determined up to unique natural isomorphism.
\subsection{Parameterization}

Let $\Irr_\epsilon^\circ(\alg{\tilde G}) = \Irr_\epsilon(\alg{\tilde G} / \OO)$ denote the set of equivalence classes of $\epsilon$-genuine irreducible $G^\circ$-spherical representations of $\tilde G$.  Write $\WP_\epsilon^\circ(\alg{\tilde G}) = \WP_\epsilon(\alg{\tilde G} / \OO)$ for the set of equivalence classes (i.e., $\tilde G^\vee$-orbits) of unramified Weil parameters $\phi \From \Weil_\OO = \langle \Fr \rangle \To {}^\EL \tilde G$.  In five steps below, we define a bijection
$$\Lang_\epsilon(\alg{\tilde G}) \From \Irr_\epsilon^\circ(\alg{\tilde G}) \xrightarrow{\sim} \WP_\epsilon^\circ(\alg{\tilde G}).$$
Thus the $\epsilon$-genuine irreducible spherical representations of $\tilde G$ are parameterized (bijectively) by Weil parameters.

For this parameterization, let $\alg{A}$ be a maximal $\OO$-split torus in $\alg{G}$, and let $\alg{T}$ be the centralizer of $\alg{A}$ in $\alg{G}$.  Then $\alg{T}$ is a maximally split maximal torus in $\alg{G}$.  Let $\alg{B}$ be a Borel subgroup of $\alg{G}$ containing $\alg{T}$.  Let $\sheaf{W}$ be the Weyl group of the pair $(\alg{G}, \alg{T})$, viewed as a sheaf of finite groups on $\OO_{\et}$.  Write $W = \sheaf{W}(\bar F)$ and $W^\circ = W^{\Fr}$.  

Let $\alg{A}_{Q,n}$ be the $\OO$-split torus with cocharacter lattice $Y_{Q,n}^{\Fr}$.  Then the inclusion $Y_{Q,n}^{\Fr} \Into Y$ defines a $W^\circ$-equivariant homomorphism of $\OO$-tori $\alg{A}_{Q,n} \To \alg{T}$.  Let $\alg{\tilde A}_{Q,n}$ be the sharp cover of $\alg{A}_{Q,n}$ obtained by pulling back the cover $\alg{\tilde T}$.  

\subsubsection{Satake step}
From Corollary \ref{SatBij}, the Satake isomorphism $\Sat \From \hecke_\epsilon(\tilde G, G^\circ) \xrightarrow{\sim} \hecke_\epsilon(\tilde T, T^\circ)^{W^\circ}$ gives a bijection
$$\Sat^\ast \From \Irr_\epsilon^\circ(\alg{\tilde G}) \xrightarrow{\sim} \frac{\Irr_\epsilon^\circ(\alg{\tilde T})}{W^\circ}.$$
\subsubsection{Support step}
From Proposition \ref{HeckeSupport}, restriction gives a $W^\circ$-equivariant isomorphism of Hecke algebras, $\hecke_\epsilon(\tilde T, T^\circ) \xrightarrow{\sim} \hecke_\epsilon(\tilde A_{Q,n}, A_{Q,n}^\circ)$.  See also \cite[\S 3.3]{MWToriNew} for a description of this Hecke algebra.  This gives a bijection,
$$\frac{\Irr_\epsilon^\circ(\alg{\tilde T})}{W^\circ} \xrightarrow{\sim} \frac{\Irr_\epsilon^\circ(\alg{\tilde A}_{Q,n})}{W^\circ}.$$
\subsubsection{Parameterization for sharp covers of split tori}
Our parameterization of Theorem \ref{UnramToriParam} gives a bijection,
$$\Lang_\epsilon^\circ(\alg{\tilde A}_{Q,n}) \From \Irr_\epsilon^\circ(\alg{\tilde A}_{Q,n}) \To \WP_\epsilon^\circ(\alg{\tilde A}_{Q,n}).$$
Each $w \in W^\circ$ can be represented by an element $\dot w \in \alg{G}(\OO)$.  The action of $W^\circ$ on the sets above can be untangled through the following commutative diagram with exact rows.
\begin{equation}
\label{ThreeRowW}
\begin{tikzcd}[column sep = 7em]
\alg{K}_2 \inarrow{r} \arrow{d}{=} & \alg{A}_{Q,n}' \onarrow{r} \arrow{d}{\Int(\dot w)} & \alg{A}_{Q,n} \arrow{d}{=} \\
\alg{K}_2 \arrow{d}{=} \inarrow{r} & \Int(w)^\ast \alg{A}_{Q,n}' \onarrow{r}{\Int(w)^{-1} \circ p} \arrow{d}{=} & \alg{A}_{Q,n} \arrow{d}{\Int(w)} \\
\alg{K}_2 \inarrow{r} & \alg{A}_{Q,n}' \onarrow{r}{p} & \alg{A}_{Q,n}
\end{tikzcd}
\end{equation}

The top two rows of \eqref{ThreeRowW} give a morphism in the category of covers $\Cat{Cov}_n^\sharp(\alg{A}_{Q,n})$,
$$\Int(\dot w) \From \alg{\tilde A}_{Q,n} \To \Int(w)^\ast \alg{\tilde A}_{Q,n}.$$
The functoriality of $\Lang_\epsilon$ for such morphisms gives a commutative diagram.
$$\begin{tikzcd}[column sep = 7em]
\Irr_\epsilon^\circ(\alg{\tilde A}_{Q,n}) \arrow{r}{\Lang_\epsilon(\alg{\tilde A}_{Q,n})} \arrow{d}{\Lang_\epsilon(\Int(\dot w))} & \WP_\epsilon^\circ(\alg{\tilde A}_{Q,n}) \arrow{d}{\Lang_\epsilon(\Int(\dot w))}  \\
\Irr_\epsilon^\circ(\Int(w)^\ast \alg{\tilde A}_{Q,n}) \arrow{r}{\Lang_\epsilon(\Int(w)^\ast \alg{\tilde A}_{Q,n})} & \WP_\epsilon^\circ(\Int(w)^\ast \alg{\tilde A}_{Q,n})
\end{tikzcd}$$

Compatibility of $\Lang_\epsilon$ with pullbacks, combined with the bottom two rows of \eqref{ThreeRowW}, gives a commutative diagram,
$$\begin{tikzcd}[column sep = 7em]
\Irr_\epsilon^\circ(\Int(w)^\ast \alg{\tilde A}_{Q,n}) \arrow{r}{\Lang_\epsilon(\Int(w)^\ast \alg{\tilde A}_{Q,n})} \arrow{d} & \WP_\epsilon^\circ(\Int(w)^\ast \alg{\tilde A}_{Q,n})  \arrow{d} \\
\Irr_\epsilon^\circ(\alg{\tilde A}_{Q,n}) \arrow{r}{\Lang_\epsilon(\alg{\tilde A}_{Q,n})}  & \WP_\epsilon^\circ(\alg{\tilde A}_{Q,n})
\end{tikzcd}$$

Combining these, we find that the parameterization $\Lang_\epsilon \From \Irr_\epsilon^\circ(\alg{\tilde A}_{Q,n}) \To \WP_\epsilon^\circ(\alg{\tilde A}_{Q,n})$ is equivariant for the action of $W^\circ$.  This gives a bijective parameterization
$$\Lang_\epsilon \From \frac{\Irr_\epsilon^\circ(\alg{\tilde A}_{Q,n})}{W^\circ} \xrightarrow{\sim} \frac{\WP_\epsilon^\circ(\alg{\tilde A}_{Q,n})}{W^\circ}.$$
\subsubsection{Split and unramified parameters}
The inclusion $\alg{\tilde A}_{Q,n} \To \alg{\tilde T}$ is a well-aligned morphism of covers, giving an L-morphism.
$$\begin{tikzcd}
\Hom(Y_{Q,n}, \CC^\times) = \tilde T^\vee \inarrow{r} \onarrow{d}{\nu} &  {}^\EL \tilde T \onarrow{r} \onarrow{d}{\nu} & \Gal_\OO \arrow{d}{=} \\
\Hom(Y_{Q,n}^{\Fr}, \CC^\times) = \tilde A_{Q,n}^\vee \inarrow{r} & {}^\EL \tilde A_{Q,n} \onarrow{r} & \Gal_\OO
\end{tikzcd}$$
Since the inclusion $\alg{\tilde A}_{Q,n} \Into \alg{\tilde T}$ is $W^\circ = W^{\Fr}$-equivariant, we find that the L-morphism $\nu \From {}^\EL \tilde T \To {}^\EL \tilde A_{Q,n}$ satisfies
$$\nu \circ {}^\EL \Int(\dot w) \text{ is naturally isomorphic to } {}^\EL \Int(\dot w) \circ \nu,$$
for all $w \in W^{\Fr}$ and all representatives $\dot w \in G^\circ$ for $w$.

The method below follows \cite[\S 6]{BorelCorvallis} very closely; minimal changes are required.  Giving an equivalence class of Weil parameters $\WP_\epsilon^\circ(\alg{\tilde A}_{Q,n})$ is the same as giving an $\tilde A_{Q,n}^\vee$-conjugacy class of elements ${}^\EL \tilde A_{Q,n}$ lying over $\Fr \in \Gal_\OO$.  Similarly, giving an equivalence class of Weil parameters in $\WP_\epsilon^\circ(\alg{\tilde T})$ is the same as giving a $\tilde T^\vee$-conjugacy class of elements in ${}^\EL \tilde T$ lying over $\Fr \in \Gal_\OO$.

As $\nu \From {}^\EL \tilde T \To {}^\EL \tilde A_{Q,n}$ is surjective and $W^{\Fr}$-equivariant, we find that $\nu$ gives a surjective map,
$$\bar \nu \From \frac{\WP_\epsilon^\circ(\alg{\tilde T})}{W^{\Fr}} \Onto \frac{\WP_\epsilon^\circ(\alg{\tilde A}_{Q,n})}{W^{\Fr}}.$$
If $[\phi], [\phi'] \in \WP_\epsilon^\circ(\alg{\tilde T})$ and $\nu([\phi])$ and $\nu([\phi'])$ are in the same $W^{\Fr}$-oribt, then there exist $t,t' \in {}^\EL \tilde T$ and $w \in W^{\Fr}$ represented by $\dot w \in G^\circ$, such that $\phi(\Fr) = t$ and $\phi'(\Fr) = t'$, and $\nu(t) = {}^\EL \Int(\dot w) (\nu(t'))$.  Thus $\nu(t) = \nu(\Int(\dot w) t')$, and so
$$t^{-1} \cdot \Int(\dot w) t' \in \Ker(\tilde T^\vee \Onto \tilde A_{Q,n}^\vee).$$
The following sequence is exact, by the same arguments as in \cite[\S 6.3]{BorelCorvallis}.  
$$\Hom(Y_{Q,n}, \CC^\times) \xrightarrow{\Fr - 1} \Hom(Y_{Q,n}, \CC^\times) \xrightarrow{\Res} \Hom(Y_{Q,n}^{\Fr}, \CC^\times).$$
It follows that $t^{-1} \cdot \Int(\dot w) t' = \tau / \Fr(\tau)$ for some $\tau \in \tilde T^\vee$.  Hence
$$\Int(\dot w) t' = \tau t \tau^{-1}.$$
Hence if $\nu([\phi])$ and $\nu([\phi'])$ are in the same $W^{\Fr}$-orbit, then $[\phi]$ and $[\phi']$ are in the same $W^{\Fr}$-orbit.  Therefore, $\bar \nu$ gives a bijective map,
$$\bar \nu \From \frac{\WP_\epsilon^\circ(\alg{\tilde T})}{W^{\Fr}} \xrightarrow{\sim} \frac{\WP_\epsilon^\circ(\alg{\tilde A}_{Q,n})}{W^{\Fr}}.$$
\subsubsection{Semisimple twisted conjugacy classes}
The Weyl group $W$ is canonically isomorphic to the Weyl group of $\tilde G^\vee$ with respect to $\tilde T^\vee$, as finite groups with $\Gal_S$-action.  In this way, any element $w \in W^{\Fr}$ corresponds to an element $w^\vee \in (W^\vee)^{\Fr}$, which may be represented by a $\Gal_S$-invariant element $n^\vee \in \tilde N^\vee \subset \tilde G^\vee$, where $\tilde N^\vee$ is the normalizer of $\tilde T^\vee$ in $\tilde G^\vee$ (see \cite[Lemma 6.2]{BorelCorvallis}).  From Theorem \ref{BigWTheorem}, the L-morphism ${}^\EL \Int(\dot w)$ of ${}^\EL \tilde T$ is naturally isomorphic to the L-morphism $\Int(n^\vee)$.  This gives a bijection
$$\frac{\WP_\epsilon^\circ(\alg{\tilde T})}{W^{\Fr}} \leftrightarrow \frac{\WP_\epsilon^\circ(\alg{\tilde T})}{(\tilde N^\vee)^{\Fr}}.$$
Finally, from \cite[Lemma 6.5]{BorelCorvallis}, inclusion provides a bijection,
$$\frac{\WP_\epsilon^\circ(\alg{\tilde T})}{(\tilde N^\vee)^{\Fr}} \xrightarrow{\sim} \WP_\epsilon^\circ(\alg{\tilde G}).$$

\begin{thm}
\label{UnramifiedParameterization}
Assembling the bijections in the five steps above, we have constructed a bijection
$$\Lang_\epsilon(\alg{\tilde G}) \From \Irr_\epsilon^\circ(\alg{\tilde G}) \xrightarrow{\sim} \WP_\epsilon^\circ(\alg{\tilde G}).$$
\end{thm}

\subsection{Automorphic L-functions}

The spherical/unramified parameterization provides a definition of automorphic L-functions, almost everywhere.  Let $F$ be a global field.  Let $\VV$ be the set of places of $F$, and let $\SS$ be a finite subset of $\VV$ containing all archimedean places.  Suppose that $\alg{\tilde G} = (\alg{G}', n)$ is a degree $n$ cover of a quasisplit reductive group $\alg{G}$ over the ring of $\SS$-integers $\OO_\SS$ (see \cite[\S 3.2]{MWIntegral} for a classification).  This defines a central extension $\mu_n \Into \tilde G_\AA \Onto G_\AA$, endowed with splittings over $G_F$ and over $G_v^\circ = \alg{G}(\OO_v)$ for all $v \not \in \SS$.

Let $\pi$ be an $\epsilon$-genuine automorphic representation of $\tilde G_\AA$.  Factorization yields $\epsilon$-genuine spherical irreducible representations $[\pi_v]$, for almost all nonarchmedean places $v \in \VV - \SS$.  Each equivalence class $[\pi_v]$ yields an equivalence class of unramified parameters $[\sigma_v] \in \WP_\epsilon^\circ(\alg{\tilde G})$ by our parameterization $\Lang_\epsilon(\alg{\tilde G}_v)$ of Theorem \ref{UnramifiedParameterization}.  Each unramified parameter yields a well-defined semisimple $\tilde G^\vee$-conjugacy class $g_v = \sigma_v(\Fr) \in {}^\EL \tilde G_v$.  Local-global compatibility in the construction of the L-group identify $g_v$ with semisimple conjugacy classes in the globally-defined L-group ${}^\EL G$.

Let $\rho \From {}^\EL \tilde G \To GL(V)$ be a representation of the L-group on a finite-dimensional complex vector space $V$.  Then, for almost all nonarchimedean places $v$ (with residue fields of cardinality $q_v$), we have a local L-function
$$L_v(\pi, \rho, s) = \det(\Id_V - q_v^{-s} \rho(g_v) \vert V)^{-1}.$$
We define the \defined{automorphic L-function} away from $\SS$,
$$L_\SS(\pi, \rho, s) = \prod_{v \not \in S} L_v(\pi, \rho, s).$$

From the remarks of \cite[\S 10.5]{B-D}, another choice $\SS'$ of places and an $\OO_{\SS'}$-model of $\alg{\tilde G}$ will yield an automorphic L-function $L_{\SS'}(\pi, \rho, s)$ which agrees with $L_\SS(\pi, \rho, s)$ at almost all places.  In this way the automorphic L-function is well-defined from $\pi$ and $\rho$, up to factors at finitely many places.  In some cases, these automorphic L-functions also appear in the recent work of Gan and Gao \cite[\S 12.4]{GanGao}, \cite[\S 5]{GaoThesis}.

\begin{remark}
Langlands' argument (see \cite[\S 13]{BorelCorvallis} or \cite[\S 2.5]{ShahidiEisenstein} for treatments), together with the spectral decomposition of automorphic forms by Moeglin and Waldspurger \cite[\S III.2.6]{MoeglinWaldspurger}, should imply the absolute convergence of $L_\SS(\pi, \rho, s)$ in some right half-plane.  The only missing ingredient at the moment is the Macdonald formula for the spherical function.  For split groups, progress has been made in the work of McNamara \cite{McNamara} and the recent thesis of Fan Gao \cite[\S 4.3]{GaoThesis}.

But, as I learned from Wee Teck Gan, this is overkill -- results of Tadic \cite[Theorem 2.5]{Tadic} imply that the Satake parameters of a unitary spherical irreducible representation lie in a compact subset of the space of all spherical irreducible representations (and Tadic's argument applies just as well to covering groups); this in turn implies the convergence of the global L-function in {\em some} right half-plane, without need for the Macdonald formula.  
\end{remark}

\section{Sharp covers of anisotropic real tori}

Let $\alg{T}$ be an \textbf{anisotropic} torus over $\RR$, and let $\alg{\tilde T}$ be a degree $2$ \textbf{sharp} cover of $\alg{T}$ over $\RR$.  Write $Y = \sheaf{Y}[\CC]$ for the cocharacter lattice, and $\sigma$ for complex conjugation.  Thus $\sigma(y) = -y$ for all $y \in Y$.  The cover $\alg{\tilde T}$ yields an extension,
$$\mu_2 \Into \tilde T \Onto T.$$
There is a natural identification of Lie groups $T = Y \otimes U(1)$, from which we identify the topological fundamental group $\pi_1(T) \ident Y$.  Connectedness of $T$ implies that $\tilde T$ is abelian, and that this extension is rigid; there are no nontrivial automorphisms of the extension $\tilde T \in \Cat{Ext}(T, \mu_2)$.  The extension $\tilde T$ is determined, up to unique isomorphism, by an associated homomorphism $\kappa \From \pi_1(T) = Y \To \half \ZZ / \ZZ$. 

The homomorphism $\kappa$ may be constructed via Pontrjagin duality as follows.  The Pontrjagin dual of $T$, $\Irr(\alg{T}) = \Hom(T, U(1))$, is naturally identified with the character lattice $X = \sheaf{X}[\CC]$.  Write $X_\kappa$ for the Pontrjagin dual of $\tilde T$, fitting into a short exact sequence
$$X \Into X_\kappa \Onto \ZZ_{/2}.$$
The genuine characters of $\tilde T$ are then identified with the elements of $X_\kappa^- = X_\kappa - X$,
$$\Irr_\epsilon(\alg{\tilde T}) \ident X_\kappa^- = X_\kappa - X.$$

For each such genuine character $x \in X_\kappa^-$, $2x \in X$; thus there exists a unique $\xi \in \half X$ such that $2x = 2 \xi$.  The map $x \mapsto \xi$ provides an embedding, equivariant for translation by $X$ (i.e., for twisting by $\Irr(\alg{T})$), $X_\kappa^- \Into \half X$.  As $X_\kappa^-$ is an $X$-torsor, there exists a unique $\kappa \in \half X / X$ such that
$$X_\kappa^- \ident \{ \xi \in \half X : \xi = \kappa \text{ mod } X \}.$$
The element $\kappa \in \half X / X$ can also be viewed as a homomorphism
$$\kappa \From \pi_1(T) = Y \To \half \ZZ / \ZZ.$$
This $\kappa$ is the homomorphism which determines the double cover $\tilde T$ uniquely up to unique isomorphism.

The Brylinski-Deligne invariants of $\alg{\tilde T}$ consist of a quadratic form $Q \From Y \To \ZZ$ and a $\Gal_\RR$-equivariant extension
$$\CC^\times \Into D \Onto Y.$$
Recall that $\alg{\tilde T}$ is a sharp cover; thus $Y = Y_{Q,2}$ and $D$ is commutative.  For each $y \in Y$, write $\iota_y \From \ZZ \Into Y$ for the homomorphism $n \mapsto n \cdot y$.  We get a $\Gal_\RR$-equivariant extension
$$\CC^\times \Into \iota_y^\ast D \Onto \ZZ,$$
where $\sigma \in \Gal_\RR$ acts by $n \mapsto -n$ on $\ZZ$, and by complex conjugation on $\CC^\times$.  Following \cite[\S 12.6]{B-D}, such $\Gal_\RR$-equivariant extensions of $\ZZ$ by $\CC^\times$ are classified up to isomorphism by elements of $\half \ZZ / \ZZ$.  (Brylinski and Deligne use $\ZZ_{/2}$ instead, but our normalization has advantages to follow).  Indeed, take $d \in \iota_y^\ast D$ lying over $1 \in \ZZ$.  Then $d \cdot \sigma(d)$ lies over $1 + (-1) = 0 \in \ZZ$, and so $d \sigma(d) \in \CC^\times$.  Since $d \sigma(d) = \sigma(d \sigma(d))$, we have $d \sigma(d) \in \RR^\times$.  Define $\eta(y) \in \frac{1}{2} \ZZ / \ZZ$ so that
$$\sgn(d \cdot \sigma(d)) = e^{2 \pi i \eta(y)}.$$
The function $\eta$ is independent of choices along the way, and defines a homomorphism 
$$\eta \From Y \To \half \ZZ / \ZZ.$$

The next proposition is a direct consequence of \cite[\S 12.6, Proposition 12.7]{B-D}.  It gives a practical recipe for determining the topological cover $\tilde T$ from the Brylinski-Deligne invariants of $\alg{\tilde T}$.
\begin{proposition}
\label{KappaEtaQ}
The homomorphism $\kappa \From Y \To \half \ZZ / \ZZ$ is given by
$$\kappa(y) = \eta(y) + \half Q(y) \quad (\text{mod } \ZZ).$$
\end{proposition}

To see that the Brylinski-Deligne framework is sufficiently strong for work on Lie groups, the following may be of interest.
\begin{corollary}
Every topological double cover $\mu_2 \Into \tilde T \Onto T$ arises from a sharp cover $\alg{\tilde T} \in \Cat{Cov}_2^\sharp(\alg{T})$ over $\RR$.  
\end{corollary}
\proof
A topological double cover $\tilde T$ arises from a homomorphism $\kappa \From Y \To \half \ZZ / \ZZ$.  Choose a basis $y_1, \ldots, y_r \in Y$, and define a quadratic form $Q \From Y \To \ZZ$ by
$$Q(a_1 y_1 + \cdots + a_r y_r) = 2 \kappa(y_1) a_1^2 + \cdots + 2 \kappa(y_r) a_r^2.$$
Then $\half Q(y) = \kappa(y)$, mod $\ZZ$.  Moreover, $Y = Y_{Q,2}$ for this quadratic form.  By the classification of \cite{B-D}, there exists a central extension $\alg{K}_2 \Into \alg{T}' \Onto \alg{T}$ over $\RR$ with first invariant $Q$ and trivial second invariant.  The previous proposition gives a unique isomorphism from $\tilde T$ to the topological double cover associated to $\alg{T}'$.  
\qed

Now, we work on the side of the L-group.  The dual group $T^\vee = \tilde T^\vee = X \otimes \CC^\times$ is endowed with the $\Gal_\RR$-action $\sigma(t^\vee) = (t^\vee)^{-1}$ for all $t^\vee \in T^\vee$.  The L-group is a short exact sequence of complex groups,
$$T^\vee \Into {}^\EL \tilde T \Onto \Gal_\RR.$$

As $\kappa \in \half X / X$, we may view $e^{2 \pi i \kappa}$ as an element of $T_{[2]}^\vee = X \otimes \mu_2$.  Explicitly,
$$y(e^{2 \pi i \kappa}) = e^{2 \pi i \kappa(y)}, \text{ for all } y \in Y.$$
\begin{lemma}
\label{gammasquare}
The lifts $\gamma \in {}^\EL \tilde T$ of $\sigma \in \Gal_\RR$ form a single orbit under $T^\vee$-conjugation, and for any such lift $\gamma$,
$$\gamma^2 = e^{2 \pi i \kappa}.$$
\end{lemma}
\proof
If $\gamma_1, \gamma_2$ are two lifts of $\sigma$, then $\gamma_1 = t^\vee \gamma_2$ for some $t^\vee \in T^\vee$.  Since squaring is surjective on $T^\vee$, there exists $\tau^\vee \in T^\vee$ such that $\tau^\vee / \sigma(\tau^\vee) = (\tau^\vee)^2 = t^\vee$.  Hence
$$\tau^\vee \cdot \gamma_1 \cdot (\tau^\vee)^{-1} = \frac{\tau^\vee}{\sigma(\tau^\vee)} \cdot \gamma_1 = t^\vee \cdot \gamma_1 = \gamma_2.$$
It follows that $\gamma_2^2 = \gamma_1^2$ and the lifts of $\sigma$ form a single $T^\vee$ orbit.  Therefore, to prove the theorem, it suffices to prove that $\gamma^2 = e^{2 \pi i \kappa}$ for a single lift $\gamma$ of $\sigma$.

Recall that ${}^\EL \tilde T$ is the Baer sum of two extensions of $\Gal_\RR$ by $T^\vee$,
$${}^\EL \tilde T = (\tau_Q)_\ast \mGal_\RR \Baer \pi_1^{\et}(\gerb{E}_\epsilon(\alg{\tilde T}), \bar z).$$
Let $\gamma_Q$ be a lift of $\sigma$ in the extension $(\tau_Q)_\ast \mGal_\RR$.  In the metaGalois group $\mu_2 \Into \mGal_\RR \Onto \Gal_\RR$, the square of any lift of $\sigma$ equals $-1$ (since $\Hilb_2(-1, -1) = -1$).  Hence we find $\gamma_Q^2 = \tau_Q(-1)$.  In other words,
$$y( \gamma_Q^2) = (-1)^{Q(y)} = e^{\pi i Q(y)} \text{ for all } y \in Y.$$

Next we construct an object $\bar z$ of $\gerb{E}_\epsilon(\alg{\tilde T})[\CC]$ and a lift of $\sigma$ in $\pi_1^{\et}(\gerb{E}_\epsilon(\alg{\tilde T}), \bar z)$.  Let $\hat T = \Hom(Y, \CC^\times)$, with $\Gal_\RR$-action arising from $y \mapsto -y$ on $Y$, and complex conjugation on $\CC^\times$.  Then $\Spl(D)$ is the $\hat T$-torsor of splittings of $D$.  Choose any $s \in \Spl(D)$.  The $\Gal_\RR$-action on $D$ gives an action on splittings; we have
$${}^\sigma s(y) = \sigma(s(\sigma(y))) = \sigma(s(-y)) = \sigma(s(y))^{-1} \text{ for all } y \in Y.$$
Define a map $h \From \hat T \To Spl(D)$ by
$$h(\hat t) = \hat t^2 \ast s.$$
Then $\bar z = (\hat T, h) \in \sqrt{\Spl(D)}$ is a square root of the $\hat T$-torsor $\Spl(D)$, i.e. $\bar z$ is an object of $\gerb{E}_\epsilon(\alg{\tilde T})[\CC]$.  Its complex conjugate ${}^\sigma \bar z$ is given by ${}^\sigma \bar z = (\hat T, {}^\sigma h)$, where
$${}^\sigma h(\hat t) = \hat t^2 \ast {}^\sigma s.$$

Choose $f \in \Hom(\bar z, {}^\sigma \bar z) \subset \pi_1^{\et}(\gerb{E}_\epsilon(\alg{\tilde T}), \bar z)$.  Then, as a morphism of $\hat T$ torsors, there exists $\hat \tau \in \hat T$ such that $f \From \hat T \To \hat T$ satisfies $f(\hat a) = \hat \tau \ast \hat a$ for all $\hat a$.   The condition ${}^\sigma h \circ f = h$ implies that 
$$\hat a^2(y) s(y) = \hat a^2(y) \hat \tau^2(y) \cdot \sigma(s(y))^{-1}$$
for all $y \in Y$.  Hence, if $y \in Y$ and $d = s(y)$, then
$$\hat \tau(y) = \sqrt{ d \sigma(d)} \in \CC^\times \quad \text{(for some square root).}$$

The square of $f \in \Hom(\bar z, {}^\sigma \bar z) \subset \pi_1^{\et}(\gerb{E}_\epsilon(\alg{\tilde T}), \bar z)$ is the element of $\Hom(\bar z, \bar z)$ given by composing morphisms
$$\bar z \xrightarrow{f} {}^\sigma \bar z \xrightarrow{ {}^\sigma f} {}^\sigma ({}^\sigma \bar z) = \bar z.$$
As a function from $\hat T$ to $\hat T$, we find that $f^2(\hat a) = {}^\sigma \hat \tau \cdot \hat \tau \cdot \hat a$, for all $\hat a \in \hat T$.  In other words, $f^2 = {}^\sigma \hat \tau \cdot \hat \tau$ as an element of $\hat T_{[2]} \subset T^\vee$.
We compute
$$y(f^2) = [{}^\sigma \hat \tau \cdot \hat \tau](y) = \sigma(\hat \tau(-y)) \cdot \tau(y) = \frac{ \sqrt{d \sigma(d)}}{\sigma \sqrt{d \sigma(d)}} = \sgn(d \sigma(d)) = e^{2 \pi i \eta(y)}.$$

Finally, define $\gamma = \gamma_Q \Baer f \in {}^\EL \tilde T$.  From Proposition \ref{KappaEtaQ}, we find
$$y(\gamma^2) = y(\gamma_Q^2) \cdot y(f^2) = e^{\pi i Q(y) + 2 \pi i \eta(y)} = e^{2 \pi i \kappa(y)}.$$
\qed

The Weil group of $\RR$ is $\Weil_\RR = \CC^\times \sqcup j \CC^\times$; the map from $\Weil_\RR$ to $\Gal_\RR$ sends $\CC^\times$ to $\Id$ and $j$ to $\sigma$.  Recall that $j z j^{-1} = \bar z$ for all $z \in \CC^\times$, and $j^2 = -1$.  A continuous homomorphism from $\CC^\times$ to $T^\vee = \Hom(Y, \CC^\times)$ has the form $z \mapsto z^{x_1 - x_2} (z \bar z)^{x_2}$ for some $x_1, x_2 \in X \otimes \CC$ satisfying $x_1 - x_2 \in X$.  We abuse notation slightly, and write $z^{x_1} \bar z^{x_2}$ for such a homomorphism, keeping in mind that $x_1 - x_2 \in X$ is necessary for this expression to make sense.

Recall that $\kappa \in \half X / X = \Hom(Y, \half \ZZ / \ZZ)$, and
$$X_\kappa^- = \{ \xi \in \half X : \xi = \kappa \text{ mod } X \}.$$
\begin{thm}
For every Weil parameter $\phi \From \Weil_\RR \To {}^\EL \tilde T$, there exists a unique element $\xi \in X_\kappa^-$ such that
$$\phi(z) = (z / \bar z)^\xi \defeq z^{2 \xi} (z \bar z)^{-\xi}.$$
The map $\phi \mapsto \xi$ defines a bijection,
$$\WP_\epsilon(\alg{\tilde T}) = \frac{\Par(\Weil_\RR, {}^\EL \tilde T) }{T^\vee-\text{conjugation} } \xrightarrow{\sim} X_\kappa^- \ident \Irr_\epsilon(\alg{\tilde T}).$$
\end{thm}
\proof
Giving a Weil parameter $\phi \From \Weil_\RR \To {}^\EL \tilde T$ is the same as giving a homomorphism $\phi(z) = z^{x_1} \bar z^{x_2}$ for some $x_1, x_2 \in X \otimes \CC$, $x_1 - x_2 \in X$, and an element $\gamma = \phi(j) \in {}^\EL \tilde T$ lying over $\sigma \in \Gal_\RR$, which satisfy the following conditions:
\begin{itemize}
\item
$\gamma^2 = \phi(-1) = (-1)^{x_1 - x_2}$;
\item
$\gamma \phi(z) \gamma^{-1} = \phi(\bar z)$ for all $z \in \CC^\times$.
\end{itemize}

From the previous lemma, $T^\vee$ acts transitively on the lifts of $\sigma$ in ${}^\EL \tilde T$, and all lifts $\gamma$ satisfy $\gamma^2 = e^{2 \pi i \kappa}$.  We find that giving a $T^\vee$-orbit on $\Par(\Weil_\RR, {}^\EL \tilde T)$ is the same as giving a homomorphism $\phi(z) = z^{x_1} \bar z^{x_2}$ for some $x_1, x_2 \in X \otimes \CC$, $x_1 - x_2 \in X$, which satisfies the following conditions:
\begin{itemize}
\item
$\phi(-1) = (-1)^{x_1 - x_2} = e^{2 \pi i \kappa}$;
\item
$\phi(\bar z) = \bar z^{x_1} z^{x_2} = \sigma(\phi(z)) = z^{-x_1} \bar z^{- x_2}$.
\end{itemize}

The second condition is equivalent to the condition $x_1 = - x_2$; if $\xi = x_1 = -x_2$, note that $2 \xi = x_1 - x_2 \in X$ and so $\xi \in \half X$.  In this case, the first condition is equivalent to $(-1)^{2 \xi} = e^{2 \pi i \xi} = e^{2 \pi i \kappa}$.  We find that giving a $T^\vee$-orbit on $\Par(\Weil_\RR, {}^\EL \tilde T)$ is the same as giving an element $\xi \in \half X$ such that $\xi = \kappa$ mod $X$.  This defines the bijection.
\qed

\section{Discrete series for covers of real semisimple groups}

\subsection{Harish-Chandra classification}

Suppose here that $\alg{G}$ is a \textbf{semisimple} quasisplit group over $\RR$, with Borel subgroup $\alg{B}$ containing a maximally split maximal torus $\alg{T}$.  Let $\alg{\tilde G} = (\alg{G}', 2)$ be a double cover of $\alg{G}$ defined over $\RR$.  Then we find a topological double cover over Lie groups,
$$\mu_2 \Into \tilde G \Onto G.$$
\begin{remark}
If, in addition, $\alg{G}$ is simply-connected and absolutely simple, then there is a unique up to unique isomorphism double cover $\alg{\tilde G}$ whose first Brylinski-Deligne invariant takes the value $1$ on all short coroots.  The resulting double cover $\mu_2 \Into \tilde G \Onto G$ coincides with Deligne's canonical central extension, after \cite[Construction 9.3, 10.3]{B-D}.  In this case, Prasad and Rapinchuk prove that $\tilde G$ is the unique nontrivial 2-fold cover of $G$ \cite[Theorem 8.4 and \S 8.5]{PrasadRapinchuk}.  Thus the framework of Brylinski-Deligne is sufficient to work with the most interesting nonlinear double covers of semisimple Lie groups.
\end{remark}

Let $G^\circ$ be the identity component of the Lie group $G$, and let $\tilde G^\circ$ be the identity component of the Lie group $\tilde G$.  Then $G^\circ$ and $\tilde G^\circ$ are connected real semisimple Lie groups with finite center.  The Harish-Chandra classification and the description of discrete series by Atiyah and Schmid (see \cite{AtiyahSchmid}) apply to $G^\circ$ and $\tilde G^\circ$. 

Let $\Irr_\epsilon^{\dis}(\alg{\tilde G})$ be the set of equivalence classes of irreducible $\epsilon$-genuine \textbf{discrete series} representations of $\tilde G$.  Here $\epsilon \From \alg{\mu}_2(\RR) \To \alg{\mu}_2(\CC)$ is the identity map, so we often refer to ``genuine representations'' without mention of $\epsilon$.  Any irreducible discrete series representation of $\tilde G$ restricts to a finite direct sum of irreducible discrete series representations of $\tilde G^\circ$, all with the same infinitesimal character.  By Frobenius reciprocity, all irreducible discrete series representations of $\tilde G^\circ$ occur in such a restriction.

The set $\Irr_\epsilon^{\dis}(\alg{\tilde G})$ is nonempty if and only if there exists a compact Cartan subgroup of $\tilde G^\circ$; such a compact Cartan subgroup exists if and only if there exists a maximal torus $\alg{S} \subset \alg{G}$, defined and anisotropic over $\RR$.  Suppose that $\alg{S}$ is such a torus.

Let $S = \alg{S}(\RR)$ and $\tilde S$ its preimage in $\tilde G$.  Let $K$ be a maximal compact subgroup of $G$ containing $S$, and $K^\circ$ its identity component.  Write $\tilde K$ for the preimage of $K$ in $\tilde G$, and $\tilde K^\circ$ for the identity component of $\tilde K$.

Consider a genuine irreducible discrete series representation $(\pi,V)$ of $\tilde G$; the $\tilde K$-finite vectors therein form an irreducible admissible $(\Lie{g}, \tilde K)$-module.  As such, it has an infinitesimal character $\chi \From \Lie{Z} \To \CC$, where $\Lie{Z}$ denotes the center of the universal enveloping algebra of $\Lie{g}$.  The Harish-Chandra isomorphism gives an isomorphism of $\CC$-algebras,
$$\Lie{Z} \xrightarrow{\sim} \CC[\Lie{s}]^W,$$
where $\Lie{s}$ is the complexified Lie algebra of $\alg{S}$, and $W$ is the Weyl group of $\alg{G}$ with respect to $\alg{S}$.  Let $X$ be the cocharacter lattice of $\alg{S}$.  Then $\Hom(\Lie{s}, \CC)$ is naturally isomorphic to $X \otimes \CC$.  Thus the infinitesimal character defines a map,
$$\inf \From \Irr_\epsilon^{\dis}(\alg{\tilde G}) \To \frac{X \otimes \CC}{W-\text{conjugation}}.$$

As $S = \alg{S}(\RR) \ident Y \otimes U(1)$, the Pontrjagin dual of $S$ is identified with $X$.  The set of genuine characters $\Irr_\epsilon(\alg{\tilde S})$ of $\tilde S$ is identified with a $X$-coset $\kappa + X \subset \half X$.  Here, the element 
$$\kappa \in \half X / X = \Hom(Y, \half \ZZ / \ZZ) = \Hom(\pi_1(S), \half \ZZ / \ZZ)$$
determines the double cover $\tilde S$ up to unique isomorphism.
  
\begin{thm}
Let $\rho$ be the half-sum of the positive roots of $\alg{G}$ (with respect to some Borel subgroup containing $\alg{S}$).  Let $X_\QQ^{\reg}$ denote the regular (=nonsingular) locus in $X_\QQ = X \otimes \QQ$.  The infinitesimal character provides a finite-to-one surjective map,
$$\inf \From \Irr_\epsilon^{\dis}(\alg{\tilde G}) \Onto \frac{(\kappa + \rho + X) \cap X_\QQ^{\reg}}{W-\text{conjugation}}.$$
\end{thm}
\proof
This follows directly from \cite[Corollary 6.13]{AtiyahSchmid}, together with our remarks on the restriction from $\tilde G$ to $\tilde G^\circ$.
\qed

In the uncovered case, the fibres of the infinitesimal character map are precisely the L-packets (under our semisimplicity assumption).  In what follows, we parameterize the genuine discrete series representations by discrete series Weil parameters, valued in the L-group of the cover $\alg{\tilde G}$.

\subsection{Discrete series parameters}

Our treatment here is largely based on the work of Langlands \cite[\S 3]{LangIrred}.  We have reproduced the structure and arguments of Langlands, making adaptations and additions where necessary.  As before, $\alg{G}$ is a \textbf{semisimple} quasisplit group over $\RR$ and $\alg{\tilde G}$ is a double cover of $\alg{G}$ over $\RR$.  The L-group of $\alg{\tilde G}$ fits into a short exact sequence,
$$\tilde G^\vee \Into {}^\EL \tilde G \Onto \Gal_\RR = \{ 1, \sigma \}.$$
Recall that the dual group $\tilde G^\vee$ is a pinned complex reductive group, with $\tilde G^\vee \supset \tilde B^\vee \supset \tilde T^\vee$, associated to the root datum
$$\tilde \Psi^\vee = (Y_{Q,n}, \tilde \Phi^\vee, \tilde \Delta^\vee, X_{Q,n}, \tilde \Phi, \tilde \Delta).$$

A \textbf{discrete series Weil parameter} is a Weil parameter $\phi \From \Weil_\RR \To {}^\EL \tilde G$ whose image is not contained in any proper parabolic subgroup ${}^\EL \tilde P$.  (As in \cite{LangIrred}, a parabolic subgroup of ${}^\EL \tilde G$ is a subgroup ${}^\EL \tilde P$ whose intersection with $\tilde G^\vee$ is a parabolic subgroup and whose projection to $\Gal_\RR$ is surjective).  Let $\WP_\epsilon^{\dis}(\alg{\tilde G})$ be the set of equivalence classes of discrete series Weil parameters.

Suppose that $\phi_0$ is a discrete series Weil parameter.  Recall that $\Weil_\RR = \CC^\times \sqcup \CC^\times j$, where $j^2 = -1$ and $j z j^{-1} = \sigma(z)$ for all $z \in \CC^\times$.
\begin{lemma}
$\phi_0$ is equivalent to a parameter $\phi$ such that $\phi(\CC^\times) \subset \tilde T^\vee$ and $\phi(j)$ normalizes $\tilde T^\vee$.
\end{lemma}
\proof
The reasoning of \cite{LangIrred} based on \cite[\S II, Theorem 5.16]{SpringerSteinberg} demonstrates that $\phi(\CC^\times) \subset \tilde G^\vee$ is contained in a maximal torus which is normalized by $\phi(\Weil_\RR)$.  Replacing $\phi$ by a $\tilde G^\vee$-conjugate parameter if necessary, we may assume that $\phi(\CC^\times) \subset \tilde T^\vee$ and $\phi(\Weil_\RR)$ normalizes $\tilde T^\vee$.
\qed

Let $\phi$ be a discrete series Weil parameter, with $\phi(\CC^\times) \subset \tilde T^\vee$ and $\gamma = \phi(j)$ normalizing $\tilde T^\vee$.  Then conjugation by $\gamma$ is a $\CC$-algebraic automorphism of $\tilde T^\vee$, and so it arises from
$$\gamma \From Y_{Q,n} \To Y_{Q,n}.$$
\begin{lemma}
For all $y \in Y_{Q,n}$, $\gamma(y) = -y$.  Equivalently, for all $t^\vee \in \tilde T^\vee$,
$$\gamma t^\vee \gamma^{-1} = (t^\vee)^{-1}.$$
\end{lemma}
\proof
Since $j^2 = -1$, $\gamma$ acts as an automorphism of order $2$ of $Y_{Q,n}$.  Thus to demonstrate that $\gamma$ acts as $-1$ on $Y_{Q,n}$, it suffices to demonstrate that $\gamma$ has no fixed points in $Y_{Q,n} \otimes \QQ$.

If $\lambda \in Y_{Q,n} \otimes \QQ$ and $\gamma \lambda = \lambda$, then define $P_\lambda^\vee$ to be the parabolic subgroup of $\tilde G^\vee$ containg $\tilde T^\vee$, whose roots are those $\tilde \alpha^\vee \in \tilde \Phi^\vee$ for which
$$\tilde \alpha(\lambda) \geq 0.$$

The conjugation action of $\gamma \To \Aut(\tilde G^\vee)$ stabilizes $\tilde P^\vee$; indeed, it stabilizes $\tilde T^\vee$ and the set of roots of $\tilde P^\vee$.  Hence the group ${}^\EL \tilde P$ generated by $\tilde P^\vee$ and $\gamma$ is a parabolic subgroup of ${}^\EL \tilde G$.  This contradicts the assumption that $\phi$ is a discrete series Weil parameter.
\qed

Now we find that the discrete series Weil parameter $\phi$, with $\phi(\CC^\times) \subset \tilde T^\vee$, satisfies
$$\phi(z) = z^{x_1} \bar z^{x_2} \in \tilde T^\vee,$$
for some $x_1, x_2 \in X_{Q,n} \otimes \CC$ satisfying $x_1 - x_2 \in X_{Q,n}$.  Since $j z j^{-1} = \bar z$, the previous lemma implies 
$$z^{x_2} \bar z^{x_1} = \phi(\bar z) = \gamma \phi(z) \gamma^{-1} = \phi(z)^{-1} = z^{-x_1} \bar z^{-x_2}.$$
Hence $x_2 = - x_1$.  Since $x_1 - x_2  \in X_{Q,n}$, we find that $2 x_1 \in X$.  Defining $\xi = x_1 \in \half X_{Q,n}$, we find that
$$\phi(z) = (z / \bar z)^{\xi} = \arg(z)^{2 \xi} \text{ for some } \xi \in \half X_{Q,n}.$$
In particular, we find
$$\gamma^2 = \phi(j^2) = \phi(-1) = e^{2 \pi i \xi}.$$
\begin{lemma}
The element $\xi \in \half X_{Q,n}$ is regular.  In particular, the centralizer of $\phi(\CC^\times)$ is the torus $\tilde T^\vee$.
\end{lemma}
\proof
This proof follows \cite[\S 10.5]{BorelCorvallis}.  If $\xi$ were singular, then the centralizer of $\phi(\CC^\times)$ would contain a nontrivial semisimple subgroup $H^\vee \subset \tilde G^\vee$, normalized by $\gamma$.  But then $\Int(\gamma)$, as an involution of $H^\vee$, would fix (pointwise) a nontrivial torus $S^\vee$ in $H^\vee$.  Thus $\Im(\phi)$ would be contained in the centralizer of $S^\vee$, which is contained in a proper parabolic subgroup of $\tilde G^\vee$, a contradiction.
\qed

\begin{proposition}
Every discrete series Weil parameter for $\alg{\tilde G}$ is equivalent to a discrete series Weil parameter $\phi$ satisfying $\phi(z) = (z / \bar z)^\xi \in \tilde T^\vee$ and $\phi(j) = \gamma \in {}^\EL \tilde G$, where $\xi$ and $\gamma$ satisfy the following conditions.
\begin{enumerate}
\item
$\xi \in \half X_{Q,n}$ is regular.
\item
$\gamma$ lies over $\sigma$ and $\gamma^2 = e^{2 \pi i \xi}$
\item
$\gamma t^\vee \gamma^{-1} = (t^\vee)^{-1}$ for all $t^\vee \in \tilde T^\vee$.
\end{enumerate}
This gives a bijection, from the set of equivalence classes of discrete series Weil parameters $\WP_\epsilon^{\dis}(\alg{\tilde G})$ to the set of $\tilde N^\vee$-orbits on the set of pairs $(\xi, \gamma)$ satisfying (1), (2), (3) above.
\end{proposition}
\proof
The previous lemmata demonstrate that every discrete series Weil parameter for $\alg{\tilde G}$ is equivalent to a parameter $\phi(z) = (z / \bar z)^\xi$, $\phi(j) = \gamma$, with $\xi$ and $\gamma$ satisfying the conditions above.  What remains is to trace through equivalence of parameters.

If $\phi_0$ is a discrete series Weil parameter for $\alg{\tilde G}$, then there exists $g^\vee \in \tilde G^\vee$ such that $\Int(g^\vee) \phi_0(z) \in \tilde T^\vee$ for all $z \in \CC^\times$.  Moreover, the regularity of $\xi$ in the previous lemma implies that $g^\vee$ is uniquely determined up to $\tilde N^\vee$.

Hence, to understand the $\tilde G^\vee$-orbits on $\WP_\epsilon^{\dis}(\alg{\tilde G})$, it suffices to verify that $\tilde N^\vee$ acts on the set of pairs $(\xi, \gamma)$ satisfying (1), (2), and (3).  If $\xi \in \half X_{Q,n}$ is regular, and $n^\vee \in \tilde N^\vee$, then $\Int(n^\vee) \xi \in \half X_{Q,n}$ is regular, so (1) is stable under the action of $\tilde N^\vee$.  If $\gamma \in {}^\EL \tilde G$ lies over $\sigma \in \Gal_\RR$, and $\gamma^2 = e^{2 \pi i \xi}$, then $\Int(n^\vee) \gamma$ lies over $\sigma \in \Gal_\RR$ and 
$$(\Int(n^\vee) \gamma )^2 = \Int(n^\vee) \gamma^2 = \Int(n^\vee) e^{2 \pi i \xi} = e^{2 \pi i \cdot \Int(n^\vee) \xi}.$$
Hence (1) and (2) are stable under the action of $\tilde N^\vee$.  Finally, if $\gamma t^\vee \gamma^{-1} = (t^\vee)^{-1}$ for all $t^\vee \in \tilde T^\vee$, then
$$\Int(n^\vee) \gamma \cdot t^\vee \cdot \Int(n^\vee) \gamma^{-1} = \Int(n^\vee) \left( \gamma \left( \Int(n^\vee)^{-1} t^\vee \right) \gamma^{-1} \right) = (t^\vee)^{-1}.$$
Hence (1), (2), and (3) are stable under the action of $\tilde N^\vee$.
\qed

Now we examine a pair $(\xi, \gamma)$ satisfying the conditions (1), (2), (3) of the above proposition.  Recall that the L-group ${}^\EL \tilde G$ is constructed as
$${}^\EL \tilde G = \frac{{}^\EL \tilde Z \ltimes \tilde G^\vee}{ \langle \zeta, \zeta^{-1} : \zeta \in \tilde Z^\vee \rangle}.$$
Here $\tilde Z^\vee$ is the center of $\tilde G^\vee$ and ${}^\EL \tilde Z$ fits into a short exact sequence
$$\tilde Z^\vee \Into {}^\EL \tilde Z \Onto \Gal_\RR.$$
In this way, ${}^\EL \tilde T$ can be viewed as a subgroup of ${}^\EL \tilde G$,
$${}^\EL \tilde T = \frac{{}^\EL \tilde Z \ltimes \tilde T^\vee}{ \langle \zeta, \zeta^{-1} : \zeta \in \tilde Z^\vee \rangle}.$$

Hence $\gamma = \zeta n^\vee$ for some $n^\vee \in \tilde G^\vee$ and some $\zeta \in {}^\EL \tilde Z$ lying over $\sigma \in \Gal_\RR$.  
\begin{lemma}
For any pair $(\xi, \gamma)$ satisfying conditions (1), (2), (3), we have $\gamma = \zeta n^\vee$ for some $\zeta \in {}^\EL \tilde Z$ and $n^\vee \in \tilde G^\vee$.  The element $n^\vee$ is contained in $\tilde N^\vee$, the normalizer of $\tilde T^\vee$, and represents a $\sigma$-fixed involution $w^\vee \in \tilde W^\vee$.  This involution acts on $Y_{Q,n}$ by $w^\vee(y) = - \sigma(y)$.
\end{lemma}
\proof
Conjugation by $\zeta$ and conjugation by $\gamma = \zeta n^\vee$ stabilize $\tilde T^\vee$.  Hence conjugation by $\zeta^{-1} \cdot \gamma = n^\vee$ stabilizes $\tilde T^\vee$.  Thus $n^\vee \in \tilde N^\vee$, the normalizer of $\tilde T^\vee$ in $\tilde G^\vee$.  Write $w^\vee$ for the associated element of the Weyl group $\tilde W^\vee = \tilde N^\vee / \tilde T^\vee$.

For all $t^\vee \in T^\vee$, we compute 
\begin{align*}
(t^\vee)^{-1} = \gamma t^\vee \gamma^{-1} & = \zeta n^\vee \cdot  t^\vee \cdot (n^\vee)^{-1} \zeta^{-1} \\
&= \zeta n^\vee \zeta^{-1} \cdot \zeta t^\vee \zeta^{-1} \cdot \zeta (n^\vee)^{-1} \zeta^{-1} \\
&= \sigma( {}^{w^\vee} t^\vee).
\end{align*}
Hence $w^\vee = - \sigma$ as automorphisms of $Y_{Q,n}$.  It follows that $\sigma \circ w^\vee = w^\vee \circ \sigma$, and $w^\vee \in (\tilde W^\vee)^\sigma$ is a Galois-fixed element of the Weyl group.

Note that $\gamma^2 = \phi(-1) \in \tilde T^\vee$.  Hence
$$\zeta n^\vee \zeta n^\vee = \zeta \zeta \cdot  \zeta^{-1} n^\vee \zeta \cdot n^\vee = \zeta^2 \sigma(n^\vee) n^\vee = \phi(-1) \in \tilde T^\vee.$$
Since $\zeta^2 \in \tilde T^\vee$, we find that ${}^\sigma n^\vee n^\vee \in \tilde T^\vee$ too.  Since ${}^\sigma w^\vee = w^\vee$, we find that $w^\vee$ is an involution in $\tilde W^\vee$.  
\qed

The function $(\xi, \gamma) \mapsto \xi$ descends to a function
$$\inf \From \WP_\epsilon^{\dis}(\alg{\tilde G}) \To \frac{\half X_{Q,n}}{W-\text{conjugation}}.$$
\begin{proposition}
The fibres of the function $\inf$ have cardinality zero or one.
\end{proposition}
\proof
We demonstrate that, if $(\xi, \gamma_1)$ and $(\xi, \gamma_2)$ are two pairs (with the same $\xi$) satisfying the conditions (1), (2), (3), then there exists $t^\vee \in \tilde T^\vee$ such that $\gamma_2 = t^\vee \gamma_1 (t^\vee)^{-1}$.  Indeed, there exist $\zeta_1, n_1^\vee, \zeta_2, n_2^\vee$ such that
$$\gamma_1 = \zeta_1 n_1^\vee, \quad \gamma_2 = \zeta_2 n_2^\vee.$$
The elements $n_1^\vee$ and $n_2^\vee$ represent the same element $w^\vee \in \tilde W^\vee$ by the previous Lemma.  Thus there exists $s^\vee \in \tilde T^\vee$ such that $n_2^\vee = n_1^\vee s^\vee$.  The elements $\zeta_1, \zeta_2 \in {}^\EL \tilde Z$ lie over $\sigma \in \Gal_\RR$, and so there exists $z^\vee \in \tilde Z^\vee$ satisfying $\zeta_2 = \zeta_1 z^\vee$.  We have
$$\gamma_1 s^\vee z^\vee = \zeta_1 n_1^\vee s^\vee z^\vee = \zeta_1 z^\vee \cdot n_1^\vee s^\vee = \zeta_2 n_2^\vee = \gamma_2.$$

Since squaring is surjective on $\tilde T^\vee$, let $t^\vee \in \tilde T^\vee$ be such that $(t^\vee)^{-2} = s^\vee z^\vee$.  The previous Lemma describes the action of $w^\vee$ on $\tilde T^\vee$, yielding
$$t^\vee \gamma_1 (t^\vee)^{-1}  = t^\vee \zeta_1 n_1^\vee (t^\vee)^{-1} = \zeta_1 n_1^\vee \cdot  {}^{\sigma w^\vee} t^\vee (t^\vee)^{-1} =  \gamma_1 (t^\vee)^{-2} = \gamma_1 s^\vee z^\vee= \gamma_2.$$
\qed 

To understand the image of $\inf$, we must understand which pairs $(\xi, \gamma)$ may occur.  
\begin{lemma}
\label{Lemmaxirhozeta}
Suppose that there exists $w^\vee \in \tilde W^\vee$ such that $w^\vee$ acts on $Y$ via $-\sigma$.  Suppose that $\xi$ is a regular element of $\half X_{Q,n}$.  Let $\rho$ be the half-sum of the positive roots of $\alg{G}$ with respect to $\alg{T}$.  Then $\xi$ occurs in the image of $\inf$ if and only if there exists $\zeta \in {}^\EL \tilde Z$ lying over $\sigma \in \Gal_\RR$ such that 
\begin{equation}
\label{XiRhoZeta}
e^{2 \pi i \xi} = e^{2 \pi i \rho} \cdot \zeta^2.
\end{equation}
\end{lemma}
\proof
For one direction, suppose that $(\xi, \gamma)$ satisfies conditions (1), (2), and (3).  Condition (2) states that $\gamma^2 = e^{2 \pi i \xi}$.  We also know that $\gamma = \zeta n^\vee$ for some $\zeta \in {}^\EL \tilde Z$ and $n^\vee \in \tilde N^\vee$ representing $w^\vee$.  Hence 
\begin{equation}
\label{Constraint2}
e^{2 \pi i \xi} = \gamma^2 = (\zeta n^\vee)^2 = \zeta^2 \cdot {}^\sigma n^\vee n^\vee.
\end{equation}
The computation of ${}^\sigma n^\vee \cdot n^\vee$ follows from a combinatorial lemma of Langlands, \cite[Lemma 3.2]{LangIrred}.
\begin{equation}
\label{LangLemma}
{}^\sigma n^\vee n^\vee = e^{2 \pi i \rho}.
\end{equation}
To compare to \cite{LangIrred}, our $n^\vee$ corresponds to Langlands' $a$, our $\rho$ corresponds to Langlands' $\delta$, and our assumption that $\alg{G}$ is semisimple implies that Langlands' $\mu^\wedge(a)$ equals 1.  Langlands' Lemma \eqref{LangLemma} and \eqref{Constraint2} imply that \eqref{XiRhoZeta} holds.

Conversely, suppose that that $\xi \in \half X_{Q,n}$ is regular and \eqref{XiRhoZeta} holds for some $\zeta \in {}^\EL \tilde Z$ lying over $\sigma$.  Define $\gamma = \zeta n^\vee$, where $n^\vee \in \tilde N^\vee$ represents $w^\vee \in \tilde W^\vee$ (acting by $-\sigma$ on $Y$).  Suppose $\xi$ and $\zeta$ satisfy \eqref{XiRhoZeta}.  Define a discrete series Weil parameter $\phi$ by setting
$$\phi(z) = (z / \bar z)^\xi, \quad \phi(j) = \gamma.$$
To see that $\phi$ defines a discrete series Weil parameter, it suffices to show that $\gamma^2 = (-1)^{2 \xi}$ and $\gamma \phi(z) \gamma^{-1} = \phi(\bar z)$.  These follow from \eqref{XiRhoZeta} and the fact that $n^\vee$ represents the appropriate involution $w^\vee$ in $\tilde W^\vee$.
\qed

\begin{remark}
If $s^\vee \in \tilde Z^\vee$, then we find
$$s^\vee = n^\vee s^\vee (n^\vee)^{-1} = {}^\sigma (s^\vee)^{-1},$$
since $w^\vee$ acts on $\tilde T^\vee$ by $t^\vee \mapsto {}^\sigma (t^\vee)^{-1}$.  Hence 
$$(\zeta s^\vee)^2 = \zeta^2 s^\vee {}^\sigma s^\vee = \zeta^2.$$
It follows that if \eqref{XiRhoZeta} holds for some $\zeta \in {}^\EL \tilde Z$ lying over $\sigma$, then \eqref{XiRhoZeta} holds for \textbf{all} $\zeta \in {}^\EL \tilde Z$ lying over $\sigma$.
\end{remark}

The previous lemma reduces the study of discrete series Weil parameters to the study of $\zeta^2$ for $\zeta \in {}^\EL \tilde Z$.  This can be related, in turn, to the L-group of an anisotropic maximal torus in $\alg{G}$.  The following lemma guarantees the existence of such a torus, and is essentially the same as \cite[Lemma 3.1]{LangIrred}
\begin{lemma}
If there exists a discrete series Weil parameter $\phi \in \WP_\epsilon^{\dis}(\alg{\tilde G})$, then there exists a maximal $\RR$-torus $\alg{S} \subset \alg{G}$ which is anisotropic, i.e., $S = \alg{S}(\RR)$ is compact.
\end{lemma}
\proof
Let $w \in W^\sigma$ be the $\Gal_\RR$-fixed element of $W$ corresponding to $w^\vee$ above.  Then $w^\vee$ acts on $Y$ as $-\sigma$, since $Y$ contains $Y_{Q,n}$ as a finite-index subgroup.  Since $w = {}^\sigma w$ and $w^2 = 1$, we find an element
$$[\eta] \in H^1(\Gal_\RR, W), \quad \eta^\vee(\Id) = 1, \eta^\vee(\sigma) = w.$$
Since we assume $\alg{G}$ is quasisplit, a result of Raghunathan \cite{Raghunathan} (proven earlier by Gill\'e \cite{Gille}) implies that $\eta$ occurs as the ``type'' of a maximal $\RR$-torus $\alg{S} \subset \alg{G}$.  

The character lattice of $\alg{S}$ is isomorphic to $Y$, but with the Galois action twisted by the cocycle $\eta$
$$\sigma_S(y) = w(\sigma(y)).$$
But $w = - \sigma$ on $Y$ and so $w(\sigma(y)) = -y$.  Therefore $\sigma_S$ has no fixed points in $Y$, and so $\alg{S}$ is anisotropic.
\qed

Now fix such an anisotropic maximal $\RR$-torus $\alg{S} \subset \alg{G}$.  We elaborate on the connection between $\alg{S}$ and the $\sigma$-fixed involution in $W$ here.  Since all maximal tori are conjugate over $\CC$, there exists $g \in G_\CC = \alg{G}(\CC)$ such that $\Int(g) \alg{S} = \alg{T}$.  Since $\alg{S}$ and $\alg{T}$ are defined over $\RR$, we find that $\Int({}^\sigma g) \alg{S} = \alg{T}$ as well.  Define $\dot w = {}^\sigma g \cdot g^{-1}$.  Then we find that $\Int(\dot w) \alg{T} = \alg{T}$, and so $\dot w \in \alg{N}_{\alg{G}}(\alg{T})(\CC)$.

Write $Y_S$ for the cocharacter lattice of $\alg{S}$, so $\Int(g) \From Y_S \To Y$ is an isomorphism of groups.  If $\dot w$ represents $w \in W$, then the following diagram commutes.
$$\begin{tikzcd}
Y_S \arrow{r}{\Int(g)} \arrow{d}{\sigma} & Y \arrow{d}{w \circ \sigma} \\
Y_S \arrow{r}{\Int(g)} & Y
\end{tikzcd}$$
Since $\alg{S}$ is anisotropic, $\sigma \From Y_S \To Y_S$ is multiplication by $-1$.  Hence $\Int(g)^{-1} \circ w \circ \sigma \circ \Int(g) = - \Id$, as automorphisms of $Y$, and so $w = - \sigma$ as an automorphism of $Y$.  Therefore ${}^\sigma w = w$ and $w \in W^{\Gal_\RR}$.  Since $\alg{G}$ is quasisplit, the $\Gal_\RR$-fixed element $w$ can be represented by a $\Gal_\RR$-fixed element of $\alg{N}_{\alg{G}}(\alg{T})(\CC)$.  Hence we may assume, without loss of generality, that 
$${}^\sigma \dot w = \dot w = {}^\sigma g \cdot g^{-1} \in \alg{G}(\RR).$$
Since ${}^\sigma \dot w \cdot \dot w = 1$ and ${}^\sigma \dot w = \dot w$, we find that $\dot w^2 = 1$ and $w^2 = 1$.  Note that $w \in W^{\Gal_\RR}$ corresponds to $w^\vee \in \tilde W^\vee$ discussed earlier.

The inner automorphism $\Int(g) \From \alg{G}_\CC \To \alg{G}_\CC$ lifts canonically to an automorphism $\Int(g) \From \alg{G}_\CC' \To \alg{G}'_\CC$ (central extensions of $\alg{G}_\CC$ by $\alg{K}_2$ on the big Zariski site over $\CC$).  Writing $\alg{S}'$ and $\alg{T}'$ for the pullbacks of $\alg{G}'$, and extending scalars to $\CC$, we find a commutative diagram of sheaves of groups on $\CC_{\Zar}$ with exact rows.
\begin{equation}
\label{STDiagram}
\begin{tikzcd}
\alg{K}_2 \inarrow{r} \arrow{d}{=} & \alg{S}_\CC' \onarrow{r} \arrow{d}{\Int(g)} & \alg{S}_\CC \arrow{d}{\Int(g)} \\
\alg{K}_2 \inarrow{r} & \alg{T}_\CC' \onarrow{r} & \alg{T}_\CC
\end{tikzcd}
\end{equation}

From the diagram \eqref{STDiagram}, the first Brylinski-Degline invariant of $\alg{S}'$ is the quadratic form
$$Q_S(y) = Q(\Int(g) y) \text{ for all } y \in Y_S.$$
We abbreviate $Y^\sharp = Y_{Q,n}$ and $Y_S^\sharp = Y_{S,Q_S,2}$ in what follows.  Then $\Int(g)$ restricts to an isomorphism from $Y_S^\sharp$ to $Y^\sharp$.

Write $\hat T = \sheaf{\hat T}[\CC] = \Hom(Y^\sharp, \CC^\times)$, and similarly $\hat S = \Hom(Y_S^\sharp, \CC^\times)$.  The underlying groups of $\hat T$ and $\hat S$ are the same as those of $T^\vee$ and $S^\vee$, but the $\Gal_\RR$-actions are different:  in $\hat T$ and $\hat S$, we consider $\CC^\times$ as a $\Gal_\RR$-module by complex conjugation, while in $T^\vee$ and $S^\vee$, we consider $\CC^\times$ as a trivial $\Gal_\RR$-module.

The action of $W$ on $Y^\sharp$ yields a homomorphism $\hat w \From \hat T \To \hat T$.  The isomorphism $\Int(g) \From Y_S^\sharp \To Y^\sharp$ yields a pullback isomorphism $\hat g \From \hat T \To \hat S$, fitting into a commutative diagram.
$$\begin{tikzcd}
\hat T \arrow{r}{\hat g} \arrow{d}{\sigma \circ \hat w} & \hat S \arrow{d}{\sigma} \\
\hat T \arrow{r}{\hat g} & \hat S
\end{tikzcd}$$

Taking the $\CC(\!(\form)\!)$-points in \eqref{STDiagram}, pushing out via $\alg{K}_2(\CC(\!(\form)\!)) \Onto \CC^\times$, and pulling back via $Y_S^\sharp \To \alg{S}(\CC(\!(\form)\!))$ and $Y^\sharp \Into \alg{T}(\CC(\!(\form)\!))$ yields a commutative diagram of abelian groups, relating the second Brylinski-Deligne invariants of $\alg{\tilde S}$ and $\alg{\tilde T}$.
\begin{equation}
\label{STBD2}
\begin{tikzcd}
\CC^\times \inarrow{r} \arrow{d}{=} & D_S^\sharp \onarrow{r} \arrow{d}{\Int(g)}  & Y_S^\sharp  \arrow{d}{\Int(g)}  \\
\CC^\times \inarrow{r} & D^\sharp  \onarrow{r} & Y^\sharp
\end{tikzcd}
\end{equation}

Recall that $\dot w \in G = \alg{G}(\RR)$ and $\dot w \cdot g = {}^\sigma g$.  It follows that, while $\Int(g) \From D_S^\sharp \To D^\sharp$ is typically not $\Gal_\RR$-equivariant, the following diagram commutes.
$$\begin{tikzcd}
D_S^\sharp \arrow{r}{\Int(g)} \arrow{d}{\sigma} & D^\sharp \arrow{d}{\Int(\dot w) \circ \sigma} \\
D_S^\sharp \arrow{r}{\Int(g)} & D^\sharp
\end{tikzcd}$$

Let $\Spl(D^\sharp)$ be the $\hat T$-torsor of splittings of $D^\sharp$, and $\Spl(D_S^\sharp)$ the $\hat S$-torsor of splittings of $D_S^\sharp$.  The $\Gal_\RR$ actions on $D^\sharp$ and $D_S^\sharp$ yield $\Gal_\RR$-actions on $\Spl(D^\sharp)$ and $\Spl(D_S^\sharp)$, compatible with those on $\hat T$ and $\hat S$.  Moreover, if $s \in \Spl(D^\sharp)$ then we find a splitting 
$$g^\ast s = \Int(g)^{-1} \circ s \circ \Int(g) \in \Spl(D_S^\sharp).$$
Similarly, we may construct a splitting $\dot w^\ast s = \Int(\dot w)^{-1} \circ s \circ \Int(\dot w) \in \Spl(D^\sharp)$.  The following diagrams commute, exhibiting the interactions of Galois actions, torsor structure, $\Int(\dot w)$, and $\Int(g)$.
$$\begin{tikzcd}
\Spl(D^\sharp) \arrow{r}{g^\ast} \arrow{d}{\dot w^\ast \circ \sigma} & \Spl(D_S^\sharp) \arrow{d}{\sigma} \\
\Spl(D^\sharp) \arrow{r}{g^\ast} & \Spl(D_S^\sharp)
\end{tikzcd}
\hfill
\begin{tikzcd}
\hat T \times \Spl(D^\sharp) \arrow{r}{\ast} \arrow{d}{\hat g \times g^\ast} & \Spl(D^\sharp) \arrow{d}{g^\ast} \\
\hat S \times \Spl(D_S^\sharp) \arrow{r}{\ast} & \Spl(D_S^\sharp)
\end{tikzcd}$$

Let $\hat g_\ast \Spl(D^\sharp)$ denote the pushout of the torsor $\Spl(D^\sharp)$ via $\hat g \From \hat T \To \hat S$.  Thus $\hat g_\ast \Spl(D^\sharp)$ has the same underlying set $\Spl(D^\sharp)$, but is now viewed as a $\hat S$-torsor.  The commutativity of the diagram above (on the right) demonstrates that 
$$g^\ast \From \hat g_\ast \Spl(D^\sharp) \xrightarrow{\sim} \Spl(D_S^\sharp)$$
is an isomorphism of $\hat S$-torsors.

Write $\sqrt{\Spl(D^\sharp)}$ for the groupoid of square roots of the $\hat T$-torsor $\Spl(D^\sharp)$.  Similarly, write $\sqrt{\Spl(D_S^\sharp)}$ for the groupoid of square roots of the $\hat S$-torsor $\Spl(S)$.  If $\bar x = (H,h) \in \sqrt{\Spl(D^\sharp)}$, then pushing out yields a square root $(\hat g_\ast H, \hat g_\ast h) \in \sqrt{\hat g_\ast \Spl(D^\sharp)}$.  Composing with the isomorphism of $\hat S$-torsors $g^\ast \From \hat g_\ast \Spl(D^\sharp) \To \Spl(D_S^\sharp)$, we find an object $g_\ast \bar x = (\hat g_\ast H, g^\ast \circ \hat g_\ast h) \in \sqrt{\Spl(D_S^\sharp)}$.  

The diagrams and discussions above define functors of groupoids,
\begin{align*}
\sigma_\ast & \From \sqrt{\Spl(D^\sharp)} \To \sqrt{\Spl(D^\sharp)} \\
\sigma_{S,\ast} & \From \sqrt{\Spl(D_S^\sharp)} \To \sqrt{\Spl(D_S^\sharp)} \\
\dot w_\ast & \From \sqrt{\Spl(D^\sharp)} \To \sqrt{\Spl(D^\sharp)} \\
g_\ast & \From \sqrt{\Spl(D^\sharp)} \To \sqrt{\Spl(D_S^\sharp)}.
\end{align*}

\begin{lemma}
Let $\dot w$ be any element of $\alg{N}_{\alg{G}}(\alg{T})(\CC)$ satisfying $\dot w^2 = 1$ (e.g, the element $\dot w = {}^\sigma g \cdot g^{-1}$ from before).  Let $\bar x = (H, h)$ be any object of $\sqrt{\Spl(D^\sharp)}$.  Then there exists an isomorphism $\dot \rho \From \bar x \To \dot w_\ast \bar x$ in the groupoid $\sqrt{\Spl(D^\sharp)}$ such that the morphism
$$\dot \rho^2 \From \bar x \xrightarrow{\dot \rho} \dot w_\ast \bar x \xrightarrow{\dot w_\ast \dot \rho} \dot w_\ast \dot w_\ast \bar x = \bar x$$
equals the identity.
\end{lemma}
\proof
Since $w$ is an involution in $W$, there exist strongly orthogonal roots $\alpha_1, \ldots, \alpha_\ell$ with associated root reflections $w_1, \ldots, w_\ell \in W$, such that $w = \prod_{i=1}^\ell w_i$ (see, for example, \cite[Proposition 1.1]{Michelson}).  Let $q_i = Q(\alpha_i^\vee)$ and let $n_i = n_{\alpha_i}$.  Thus $n_i = 1$ if $q_i$ is even, and $n_i = 2$ if $q_i$ is odd.  Define
$$\beta = \sum_{i=1}^\ell q_i \tilde \alpha_i \in X_{Q,n}.$$
Orthogonality implies that $w(\beta) = -\beta$.

Also, orthogonality implies that the coroots $\tilde \alpha_i^\vee$ are $\QQ$-linearly independently, and so there exists a splitting $s \in \Spl(D^\sharp)$ which is aligned in the sense that
\begin{equation}
\label{aligns}
s(\tilde \alpha_i^\vee) = [e_i]^{n_i}, \text{ for all } 1 \leq i \leq \ell.
\end{equation}
As in Lemma \ref{WDLemma}, if $d \in D^\sharp$ lies over $y \in Y^\sharp$, we have
$$\Int(\dot w) d = d \cdot \prod_{i = 1}^\ell \left( [e_i]^{-n_i \langle \tilde \alpha_i, y \rangle} \cdot (-1)^{q_i \varepsilon(-n_i \langle \tilde \alpha_i, y \rangle)} \right).$$
Here and below, $\varepsilon(N) = N(N+1)/2$ for any integer $N$.  As in Theorem \ref{WeylActsOnSplittings}, recall that $\varepsilon(2N) = N$ modulo $2$, and $\varepsilon(-N) - N = \varepsilon(N)$ modulo $2$.  

When $s$ is aligned as in \eqref{aligns}, we compute
\begin{align*}
[\dot w^\ast(s)](y) &= \Int(\dot w) s \left( y - \sum_{j = 1}^\ell \langle \tilde \alpha_j, y \rangle \tilde \alpha_j^\vee \right), \\
&= \Int(\dot w) s(y) \cdot \prod_j \Int(\dot w) s(\tilde \alpha_j^\vee)^{- \langle \tilde \alpha_j, y \rangle}, \\
&= s(y) \cdot \prod_i [e_i]^{- n_i \langle \tilde \alpha_i, y \rangle} \cdot (-1)^{q_i \varepsilon(-n_i \langle \tilde \alpha_i, y \rangle)}  \\
& \cdot \prod_j \left( s(\tilde \alpha_j^\vee)^{- \langle \tilde \alpha_j, y \rangle} \cdot \prod_i [e_i]^{n_i \langle \tilde \alpha_j, y \rangle \langle \tilde \alpha_i, \tilde \alpha_j^\vee \rangle} \cdot (-1)^{- q_i \langle \tilde \alpha_j, y \rangle \varepsilon(- n_i \langle \tilde \alpha_i, \tilde \alpha_j^\vee \rangle)} \right), \\
&= s(y) \cdot (-1)^{q_i \varepsilon(-n_i \langle \tilde \alpha_i, y \rangle) + q_i \langle \tilde \alpha_i, y \rangle \varepsilon(-2 n_i)}, \quad \text{(by \eqref{aligns} and orthogonality)} \\
&= s(y) \cdot (-1)^{q_i \varepsilon(-n_i \langle \tilde \alpha_i, y \rangle) - q_i n_i \langle \tilde \alpha_i, y \rangle},\\
&= s(y) \cdot \prod_{i=1}^\ell (-1)^{q_i \varepsilon( n_i \langle \tilde \alpha_i, y \rangle)}.
\end{align*}

If $n_i = 1$, then $q_i$ is even and the exponent $q_i \varepsilon(n_i \langle \tilde \alpha, y \rangle)$ is even.  If $n_i = 2$ then $\varepsilon(n_i \langle \tilde \alpha, y \rangle)$ has the same parity as $\langle \tilde \alpha, y \rangle$.  In both cases, we find
$$\prod_{i=1}^\ell (-1)^{q_i \varepsilon( n_i \langle \tilde \alpha_i, y \rangle)} = \prod_{i=1}^\ell (-1)^{q_i \langle \tilde \alpha_i, y \rangle}.$$
Hence 
$$\Int(\dot w) s = (-1)^{\beta} \ast s = e^{i \pi \beta} \ast s.$$

Since $h \From H \To \Spl(D^\sharp)$ is surjective, let $a \in H$ be such that $h(a) = s$, with $s$ aligned as above.  Let $\rho \From \bar x \To \dot w_\ast \bar x$ be any isomorphism in the connected groupoid $\sqrt{\Spl(D^\sharp)}$.  Then, as a function from $H$ to $H$, we have $\rho(a) = \hat r \ast a$ for some $\hat r \in \hat T$.    We have
$$s = h(a) = \Int(\dot w)(h(\rho(a))) = \Int(\dot w)(h(\hat r \ast a)) = \hat w(\hat r^2) \ast \Int(\dot w)(s),$$
and so
$$\hat w(\hat r)^2 = e^{i \pi \beta} \in \hat T.$$
Since $w(\beta) = - \beta$, $\hat w(e^{i \pi \beta}) = e^{- i \pi \beta} = e^{i \pi \beta}$, and so $\hat r^2 = e^{i \pi \beta}$ too.  Define
$$\dot \rho = \hat r^{-1} e^{i \pi \beta / 2} \circ \rho, \quad \dot a = e^{i \pi \beta / 4} \ast a \in H.$$
Note that $\hat r^{-1} e^{i \pi \beta / 2} \in \hat T_{[2]}$, and so $\dot \rho \From \bar x \to \dot w_\ast \bar x$ is again a morphism in $\sqrt{\Spl(D^\sharp)}$.

As in Lemma \ref{RhoAA}, we compute
\begin{align*}
\dot \rho(\dot a) &= \hat r^{-1} e^{i \pi \beta/2} \ast \rho(e^{i \pi \beta / 4} \ast a), \\
&= \hat r^{-1} e^{i \pi \beta/2} e^{-i \pi \beta / 4} \ast \rho(a), \quad (\text{since } \hat w(e^{i \pi \beta / 4}) = e^{- i \pi \beta/4} )\\
&= \hat r^{-1} e^{i \pi \beta/2} e^{-i \pi \beta / 4} \hat r \ast a, \\
&= \hat r^{-1} e^{i \pi \beta/2} e^{-i \pi \beta / 4} \hat r e^{-i \pi \beta / 4} \ast \dot a = \dot a.
\end{align*}
Hence $\dot \rho(\dot \rho(\dot a)) = \dot a$.  But $\dot \rho \circ \dot \rho$ is a $\hat T$-torsor isomorphism from $H$ to itself, and so $\dot \rho^2 = \Id$.
\qed

Now fix an object $\bar z = (H, h, j)$ of the groupoid $\gerb{E}_\epsilon(\alg{\tilde G})[\CC]$.  Define $\bar x = (H, h)$, an object of $\sqrt{\Spl(D^\sharp)}$.  Pushing out yields an object $g_\ast \bar x \in \sqrt{\Spl(D_S^\sharp)}$.

\begin{thm}
There exists an L-morphism ${}^\EL g$ making the following diagram commute.
$$\begin{tikzcd}
\tilde Z^\vee \inarrow{r} \arrow{d}{g^\vee} & {}^\EL \tilde Z \onarrow{r} \arrow{d}{{}^\EL g} & \Gal_\RR \arrow{d}{=} \\
\tilde S^\vee \inarrow{r} & {}^\EL \tilde S \onarrow{r} & \Gal_\RR
\end{tikzcd}$$
\end{thm}  
\proof
We begin by producing an L-morphism $I(g)$ making the following diagram commute.
$$\begin{tikzcd}
\tilde Z^\vee \inarrow{r} \arrow{d}{g^\vee} & \pi_1^{\et}(\gerb{E}_\epsilon(\alg{\tilde G}), \bar z) \onarrow{r} \arrow{d}{I(g)} & \Gal_\RR \arrow{d}{=} \\
\tilde S^\vee \inarrow{r} &\pi_1^{\et}(\gerb{E}_\epsilon(\alg{\tilde S}), g_\ast \bar x) \onarrow{r} & \Gal_\RR
\end{tikzcd}$$

Consider a morphism $f \From \bar z \To \sigma_\ast \bar z$ in the groupoid $\gerb{E}_\epsilon(\alg{\tilde G})[\CC]$.  Thus $f$ is an element of $\pi_1^{\et}(\gerb{E}_\epsilon(\alg{\tilde G}), \bar z)$ lying over $\sigma$.  $f$ restricts to an isomorphism from $\bar x \To \sigma_\ast \bar x$ in the groupoid $\gerb{E}_\epsilon(\alg{\tilde T})[\CC]$.  

Choose an isomorphism $\rho \From \bar x \To \dot w_\ast \bar x$ in the groupoid $\sqrt{\Spl(D^\sharp)}$ such that $\rho^2 = \Id$, as in the previous lemma.  From the proof of Theorem \ref{BigWTheorem}, the following diagram in the groupoid $\gerb{E}_\epsilon(\alg{\tilde T})[\CC]$ commutes.
$$\begin{tikzcd}
\bar x \arrow{r}{f} \arrow{d}{\rho} & \sigma_\ast \bar x \arrow{d}{\rho} \\
\dot w_\ast \bar x \arrow{r}{f} & \dot w_\ast \sigma_\ast \bar x = \sigma_\ast \dot w_\ast \bar x
\end{tikzcd}$$
The equality $\dot w_\ast \sigma_\ast \bar x = \sigma_\ast \dot w_\ast \bar x$ follows from the fact that $\dot w$ is $\Gal_\RR$-invariant.  It follows that $f^2 = (\rho f)^2 \in \tilde Z^\vee$.
$$\begin{tikzcd}
\bar x \arrow{r}{f} \arrow[bend right=10]{rrrr}[swap]{f \circ f} & \sigma_\ast \bar x \arrow{r}{\rho} & \sigma_\ast \dot w_\ast \bar x = \dot w_\ast \sigma_\ast \bar x \arrow{r}{f} & \dot w_\ast \sigma_\ast \sigma_\ast \bar x = \dot w_\ast \bar x \arrow{r}{\rho = \rho^{-1}} & \bar x.
\end{tikzcd}$$

To compute $g^\vee(f^2)$, we apply the functor $g_\ast$ throughout.
$$\begin{tikzcd}
g_\ast \bar x \arrow{r}{g_\ast f} \arrow[bend right=10]{rrrr}[swap]{g^\vee(f \circ f)}& g_\ast \sigma_\ast \bar x \arrow{r}{g_\ast \rho} & g_\ast \sigma_\ast \dot w_\ast \bar x = g_\ast \dot w_\ast \sigma_\ast \bar x \arrow{r}{g_\ast f} & g_\ast \dot w_\ast \sigma_\ast \sigma_\ast \bar x = g_\ast \dot w_\ast \bar x \arrow{r}{g_\ast \rho} & g_\ast \bar x.
\end{tikzcd}$$

But $g_\ast \dot w_\ast \sigma_\ast \bar x = \sigma_\ast g_\ast \bar x$, so the diagram above yields a commutative diagram in $\gerb{E}_\epsilon(\alg{\tilde S})[\CC]$,
\begin{equation}
\label{gfsquare}
\begin{tikzcd}[column sep = 7em]
g_\ast \bar x \arrow{r}{g_\ast (\rho f)} \arrow[bend right=10]{rr}[swap]{g^\vee(f \circ f)} & \sigma_\ast g_\ast \bar x \arrow{r}{g_\ast (\rho f)}& g_\ast \bar x
\end{tikzcd}
\end{equation}

Now we define the L-morphism $I(g)$ as follows.  Since $f$ lies over $\sigma \in \Gal_\RR$, 
$$\pi_1^{\et}(\gerb{E}_\epsilon(\alg{\tilde G}), \bar z) = \tilde Z^\vee \sqcup f \cdot \tilde Z^\vee.$$
Define $I(g) \zeta = g^\vee(\zeta)$ for all $\zeta \in \tilde Z^\vee$.  Define $I(g) f = g_\ast(\rho f)$.  To see that $I(g)$ extends uniquely to a homomorphism from $\pi_1^{\et}(\gerb{E}_\epsilon(\alg{\tilde G}), \bar z)$ to $\pi_1^{\et}(\gerb{E}_\epsilon(\alg{\tilde S}), g_\ast \bar x)$, it suffices to observe two facts.
\begin{itemize}
\item
For all $\zeta \in \tilde Z^\vee$, we have $g^\vee \sigma(\zeta) = g^\vee(\sigma(w^\vee(\zeta))) = \sigma(g^\vee(\zeta))$, and so $g^\vee$ is $\Gal_\RR$-equivariant.
\item
Diagram \eqref{gfsquare} implies that $I(g) f^2 = g^\vee(f^2) = (g_\ast (\rho f) )^2 = (I(g) f)^2$.
\end{itemize}

To finish the proof, we recall that the first Brylinski-Deligne invariant of $\alg{\tilde S}$ satisfies $Q_S(y) = Q(\Int(g) y)$ for all $y \in Y_S$.  It follows quickly that $\tau_{Q_S} \From \mu_2 \To \tilde S^\vee$ satisfies $\tau_{Q_S}(\pm 1) = g^\vee(\tau_Q(\pm 1))$.  Hence we find a L-morphism $\tau_g$ making the following diagram commute.
$$\begin{tikzcd}
\tilde Z^\vee \inarrow{r} \arrow{d}{g^\vee} & (\tau_Q)_\ast \mGal_\RR \onarrow{r} \arrow{d}{\tau_g} & \Gal_\RR \arrow{d}{=} \\
\tilde S^\vee \inarrow{r} &(\tau_{Q_S})_\ast \mGal_\RR \onarrow{r} & \Gal_\RR
\end{tikzcd}$$

The Baer sum yields the desired L-morphism ${}^\EL g = I(g) \Baer \tau_g$.
\qed

\begin{corollary}
Let $\rho \in \half X$ be the half-sum of the positive roots for $\alg{G}$ with respect to $\alg{B} \supset \alg{T}$.  Let $\kappa \in \half X_{Q,n} / X_{Q,n}$ be the element defining the double cover $\tilde S_{Q,n} \To S_{Q,n}$.  Then the map $\inf$ gives a bijection,
$$\WP_\epsilon^{\dis}(\alg{\tilde G}) \rightarrow \frac{(\kappa + \rho + X_{Q,n}) \cap X_\QQ^{\reg}}{W-\text{conjugation}}.$$
\end{corollary}
\proof
From Lemma \ref{Lemmaxirhozeta}, we find that $\xi \in \half X_{Q,n}$ is in the image of $\inf$ if and only if  $\xi$ is regular and $e^{2 \pi i \xi} = e^{2 \pi i \rho} \cdot \zeta^2$ for some $\zeta \in {}^\EL \tilde Z$ lying over $\sigma \in \Gal_\RR$.  The previous Theorem, together with Lemma \ref{gammasquare} implies that $g^\vee(\zeta^2) = (I(g) \zeta)^2 = e^{2 \pi i \kappa} \in \tilde S^\vee$.  The result follows immediately.
\qed

Recall that the infinitesimal character provided a finite-to-one surjective map,
$$\inf \From \Irr_\epsilon^{\dis}(\alg{\tilde G}) \Onto \frac{(\kappa + \rho + X) \cap X_\QQ^{\reg}}{W-\text{conjugation}}.$$
The $W$-equivariant inclusion $X \Into X_{Q,n}$ uniquely determines a parameterization of discrete series $\Lang_\epsilon^{\dis}$.
\begin{thm}
There is a unique finite-to-one function $\Lang_\epsilon^{\dis} \From \Irr_\epsilon^{\dis}(\alg{\tilde G})  \To \WP_\epsilon^{\dis}(\alg{\tilde G})$, making the following diagram commute.
$$\begin{tikzcd}
\Irr_\epsilon^{\dis}(\alg{\tilde G})  \onarrow{r}{\inf} \arrow{d}{\Lang_\epsilon^{\dis}} & \frac{(\kappa + \rho + X) \cap X_\QQ^{\reg}}{W-\text{conjugation}} \inarrow{d} \\
\WP_\epsilon^{\dis}(\alg{\tilde G}) \arrow{r}{\inf} & \frac{(\kappa + \rho + X_{Q,n}) \cap X_\QQ^{\reg}}{W-\text{conjugation}} 
\end{tikzcd}$$
\end{thm}
\proof
This follows from the bijectivity of the bottom row, and the finite-to-one nature of the top row.
\qed

This theorem provides a finite-to-one parameterization of the genuine discrete series representations of $\tilde G$, by discrete series Weil parameters valued in ${}^\EL \tilde G$.

\part*{Appendix}

\section{Torsors, gerbes, and fundamental groups}
\label{GroupsTorsors}
Let $S$ be a connected scheme, and $S_{\et}$ the \'etale site.  Our treatment of sheaves on $S_{\et}$ follows \cite[\S II]{SGA4half}.  Recall that a geometric point of $S$ is a morphism of schemes $\bar s \From \Spec(\bar F) \To S$, where $\bar F$ is a separably closed field.  

An \defined{open \'etale neighborhood} of $\bar s$ is an \'etale morphism $U \To S$ endowed with a lift $\bar u \From \Spec(\bar F) \To U$ of the geometric point $\bar s$.  If $\bar s$ is a geometric point, we write $\pi_1^{\et}(S, \bar s)$ for the \'etale fundamental group.  When $\bar s$ is fixed, we define $\Gal_S = \pi_1^{\et}(S, \bar s)$.

\subsection{Local systems on $S_{\et}$}

\begin{definition}
A \defined{local system on $S_{\et}$} is a locally constant sheaf $\sheaf{J}$ of sets on $S_{\et}$.
\end{definition}
When $\sheaf{J}$ is a local system on $S_{\et}$ and $U \To S$ is \'etale, we write $\sheaf{J}[U]$ for the set of sections over $U$ and we write $\sheaf{J}_U$ for the local system on $U_{\et}$ obtained by restriction.  If $\bar s$ is a geometric point of $S$, then the fibre $\sheaf{J}_{\bar s}$ is the inductive limit  $\limdir_U \sheaf{J}[U]$, over open \'etale neighborhoods of $\bar s$.  By local constancy, $\sheaf{J}_{\bar s} = \sheaf{J}[U]$ for some such open \'etale neighborhood.  Often in this paper we work locally on $S_{\et}$ and abuse notation a bit by writing $j \in \sheaf{J}$ rather than $j \in \sheaf{J}[U]$ (for an  \'etale $U \To S$).  

More generally, if $\Cat{C}$ is a category, then one may work with local systems on $S_{\et}$ of objects of $\Cat{C}$, or ``$\Cat{C}$-valued local systems.''  If $\Cat{C} \To \Cat{D}$ is a functor, then one finds a corresponding functor from the category of $\Cat{C}$-valued local systems to the category of $\Cat{D}$-valued local systems. Fibres of such $\Cat{C}$-valued local systems over geometric points make sense in this generality, by local constancy.

\begin{example}
\label{LocSpec}
Let $\sheaf{M}$ be a local system on $S_{\et}$ of finitely-generated abelian groups, and let $R$ be a commutative ring.  Then $\Spec(R[\sheaf{M}])$ will denote the local system on $S_{\et}$ of affine group schemes over $R$ given by
$$\Spec(R[\sheaf{M}])[U] = \Spec(R[\sheaf{M}[U]]).$$
\end{example}
We will work with local systems of groups, local systems of affine group schemes over $\ZZ$, local systems of root data, etc..

\subsection{Torsors on $S_{\et}$}
\begin{definition}
Let $\sheaf{G}$ be a sheaf of groups on $S_{\et}$.  A \defined{$\sheaf{G}$-torsor} is a locally nonempty sheaf of sets $\sheaf{V}$ on $S_{\et}$, endowed with an action $\ast \From \sheaf{G} \times \sheaf{V} \To \sheaf{V}$, such that
$$\sheaf{G} \times \sheaf{V} \To \sheaf{V} \times \sheaf{V}, \quad (g,v) \mapsto (g \ast v, v)$$
is an isomorphism of sheaves of sets on $S_{\et}$.  Morphisms of $\sheaf{G}$-torsors are morphisms of sheaves on $S_{\et}$ which intertwine the $\sheaf{G}$-action.  The \defined{category of $\sheaf{G}$-torsors} will be denoted $\Cat{Tors}(\sheaf{G})$.
\end{definition}

If $\sheaf{V}$ is a $\sheaf{G}$-torsor, we write $[ \sheaf{V} ]$ for its isomorphism class.  The isomorphism classes of $\sheaf{G}$-torsors form a pointed set denoted $H_{\et}^1(S, \sheaf{G})$.  The set is pointed by the isomorphism class of the \defined{neutral} $\sheaf{G}$-torsor:  $\sheaf{G}$ itself, as a $\sheaf{G}$-torsor by left-multiplication.  If $\sheaf{V}$ is a $\sheaf{G}$-torsor, $U \To S$ is \'etale, and $v \in \sheaf{V}[U]$, then write $v_{U'}$ for the restriction of $v$ to any further \'etale $U' \To U$.  For any such $U'$, and any $w \in \sheaf{V}[U']$, there exists a unique $g \in \sheaf{G}[U']$ satisfying $w = g \ast v_{U'}$.  Allowing $U'$ to vary, the map $w \mapsto g$ gives a $\sheaf{G}_U$-torsor isomorphism from $\sheaf{V}_{U}$ to $\sheaf{G}_U$.  Thus we say that the point $v \in \sheaf{V}[U]$ \defined{neutralizes} the torsor $\sheaf{V}$ over $U$.


The category of torsors has more structure in the abelian case.  If $\sheaf{A}$ is a sheaf of {\em abelian} groups on $S_{\et}$, then the category $\Cat{Tors}(\sheaf{A})$ inherits a monoidal structure.  Namely, if $\sheaf{V}_1$ and $\sheaf{V}_2$ are two $\sheaf{A}$-torsors, define
$$\sheaf{V}_1 \Baer \sheaf{V}_2 = \frac{\sheaf{V}_1 \times \sheaf{V}_2}{ (a \ast v_1, v_2) \sim (v_1, a \ast v_2)}.$$
With this monoidal structure, the trivial torsor $\sheaf{A}$ as zero object, and obvious isomorphisms for commutativity and associativity and unit, the category $\Cat{Tors}(\sheaf{A})$ becomes a Picard groupoid (i.e., a strictly commutative Picard category, in the terminology of \cite[Expos\'e XVIII \S 1.4]{SGA4T3}.  The pointed set of isomorphism classes $H_{\et}^1(S, \sheaf{A})$ becomes an abelian group, with
$$[\sheaf{V}_1] + [\sheaf{V}_2] \defeq \left[ \sheaf{V}_1 \Baer \sheaf{V}_2 \right].$$
The group $H_{\et}^1(S, \sheaf{A})$ is identified with the \'etale cohomology with coefficients $\sheaf{A}$.

Suppose that $f \From \sheaf{A} \To \sheaf{G}$ is a homomorphism of sheaves of groups on $S_{\et}$, with $\sheaf{A}$ abelian, and $f$ central (i.e., $f$ factors through the inclusion of the center $\sheaf{Z} \Into \sheaf{G}$).  If $\sheaf{V}$ is an $\sheaf{A}$-torsor, then we write $f_\ast \sheaf{V}$ for the pushout,
$$f_\ast \sheaf{V} \defeq \frac{ \sheaf{G} \times \sheaf{V} }{ (g,a \ast v) \sim (f(a) \cdot g,v) }.$$
This operation of \defined{pushing out torsors} defines a functor,
$$f_\ast \From \Cat{Tors}(\sheaf{A}) \To \Cat{Tors}(\sheaf{G}).$$
If $g \in \sheaf{G}$ and $v \in \sheaf{V}$, we write $g \wedge v$ for its image in $f_\ast \sheaf{V}$.  Then 
$$g \wedge (a \ast v) = g f(a) \wedge v = f(a)g \wedge v.$$
The $\sheaf{G}$-torsor structure is given by
$$\gamma \ast (g \wedge v) = (\gamma g) \wedge v, \text{ for all } \gamma, g \in \sheaf{G}, v \in \sheaf{V}.$$

Suppose that $c \From \sheaf{A}_1 \To \sheaf{A}_2$ is a homomorphism of sheaves of abelian groups on $S_{\et}$.  Let $\sheaf{V}_1$ be an $\sheaf{A}_1$-torsor and $\sheaf{V}_2$ an $\sheaf{A}_2$-torsor.  A map of torsors $\pi \From \sheaf{V}_1 \To \sheaf{V}_2$ \defined{lying over} $c$ means a morphism of sheaves of sets on $S_{\et}$ satisfying
$$\pi(a_1 \ast v_1) = c(a_1) \ast \pi(v_1), \text{ for all } a_1 \in \sheaf{A}_1, v_1 \in \sheaf{V}_1.$$
Such a map factors uniquely through $c_\ast \sheaf{V}_1$.  

A short exact sequence of sheaves of abelian groups on $S_{\et}$,
\begin{equation}
\label{ABC}
\sheaf{A} \xhookrightarrow{\alpha}\sheaf{B} \xtwoheadrightarrow{\beta} \sheaf{C},
\end{equation}
yields two more constructions of torsors.

First, the sequence yields a boundary map in cohomology, $\partial \From H_{\et}^0(S, \sheaf{C}) \To H_{\et}^1(S, \sheaf{A})$.  There is a corresponding map from global sections of $\sheaf{C}$ to objects of the category of $\sheaf{A}$-torsors as follows.

Begin with $c \in \sheaf{C}[S]$ and write $[c]$ to consider it as an element of $H_{\et}^0(S, \sheaf{C})$.  For any \'etale $U \To S$, write $c_U \in \sheaf{C}[U]$ for the restriction of $c$ to $U$.  Define $\partial c$ to be the sheaf on $S_{\et}$ whose sections are given by
$$\partial c[U] = \{ b \in \sheaf{B}[U] :  \beta(b) = c_U \}.$$
The sheaf $\partial c$ is naturally an $\sheaf{A}$-torsor; the equivalence class $[\partial c] \in H_{\et}^1(S, \sheaf{A})$ coincides with $\partial [c]$.  The sheaf $\partial c$ is called the \defined{torsor of liftings} of $c$ via $\beta$.

Next, write $\shom(\sheaf{C}, \sheaf{B})$ for the sheaf of homomorphisms (``sheaf-hom'') from $\sheaf{C}$ to $\sheaf{B}$.  This is a sheaf of abelian groups on $S_{\et}$, and there is a subsheaf of sets $\sspl(\sheaf{B})$ consisting of those homomorphisms which split the extension \eqref{ABC}.  This subsheaf $\sspl(\sheaf{B})$ is naturally a $\shom(\sheaf{C}, \sheaf{A})$-torsor, called the \defined{torsor of splittings}.  This corresponds to a familiar map in the cohomology of sheaves of abelian groups,
$$\Ext(\sheaf{C}, \sheaf{A}) \To H_{\et}^1(S, \shom(\sheaf{C}, \sheaf{A})).$$

\subsection{Gerbes on $S_{\et}$}

Here we introduce a class of gerbes on $S_{\et}$.  In what follows, let $\sheaf{A}$ be a sheaf of {\em abelian} groups on $S_{\et}$.
\begin{definition}
A \defined{gerbe on $S_{\et}$ banded by $\sheaf{A}$} is a (strict) stack $\gerb{E}$ on $S_{\et}$ of groupoids such that $\gerb{E}$ is locally nonempty, locally connected, and banded by $\sheaf{A}$.
\end{definition}
We unravel this definition here, beginning with the data.
\begin{description}
\item[a (strict) stack $\gerb{E}$ on $S_{\et}$ of groupoids]  
For each \'etale $U \To S$, we have a (possibly empty) groupoid $\gerb{E}[U]$.  For $\gamma \From U' \To U$, a morphism of schemes \'etale over $S$, we have a pullback functor $\gamma^\ast \From \gerb{E}[U] \To \gerb{E}[U']$.  
\item[banded by $\sheaf{A}$]
For every object $x$ of $\gerb{E}[U]$, there is given an isomorphism $\sheaf{A}[U] \To \Aut(x)$ (written $\alpha \mapsto \alpha_x$).
\end{description}
This data satisfies additional axioms:
\begin{description}
\item[(strict) stack axioms]
For each pair $\gamma \From U' \To U$ and $\delta \From U'' \To U'$, we require {\em equality} of functors $\delta^\ast \circ \gamma^\ast = (\gamma \delta)^\ast$.  (``Strictness'' refers to the requirement of equality rather than extra data of natural isomorphisms).  Descent for objects and morphisms is effective.  
\item[locally nonempty]
There exists a finite \'etale $U \To S$ such that $\gerb{E}[U]$ is nonempty.
\item[locally connected]
For any \'etale $U \To S$ and pair of objects $x,y \in \gerb{E}[U]$, there exists a finite \'etale $\gamma \From U' \To U$ such that $\gamma^\ast x$ is isomorphic to $\gamma^\ast y$ in $\gerb{E}[U']$.
\item[banding]
Given a morphism $f \From x \To y$ in $\gerb{E}[U]$, and $\alpha \in \sheaf{A}[U]$, $\alpha_y \circ f = f \circ \alpha_x$.  Also, given $\gamma \From U' \To U$ \'etale, $\alpha_{\gamma^\ast x} = \gamma^\ast \alpha_x \in \Aut(\gamma^\ast x)$.   
\end{description}

\begin{remark}
We will not require gerbes banded by nonabelian groups -- commutativity greatly simplifies the theory.  For a fuller treatment of gerbes, one can consult the original book of Giraud \cite{Giraud}, work of Breen \cite{BreenNotes}, \cite{BreenOriginal}, the book of Brylinski \cite[Chapter V]{BrylBook}, the article of Deligne \cite{DelModere}, the introduction of Garland and Patnaik \cite{Garl}, and the Stacks Project \cite{Stacks}, among others.  We don't keep track of universes along the way, while Giraud \cite{Giraud} is careful about set-theoretic subtleties.  Our ``strictness'' assumption is typically referred to as an assumption that the fibred category $\gerb{E} \To S_{\et}$ is ``split'' (see \cite{Vistoli}).  The strictness assumption is not so restrictive, since every stack is equivalent to a strict stack (cf. \cite[Theorem 3.45]{Vistoli}).  See also \cite[\S I.1]{Giraud}.
\end{remark}

If $\gerb{E}$ is a gerbe on $S_{\et}$ banded by $\sheaf{A}$, and $U \To S$ is \'etale, then we write $\gerb{E}_U$ for its restriction to $U_{\et}$; this is a gerbe on $U_{\et}$ banded by $\sheaf{A}_U$.

Our strictness assumption allows us to easily define the \defined{fibre} of a gerbe $\gerb{E}$ at a geometric point $\bar s$.  This is the category $\gerb{E}_{\bar s}$ whose object set is the direct limit $\limdir_{U} \gerb{E}[U]$ of object sets, indexed by open \'etale neighborhoods of $\bar s$.  Write $A = \sheaf{A}_{\bar s}$ for the fibre of $\sheaf{A}$ over $\bar s$.  If $\bar z_1, \bar z_2$ are objects of $\gerb{E}_{\bar s}$, then $\Hom(\bar z_1, \bar z_2)$ is naturally an $A$-torsor.

If $x,y \in \gerb{E}[U]$, then for $\gamma \From U' \To U$ define
$$\shom(x,y)[U'] = \Hom(\gamma^\ast x, \gamma^\ast y).$$
In this way, we construct a sheaf $\shom(x,y)$ of sets on $U_{\et}$.  In particular, we find a sheaf of groups $\sAut(x)$ on $U_{\et}$.  The banding provides an isomorphism of sheaves of groups, $\sheaf{A}_U \xrightarrow{\sim} \sAut(x)$, and $\shom(x,y)$ becomes an $\sheaf{A}_U$-torsor.

\subsubsection{Functors of gerbes}

Suppose that $c \From \sheaf{A}_1 \To \sheaf{A}_2$ is a homomorphism of sheaves of abelian groups on $S_{\et}$, and $\gerb{E}_1, \gerb{E}_2$ are gerbes on $S_{\et}$ banded by $\sheaf{A}_1, \sheaf{A}_2$, respectively.  A \defined{functor of gerbes} $\phi \From \gerb{E}_1 \To \gerb{E}_2$, \defined{lying over $c$}, is a (strict) functor of stacks lying over $c$.  This entails the following.
\begin{description}
\item[(strict) functor of stacks]
For each \'etale $U \To S$, a functor of categories $\phi[U] \From \gerb{E}_1[U] \To \gerb{E}_2[U]$.  For every $\gamma \From U' \To U$, with pullback functors $\gamma_1^\ast$ in $\gerb{E}_1$ and $\gamma_2^\ast$ in $\gerb{E}_2$, the ``strictness'' condition requires an {\em equality} of functors, $\gamma_2^\ast \circ \phi[U] = \phi[U'] \circ \gamma_1^\ast$.
\item[lying over $c$]
The ``lying over $c$'' condition requires that, for each $\alpha_1 \in \sheaf{A}_1[U]$ with $\alpha_2 = c(\alpha_1)$, and object $x_1 \in \gerb{E}_1[U]$ with $x_2 = \phi[U] x_1$, we have
$$(\alpha_2)_{x_2} = \phi[U] \left( (\alpha_1)_{x_1} \right) \in \Aut(x_2).$$
\end{description} 

Gerbes on $S_{\et}$ banded by a fixed sheaf of abelian groups $\sheaf{A}$ form a 2-category (in the sense of \cite[\S I.1.8]{Giraud}); if $\gerb{E}_1$ and $\gerb{E}_2$ are two such gerbes banded by the same $\sheaf{A}$, an \defined{equivalence of gerbes} $\phi \From \gerb{E}_1 \To \gerb{E}_2$ is a functor of gerbes lying over $\Id \From \sheaf{A} \To \sheaf{A}$.  Given two such equivalences of gerbes $\phi, \phi' \From \gerb{E}_1 \To \gerb{E}_2$, a natural isomorphism $\phi \xRightarrow{\sim} \phi'$ consists of natural isomorphisms of functors $\phi[U] \Rightarrow \phi'[U]$ for each $U$, compatible with pullback.  This defines a 2-category of gerbes banded by $\sheaf{A}$, equivalences, and natural isomorphisms of equivalences. 


Given two gerbes $\gerb{E}_1, \gerb{E}_2$ banded by $\sheaf{A}$, one may ``contract'' them to form another gerbe $\gerb{E}_1 \Baer \gerb{E}_2$ banded by $\sheaf{A}$.    The family of categories of torsors, $\stors(\sheaf{A})$, given by $\stors(\sheaf{A})[U] = \Cat{Tors}(\sheaf{A}_U)$ (for each \'etale $U \To S$), with pullbacks given by restriction of sheaves, forms the \defined{neutral} $\sheaf{A}$-gerbe on $S_{\et}$.  $\gerb{E} \Baer \stors(\sheaf{A})$ is equivalent to $\gerb{E}$, for any gerbe $\gerb{E}$ banded by $\sheaf{A}$ (and the equivalence is determined up to unique natural isomorphism).

Suppose that $\gerb{E}$ is a gerbe on $S_{\et}$ banded by $\sheaf{A}$, and $x$ is an object of $\gerb{E}[U]$ for some \'etale $U \To S$.  Then, for $\gamma \From U' \To U$ \'etale, and $y \in \gerb{E}[U']$, we have a $\sheaf{A}_{U'}$-torsor $\shom(\gamma^\ast x, y)$.  This map $y \mapsto \shom(\gamma^\ast x, y)$ extends to an equivalence of gerbes from $\gerb{E}_U$ to $\stors(\sheaf{A})_U$.  In this way, we say that $x$ \defined{neutralizes} the gerbe $\gerb{E}$ over $U$.

If $\gerb{E}$ is a gerbe banded by $\sheaf{A}$, we write $\left[ \gerb{E} \right]$ for its equivalence class.  The set of such equivalence classes is denoted $H_{\et}^2(S, \sheaf{A})$.  This forms an abelian group, with zero corresponding to the neutral gerbe of $\sheaf{A}$-torsors, and addition arising from contraction.  From \cite{Giraud}, we identify $H_{\et}^2(S, \sheaf{A})$ with the \'etale cohomology of $S$ with coeffiencts $\sheaf{A}$.

\subsubsection{Pushouts}
\label{PushoutGerbe}
Given a gerbe $\gerb{E}_1$ banded by $\sheaf{A}_1$, and $c \From \sheaf{A}_1 \To \sheaf{A}_2$ as above, one may construct a gerbe $c_\ast \gerb{E}_1$ banded by $\sheaf{A}_2$ called the \defined{pushout} of $\gerb{E}_1$ by $c$.  Any functor of gerbes $\phi \From \gerb{E}_1 \To \gerb{E}_2$ lying over $c$ factors through a functor $c_\ast \gerb{E}_1 \To \gerb{E}_2$ (the functor being determined uniquely up to unique natural isomorphism, see \cite[\S 5.3]{DelModere}).  The objects of $c_\ast \gerb{E}_1$ are the same as those of $\gerb{E}_1$.  But, given two such objects $x,y \in \gerb{E}_1[U]$, the morphism set $\Hom_{c_\ast \gerb{E}_1}(x,y)$ is defined as the pushout of torsors,
$$\Hom_{c_\ast \gerb{E}_1}(x,y) = c_\ast \Hom_{\gerb{E}_1}(x,y).$$
If $x,y \in \gerb{E}_1[U]$, $f \in \Hom_{\gerb{E}_1}(x,y)$, and $\alpha_2 \in \sheaf{A}_2[U]$, we write $\alpha_2 \wedge f $ for the resulting morphism from $x$ to $y$ in $c_\ast \gerb{E}_1$.  If $\alpha_1 \in \sheaf{A}_1[U]$, then we have
$$\alpha_2 \wedge (\alpha_{1,y} \circ f) = (\alpha_2 \cdot c(\alpha_1)) \wedge f.$$

The pushout of gerbes corresponds to the map in cohomology,
$$H_{\et}^2(S, \sheaf{A}_1) \To H_{\et}^2(S, \sheaf{A}_2), \quad [\gerb{E}_1] \mapsto c_\ast [\gerb{E}_1].$$
See \cite[Chapitre IV, \S 3.3, 3.4]{Giraud} for details.

\subsubsection{The gerbe of liftings}

\label{AppendixGerbeLiftings}

If $\sheaf{A} \xhookrightarrow{\alpha}\sheaf{B} \xtwoheadrightarrow{\beta} \sheaf{C}$ is a short exact sequence of sheaves of abelian groups on $S_{\et}$, then the sequence of cohomology groups,
$$H_{\et}^2(S, \sheaf{A}) \xrightarrow{\alpha} H_{\et}^2(S, \sheaf{B}) \xrightarrow{\beta} H_{\et}^2(S, \sheaf{C})$$
is also exact.  The analogous construction with gerbes is the following:  suppose that $\gerb{E}$ is a gerbe banded by $\sheaf{B}$, $\gerb{F}$ is a gerbe banded by $\sheaf{C}$, and $\gerb{p} \From \gerb{E} \To \gerb{F}$ is a functor of gerbes lying over $\sheaf{B} \xrightarrow{\beta} \sheaf{C}$.  If $z$ is an $S$-object of $\gerb{F}$ (neutralizing $\gerb{F}$, so that $0 = [\gerb{F}] \in H_{\et}^2(S, \sheaf{C})$), then cohomology suggests that $\gerb{E}$ arises as the pushout of a gerbe banded by $\sheaf{A}$.

Indeed, we define the gerbe $\gerb{p}^{-1}(z)$ as follows:  the objects of $\gerb{p}^{-1}(z)[U]$ are pairs $(y, j)$ where $y$ is an object of $\gerb{E}$[U], and $j \From \gerb{p}(y) \To z$ is an isomorphism in $\gerb{F}[U]$.  The morphisms in $\gerb{p}^{-1}(z)$ are those in $\gerb{E}$ which are compatible with the isomorphisms to $z$.  The gerbe $\gerb{p}^{-1}(z)$ will be called the \defined{gerbe of liftings} of $z$ via $p$.  It is a gerbe banded by $\sheaf{A}$, and there is a natural equivalence from $\alpha_\ast \gerb{p}^{-1}(z)$ to $\gerb{E}$, given by ``forgetting $j$.''

\subsubsection{The gerbe of $n^{\th}$ roots}
\label{AppendixGerbeRoots}
Suppose that $\sheaf{C}$ is a sheaf of abelian groups on $S_{\et}$, and the homomorphism $\sheaf{C} \xrightarrow{n} \sheaf{C}$ is surjective.  An important example of a gerbe of liftings arises from the Kummer sequence $\sheaf{C}_{[n]} \Into \sheaf{C} \xtwoheadrightarrow{n} \sheaf{C}$.

Pushing out gives to a functor of gerbes, $n_\ast \From \stors(\sheaf{C}) \To \stors(\sheaf{C})$, lying over $\sheaf{C} \xtwoheadrightarrow{n} \sheaf{C}$.  Given a $\sheaf{C}$-torsor $\sheaf{V}$, the gerbe of liftings of $\sheaf{V}$ via $n_\ast$ will be called the \defined{gerbe of $n^{\th}$ roots}, denoted $\sqrt[n]{\sheaf{V}}$.  It is a gerbe on $S_{\et}$ banded by $\sheaf{C}_{[n]}$.  The map which sends a $\sheaf{C}$-torsor to its gerbe of $n^{\th}$ roots corresponds to the Kummer coboundary map $\Kum \From H_{\et}^1(S, \sheaf{C}) \To H_{\et}^2(S, \sheaf{C}_{[n]})$.

Explicitly, an object of $\sqrt[n]{\sheaf{V}}$ is a pair $(\sheaf{H}, h)$ where $\sheaf{H}$ is a $\sheaf{C}$-torsor, and $h \From \sheaf{H} \To \sheaf{V}$ is a morphism of sheaves making the following diagram commute.
$$\begin{tikzcd}
\sheaf{C} \times \sheaf{H} \arrow{r}{\ast} \arrow{d}{n \times h} & \sheaf{H} \arrow{d}{h} \\
\sheaf{C} \times \sheaf{V} \arrow{r}{\ast} & \sheaf{V}
\end{tikzcd}$$

The construction of the gerbe of $n^{\th}$ roots is itself functorial.  Consider a homomorphism of sheaves of abelian groups, $c \From \sheaf{C}_1 \To \sheaf{C}_2$, and assume that $\sheaf{C}_1 \xrightarrow{n} \sheaf{C}_1$ and $\sheaf{C}_2 \xrightarrow{n} \sheaf{C}_2$ are surjective.  Suppose that $\sheaf{V}_1$ is a $\sheaf{C}_1$-torsor, and $\sheaf{V}_2$ is a $\sheaf{C}_2$-torsor.  Suppose that $f \From \sheaf{V}_1 \To \sheaf{V}_2$ is a morphism of torsors lying over the homomorphism $c \From \sheaf{C}_1 \To \sheaf{C}_2$.  Then we find a functor of gerbes $\sqrt[n]{f} \From \sqrt[n]{\sheaf{V}_1} \To \sqrt[n]{\sheaf{V}_2}$ lying over the homomorphism of bands $\sheaf{C}_{1,[n]} \To \sheaf{C}_{2,[n]}$.

\subsection{Fundamental group}
Now let $\sheaf{A}$ be a \textbf{local system} of abelian groups on $S_{\et}$.  Let $\gerb{E}$ be a gerbe on $S_{\et}$ banded by $\sheaf{A}$.  Let $\bar F$ be a separably closed field, let $\bar s \From \Spec(\bar F) \To S$ be a geometric point, and recall that $\Gal_S = \pi_1^{\et}(S,\bar s)$ denotes the \'etale fundamental group.  

Suppose that $U \To S$ is a Galois cover, and $\bar u \From \Spec(\bar F) \To U$ lifts the geometric point $\bar s$.  Writing $\Gal_U = \pi_1^{\et}(U, \bar u)$, we find a short exact sequence
$$\Gal_U \Into \Gal_S \Onto \Gal(U / S).$$
If $\gamma \in \Gal_S$ we write $\gamma_U$ for its image in $\Gal(U / S)$.  

Suppose moreover that $\sheaf{A}_U$ is a constant sheaf and $\gerb{E}[U]$ is a nonempty groupoid (i.e., $\gerb{E}$ is neutral over $U$).  Write $A = \sheaf{A}[U]$ for the resulting abelian group.  Then $A$ is endowed with an action of $\Gal_S$ that factors through the finite quotient $\Gal(U/S)$.  An object $z \in \gerb{E}[U]$ will be called a \defined{base point} for the gerbe $\gerb{E}$ (over $U$).  The banding identifies $A$ with the automorphism group of $z$.

Without loss of generality, pulling back to a larger Galois cover if necessary, we may assume that $\Hom(z, \gamma_U^\ast z)$ is nonempty for all $\gamma \in \Gal_S$.  In this way, the base point $z \in \gerb{E}[U]$ and $\gamma \in \Gal_S$ define an $A$-torsor,
$$\Aut_\gamma(z) \defeq \Hom(z, \gamma_U^\ast z).$$
We write $\gamma^\ast$ instead of $\gamma_U^\ast$, when there is little chance of confusion.  

Define the \defined{\'etale fundamental group of the gerbe} $\gerb{E}$, at the base point $z$, by
$$\pi_1(\gerb{E}, z) = \bigsqcup_{\gamma \in \Gal_S} \Aut_\gamma(z).$$
The group structure is given, for $\gamma_1, \gamma_2 \in \Gal_S$, by the following sequence.
\begin{align*}
\Aut_{\gamma_1}(z) \times \Aut_{\gamma_2}(z) &= \Hom(z, \gamma_1^\ast z) \times \Hom(z, \gamma_2^\ast z)  \\
& \xrightarrow{\gamma_2^\ast \times \Id} \Hom(\gamma_2^\ast z, \gamma_2^\ast \gamma_1^\ast z) \times \Hom(z, \gamma_2^\ast z) \\
& \xrightarrow{\circ} \Hom(z, \gamma_2^\ast  \gamma_1^\ast z) \\
& \xrightarrow{=} \Hom(z, (\gamma_1 \gamma_2)^\ast z) = \Aut_{\gamma_1 \gamma_2}(z)
\end{align*}
As $\Aut_{\Id}(z) = \Aut(z)$, the isomorphism $A \xrightarrow{\sim} \Aut(z)$, $\alpha \mapsto \alpha_z$, gives an extension of groups,
\begin{equation}
\label{FundGpGerbe}
A \Into \pi_1^{\et}(\gerb{E}, z) \Onto \Gal_S.
\end{equation}
If $\gamma \in \Gal_U \subset \Gal_S$, then $\gamma_U = \Id$ and so $\Aut_\gamma(z) = \Aut_{\Id}(z)$.  In this way, we find a splitting $\Gal_U \Into \pi_1(\gerb{E}, z)$.  In other words, $\pi_1^{\et}(\gerb{E}, z)$ arises as the pullback of an extension of $\Gal(U/S)$ by $A$.  The conjugation action of $\Gal_S$ on $A$, in the extension \eqref{FundGpGerbe}, coincides with the canonical action of $\Gal(U/S)$ on $A = \sheaf{A}[U]$.

The sequence \eqref{FundGpGerbe} describes the fundamental group of a gerbe (with base point) as an extension of $\Gal_S$ by $A$.  Here we analyze how this fundamental group depends on the choice of base point, and how it behaves under equivalence of gerbes.  

Consider a further Galois cover $\delta \From U' \To U$ and geometric base point $\bar u'$ lifting $\bar u$.  By constancy of $\sheaf{A}_U$, we identify $A = \sheaf{A}[U] = \sheaf{A}[U']$.  For all $\gamma \in \Gal_S$, we have $\gamma_U \circ \delta = \delta \circ \gamma_{U'}$.  This defines an isomorphism of $A$-torsors, $\delta^\ast \From \Aut_\gamma(z) \xrightarrow{\sim} \Aut_\gamma(\delta^\ast z)$, using the sequence below.
\begin{align*}
\Aut_\gamma(z) = \Hom(z, \gamma_U^\ast z) & \xrightarrow{\delta^\ast} \Hom( \delta^\ast z, \delta^\ast \gamma_U^\ast z) \\
& \xrightarrow{=} \Hom \left( \delta^\ast z, (\gamma_U \delta)^\ast z \right) \\
& \xrightarrow{=} \Hom \left( \delta^\ast z, (\delta \gamma_{U'})^\ast z \right) \\
& \xrightarrow{=} \Hom \left( \delta^\ast z, \gamma_{U'}^\ast \delta^\ast z \right)  = \Aut_\gamma(\delta^\ast z).
\end{align*}
Putting these isomorphisms together, we find an isomorphism of extensions.
$$\begin{tikzcd}
A \inarrow{r} \arrow{d}{=} & \pi_1^{\et}(\gerb{E}, z) \onarrow{r} \arrow{d}{\iota_\delta} & \Gal_S \arrow{d}{=} \\
A \inarrow{r} & \pi_1^{\et}(\gerb{E}, \delta^\ast z) \onarrow{r} & \Gal_S
\end{tikzcd}$$

A further cover $\delta' \From U'' \To U'$, with $\delta'' = \delta \circ \delta' \From U'' \To U$, gives a commutative diagram in the category of extensions of $\Gal_S$ by $A$.
$$\begin{tikzcd}
\pi_1^{\et}(\gerb{E}, z) \arrow{r}[swap]{\iota_\delta} \arrow[bend left=20]{rr}{\iota_{\delta''}} & \pi_1^{\et}(\gerb{E}, \delta^\ast z)  \arrow{r}[swap]{\iota_{\delta'}} & \pi_1^{\et}(\gerb{E}, (\delta')^\ast \delta^\ast z) = \pi_1^{\et}(\gerb{E}, (\delta'')^\ast z)  \end{tikzcd}$$

Define $\bar z \in \gerb{E}_{\bar s}$ to be the image of the base point $z$ in the direct limit.  We call $\bar z$ a \defined{geometric base point} for the gerbe $\gerb{E}$.  Define
$$\pi_1^{\et}(\gerb{E}, \bar z) = \limdir_{U'} \pi_1^{\et}(\gerb{E}, \delta^\ast z),$$
the direct limit over Galois covers $\delta \From (U', \bar u') \To (U, \bar u)$, via the isomorphisms $\iota_\delta$ described above.  This gives an extension of groups, depending (up to unique isomorphism) only on the {\em geometric} base point $\bar z \in \gerb{E}_{\bar s}$.
$$A \Into \pi_1^{\et}(\gerb{E}, \bar z) \Onto \Gal_S.$$
This extension is also endowed with a family of splittings over finite-index subgroups $\Gal_U \subset \Gal_S$, arising from base points $z \in \gerb{E}[U]$ mapping to $\bar z$.  Having such splittings is useful for topological purposes, e.g., $\pi_1^{\et}(\gerb{E}, \bar z) \To \Gal_S$ is naturally a continuous homomorphism of profinite groups when $A$ is finite.

Consider a second geometric base point $\bar z_0 \in \gerb{E}_{\bar s}$ (over the {\em same} $\bar s \To S$).  There exists an isomorphism $\bar f \From \bar z_0 \To \bar z$ in $\gerb{E}_{\bar s}$.  For a sufficiently large Galois cover $U \To S$, we may assume that $\bar f \From \bar z_0 \To \bar z$ arises from a morphism $f \From z_0 \To z$ in $\gerb{E}[U]$.  

Define $\iota_f \From \Aut_\gamma(z_0) \To \Aut_\gamma(z)$ to be the bijection
$$\iota_f(\eta) = \gamma^\ast f \circ \eta \circ f^{-1}, \text{ for all } \eta \in \Aut_\gamma(z_0),$$
making the following diagram commute (in the groupoid $\gerb{E}[U]$).
$$\begin{tikzcd}
z_0 \arrow{r}{f} \arrow{d}{\eta} & z \arrow{d}{\iota_f(\eta)} \\
\gamma^\ast z_0 \arrow{r}{\gamma^\ast f} & \gamma^\ast z 
\end{tikzcd}$$
As $\gamma$ varies over $\Gal_S$, this provides an isomorphism of extensions, $\iota_f \From  \pi_1^{\et}(\gerb{E}, z_0) \To  \pi_1^{\et}(\gerb{E}, z)$.  Passing to the direct limit, we find an isomorphism of extensions depending only on $\bar f \From \bar z_0 \To \bar z$,
$$\begin{tikzcd}
A \inarrow{r} \arrow{d}{=} & \pi_1^{\et}(\gerb{E}, \bar z_0) \onarrow{r} \arrow{d}{\iota_{\bar f}} & \Gal_S \arrow{d}{=} \\
A \inarrow{r} & \pi_1^{\et}(\gerb{E}, \bar z) \onarrow{r} & \Gal_S
\end{tikzcd}$$

Given {\em another} isomorphisms $\bar g \From \bar z_0 \To \bar z$, there exists a unique element $\alpha \in A$ such that $\bar g = \alpha_{\bar z} \circ \bar f$.  As for $f$, we may assume that $\bar g$ arises from $g \From z_0 \To z$ in $\gerb{E}[U]$.  It follows that, for all $\eta \in \Aut_\gamma(z_0)$,
\begin{align*}
\iota_g(\eta) &= \gamma^\ast g \circ \eta \circ g^{-1} \\
&= \gamma^\ast (\alpha_z \circ f) \circ \eta \circ f^{-1} \circ \alpha_z^{-1} \\
&= \gamma^\ast \alpha_z \circ \iota_f(\eta) \circ \alpha_z^{-1} \\
&= \alpha_{\gamma^\ast z} \circ \iota_f(\eta) \circ \alpha_z^{-1}
\end{align*}
In other words, we have $\iota_{\bar g} = \Int(\alpha) \circ \iota_{\bar f}$.

To summarize the relationship between gerbes banded by $\sheaf{A}$ and extensions of $\Gal_S$ by $A$, we have the following.
\begin{thm}
To each geometric base point $\bar z \in \gerb{E}_{\bar s}$, we obtain an extension
$$A \Into \pi_1^{\et}(\gerb{E}, \bar z) \Onto \Gal_S,$$
known as the fundamental group of the gerbe $\gerb{E}$ at $\bar z$.  For any two geometric base points $\bar z_0, \bar z$, we obtain a {\em family} of isomorphisms of extensions
$$\pi_1^{\et}(\gerb{E}, \bar z_0) \To \pi_1^{\et}(\gerb{E}, \bar z),$$
any two of which are related by $\Int(\alpha)$ for a uniquely determined $\alpha \in A$.
\end{thm}

This theorem may seem more natural using 2-categorical language as follows:  consider the 2-category $\Cat{OpExt}(\Gal_S, \sheaf{A})$ whose objects are extensions of $\Gal_S$ by $A$ in which the $\Gal_S$ action on $A = \sheaf{A}_{\bar s}$ coincides with that which arises from the local system $\sheaf{A}$.  The morphisms in this category are isomorphisms of extensions (giving equality on $\Gal_S$ and $A$).  Given two such morphisms $\iota, \iota'$ sharing the same source and target, a natural transformation $\iota \Rightarrow \iota'$ is an element $\alpha \in A$ such that $\iota' = \Int(\alpha) \circ \iota$.  A restatement of the above theorem is the following.
\begin{thm}  
\label{AppxFundGpDefined}
A gerbe $\gerb{E}$ banded by $\sheaf{A}$, and a geometric point $\bar s \To S$, yield an object $\pi_1^{\et}(\gerb{E}, \bar s)$ of $\Cat{OpExt}(\Gal_S, \sheaf{A})$, well-defined up to equivalence, the equivalence being uniquely determined up to unique natural isomorphism.
\end{thm}

\begin{remark}
If $\sheaf{A}$ is a constant sheaf, then the extension $A \Into \pi_1^{\et}(\gerb{E}, \bar z) \Onto \Gal_S$ is a central extension.  It follows quickly from the theorem that we may define an extension
$$A \Into \pi_1^{\et}(\gerb{E}, \bar s) \Onto \Gal_S$$
up to unique isomorphism (without choice of geometric base point $\bar z$).  In this case, the 2-category $\Cat{OpExt}(\Gal_S, A)$ is an ordinary category:  the only transformations are the identities.
\end{remark}

Next, consider a functor of gerbes $\phi \From \gerb{E}_1 \To \gerb{E}_2$ lying over a homomorphism $c \From \sheaf{A}_1 \To \sheaf{A}_2$ of local systems of abelian groups.  If $\bar z_1$ is a geometric base point for $\gerb{E}_1$ over $\bar s$, arising from $z_1 \in \gerb{E}_1[U]$, then let $z_2 = \phi[U](z_1)$.  Define $\bar z_2$ to be the resulting geometric base point of $\gerb{E}_2$.  Our strictness assumption for functors of gerbes implies that the geometric base point $\bar z_2$ depends only on the geometric base point $\bar z_1$, and not on the choice of $z_1$.

For any $\gamma \in \Gal_S$, we obtain a map $\phi_\gamma \From \Aut_\gamma(z_1) \To \Aut_\gamma(z_2)$, given by
\begin{align*}
\Aut_\gamma(z_1) = \Hom(z_1, \gamma^\ast z_1) & \xrightarrow{\phi[U]} \Hom(\phi[U](z_1), \phi[U](\gamma^\ast z_1) ) \\
& = \Hom( z_2, \gamma^\ast \phi[U](z_1) ) \\
& = \Hom(z_2, \gamma^\ast z_2) = \Aut_\gamma(z_2).
\end{align*}
Putting these together yields a homomorphism of extensions,
$$\begin{tikzcd}
A_1 \inarrow{r} \arrow{d}{c} & \pi_1^{\et}(\gerb{E}_1, \bar z_1) \onarrow{r} \arrow{d}{\phi} & \Gal_S \arrow{d}{=} \\
A_2 \inarrow{r} & \pi_1^{\et}(\gerb{E}_2, \bar z_2) \onarrow{r} & \Gal_S
\end{tikzcd}$$

If $\phi, \phi' \From \gerb{E}_1 \To \gerb{E}_2$ are two functors of gerbes lying over $c \From \sheaf{A}_1 \To \sheaf{A}_2$, and $\phi(\bar z_1) = \phi'(\bar z_1) = \bar z_2$, then we find two such homomorphisms of extensions,
$$\phi, \phi' \From \pi_1^{\et}(\gerb{E}_1, \bar z_1) \To \pi_1^{\et}(\gerb{E}_2, \bar z_2),$$
lying over $c \From A_1 \To A_2$.  

If $N \From \phi \Rightarrow \phi'$ is a natural isomorphism of functors, then $N$ determines an isomorphism $\phi(\bar z_1) \To \phi'(\bar z_1)$, whence an isomorphism $\bar z_2 \To \bar z_2$.  Such an isomorphism is given by an element $\alpha_2 \in A_2 = \Aut(\bar z_2)$, and one may check that
$$\phi' = \Int(\alpha_2) \circ \phi \From \pi_1^{\et}(\gerb{E}_1, \bar z_1) \To \pi_1^{\et}(\gerb{E}_2, \bar z_2).$$

\printbibliography

\end{document}